\documentclass{article}
\usepackage{amssymb}
\usepackage{amsmath}

\input{tcilatex}

\begin{document}

\bigskip 

\bigskip

\bigskip

\bigskip

\section{Non-archimedean analysis on the extended \ \ \ \ \ \ \ \ \ \ \ \ \
\ \ \ \ \ \ \ hyperreal line $^{\ast }%
\mathbb{R}
_{\mathbf{d}}$ and some transcendence conjectures over field $%
\mathbb{Q}
$ and $^{\ast }%
\mathbb{Q}
_{\protect\omega }.$}

\bigskip\ \ \ \ \ \ \ \ \ \ \ \ \ \ \ \ \ \ \ \ \ \ \ \ \ \ \ \ \ 

\ \ \ \ \ \ \ \ \ \ \ \ \ \ \ \ \ \ \ \ \ \ \ \ \ \ \ \ \ \ \ \ \ \ \ \ \ \
\ \ \ \ \ \ \ \ Jaykov Foukzon

\ \ \ \ \ \ \ \ \ \ \ \ \ \ \ \ \ \ \ \ \ \ \ \ \ \ \ \ \ \ \ \ \ \ \ \ \
Israel Institute of Technology\ 

\ \ \ \ \ \ \ \ \ \ \ \ \ \ \ \ \ \ \ \ \ \ \ \ \ \ \ \ \ \ \ \ \ 

\ \ \ \ \ \ \ \ \ \ \ \ \ \ \ \ \ \ \ \ \ \ \ \ \ \ \ \ \ \ \ \ \ \ \ \ \ \
\ \ \ \ jaykovfoukzon@list.ru

\bigskip

\textbf{Abstract. }In this paper\textbf{\ }possible completion of the
Robinson non-archimedean field $^{\ast }%
\mathbb{R}
$ constructed by Dedekind sections. As interesting example I show how, a few
simple ideas from non-archimedean analysis on the pseudo-ring $^{\ast }%
\mathbb{R}
_{\mathbf{d}}$ gives a short clear nonstandard reconstruction for the
Euler's original proof of the Goldbach-Euler theorem. Given an analytic
function of one complex variable $f\in 
\mathbb{Q}
\left[ z\right] ,$we investigate the arithmetic nature of the values of $f$
at transcendental points.

\textbf{Contents}

Introduction.

\textbf{1.}Some transcendence conjectures over field $%
\mathbb{Q}
.$

\textbf{2.}Modern nonstandard analysis and non-archimedean analysis on

the extended hyperreal line $^{\ast }%
\mathbb{R}
_{\mathbf{d}}.$

\textbf{Chapter I.}The classical hyperreals numbers.

\textbf{I.1.1. \ }The construction non-archimedean field $^{\ast }%
\mathbb{R}
$.

\textbf{I.1.2. \ }The brief nonstandard vocabulary.

\textbf{I.2. \ \ \ \ }The higher orders of hyper-method.Second\textbf{\ }%
order\textbf{\ }transfer principle.

\textbf{I.2.1. \ }What are the higher orders of hyper-method?

\textbf{I.2.2. \ }The higher orders of hyper-method by using countable
universes.

\textbf{I.2.3. \ }Divisibility of hyperintegers.

\textbf{I.3. \ \ \ \ }The construction non-archimedean pseudo-ring $^{\ast }%
\mathbb{R}
_{\mathbf{d}}.$

\textbf{I.3.1. \ }Generalized pseudo-ring of Wattenberg-Dedekind hyperreals $%
^{\ast }%
\mathbb{R}
_{\mathbf{d}}$\textbf{\ }

\ \ \ \ \ \ \ \ \ and hyperintegers $^{\ast }%
\mathbb{Z}
_{\mathbf{d}}.$

\textbf{I.3.1.1. }Strong and weak Dedekind cuts.Wattenberg-Dedekind
hyperreals

\ \ \ \ \ \ \ \ \ \ and hyperintegers.

\textbf{I.3.2. \ \ }The topology of $^{\ast }%
\mathbb{R}
_{\mathbf{d}}.$

\textbf{I.3.3. \ \ \ }Absorption numbers in $^{\ast }%
\mathbb{R}
_{\mathbf{d}}.$

\textbf{I.3.3.1.} Absorption function and numbers in $^{\ast }%
\mathbb{R}
_{\mathbf{d}}$.

\textbf{I.3.3.2. }Special Kinds of Idempotents in $^{\ast }%
\mathbb{R}
_{\mathbf{d}}$.

\textbf{I.3.3.3.} Types of $\alpha $ with a given $\mathbf{ab.p.}(\alpha ).$

\textbf{I.3.3.4}. $\varepsilon $-Part of $\alpha $ with $\mathbf{ab.p.}%
(\alpha )\neq 0.$

\textbf{I.3.3.5. }Multiplicative idempotents.

\textbf{I.3.3.6. }Additive monoid of Dedekind hyperreal integers $^{\ast }%
\breve{%
\mathbb{Z}%
}_{\mathbf{d}}.$

\textbf{I.3.5. \ \ }Pseudo-ring of Wattenberg hyperreal integers $^{\ast }%
\mathbb{Z}
_{\mathbf{d}}.$

\textbf{I.3.6. \ \ }External summation of countable and hyperfinite
sequences in $^{\ast }%
\mathbb{R}
_{\mathbf{d}}.$

\textbf{I.3.7. \ \ }The construction non-archimedean field $^{\ast }%
\mathbb{R}
_{\mathbf{d}}^{\omega }$ as Dedekind

\ \ \ \ \ \ \ \ \ \ \ completion of countable non-standard models of $%
\mathbb{R}
.$\ 

\textbf{I.4. \ \ \ \ \ }The construction non-archimedean field $^{\ast }%
\mathbb{R}
_{\mathbf{c}}.$

\textbf{I.4.1. \ \ }Completion of ordered group and fields in general by
using

\ \ \ \ \ \ \ \ \ \ \ "Cauchy pregaps".

\textbf{I.4.1.1. }Totally ordered group and fields

\textbf{I.4.1.2. }Cauchy completion of ordered group and fields.

\textbf{I.4.2.1. }The construction non-archimedean field $^{\ast }%
\mathbb{R}
_{\mathbf{c}}$ by using Cauchy

\ \ \ \ \ \ \ \ \ \ \ \ hypersequence in ancountable field $^{\ast }%
\mathbb{Q}
.$

\textbf{I.4.2.2. }The construction non-archimedean field $^{\ast }%
\mathbb{R}
_{\mathbf{c}}^{\omega }$ as \ Cauchy

\ \ \ \ \ \ \ \ \ \ \ \ completion of countable non-standard models of $%
\mathbb{Q}
.$\ 

\textbf{Chapter II.}Euler's proofs by using non-archimedean analysis on the
\ \ \ \ \ \ \ \ \ \ \ \ \ \ \ \ \ \ \ \ \ \ \ \ \ \ \ \ \ \ \ \ \ \ \ \ \ \
\ \ \ 

\ \ \ \ \ \ \ \ \ \ \ \ \ \ \ \ \ pseudo-ring $^{\ast }%
\mathbb{R}
_{\mathbf{d}}$ revisited.

\textbf{II.1.}Euler's original proof of the Goldbach-Euler Theorem revisited.

\textbf{III. \ }Non-archimedean analysis on the extended hyperreal line $%
^{\ast }%
\mathbb{R}
_{\mathbf{d}}$ and

transcendence conjectures over field $%
\mathbb{Q}
.$Proof that $\ e+\pi $ and $e\cdot \pi $ is

irrational.\bigskip

\textbf{Chapter III.}Non-archimedean analysis on the extended hyperreal line 
$^{\ast }%
\mathbb{R}
_{\mathbf{d}}$

\ \ \ \ \ \ \ \ \ \ \ \ \ \ \ \ \ and transcendence conjectures over field $%
\mathbb{Q}
.$

\textbf{III.1. }Proof that $e$ is $\#$-transcendental and that $\ e+\pi $
and $e\cdot \pi $ is irrational.

\bigskip \textbf{III.2.} Nonstandard generalization of the Lindeman Theorem.

\bigskip\ \ \ \ \ \ \ \ \ \ \ \ \ \ \ \ \ \ \ \ \ \ \ \ \ \ \ \ \ \ \ \ \ \
\ \ \ \ \ \ \ \ \ \ \ \ \ \ \ \ \ \ \ \ \ \ \ \ \ \ \ \ \ \ \ \ \ \ \ \ \ \
\ \ \ \ \ \ \ 

\bigskip $\ \ \ \ \ \ \ \ \ \ \ \ \ \ \ \ \ \ \ \ \ \ \ $

\section{List of Notation.}

\ \ 

$^{\ast }%
\mathbb{N}
_{\infty }$..................................................... the set of
infinite natural numbers

$^{\ast }%
\mathbb{R}
_{\infty }$.................................................the set of
infinite hyper real numbers

$\mathbf{L}_{\ast }\mathbf{=L}\left( ^{\ast }%
\mathbb{R}
\right) $.........................................the set of the limited
members of $^{\ast }%
\mathbb{R}
$

$\mathbf{I}_{\ast }=\mathbf{I}\left( ^{\ast }%
\mathbb{R}
\right) $..................................the set of the infinitesimal
members of $^{\ast }%
\mathbb{R}
$

\textbf{halo}$\left( x\right) =\mu \left( x\right) =x+\mathbf{I}_{\ast }$%
............................................halo (monad) of $x\in 
\mathbb{R}
$

$\mathbf{st}\left( a\right) $%
......................................................Robinson standard part
of $a\in $ $^{\ast }%
\mathbb{R}
$

$^{\ast }%
\mathbb{Z}
\left[ z_{1},...,z_{\mathbf{n}}\right] $.....................internal
polynomials over $^{\ast }%
\mathbb{Z}
$ in $\ \mathbf{n}\in $ $^{\ast }%
\mathbb{N}
$ variables

\ 

$^{\ast }%
\mathbb{R}
\left[ z_{1},...,z_{\mathbf{n}}\right] $.....................internal
polynomials over $^{\ast }%
\mathbb{R}
$ in $\mathbf{n}\in $ $^{\ast }%
\mathbb{N}
$ variables

$^{\ast }%
\mathbb{C}
\left[ z_{1},...,z_{\mathbf{n}}\right] $.....................internal
polynomials over $^{\ast }%
\mathbb{C}
$ in $\mathbf{n}\in $ $^{\ast }%
\mathbb{N}
$ variables

$^{\ast }%
\mathbb{Z}
_{\mathbf{d}}$......................................................Dedekind
completion of the ring $^{\ast }%
\mathbb{Z}
$

$^{\ast }%
\mathbb{R}
_{\mathbf{d}}$......................................................Dedekind
completion of the field $^{\ast }%
\mathbb{R}
$

$^{\ast }%
\mathbb{R}
_{\mathbf{c}}$%
.........................................................Cauchy completion
of the field $^{\ast }%
\mathbb{R}
$

$\varepsilon _{\mathbf{d}}=\sup \left[ x|x\in \mu \left( 0\right) \right]
=\inf \left[ x|x\in 
\mathbb{R}
_{+}\right] $..............................$\mu \left( 0\right) $ $\subset $ 
$^{\ast }%
\mathbb{R}
_{\mathbf{d}},%
\mathbb{R}
_{+}$ $\subset $ $^{\ast }%
\mathbb{R}
_{\mathbf{d}}$

$\Delta _{\mathbf{d}}=\sup \left( 
\mathbb{R}
_{+}\right) =\inf \left( ^{\ast }%
\mathbb{R}
_{+\infty }\right) $%
.........................................................................

$WST\left( \alpha \right) $.........Wattenberg standard part of $\ \alpha
\in \left( -\Delta _{\mathbf{d}},\Delta _{\mathbf{d}}\right) _{^{\ast }%
\mathbb{R}
_{\mathbf{d}}}$ \ \ \ (Def.1.3.2.3)

\bigskip $\mathbf{ab.p.}\left( \alpha \right) $%
....................................absorption part of $\alpha \in $ $^{\ast
}%
\mathbb{R}
_{\mathbf{d}}$ \ (Def 1.3.3.1.1)

$\left[ \alpha \right] _{\varepsilon }$%
.........................................................................$%
\varepsilon $-part of $-\Delta _{\mathbf{d}}<\alpha <\Delta _{\mathbf{d}}$

$\left[ \alpha |b^{\#}\right] _{\varepsilon }$%
...................................................$\varepsilon $-part of $%
\alpha \in $ $^{\ast }%
\mathbb{R}
_{\mathbf{d}}$ for a given $b\in $ $^{\ast }%
\mathbb{R}
$

\ \ \ \ \ \ \ \ \ \ \ \ \ \ \ \ \ \ \ \ \ \ \ \ \ \ \ \ \ \ \ \ \ \ \ \ \ \
\ \ \ \ \ 

\bigskip

"Arthur stopped at the steep descent into the quarry, froze in his
steps,straining to look down and into the distance, extending his long
neck.Redrick joined him. But he did not look where Arthur was looking. Right
at their feet the road into the quarry began, torn up many years ago by the
treads and wheels of heavy vehicles.To the right was a white steep
slope,cracked by the heat; the next slope was half excavated, and among the
rocks and rubble stood a dredge, its lowered bucket jammed impotently
against the side of the road. And,as was to be expected, there was nothing
else to be seen on the road..."

\ \ \ \ \ \ \ \ \ \ \ \ \ \ \ \ \ \ \ \ \ \ \ \ \ \ \ \ \ \ \ \ \ \ \ \ \ \
\ \ \ \ \ \ \ \ \ \ \ \ \ \ \ \ \ \ \ \ \ \ \ \ \ \ \ \ \ \ \ \ \ \ \ \ \ 

\ \ \ \ \ \ \ \ \ \ \ \ \ \ \ \ \ \ \ \ \ \ \ \ \ \ \ \ \ \ \ \ \ \ \ \ \ \
\ \ \ \ \ \ \ \ \ \ \ \ \ \ \ \ \ \ \ \ \ \ \ \ \ \ \ Arkady and Boris
Strugatsky

\bigskip\ \ \ \ \ \ \ \ \ \ \ \ \ \ \ \ \ \ \ \ \ \ \ \ \ \ \ \ \ \ \ \ \ \
\ \ \ \ \ \ \ \ \ \ \ \ \ \ \ \ \ \ \ \ \ \ \ \ \ \ \ \ \ \ \ \ \ \ \ \ \ \
"Roadside Picnic"

\bigskip $\ \ \ \ \ \ \ \ \ \ \ \ \ \ \ \ \ \ \ \ \ \ \ $

\section{Introduction.}

\textbf{1.Some transcendence conjectures over field} $%
\mathbb{Q}
.$In 1873 French mathematician, Charles Hermite, proved that $e$ is
transcendental. Coming as it did 100 years after Euler had established the
significance of $e,$ this meant that the issue of transcendence was one
mathematicians could not afford to ignore.Within 10 years of Hermite's
breakthrough,his techniques had been extended by Lindemann and used to add $%
\pi $ to the list of known transcendental numbers. Mathematician then tried
to prove that other numbers such as $e+\pi $ and $e\times \pi $ are
transcendental too,but these questions were too difficult and so no further
examples emerged till today's time. The transcendence of $e^{\pi }$had been
proved in1929 by A.O.Gel'fond.

\textbf{Conjecture} \textbf{1. }The numbers $e+\pi $ and $e\times \pi $ are
irrational.

\bigskip \textbf{Conjecture} \textbf{2.} The numbers $e$ and $\pi $ are
algebraically independent.

However, the same question with $e^{\pi }$ and $\pi $ has been answered:

\textbf{Theorem.1.}(Nesterenko,1996 [22]) The numbers $e^{\pi }$ and $\pi $
are algebraically

independent.

During of XX th century,a typical question: is whether $f(\alpha )$ is a
transcen-

dental number for each algebraic number $\alpha $ has been investigated and

answered many authors.Modern result in the case of entire functions

satisfying a linear differential equation provides the strongest results,

related with Siegel's $E$-functions [22],[27].Ref. [22] contains references

to the subject before 1998, including Siegel $E$ and $G$ functions.

\textbf{Theorem.2.}(Siegel C.L.) Suppose that $\lambda \in 
\mathbb{Q}
,\lambda \neq -1,-2,...,\alpha \neq 0.$

\bigskip

\bigskip\ \ \ \ $%
\begin{array}{cc}
\begin{array}{c}
\\ 
\varphi _{\lambda }\left( z\right) =\sum_{n=0}^{\infty }\dfrac{z^{n}}{\left(
\lambda +1\right) \left( \lambda +2\right) \cdot \cdot \cdot \left( \lambda
+n\right) }. \\ 
\end{array}
& \text{ }\left( 1.1\right) \text{\ \ }%
\end{array}%
$ \ \ \ \ \ \ \ \ \ \ \ \ \ \ \ \ \ \ \ \ \ \ \ \ \ \ \ \ \ \ \ \ \ \ 

\bigskip

Then $\varphi _{\lambda }\left( \alpha \right) $ is a transcendental number
for each algebraic number $\alpha \neq 0.$

\bigskip

Given an analytic function of one complex variable $f\left( z\right) \in 
\mathbb{Q}
\left[ z\right] ,$we

investigate the arithmetic nature of the values of $f\left( z\right) $ at
transcendental

points.

\bigskip

\textbf{Conjecture} \textbf{3.}Is whether $f(\alpha )$ is a irrational
number for given

transcendental number $\alpha .$

\textbf{Conjecture} \textbf{4.}Is whether $f(\alpha )$ is a transcendental
number for given

transcendental number $\alpha .$

In particular we investigate the arithmetic nature of the values of
classical polylogarithms $Li_{s}\left( z\right) $ at transcendental
points.The classical polylogarithms

\bigskip

\bigskip\ \ \ \ \ \ \ \ \ \ \ \ \ \ \ \ \ \ \ \ \ \ \ \ \ \ \ $\ \ 
\begin{array}{cc}
\begin{array}{c}
\\ 
Li_{s}\left( z\right) =\sum_{n\geq 1}\dfrac{z^{n}}{n^{s}} \\ 
\end{array}
& \left( 1.2\right)%
\end{array}%
$\ \ \ \ \ \ 

\ \ \ \ \ \ \ \ \ \ \ \ \ \ \ \ \ \ \ \ \ \ \ \ \ \ \ \ \ \ \ \ \ \ \ \ \ \
\ \ \ \ \ \ \ \ \ \ \ \ \ \ \ \ \ \ \ \ \ \ 

for $s=1,2,...$ and $|z|\leq 1$ with $(s;z)=(1;1),$ are ubiquitous. The study

of the arithmetic nature of their special values is a fascinating subject

[35] very few is known.Several recent investigations concern the values

of these functions at $z=1:$ these are the values at the positive integers

of Riemann zeta function

\bigskip

\bigskip $\ \ \ \ \ \ \ \ \ \ \ \ \ \ \ \ \ \ \ \ \ \ \ \ \ \ \ \ \ \ \ \ \
\ \ \ 
\begin{array}{cc}
\begin{array}{c}
\\ 
\zeta \left( s\right) =Li_{s}\left( z=1\right) =\sum_{n\geq 1}\dfrac{1}{n^{s}%
} \\ 
\end{array}
& \left( 1.3\right)%
\end{array}%
$

\bigskip

One knows that $\zeta (3)$ is irrational [36],and that inInitely many values

$\zeta (2n+1)$ of the zeta function at odd integers are irrational.

\textbf{Conjecture} \textbf{4.}Is whether $Li_{s}\left( \alpha \right) $ is
a irrational number for given

transcendental number $\alpha .$

\textbf{Conjecture} \textbf{5.}Is whether $Li_{s}\left( \alpha \right) $ is
a transcendental number for given

transcendental number $\alpha .$

\textbf{2.Modern nonstandard analysis and} \textbf{non-archimedean analysis
on the extended hyperreal line} $^{\ast }%
\mathbb{R}
_{\mathbf{d}}.$Nonstandard analysis, in its early period of development,
shortly after having been established by A. Robinson [1],[4],[5] dealt
mainly with nonstandard extensions of some traditional mathematical
structures. The system of its foundations, referred to as "model-theoretic
foundations" was proposed by Robinson and E. Zakon [12]. Their approach was
based on the type-theoretic concept of superstructure $V(S)$ over some set
of individuals $S$ and its nonstandard extension (enlargement) $^{\ast
}V(S), $ usually constructed as a (bounded) ultrapower of the "standard"
superstructure $V(S).$They formulated few principles concerning the
elementary embedding $V(S)\longmapsto $ $^{\ast }V(S),$ enabling the use of
methods of nonstandard analysis without paying much attention to details of
construction of the particular nonstandard extension.

In classical Robinsonian nonstandard analysis we usualy deal only with
completely internal objects which can defined by internal set theory $%
\mathbf{IST}$ introduced by E.Nelson [11]. It is known that $\mathbf{IST}$
is a conservative extension of $ZFC.$ In $\mathbf{IST}$ all the classical
infinite sets, e.g., $%
\mathbb{N}
,%
\mathbb{Z}
,%
\mathbb{Q}
$ or $%
\mathbb{R}
,$ acquire new, nonstandard elements (like \textit{"infinite" natural numbers%
} or \textit{"infinitesimal" reals}). At the same time, the families $%
^{\sigma }%
\mathbb{N}
$ $=$ $\left\{ x\in 
\mathbb{N}
:\mathbf{st}\left( x\right) \right\} $ or $^{\sigma }%
\mathbb{R}
$ $=$ $\left\{ x\in 
\mathbb{R}
:\mathbf{st}\left( x\right) \right\} $ of all standard,i.e., "true," natural
numbers or reals, respectively, are not sets in $\mathbf{IST}$ at all. Thus,
for a traditional mathematician inclined to ascribe to mathematical objects
a certain kind of objective existence or reality, accepting $\mathbf{IST}$
would mean confessing that everybody has lived in confusion, mistakenly
having regarded as, e.g., the set $%
\mathbb{N}
$ just its tiny part $^{\sigma }%
\mathbb{N}
$ (which is not even a set) and overlooked the rest. Edvard Nelson and Karel
Hrb\u{c}ek have\ improved this lack by introducing several "nonstandard" set
theories dealing with standard, internal and \textit{external sets }[13].
Note that in contrast with early period of development of the nonstandard
analysis in latest period many mathematicians dealing with external and
internal set simultaneously,for example see [14],[15],[16],[17].

Many properties of the standard reals $x\in $ $%
\mathbb{R}
$ suitably reinterpreted, can be transfered to the internal hyperreal number
system. For example, we have seen that $^{\ast }%
\mathbb{R}
,$like $%
\mathbb{R}
,$is a totally ordered field. Also, jast $%
\mathbb{R}
$ contain the natural number $%
\mathbb{N}
$ as a discrete subset with its own characteristic properties,$^{\ast }%
\mathbb{R}
$ contains the hypernaturals $^{\ast }%
\mathbb{N}
$ as the corresponding discrete subset with analogous properties.For
example, the standard archimedean property

\bigskip

$\ \ \ 
\begin{array}{cc}
\begin{array}{c}
\\ 
\forall x_{x\in 
\mathbb{R}
}\forall y_{y\in 
\mathbb{R}
}\exists n_{n\in 
\mathbb{N}
}\left[ \left( \left\vert x\right\vert <\left\vert y\right\vert \right)
\dashrightarrow n\left\vert x\right\vert \geq \left\vert y\right\vert \right]
\\ 
\end{array}
& \left( 1.4\right)%
\end{array}%
$

$\ \ \ \ \ \ \ \ \ \ \ \ \ \ \ \ \ \ \ \ \ \ \ \ \ \ \ \ \ \ \ \ \ $

\bigskip

is preserved in non-archimedean field $^{\ast }%
\mathbb{R}
$ in respect hypernaturals $^{\ast }%
\mathbb{N}
,$i.e. the next property is satisfied

$\ \ \ \ \ 
\begin{array}{cc}
\begin{array}{c}
\\ 
\ \forall x_{x\in 
\mathbb{R}
}\forall y_{y\in 
\mathbb{R}
}\exists n_{n\in 
\mathbb{N}
}\left[ \left( \left\vert x\right\vert <\left\vert y\right\vert \right)
\dashrightarrow n\left\vert x\right\vert \geq \left\vert y\right\vert \right]
. \\ 
\end{array}
& \left( 1.5\right)%
\end{array}%
\ \ \ \ $

$\bigskip $

$\ \ \ \ \ \ \ \ \ \ \ \ \ \ \ \ \ \ \ \ \ \ \ \ $

However, there are many fundamental properties of $%
\mathbb{R}
$ do not transfered to $^{\ast }%
\mathbb{R}
.$

\textbf{I}. This is the case one of the fundamental \textit{supremum} 
\textit{property} of the standard totally ordered field $%
\mathbb{R}
.$It is easy to see that it apper bound property does not necesarily holds
by considering, for example, the (external) set $%
\mathbb{R}
$ itself which we ragard as canonically imbedded into hyperreals $^{\ast }%
\mathbb{R}
.$ This is a non-empty set which is bounded above (by any of the infinite
member in $^{\ast }%
\mathbb{R}
$) but does not have a least apper bound in $^{\ast }%
\mathbb{R}
.$ However by using transfer one obtain the next statement [18] :

\bigskip

\textbf{Weak supremum property for }$^{\ast }%
\mathbb{R}
:$Every non-empty \textit{internal} subset

$A\subsetneqq $ $^{\ast }%
\mathbb{R}
$ which has an apper bound in $^{\ast }%
\mathbb{R}
$ has a least apper bound in $^{\ast }%
\mathbb{R}
.$

This is a problem, because any advanced variant of the analysis on the field 
$^{\ast }%
\mathbb{R}
$ is needed more strong fundamental supremum property. At first sight one
can improve this lack by using corresponding external constructions which
known as Dedekind sections and Dedekind completion (see section \textbf{I.3.}%
).We denote corresponding Dedekind completion by symbol $^{\ast }%
\mathbb{R}
_{\mathbf{d}}$. It is clear that $^{\ast }%
\mathbb{R}
_{\mathbf{d}}$ is completely external object. But unfortunately $^{\ast }%
\mathbb{R}
_{\mathbf{d}}$ is not iven a non-archimedean ring but non-archimedean 
\textit{pseudo-ring }only. However this lack does not make greater
difficulties because non-archimedean pseudo-ring $^{\ast }%
\mathbb{R}
_{\mathbf{d}}$ contains non-archimedean subfield $\mathbf{\Re }_{\mathbf{c}%
}\subset $ $^{\ast }%
\mathbb{R}
_{\mathbf{d}}$ such that $\mathbf{\Re }_{\mathbf{c}}\approx $ $^{\ast }%
\mathbb{R}
_{\mathbf{c}}.$Here $^{\ast }%
\mathbb{R}
_{\mathbf{c}}$\textit{\ }this is a Cauchy completion of the non-archimedean
field $^{\ast }%
\mathbb{R}
$ (see section \textbf{I.4.}).

\textbf{II. }This is the case two of the fundamental \textit{Peano's
induction} \textit{property:}

$\ \ \ \ \ \ \ \ \ \ \ \ \ \ \ \ \ \ 
\begin{array}{cc}
\begin{array}{c}
\\ 
\forall B\left[ \left[ \left( 1\in B\right) \wedge \forall x\left( x\in
B\implies x+1\in B\right) \right] \implies B=%
\mathbb{N}
\right] \\ 
\end{array}
& \text{ \ \ \ \ \ }\left( 2.1\right)%
\end{array}%
$

\bigskip

does not necesarily holds for arbitrary subset $B\subset $ $^{\ast }%
\mathbb{N}
.$ Therefore (2.1)

is true for $^{\ast }%
\mathbb{N}
$ when interpreted in $^{\ast }%
\mathbb{N}
$ i.e.,

\bigskip

$\ \ \ \ \ \ \ \ \ \ 
\begin{array}{cc}
\begin{array}{c}
\\ 
\forall ^{\mathbf{int}}B\left[ \left[ \left( 1\in B\right) \wedge \forall
x\left( x\in B\implies x+1\in B\right) \right] \implies B=\text{ }^{\ast }%
\mathbb{N}
\right] \\ 
\end{array}
& \text{ \ \ \ \ \ \ }\left( 2.2\right)%
\end{array}%
$

\bigskip

is true for $^{\ast }%
\mathbb{N}
$ provided that we read "$\forall B$" as "for each internal subset $B$ of $%
^{\ast }%
\mathbb{N}
$", i.e. as $\forall ^{\mathbf{int}}B.$\ In general the importance of
internal versus external entities rests on the fact that each statement that
is true for $%
\mathbb{R}
$ is true for $^{\ast }%
\mathbb{R}
$ provided its quantifiers are restricted to the internal entities (subset)
of $^{\ast }%
\mathbb{R}
$ only [18].This is a problem, because any advanced variant of the analysis
on the field $^{\ast }%
\mathbb{R}
$ is needed more strong induction property\textit{\ }than property\textit{\ }%
(2).In this paper I have improved this lack by using external construction
two different types for operation of exteral summation:

$\bigskip $

$\ \ \ \ \ 
\begin{array}{cc}
\begin{array}{c}
\\ 
\ Ext-\dsum\limits_{n\in 
\mathbb{N}
}q_{n}, \\ 
\\ 
\#Ext-\dsum\limits_{n\in 
\mathbb{N}
}q_{n}^{\#} \\ 
\end{array}
& \left( 2.3\right)%
\end{array}%
\ \ \ \ \ \ \ \ \ \ \ $

$\ \ \ \ \ \ \ \ \ \ \ \ \ \ \ \ \ \ \ \ \ \ \ \ \ \ \ \ \ $

and two different types for operation of exteral multiplication:

$\bigskip \ \ 
\begin{array}{cc}
\begin{array}{c}
\\ 
Ext-\dprod\limits_{n\in 
\mathbb{N}
}q_{n}, \\ 
\\ 
\#Ext-\dprod\limits_{n\in 
\mathbb{N}
}q_{n}^{\#} \\ 
\end{array}
& \text{ \ \ \ \ \ \ }\left( 2.4\right)%
\end{array}%
$

for arbitrary countable sequences \ such as $q_{n}:%
\mathbb{N}
\rightarrow $ $%
\mathbb{R}
$ and $q_{n}^{\#}:%
\mathbb{N}
\rightarrow $ $^{\ast }%
\mathbb{R}
_{\mathbf{d}}.$

As interesting example I show how, this external constructions from
non-archimedean analysis on the pseudo-ring $^{\ast }%
\mathbb{R}
_{\mathbf{d}}$ gives a short and clear nonstandard reconstruction for the
Euler's original proof of the Goldbach-Euler theorem.

\bigskip

\section{I.The classical hyperreals numbers.}

\section{I.1.1.The construction non-archimedean field $^{\ast }%
\mathbb{R}
$.}

Let $\Re $ denote the ring of real valued sequences with the usual pointwise
operations.If $x$ is a real number we let $s_{x}$ denote the constant
sequence,$\mathbf{s}_{x}=x$ for all $n.$ The function sending $x$ to $%
\mathbf{s}_{x}$ is a one-to-one ring homomorphism,providing an embedding of $%
\mathbb{R}
$ into $\Re $. In the following, wherever it is not too confusing we will
not distinguish between $x\in $ $%
\mathbb{R}
$ and the constant function $\mathbf{s}_{x},$leaving the reader to derive
intent from context. The ring $\Re $ has additive identity $0$ and
multiplicative identity $1.$ $\Re $ is not a field because if $r$ is any
sequence having $0$ in its range it can have no multiplicative inverse.
There are lots of zero divisors in $\Re $.

We need several definitions now. Generally, for any set $S,$ $\mathbf{P}(S)$
denotes the set of all subsets of $S.$ It is called the power set of $S.$
Also, a subset of $%
\mathbb{N}
$ will be called cofinite if it contains all but finitely many members of $%
\mathbb{N}
$. The symbol $\varnothing $ denotes the empty set. A partition of a set $S$
is a decomposition of $S$ into a union of sets, any pair of which have no
elements in common.

\textbf{Definition.1.1.1.} An ultrafilter $\mathbf{H}$ over $%
\mathbb{N}
$ is a family of sets for which:

\textbf{(i)} $\varnothing \notin \mathbf{H\subset P(}%
\mathbb{N}
\mathbf{),%
\mathbb{N}
\in H.}$

\textbf{(ii)} Any intersection of finitely many members of $\mathbf{H}$ is
in $\mathbf{H.}$

\textbf{(iii)} $A\subset 
\mathbb{N}
,B\in \mathbf{H}\Rightarrow A\cup B\in \mathbf{H.}$

\textbf{(iv)} If $V_{1},...,V_{n}$ is any finite partition of $%
\mathbb{N}
$ then $\mathbf{H}$ contains exactly

one of the $V_{i}.$

If, further,

\textbf{(v)} $\mathbf{H}$ contains every cofinite subset of $%
\mathbb{N}
$.

the ultrafilter is called \textit{free.}

If an ultrafilter on $%
\mathbb{N}
$ contains a finite set then it contains a one-point set,

and is nothing more than the family of all subsets of $%
\mathbb{N}
$ containing that

point. So if an ultrafilter is not free it must be of this type, and is
called

a \textit{principal} \textit{ultrafilter.}

The existence of a free ultrafilter containing any given infinite subset of $%
\mathbb{N}
$

is implied by the Axiom of Choice.

\textbf{Remark 1.1.1. }Suppose that $x\in X.$ An ultrafilter denoted

$\mathbf{prin}_{X}\left( x\right) \subseteq X$ \ consisting of all subsets $%
S\subseteq X$ which contain $x,$ and called

the \textit{principal ultrafilter} generated by $x.$

\textbf{Proposition} \textbf{1.1.1.} If an ultrafilter $\mathbf{\tciFourier }
$ on $X$ contains a finite set $S\subseteq X,$

then $\mathbf{\tciFourier }$ is \ principal.

\textbf{Proof:} It is enough to show $\mathbf{\tciFourier }$ contains $%
\left\{ x\right\} $ for some $x\in S.$ If not, then

$\mathbf{\tciFourier }$ contains the complement $X\backslash \left\{
x\right\} $ for every $x\in S$, and therefore also

the finite intersection $\mathbf{\tciFourier }\ni \dbigcap\limits_{x\in
S}X\backslash \left\{ x\right\} =X\backslash S,$ which contradicts the fact

that $S\in \tciFourier .$It follows that nonprincipal ultrafilters can exist
only on infinite

sets $X$, and that every cofinite subset of $X$ (complement of a finite set)

belongs to such an ultrafilter.

$\mathbf{Remark}$ $\mathbf{1.1.2.}$Our construction below depends on the use
of a free-not a

principal-ultrafilter.

We are going to be using conditions on sequences and sets to define subsets
of $%
\mathbb{N}
$. We introduce a convenient shorthand for the usual \textquotedblleft set
builder\textquotedblright\ notation. If $P$ is a property that can be true
or false for natural numbers we use $[[P]]$ to denote $\{n\in 
\mathbb{N}
|P(n)$ is true $\}.$ This notation will only be employed during a discussion
to decide if the set of natural numbers defined by $P$ is in $\mathbf{H,}$
or not. For example, if $s,t$ is a pair of sequences in $\Re $ we define
three sets of integers For example, if $s,t$ is a pair of sequences in $S$
we define three sets of integers

$\ \ \ \ \ \ \ \ \ \ \ 
\begin{array}{cc}
\begin{array}{c}
\\ 
\ [[s<t]],\ \ [[s=t]],\ [[s>t]]. \\ 
\end{array}
& \text{ \ \ \ \ \ \ \ \ \ }\left( 1.1.1\right)%
\end{array}%
\ \ $

Since these three sets partition $%
\mathbb{N}
$, exactly one of them is in $\mathbf{H,}$ and we \ \ \ \ \ \ \ \ \ \ \ \ \
\ \ \ \ \ \ \ \ \ 

declare $s\equiv t$ when $[[s=t]]\in \mathbf{H.}$

\textbf{Lemma 1.1.1.} $\equiv $ is an equivalence relation on $\Re $. We
denote the

equivalence class of any sequence $s$ under this relation by $[s].$

Define for each $r\in \Re $ the \ sequence $\tilde{r}$ by

\ \ \ \ \ \ \ \ \ \ \ \ \ \ \ \ \ \ \ 

$\ \ \ \ \ \ \ \ \ \ 
\begin{array}{cc}
\begin{array}{c}
\\ 
\ \ \ \tilde{r}=\left\{ 
\begin{array}{c}
0\text{ \ \ }\mathbf{iff}\text{ }r_{n}=0 \\ 
r_{n}^{-1}\text{ }\mathbf{iff}\text{ }r_{n}\neq 0%
\end{array}%
\right\} . \\ 
\end{array}
& \text{ \ \ \ \ \ \ \ \ \ \ \ \ \ \ \ \ \ \ \ \ \ \ \ \ \ }\left(
1.1.2\right) 
\end{array}%
$

\textbf{Lemma} \textbf{1.1.2.} (\textbf{a}) There is at most one constant
sequence in any class $[r].$

(\textbf{b}) $[0]$ is an ideal in $\Re $ so $\Re /[0]$ is a commutative ring
with identity [1].

(\textbf{c}) Consequently $[r]=r+[0]=\{r+t|t\in \lbrack 0]\}$ for all $r\in
\Re $.

(\textbf{d}) If $[r]\neq \lbrack 0]$ then $[\tilde{r}]\cdot \lbrack r]=[1].$
So $[r]^{-1}=[\tilde{r}].$

From Lemma 1.1.2., we conclude that $^{\ast }%
\mathbb{R}
,$ defined to be $\Re /[0],$ is a field \ \ \ \ \ \ \ \ \ \ \ \ \ \ \ \ \ \
\ \ \ \ \ \ 

containing an embedded image of as a subfield. $[0]$ is a maximal ideal in

$\Re $. \ \ \ \ \ \ \ \ \ \ \ \ \ \ \ \ \ \ \ \ \ \ \ \ \ \ \ \ \ \ \ \ \ \
\ \ \ \ \ \ \ \ \ \ 

\textbf{Definition.1.1.2.}This quotient ring is called the field $^{\ast }%
\mathbb{R}
$ of \textbf{classical }

\textbf{hyperreal numbers}.

We declare $[s]<[t]$ provided $[[s<t]]\in \mathbf{H.}$

Recall that any field with a linear order $<$ is called an ordered field

provided:

\textbf{(i)} $x+y>0$ whenever $x,y>0$

\textbf{(ii)} $x\cdot y>0$ whenever $x,y>0$

\textbf{(iii)} $x+z>y+z$ whenever $x>y$

\textbf{Theorem 1.1.3.} (\textbf{a}) The relation given above is a linear
order on $^{\ast }%
\mathbb{R}
,$ and

makes $^{\ast }%
\mathbb{R}
$ into an ordered field. As with any ordered field, we define $|x|$

for $x\in $ $^{\ast }%
\mathbb{R}
$ to be $x$ or $-x,$ whichever is nonnegative.

\textbf{(b)} If $x,y$ are real then $x\leq y$ if and only if $[x]\leq
\lbrack y].$ So the ring morphism

of $%
\mathbb{R}
$ into $^{\ast }%
\mathbb{R}
$ is also an order isomorphism onto its image in $^{\ast }%
\mathbb{R}
.$

Because of this last theorem and the essential uniqueness of the real
numbers it is common to identify the embedded image of $%
\mathbb{R}
$ in $^{\ast }%
\mathbb{R}
$ with $%
\mathbb{R}
$ itself. Though obviously circular, one does something similar when
identifying $%
\mathbb{Q}
$ with its isomorphic image in $^{\ast }%
\mathbb{R}
$, and $%
\mathbb{N}
$ itself with the corresponding subset of $%
\mathbb{Q}
$. This kind of notational simplification usually does not cause problems.

Now we get to the ideas that prompted the construction. Define the sequence $%
r$ by $r_{n}=\left( n+1\right) ^{-1}$ . For every positive integer $%
k,[[r<k^{-1}]]\in \mathbf{H.}$So $0<[r]<1/k.$ We have found a positive
hyperreal smaller than (the embedded image of) any real number. This is our
first nontrivial infinitesimal number. The sequence $\tilde{r}$ is given by $%
\tilde{r}_{n}=n+1.$So $[r]^{-1}=[\tilde{r}]>k$ for every positive integer $%
k.[r]^{-1}$ is a hyperreal larger than any real number.

\bigskip

\bigskip

\section{I.1.2.The brief nonstandard vocabulary.}

\bigskip

\textbf{Definition.1.1.2.1. }We call a member $x\in $ $^{\ast }%
\mathbb{R}
$ $%
\mathbb{R}
$-\textit{limited} if there are members $\ \ \ \ \ \ \ \ \ \ \ \ \ \ $ $\ \
\ \ \ \ a,b\in 
\mathbb{R}
$ with $a<x<b.$

We will use $\mathbf{L}_{\ast }\mathbf{=L}\left( ^{\ast }%
\mathbb{R}
\right) $ to indicate the limited members of $^{\ast }%
\mathbb{R}
$. $x$ is called

$%
\mathbb{R}
$-\textit{unlimited} if it not $%
\mathbb{R}
$-limited.

These terms are preferred to \textquotedblleft finite\textquotedblright\ and
\textquotedblleft infinite,\textquotedblright

which are reserved for concepts related to cardinality.

\textbf{Definition.1.1.2.2. }If $x,y\in $ $^{\ast }%
\mathbb{R}
$ and $x<y$ we use $^{\ast }[x,y]$ to denote

$\{t\in $ $^{\ast }%
\mathbb{R}
|x\leq t\leq y\}.$ \ \ \ \ \ \ \ \ \ 

This set is called a \textit{closed hyperinterval.} Open and half-open

hyperintervals are defined and denoted similarly.

\textbf{Definition.1.1.2.3. }A set $S\subset $ $^{\ast }%
\mathbb{R}
$ is called \textit{hyperbounded} if there are

members $x,y$ of $^{\ast }%
\mathbb{R}
$ for which $S$ is a subset of the hyperinterval

$^{\ast }[x,y].$Abusing standard vocabulary for ordered sets, $S$ is called
bounded

if $x$ and $y$ can be chosen to be limited members of $^{\ast }%
\mathbb{R}
$ $x$ and $y$ could, in

fact, be chosen to be real if $S$ is bounded.

\textit{The vocabulary of bounded or hyperbounded above and below can be
used.}

\textbf{Definition.1.1.2.3. }We call a member $x\in $ $^{\ast }%
\mathbb{R}
$ infinitesimal if $|x|<a$ for

every positive $a\in 
\mathbb{R}
.$ We write $x\approx 0$ iff $x$ is infinitesimal.

The only real infinitesimal is obviously $0.$

We will use $\mathbf{I}_{\ast }=\mathbf{I}\left( ^{\ast }%
\mathbb{R}
\right) $ to indicate the infinitesimal members of $^{\ast }%
\mathbb{R}
.$

\textbf{Definition.1.1.2.4. }A member $x\in $ $^{\ast }%
\mathbb{R}
$ is called \textit{appreciable} if it is limited

but not infinitesimal.

\textbf{Definition.1.1.2.5. }Hyperreals $x$ and $y$ are said to have

\textit{appreciable separation} if $|x-y|$ is appreciable.

We will be working with various subsets $S$ of $^{\ast }%
\mathbb{R}
$ and adopt the

following convention: $S_{\infty }=S\backslash \mathbf{L}_{\ast }=\{x\in
S|x\notin \mathbf{L}_{\ast }\}.$ These are the

unlimited members of $S,$if any.

\textbf{Definition.1.1.2.6. }(\textbf{a}) We say two hyperreals $x,y$ are
infinitesimally close

or have infinitesimal separation if $|x-y|\in \mathbf{I}_{\ast }.$

We use the notation $x\approx y$ to indicate that $x$ and $y$ are
infinitesimally close.

(\textbf{b}) They have limited separation if $|x-y|\in \mathbf{L}_{\ast }%
\mathbf{.}$ \ \ \ \ \ \ \ \ \ \ \ \ \ \ 

(\textbf{c}) Otherwise they are said to have unlimited separation.

We define the halo of $x$ by $\mathbf{halo}(x)=x+\mathbf{I}_{\ast }.$ There
can be at most one real

number in any halo. Whenever $\mathbf{halo}(x)\cap 
\mathbb{R}
$ is nonempty we define the

shadow of $x,$ denoted $\mathbf{shad}(x),$ to be that unique real number.

The galaxy of $x$ is defined to be $\mathbf{gal}(x)=x+\mathbf{L}_{\ast }%
\mathbf{.}$ $\mathbf{gal}(x)$ is the set of

hyperreal numbers $a$ limited distance away from $x.$ So if $x$ is limited

$\mathbf{gal}(x)=\mathbf{L}_{\ast }\mathbf{.}$

If n is any fixed positive integer we define $^{\ast }%
\mathbb{R}
^{n}$ to be the set of

equivalence classes of sequences in $%
\mathbb{R}
^{n}$ under the equivalence

relation $x\equiv y$ exactly when $[[x=y]]\in \mathbf{H.}$

\textbf{Definition.1.1.2.7. }We call $^{\ast }%
\mathbb{N}
$ the set of \textit{classical hypernatural or }

\textit{A. Robinson's hypernatural numbers,} $^{\ast }%
\mathbb{N}
_{\infty }$ the set of \textit{classical infinite }

\textit{hypernatural or \ A. Robinson's infinite hypernatural \ numbers,}$%
^{\ast }%
\mathbb{R}
_{\infty }$ the

set of \textit{classical infinite hyperreal \ or A. Robinson's infinite
hyperreal }

\textit{numbers,}$^{\ast }%
\mathbb{Z}
$ the set of \textit{classical hyperintegers or A. Robinson's}\ 

\textit{hyperintegers}, and $^{\ast }%
\mathbb{Q}
$ the set of \textit{classical hyperrational numbers or }

\textit{A. Robinson's} \textit{hyperrational numbers.}

\textbf{Theorem 1.1.2.1.} $^{\ast }%
\mathbb{R}
$ $\mathbf{is}$\textbf{\ }$\mathbf{not}$\textbf{\ }$\mathbf{Dedekind}$%
\textbf{\ }$\mathbf{complete}.$

(hint: $%
\mathbb{N}
$ is bounded above by the member $[\mathbf{t}]\in $ $^{\ast }%
\mathbb{N}
_{\infty },$ where $\mathbf{t}$ is the sequence

given by $t_{n}=n$ for all $n\in 
\mathbb{N}
.$ But $%
\mathbb{N}
$ can have no least upper bound: if $n\leq c$

for all $n\in 
\mathbb{N}
$ then $n\leq c-1$ for all $n\in $ $%
\mathbb{N}
$.

As another example consider $\mathbf{I}_{\ast }\mathbf{.}$ This set is
(very) bounded, but has no

least upper bound.)

\textbf{Theorem 1.1.2.2. }For every $r\in $ $^{\ast }%
\mathbb{R}
$ there is unique $n\in $ $^{\ast }%
\mathbb{N}
$ with

$n\leq r<n+1.$

\bigskip

\section{I.2.The higher orders of hyper-method.Second\textbf{\ }order\textbf{%
\ \ }transfer principle.}

\bigskip

\section{I.2.1.What are the higher orders of hyper-method?}

Usual nonstandard analysis essentially consists only of two fundamental
tools:the (\textbf{first order}) \textbf{star-map} $\ast _{1}\triangleq \ast 
$ and the (\textbf{first order}) \textbf{transfer principle}. In most
applications, a third fundamental tool is also considered, namely the
saturation property.

\textbf{Definition.1.2.1.1. }Any universe $\mathbf{U}$ is a nonempty
collection of

\textit{"standard mathematical} \textit{objects" }that is closed under
subsets,

i.e. $a\subseteqq A\in \mathbf{U}$ $\implies a\in \mathbf{U}$ and closed
under the basic mathematical

operations.Precisely, whenever $A,B\in \mathbf{U},$ we require that also the

union $A\cup B$, the intersection $A\cap B,$the set-difference $A\backslash
B $ the

ordered pair $\left\{ A,B\right\} ,$the Cartesian product $A\times B,$ the

powerset $P(A)=\left\{ a|a\subseteqq A\right\} ,$the function-set $%
B^{A}=\left\{ f\text{ }|\text{ }f:A\rightarrow B\right\} $ all

belong to $\mathbf{U.}$A universe $\mathbf{U}$ is also assumed to contain
(copies of) all

sets of numbers $%
\mathbb{N}
,$ $%
\mathbb{Z}
,%
\mathbb{Q}
,%
\mathbb{R}
,%
\mathbb{C}
$ $\in \mathbf{U,}$ and to be transitive, i.e. members

of members of $\mathbf{U}$ belong to $\mathbf{U}$ or in formulae: $a\in A\in 
\mathbf{U}$ $\implies a\in \mathbf{U}$.

\bigskip

The notion of "\textit{standard }mathematical object" includes all objects
used

in the ordinary practice of mathematics, namely: numbers, sets,

functions, relations, ordered tuples, Cartesian products, etc. It is well-

known that all these notions can be defined as sets and formalized in

the foundational framework of Zermelo-Fraenkel axiomatic set theory

$\mathbf{ZFC.}$

\bigskip

From standard assumption: $Con\left( \mathbf{ZFC}\right) $ and G\"{o}del's
completeness theorem

one obtain that $\mathbf{ZFC}$ has a model $M$. A model $M$ of set theory is
called

standard if the element relation $\in _{M}$of the model $M$ is the actual
element

relation $\in $ restricted to the model $M,$ i.e.$\in _{M}\triangleq $ $\in
|_{M}$ . A model is called

transitive when it is standard and the base class is a transitive class of

sets. A model of set theory is often assumed to be transitive unless it is

explicitly stated that it is non-standard. Inner models are transitive,

transitive models are standard, and standard models are well-founded.

The assumption that there exists a standard model of $\mathbf{ZFC}$ (in a
given

universe) is stronger than the assumption that there exists a model. In

fact, if there is a standard model, then there is a smallest standard

model called the minimal model contained in all standard models. The

minimal model contains no standard model (as it is minimal) but

(assuming the consistency of $\mathbf{ZFC}$) it contains some model of $%
\mathbf{ZFC}$ by

the G\"{o}del completeness theorem. This model is necessarily not well

founded otherwise its Mostowski collapse would be a standard model.

(It is not well founded as a relation in the universe, though it

satisfies the axiom of foundation so is "internally" well founded.

Being well founded is not an absolute property[2].) In particular in the

minimal model there is a model of $\mathbf{ZFC}$ but there is no standard
model

of $\mathbf{ZFC.}$

\bigskip

By the theorem of L\"{o}wenheim-Skolem, we can choose transitive models

$M_{\omega }$ of $\mathbf{ZFC}$ of countable cardinality.

\bigskip

\textbf{Remark 1.2.1.1.}In $\mathbf{ZFC,}$ an ordered pair $\left\{
a,b\right\} $ is defined as the

Kuratowski pair $\left\{ \left\{ a\right\} ,\left\{ a,b\right\} \right\} ;$%
an $n$-tuple is inductively defined by

$\left\{ a_{1},...,a_{n},a_{n+1}\right\} $ $=\left\{ \left\{
a_{1},...,a_{n}\right\} ,a_{n+1}\right\} ;$ an $n$-place relation $R$ on $A$

is identified with the set $R\subseteq A^{n}$ of $n$-tuples that satisfy it;
a function

$f:A\rightarrow B$ is identi\U{407}ed with its graph $\left[ \left\{
a,b\right\} \in A\times B|b=f(a)\right] $.

As for numbers, complex numbers $%
\mathbb{C}
=$ $%
\mathbb{R}
\times 
\mathbb{R}
$ are defined as ordered pairs of

real numbers, and the real numbers $%
\mathbb{R}
$ are defined as equivalence classes

of suitable sets of rational numbers namely, Dedekind cuts or Cauchy

sequences.

The rational numbers $%
\mathbb{Q}
$ are a suitable quotient $%
\mathbb{Z}
\times 
\mathbb{Z}
/_{\approx },$ and the integers

$%
\mathbb{Z}
$ are in turn a suitable quotient $%
\mathbb{N}
\times 
\mathbb{N}
/_{\approx }$. The natural numbers of $\mathbf{ZFC}$ are

defined as the set of von Neumann naturals: $0=\NEG{0}$ and $n+1=\left\{
n\right\} $

(so that each natural number $\left\{ n=0,1,...,n-1\right\} $ is identified
with the set

of its predecessors.)

\bigskip

Each countable model $M_{\omega }$ of $\mathbf{ZFC}$ contains countable
model $%
\mathbb{R}
_{\omega }$ of the

real numbers $%
\mathbb{R}
.$Every element $x\in 
\mathbb{R}
_{\omega }$ defines a Dedekind cut:

$%
\mathbb{Q}
_{x}\triangleq \left\{ q\in 
\mathbb{Q}
|q\leq x\right\} \cup \left\{ q\in 
\mathbb{Q}
|q>x\right\} .$We therefore get a order preserving

map

\ \ \ \ \ \ \ \ \ $\ \ 
\begin{array}{cc}
\begin{array}{c}
\\ 
f_{p}:%
\mathbb{R}
_{\omega }\rightarrow 
\mathbb{R}
\\ 
\end{array}
& \text{ \ \ }\left( 1.2.1\right)%
\end{array}%
$

and which respects addition and multiplication.

\bigskip We address the question what is the possible range of $f_{p}?$

\textbf{Proposition} \textbf{1.2.1.1. }Choose an arbitrary subset $\Theta
\subset 
\mathbb{R}
.$Then there is a

model $%
\mathbb{R}
\left( \Theta \right) $ such that $f_{p}\left[ 
\mathbb{R}
\left( \Theta \right) \right] \supset \Theta .$Moreover, the cardinality of $%
\mathbb{R}
\left( \Theta \right) $ can

be chosen to coincide with $\Theta $, if $\Theta $ is infinite.

\textbf{Proof. }Choose $\Theta \subset 
\mathbb{R}
.$For each $\alpha \in \Theta $ choose

$q_{1}\left( \alpha \right) <q_{2}\left( \alpha \right) <...<p_{1}\left(
\alpha \right) <p_{2}\left( \alpha \right) $ with $\underset{n\rightarrow
\infty }{\lim }q_{n}\left( \alpha \right) =$ $\underset{n\rightarrow \infty }%
{\lim }p_{n}\left( \alpha \right) =\alpha .$

\bigskip We add to the axioms of $%
\mathbb{Q}
$ the following axioms:

$\ \ \ \ \ \ \ \ \ 
\begin{array}{cc}
\begin{array}{c}
\\ 
\ \ \forall \alpha \left( \alpha \in \Theta \right) \exists e_{\alpha
}\forall k\left( k\in 
\mathbb{N}
\right) \left[ q_{k}\left( \alpha \right) <e_{\alpha }<p_{k}\left( \alpha
\right) \right] \\ 
\end{array}
& \text{ \ \ \ \ \ \ \ }\left( 1.2.2\right)%
\end{array}%
$

Again $%
\mathbb{Q}
$ is a model for each finite subset of these axioms,so that the

compactness theorem implies the existence of $%
\mathbb{R}
\left( \Theta \right) $ as required,where

the cardinality of $%
\mathbb{R}
\left( \Theta \right) $ can be chosen to be the cardinality of the set of

axioms, i.e. of $\Theta $, if $\Theta $ is infinite. Note that by
construction $f_{p}\left( e_{\alpha }\right) =\alpha .$

$\bigskip $

\textbf{Remark 1.2.1.2. }It follows by the theorem of L\"{o}%
wenheim-Skolem,that

for each countable subset $\Theta _{\omega }\subset 
\mathbb{R}
$ we can find a countable model

$%
\mathbb{R}
_{\omega }=%
\mathbb{R}
\left( \Theta _{\omega }\right) $ of $%
\mathbb{R}
$ such that the image of $f_{p}$ contains this subset.Note,

on the other hand, that the image will only be countable, so that the

different models $%
\mathbb{R}
\left( \Theta _{\omega }\right) $ will have very different ranges.

$\ \ \ \ \ \ \ \ \ \ \ \ \ \ \ \ \ \ \ \ \ \ \ \ \ $

\bigskip

\bigskip \textbf{Hyper-Tool \# 1: FIRST ORDER STAR-MAP.}

\textbf{Definition.1.2.1.2. }The \textbf{first order} star-map is a function 
$\ast _{1}:\mathbf{U\rightarrow V}_{1}$

between two universes that associates to each object $A$ $\in $ $\mathbf{U}$
its \textbf{first order}

hyper-extension (or \textbf{first order} non-standard extension) $^{\ast
_{1}}A$ $\in $ $\mathbf{V}_{1}\mathbf{.}$ It is also

assumed that $^{\ast _{1}}n=n$ for all natural numbers $n$ $\in $ $%
\mathbb{N}
$, and that the properness

condition $^{\ast _{1}}%
\mathbb{N}
\neq 
\mathbb{N}
$ holds.

\textbf{Remark 1.2.1.3. }It is customary to call standard any object $A$ $%
\in $ $\mathbf{U}$ in the

domain of the first order star-map $\ast _{1},$ and first order nonstandard
any object

$B$ $\in $ $\mathbf{V}_{1}$ in the codomain. The adjective standard is also
often used in the

literature for first order hyper-extensions $^{\ast _{1}}A\in \mathbf{V}%
_{1}. $

\textbf{Hyper-Tool \# 2: SECOND ORDER STAR-MAP.}

\textbf{Definition.1.2.1.3. }The \textbf{second order} star-map is a function

$\ast _{2}:\mathbf{V}_{1}$ $\mathbf{\rightarrow V}_{2}$ between two
universes that associates to each

object $A$ $\in $ $\mathbf{V}_{1}$ its \textbf{second order} hyper-extension
(or \textbf{second }

\textbf{order} non-standard extension) $^{\ast _{2}}A$ $\in $ $\mathbf{V}_{2}%
\mathbf{.}$

It is also assumed that $^{\ast _{2}}N=N$ for all hyper natural numbers

$N$ $\in $ $^{\ast _{1}}%
\mathbb{N}
,$and that the properness condition $^{\ast _{2}}%
\mathbb{N}
\neq $ $^{\ast _{1}}%
\mathbb{N}
$ holds.

\textbf{Hyper-Tool \# 3: FIRST ORDER TRANSFER PRINCIPLE.}

\textbf{Definition.1.2.1.4. }Let $P(a_{1},...,a_{n})$ be a property of the
standard objects

$a_{1},...,a_{n}\in $ $\mathbf{U}$ expressed as an "elementary sentence".
Then $P(a_{1},...,a_{n})$ is

true if and only if corresponding sentence $^{\ast _{1}}P\left(
c_{1},...,c_{n}\right) $ is true about the

corresponding hyper-extensions $^{\ast _{1}}a_{1},...,^{\ast _{1}}a_{n}\in 
\mathbf{V}_{1}$. That is:

\bigskip\ $\ \ \ \ \ \ \ \ \ \ \ \ \ \ \ \ \ \ \ \ \ 
\begin{array}{cc}
\begin{array}{c}
\\ 
P(a_{1},...,a_{n})\iff \text{ }^{\ast _{1}}P\left( ^{\ast
_{1}}a_{1},...,^{\ast _{1}}a_{n}\right) . \\ 
\end{array}
& \left( 1.2.3\right)%
\end{array}%
$

In particular $P(a_{1},...,a_{n})$ is true if and only if the same sentence $%
P\left( c_{1},...,c_{n}\right) $

is true about the corresponding hyper-extensions $^{\ast
_{1}}a_{1},...,^{\ast _{1}}a_{n}\in \mathbf{V}_{1}$. That

is: $P(a_{1},...,a_{n})\iff $ $P\left( ^{\ast _{1}}a_{1},...,^{\ast
_{1}}a_{n}\right) .$

\bigskip \textbf{Hyper-Tool \# 4: SECOND ORDER TRANSFER PRINCIPLE.}

\textbf{Definition.1.2.1.5.}Let $^{\ast _{1}}P\left( ^{\ast
_{1}}a_{1},...,^{\ast _{1}}a_{n}\right) $ be a property of the \textbf{first
order}

non-standard objects $^{\ast _{1}}a_{1},...,^{\ast _{1}}a_{n}\in \mathbf{V}%
_{1}$ expressed as an "elementary

sentence".

\bigskip

\bigskip

\bigskip

\section{I.2.2.The higher orders of hyper-method by \ \ \ \ \ \ \ \ \ \ \ \
\ \ \ \ \ \ \ \ \ \ using countable universes.}

\bigskip

\textbf{Definition.1.2.2.1. }Any countable universe $\mathbf{U}_{\omega }$
is a nonempty countable

collection of \textit{"standard mathematical} \textit{objects" }that is
closed under subsets,

i.e. $a\subseteqq A\in \mathbf{U}$ $\implies a\in \mathbf{U}$ and closed
under the basic mathematical

operations. Precisely, whenever

$A,B\in \mathbf{U},$ we require that also the union $A\cup B$, the
intersection $A\cap B,$

the set-difference $A\backslash B$ the ordered pair $\left\{ A,B\right\} ,$%
the Cartesian product

$A\times B,$ the powerset $P(A)=\left\{ a|a\subseteqq A\right\} ,$the
function-set

$B^{A}=\left\{ f\text{ }|\text{ }f:A\rightarrow B\right\} $ all belong to $%
\mathbf{U}_{\omega }\mathbf{.}$A countable universe $\mathbf{U}_{\omega }$
is also

assumed to contain (copies of) all sets of numbers $%
\mathbb{N}
,$ $%
\mathbb{Z}
,%
\mathbb{Q}
\in \mathbf{U}_{\omega },%
\mathbb{R}
_{\omega },$

$%
\mathbb{C}
_{\omega }$ $\in \mathbf{U}_{\omega }\mathbf{,}$ and to be transitive, i.e.
members of members of $\mathbf{U}_{\omega }$ belong

to $\mathbf{U}_{\omega }$ or in formulae: $a\in A\in \mathbf{U}_{\omega }$ $%
\implies a\in \mathbf{U}_{\omega }$.

\textbf{Remark 1.2.2.1.}In any countable model $M_{\omega }$ of $\mathbf{ZFC,%
}$an ordered pair $\left\{ a,b\right\} $

is defined as the Kuratowski pair

$\left\{ \left\{ a\right\} ,\left\{ a,b\right\} \right\} ;$an $n$-tuple is
inductively defined by $\left\{ a_{1},...,a_{n},a_{n+1}\right\} $ $=$

$\left\{ \left\{ a_{1},...,a_{n}\right\} ,a_{n+1}\right\} ;$ an $n$-place
relation $R$ on $A$ is identified with the

countable set $R\subseteq A^{n}$ of $n$-tuples that satisfy it; a function $%
f:A\rightarrow B$

is identified with its graph $\left[ \left\{ a,b\right\} \in A\times B|b=f(a)%
\right] .$

As for numbers, complex numbers $%
\mathbb{C}
_{\omega }=$ $%
\mathbb{R}
_{\omega }\times 
\mathbb{R}
_{\omega }$ are defined as ordered

pairs of real numbers, and the real numbers $%
\mathbb{R}
_{\omega }$ are defined as

countable set of countable equivalence classes of suitable sets of

rational numbers namely,Dedekind cuts or Cauchy sequences.

The rational numbers $%
\mathbb{Q}
$ are a suitable quotient $%
\mathbb{Z}
\times 
\mathbb{Z}
/_{\approx },$ and the

integers $%
\mathbb{Z}
$ are in turn a suitable quotient $%
\mathbb{N}
\times 
\mathbb{N}
/_{\approx }$. The natural

numbers of $\mathbf{ZFC}$ are defined as the set of von Neumann naturals:

$0=\NEG{0}$ and $n+1=\left\{ n\right\} $ (so that each natural number

$\left\{ n=0,1,...,n-1\right\} $ is identified with the set of its
predecessors.)

\bigskip

\bigskip

\bigskip

\section{I.2.3.Divisibility of hyperintegers.}

\bigskip

\textbf{Definition.1.2.3.1.}If $n$ and $d$ are hypernaturals,i.e. $n,d\in $ $%
^{\ast }%
\mathbb{N}
$ or hyperintegers,

i.e. $n,d\in $ $^{\ast }%
\mathbb{Z}
$ and $d\neq 0,$then $n$ is divisible by $d$ provided $n=d\cdot k$ for some

hyperinteger $k.$Alternatively, we say:

\textbf{1.}$n$ is a multiple of $d,$

\textbf{2.}$d$ is a factor of $n,$

\textbf{3.}$d$ is a divisor of $n,$

\textbf{4.}$d$ divides $n$ (denoted with $d$ $|$ $n$).

\textbf{Theorem 1.2.3.1.}Transitivity of Divisibility.

For all $a,b,c\in $ $^{\ast }%
\mathbb{Z}
,$ if $a|b$ and $b|c,$ then $a|c.$

\textbf{Theorem 1.2.3.2.}Every positive hyperinteger greater than $1$

is divisible by a hyperprime number.

\bigskip

\textbf{Definition.1.2.3.2.}Given any integer $n>1,$ the

standard factored form of $n\in $ $^{\ast }%
\mathbb{N}
$ is an expression of $n=$ $^{\ast }\prod_{k=1}^{m}p_{k}^{e_{k}},$

where $m$ is a positive hyperinteger, $p_{1},p_{2},...,p_{m}$ are

hyperprime numbers with $p_{1}<p_{2}<...<p_{m}$ and

$e_{1},e_{2},...,e_{m}$ are positive hyperintegers.

\textbf{Theorem 1.2.3.3.}Given any hyperinteger $n>1,$ there exist

positive hyperinteger $m,$hyperprime numbers $p_{1},p_{2},...,p_{m}$

\bigskip and positive hyperintegers $e_{1},e_{2},...,e_{m}$ with $n=$ $%
^{\ast }\prod_{k=1}^{m}p_{k}^{e_{k}}.$\bigskip

\textbf{Theorem 1.2.3.1. }(\textbf{i})\textbf{\ }Every pair of elements $%
m,n\in $ $^{\ast }%
\mathbb{N}
$ has a highest

common factor $d=s\times m+t\times n$ for some $s,t\in $ $^{\ast }%
\mathbb{Z}
.$

(\textbf{ii}) For every pair of elements $a,d\in $ $^{\ast }%
\mathbb{N}
$ dividend $a$ and divisor $d,$ with $d\neq 0$

there exist unique integers $q$ and $r$ such that $a=q\times d+r$ and $0\leq
r<\left\vert d\right\vert .$

\bigskip

\textbf{Definition.1.2.3.2. }Suppose\textbf{\ }that $a=q\times d+r$ and $%
0\leq r<\left\vert d\right\vert .$We call $d$

the quotient and $r$ the remainder.

\bigskip

\bigskip

Redrick, squinting his swollen eyes against the blinding light, silently
watched

him go. He was cool and calm, he knew what was about to happen, and he

knew that he would not watch,but it was still all right to watch, and he did,

feeling nothing in particular,except that deep inside a little worm started

wriggling around and twisting its sharp head in his gut.

\ \ \ \ \ \ \ \ \ \ \ \ \ \ \ \ \ \ \ \ \ \ \ \ \ \ \ \ \ \ \ \ \ \ \ \ \ \
\ \ \ \ \ \ \ \ \ \ \ \ \ \ \ \ \ \ \ \ \ \ \ \ \ \ \ \ \ \ \ \ \ \ \ \
Arkady and Boris Strugatsky

\bigskip\ \ \ \ \ \ \ \ \ \ \ \ \ \ \ \ \ \ \ \ \ \ \ \ \ \ \ \ \ \ \ \ \ \
\ \ \ \ \ \ \ \ \ \ \ \ \ \ \ \ \ \ \ \ \ \ \ \ \ \ \ \ \ \ \ \ \ \ \ \ \ \
\ \ \ \ \ \ \ \ \ \ \ \ \ "Roadside Picnic"

\section{I.3.The construction non-archimedean \ \ \ \ \ \ \ \ \ \ \ \ \ \ \
\ \ \ \ pseudo-ring$^{\ast }%
\mathbb{R}
_{\mathbf{d}}.$\ \ \ \ \ \ \ \ \ \ \ \ \ \ \ }

\bigskip

\section{I.3.1.Generalized pseudo-rings \ \ \ \ \ \ \ \ \ \ \ \ \ \ \ \ \ \
\ \ \ \ Wattenberg-Dedekind hyperreals$\ $\ \ \ \ \ \ \ \ \ \ \ \ \ \ \ \ \
\ \ \ \ \ \ \ \ \ \ \ \ \ \ \ \ \ \ \ \ \ \ $^{\ast }%
\mathbb{R}
_{\mathbf{d}}$ and hyperintegers $^{\ast }%
\mathbb{Z}
_{\mathbf{d}}$.}

\bigskip

\ \ \ \ \ \ \ \ \ \ \ \ \ \ \ \ \ \ \ \ \ \ \ \ \ \ \ \ \ \ \ \ \ \ \ \ \ \
\ \ \ \ \ \ \ \ \ \ \ \ \ \ \ \ \ \ \ \ \ \ \ \ \ \ \ \ \ \ \ \ \ \ \ \ \ \
\ \ \ \ \ \ \ \ \ 

\section{I.3.1.1.Strong and weak Dedekind cuts. Wattenberg-Dedekind
hyperreals and \ \ \ \ \ \ \ \ \ \ \ \ \ \ \ \ \ \ \ \ \ \ \ hyperintegers.}

\ \ \ \ \ \ \ \ \ \ \ \ \ \ \ \ \ \ \ \ \ \ \ \ \ \ \ 

From\textbf{\ Theorem 1.2.1.1} above we knov that: $^{\ast }%
\mathbb{R}
$ $\mathbf{is}$\textbf{\ }$\mathbf{not}$\textbf{\ }$\mathbf{Dedekind}$%
\textbf{\ }$\mathbf{complete}.$

For example, $\mu (0)$ and $%
\mathbb{R}
$ are bounded subsets of $^{\ast }%
\mathbb{R}
$ which have no

suprema or infima in $^{\ast }%
\mathbb{R}
$.

Possible standard completion of the field $^{\ast }%
\mathbb{R}
$ can be constructed by \ \ \ \ \ \ \ \ \ \ \ \ \ \ \ \ \ \ \ \ \ \ 

Dedekind sections [23],[24]. In [24] Wattenberg constructed the Dedekind

completion of a nonstandard model of the real numbers and applied the

construction to obtain certain kinds of special measures on the set of
integers.

Thus was established that the Dedekind completion $^{\ast }%
\mathbb{R}
_{\mathbf{d}}$ of the field $^{\ast }%
\mathbb{R}
$ is a

structure of interest not for its own sake only and we establish further

importent applications here. Importent concept was introduce Gonshor [23]

is that of the \textbf{absorption number} of an element $\mathbf{a\in }%
^{\ast }%
\mathbb{R}
_{\mathbf{d}}$ which, roughly

speaking,measures the degree to which the cancellation law

$\mathbf{a}+b=\mathbf{a}+c\implies b=c$ fails for $\mathbf{a}$.

\bigskip\ \ \ \ \ \ \ \ \ \ \ \ 

More general construction well known from topoi theory [10].\bigskip

\textbf{Definition 1.3.1.1.1.} A \textit{Dedekind hyperreal }$\alpha \in $ $%
^{\ast }%
\mathbb{R}
_{\mathbf{d}}$ is a pair

$(U,V)\in \mathbf{P}\left( ^{\ast }%
\mathbb{Q}
\right) \times $ $\mathbf{P}\left( ^{\ast }%
\mathbb{Q}
\right) $ satisfying the next conditions:

\textbf{1.}$\exists x\exists y\left( x\in U\wedge y\in V\right) .$

\textbf{2.} $U\cap V=\varnothing .$

\textbf{3.}$\forall x\left( x\in U\iff \exists y\left( y\in V\wedge
x<y\right) \right) .$

\textbf{4.} $\forall x\left( x\in V\iff \exists y\left( y\in V\wedge
y<x\right) \right) .$

\textbf{5.} $\forall x\forall y\left( x<y\implies x\in U\vee y\in V\right) .$%
\bigskip

\textbf{Remark.} The monad of $\alpha \in $ $^{\ast }%
\mathbb{R}
,$ the set $\left\{ x\in \text{ }^{\ast }%
\mathbb{R}
|\text{ }x\approx \alpha \right\} $ is

denoted: $\mu \left( \alpha \right) .$

Monad $\mu \left( 0\right) $ is denoted: $\mathbf{I}_{\ast }.$Supremum of $%
\mathbf{I}_{\ast }$ is denoted: $\varepsilon _{\mathbf{d}}.$

Let $A$ be a subset of $^{\ast }%
\mathbb{R}
$ is bounded or hyperbounded above

then $\sup \left( A\right) $ exists in $^{\ast }%
\mathbb{R}
_{\mathbf{d}}.$\bigskip\ 

\textbf{Example. }(\textbf{i}) $\Delta _{\mathbf{d}}=\sup \left( 
\mathbb{R}
_{+}\right) \in $ $^{\ast }%
\mathbb{R}
_{\mathbf{d}}\left\backslash ^{\ast }%
\mathbb{R}
\right. ,$(\textbf{ii}) $\varepsilon _{\mathbf{d}}=\sup \left( \ \mathbf{I}%
_{\ast }\right) \in $ $^{\ast }%
\mathbb{R}
_{\mathbf{d}}\left\backslash ^{\ast }%
\mathbb{R}
\right. .$

\textbf{Remark. }Anfortunately the set $^{\ast }%
\mathbb{R}
_{\mathbf{d}}$ inherits some but by \textbf{no means all}

of the algebraic structure on $^{\ast }%
\mathbb{R}
.$For example,$^{\ast }%
\mathbb{R}
_{\mathbf{d}}$ is not a group with

respect to addition since if $x+_{^{\ast }%
\mathbb{R}
_{\mathbf{d}}}y$ denotes the addition in $^{\ast }%
\mathbb{R}
_{\mathbf{d}}$ then:

\bigskip

\ \ \ \ \ \ \ $%
\begin{array}{cc}
\begin{array}{c}
\\ 
\varepsilon _{\mathbf{d}}+_{^{\ast }%
\mathbb{R}
_{\mathbf{d}}}\varepsilon _{\mathbf{d}}=\varepsilon _{\mathbf{d}}+_{^{\ast }%
\mathbb{R}
_{\mathbf{d}}}0_{^{\ast }%
\mathbb{R}
_{\mathbf{d}}}=\varepsilon _{\mathbf{d}} \\ 
\end{array}
& \left( 1.3.1.1.1\right)%
\end{array}%
$

\bigskip

Thus $^{\ast }%
\mathbb{R}
_{\mathbf{d}}$ is not iven a ring but pseudo-ring only. Thus, one must

proceed somewhat cautiously. In this section more details than is

customary will be included in proofs because some standard properties

which at first glance appear clear often at second glance reveal themselves

to be false in $^{\ast }%
\mathbb{R}
_{\mathbf{d}}$.

\bigskip

We shall briefly remind a way Dedekind's constructions of a pseudo-field

\bigskip\ $^{\ast }%
\mathbb{R}
_{\mathbf{d}}.$

\textbf{Definition 1.3.1.2.a.(Strong and weak Dedekind cuts) }

(\textbf{1}) Suppose $\preceq $ is a total ordering on $X.$ We write $%
x\preceq y$ if $x$ is

less than or equal\ to $y$ and we write $x\prec y$ if $x\preceq y$\ and $%
x\neq y.$

Then $\left\{ A,B\right\} $ is said to be a \textbf{strong} \textbf{Dedekind
cut} of $\left\langle X,\preceq \right\rangle ,$

if and only if:

\textbf{1.} $A$ and $B$ are nonempty subsets of $X.$

\textbf{2. }$A\cup B=X.$

\textbf{3. }For each $x$ in $A$ and each $y$ in $B,x\preceq y.$\bigskip

(\textbf{2}) Suppose $\left\{ A,B\right\} $ is a strong Dedekind cut of $%
\left\langle X,\preceq \right\rangle $

Then $\left\{ A,\widetilde{A};B,\widetilde{B}\right\} $ is said to be a 
\textbf{weak} \textbf{Dedekind cut} of $\left\langle X,\preceq \right\rangle
,$

if and only if:

\textbf{1. }$A\subsetneqq \widetilde{A},B\subsetneqq \widetilde{B}.$

\textbf{2. }For each $x$ in $A$ there is exist $\widetilde{x}_{1}\in 
\widetilde{A}$ such that $x\prec \widetilde{x}$ and

$\widetilde{x}_{2}\in \widetilde{A}$ such that $\widetilde{x}_{2}\prec x.$

\textbf{3. }For each $y$ in $B$ there is exist $\widetilde{y}_{1}\in 
\widetilde{B}$ such that $\widetilde{y}_{1}\prec y$ and

$\widetilde{y}_{2}\in \widetilde{B}$ such that $y\prec \widetilde{y}_{2}.$

\bigskip

(\textbf{3}) $A$ is the \textit{left-hand} part of the strong cut $\left\{
A,B\right\} $ and $B$ is the

\textit{riht-hand} part\ of the strong cut\ $\left\{ A,B\right\} $.

We denote the strong cut as $x=A|B$\ or simple $x=A.$

The strong cut $x=A|B$ is less than or equal\ to the strong cut

$y=C|D$\ if $A\subseteqq C.$

(\textbf{4}) $\widetilde{A}$ is the \textit{left-hand} part of the weak cut $%
\left\{ A,\widetilde{A};B,\widetilde{B}\right\} $ and $\widetilde{B}$ is the

\textit{riht-hand} part\ of the weak cut\ $\left\{ A,\widetilde{A};B,%
\widetilde{B}\right\} $.

We denote the weak cut as $x=\widetilde{A}|\widetilde{B}$\ or simple $x=%
\widetilde{A}.$

The weak cut $x=\widetilde{A}|\widetilde{B}$ is less than or equal\ to the
weak cut

$y=\widetilde{C}|\widetilde{D}$\ iff $A\subseteqq C.$

\ \ \ \ \ \ \ \ \ \ \ \ \ \ \ \ \ \ \ \ \ \ \ \ \ \ \ \ \ \ \ \ \ \ \ \ \ \
\ \ \ \ \ \ \ \ \ \ \ \ \ \ \ \ \ 

\textbf{Definition 1.3.2.b. }

(\textbf{1}) $c\in X$ is said to be a cut element of $\left\{ A,B\right\} $

if and only if either:

(\textbf{i}) $c$ is in $A$ and $x\preceq c\preceq y$ for each $x$ in $A$ and
each $y$ in $B,$ or

(\textbf{ii}) $c$ is in $B$ and $x\preceq c\preceq y$ for each $x$ in $A$
and each $y$ in $B.$

(\textbf{2}) $c\in X$ is said to be a cut element of $\left\{ A,\widetilde{A}%
;B,\widetilde{B}\right\} $

if and only if either:

(\textbf{i}) $c$ is in $\widetilde{A}$ and $x\preceq c\preceq y$ for each $x$
in $\widetilde{A}$ and each $y$ in $\widetilde{B},$ or

(\textbf{ii}) $c$ is in $\widetilde{B}$ and $x\preceq c\preceq y$ for each $%
x $ in $\widetilde{A}$ and each $y$ in $\widetilde{B}.$

\textbf{Definition 1.3.2.c.}$\left\langle X,\preceq \right\rangle $ is said
to be \textit{Dedekind complete} if

and only if each strong Dedekind cut of $\left\langle X,\preceq
\right\rangle $,has a cut

element.

Equivalently $\left\langle X,\preceq \right\rangle $ is said to be \textit{%
Dedekind complete} if and

only if each weak Dedekind cut $\widetilde{A}|\widetilde{B}$ of $%
\left\langle X,\preceq \right\rangle $,has a cut

element.

\textbf{Example.} The following theorem is well-known.

Theorem. $\left\langle 
\mathbb{R}
,\leqslant \right\rangle $ is Dedekind complete, and for each Dedekind cut

$\left\{ A,B\right\} $,of $\left\langle 
\mathbb{R}
,\leqslant \right\rangle $ if $r$ and $s$ are cut elements of $\left\{
A,B\right\} $, then $r=s.$

Making a semantic leap, we now answer the question "what is a

\textit{Wattenberg}-\textit{Dedekind hyperreal number }?"

\textbf{Definition 1.3.2.d. }A \textbf{Wattenberg-Dedekind hyperreal number }%
is

a cut in $^{\ast }%
\mathbb{Q}
.$

$^{\ast }%
\mathbb{R}
_{\mathbf{d}}$ is the class of all Dedekind hyperreal numbers $x=A|B$ ($x=A$%
).

We will show that in a natural way $^{\ast }%
\mathbb{R}
_{\mathbf{d}}$ is a \textbf{complete ordered \ \ \ \ \ \ \ \ \ \ \ \ \ \ \ \
\ \ \ \ \ \ \ \ \ \ \ \ \ \ }

\textbf{generalized} \textbf{pseudo-ring }containing $^{\ast }%
\mathbb{R}
.$

Before spelling out what this means, here are some examples of cuts.

(\textbf{i}) \ \ 

$A|B=\left. \left\{ r\in \text{ }^{\ast }%
\mathbb{Q}
\text{ }|\text{ }r<1\right\} \right\vert \left\{ r\in \text{ }^{\ast }%
\mathbb{Q}
\text{ }|\text{ }r\geq 1\right\} .$

(\textbf{ii}) $\ $

$\ A|B=\left. \left\{ r\in \text{ }^{\ast }%
\mathbb{Q}
\text{ }|\left( r\leq 0\right) \vee \left( \text{ }r^{2}<2\right) \right\}
\right\vert \left\{ r\in \text{ }^{\ast }%
\mathbb{Q}
\text{ }|\left( r>0\right) \wedge \text{ }\left( r^{2}\geq 2\right) \right\}
.$

(\textbf{iii})

$A|B=\left. \left\{ r\in \text{ }^{\ast }%
\mathbb{Q}
\text{ }|\text{ }r<\omega \right\} \right\vert \left\{ r\in \text{ }^{\ast }%
\mathbb{Q}
\text{ }|\text{ }r\geq \omega \right\} ,$where $\omega \in $ $^{\ast }%
\mathbb{Q}
_{+}\backslash 
\mathbb{Q}
_{+}.$\ 

(\textbf{iv})

$A|B=$

$\left. \left\{ r\in \text{ }^{\ast }%
\mathbb{Q}
\text{ }|\left( r\leq 0\right) \vee \left( r\in \mathbf{I}_{\ast }\right)
\vee \left( \text{ }r\in 
\mathbb{Q}
_{+}\right) \right\} \right\vert $

$\left\{ r\in \text{ }^{\ast }%
\mathbb{Q}
\text{ }|\left( r>0\right) \wedge \text{ }\left( r\in \text{ }^{\ast }%
\mathbb{Q}
_{+}\backslash \left( 
\mathbb{Q}
_{+}\cup \mathbf{I}_{\ast }\right) \right) \right\} .$

(\textbf{v})

$A|B=$

$\left. \left\{ r\in \text{ }^{\ast }%
\mathbb{Q}
\text{ }|\left( r\leq 0\right) \vee \left( r\in \mathbf{I}_{\ast }\right)
\right\} \right\vert \left\{ r\in \text{ }^{\ast }%
\mathbb{Q}
\text{ }|\left( r>0\right) \wedge \text{ }\left( r\in \text{ }^{\ast }%
\mathbb{Q}
_{+}\backslash \mathbf{I}_{\ast }\right) \right\} .$

\textbf{Remark.} It is convenient to say that $A|B\in $ $^{\ast }%
\mathbb{R}
_{\mathbf{d}}$ is a \textbf{rational} (\textbf{hyperrational})\textbf{\ \ \
\ \ \ \ \ \ \ \ \ \ \ \ \ \ \ \ \ }

\textbf{cut }in $^{\ast }%
\mathbb{Q}
$\textit{\ }if it is like the cut in examples (\textbf{i}),(\textbf{iii}):
fore some fixed rational

(hyperrational) number $c\in $ $^{\ast }%
\mathbb{Q}
,A$ is the set of all hyperrational $r$ such

that $r<c$ while $B$ is the rest of $^{\ast }%
\mathbb{Q}
.$ \ \ \ \ \ \ \ \ \ \ \ \ \ \ \ \ \ \ \ \ \ \ \ \ \ \ \ \ 

The $B$-set of a rational (hyperrational) cut contains a smollest $c\in $ $%
^{\ast }%
\mathbb{Q}
,$

and conversaly if $A|B$ is a cut in $^{\ast }%
\mathbb{Q}
$ and $B$ contains a smollest element

$c$ then $A|B$ is a rational or hyperrational cut at $c.$We write $\breve{c}$
for the rational

hyperrational cut at $c.$This lets us think of $^{\ast }%
\mathbb{Q}
\subset $ $^{\ast }%
\mathbb{R}
_{\mathbf{d}}$ by identifying $c$

with $\breve{c}.$

\textbf{Remark. }It is convenient to say that:

(\textbf{1}) $A|B\in $ $^{\ast }%
\mathbb{R}
_{\mathbf{d}}$ is an \textbf{standard cut} in $^{\ast }%
\mathbb{Q}
$\textit{\ }if it is like the cut in examples

(\textbf{i})-(\textbf{ii}):fore some cut $A^{\prime }|B^{\prime }\in 
\mathbb{R}
$ the next equality is satisfied:

$A|B=$ $^{\ast }\left( A^{\prime }\right) |^{\ast }\left( B^{\prime }\right)
,$i.e. $A$-set of a cut is an standard set.

(\textbf{2}) $A|B\in $ $^{\ast }%
\mathbb{R}
_{\mathbf{d}}$ is an \textbf{internal} \textbf{cut }or \textbf{nonstandard
cut} in $^{\ast }%
\mathbb{Q}
$\textit{\ }if it is like

the cut in example (\textbf{iii}), i.e. $A$-set of a cut is an \textit{%
internal nonstandard }

\textit{set.}

(\textbf{3}) $A|B\in $ $^{\ast }%
\mathbb{R}
_{\mathbf{d}}$ is an \textbf{external} \textbf{cut }in $^{\ast }%
\mathbb{Q}
$\textit{\ }if it is like the cut in

examples (\textbf{iv})-(\textbf{v}),i.e. $A$-set of a cut is an \textit{%
external set. }\ 

There is an order relation $\left( \cdot \leq \cdot \right) $ on cuts that
fairly cries out for

attention.

\textbf{Definition 1.3.2.e. }The cut $x=A|B$ is less than or equal\ to the
cut

$y=C|D$\ if $A\subseteqq C.$

We write $x\leq y$ if $x$ is less than or equal\ to $y$ and we write $x<y$ if

$x\leq y$\ and $x\neq y.$If\ $x=A|B$\ is less than $y=C|D$ then $A\subset C$
and

$A\neq C,$so there is some $c_{0}\in C\backslash A.$Sinse the $A$-set of a
cut contains

no largest element, there is also a $c_{1}\in C$ with $c_{0}<c_{1}.$All the

hyperrational numbers $c$ with $c_{0}\leq c\leq c_{1}$ belong to $%
C\backslash A.$

\textbf{Remark.} The property distinguishing $^{\ast }%
\mathbb{R}
_{\mathbf{d}}$ from $^{\ast }%
\mathbb{Q}
$ and from $^{\ast }%
\mathbb{R}
$ and

which is the bottom of every significant theorem about $^{\ast }%
\mathbb{R}
_{\mathbf{d}}$ involves

upper bounds and least upper bounds or equivalently,lower bounds

and gretest lower bounds.

\textbf{Definition 1.3.2.f. }$M\in $ $^{\ast }%
\mathbb{R}
_{\mathbf{d}}$ is an \textbf{upper bound }for a set $S\subset $ $^{\ast }%
\mathbb{R}
_{\mathbf{d}}$ if

each $s\in S$ satisfies $s\leq M.$ We also say that the set $S$ is \textbf{%
bounded }

\textbf{above }by $M$ iff $M\in $ $\mathbf{L}\left( ^{\ast }%
\mathbb{R}
\right) $\ We also say that the set $S$ is

\textbf{hyperbounded above }iff $M\notin $ $\mathbf{L}\left( ^{\ast }%
\mathbb{R}
\right) ,$i.e.$\left\vert M\right\vert \in $ $^{\ast }%
\mathbb{R}
_{+}\backslash 
\mathbb{R}
_{+}.$

\textbf{Definition 1.3.2.g. }An\textbf{\ }upper bound for $S$ that is less
than all other

upper bound for $S$ is a \textbf{least upper bound} for $S.$

The concept of a \textit{pseudo-ring} originally was introduced by \ \ \ \ \
\ \ \ \ \ \ \ \ \ \ \ \ \ \ \ \ \ \ \ \ \ \ \ \ \ \ \ \ \ \ \ \ \ \ \ \ \ \
\ \ \ \ \ \ \ \ 

E. M.Patterson [21].Briefly,Patterson's pseudo-ring is an algebraic \ \ \ \
\ \ \ \ \ \ \ \ \ \ \ \ \ \ \ \ \ \ \ \ \ \ \ \ \ \ \ \ 

system consisting of an additive abelian group $\mathbf{A}$, a distinguished
\ \ \ \ \ \ \ \ \ \ \ \ \ \ \ \ \ \ \ \ \ \ \ \ \ \ \ \ \ \ \ 

subgroup $\mathbf{A}^{\mathbf{\ast }}$ of $\mathbf{A,}$and a multiplication
operation $\mathbf{A}^{\ast }\mathbf{\times A}\rightarrow \mathbf{A}$ under
\ \ \ \ \ \ \ \ \ \ \ \ \ \ \ \ \ \ \ \ \ \ \ \ \ \ \ \ \ \ \ \ \ \ 

which $\mathbf{A}^{\ast }$ is a ring and $\mathbf{A}$ a left $\mathbf{A}%
^{\ast }$-module.For convenience, we

denote the pseudo-ring by $\Re =(A^{\ast },A).$

\textbf{Definition 1.3.1.2.h.Generalized} \textbf{pseudo-ring
(m-pseudo-ring)\ }

is an algebraic system consisting of\textbf{\ }an\textbf{\ }abelian
semigroup $\mathbf{A}_{\mathbf{s}}$\textbf{\ }(or

abelian monoid $\mathbf{A}_{\mathbf{m}}$),a distinguished subgroup $\mathbf{A%
}_{\mathbf{s}}^{\mathbf{\ast }}$ of $\mathbf{A}_{\mathbf{s}}$ (or a

distinguished subgroup $\mathbf{A}_{\mathbf{m}}^{\mathbf{\ast }}$ of $%
\mathbf{A}_{\mathbf{m}}$)\textbf{,}and a multiplication operation

$\mathbf{A}_{\mathbf{s}}^{\ast }\mathbf{\times A}_{\mathbf{s}}\rightarrow 
\mathbf{A}_{\mathbf{s}}$ ($\mathbf{A}_{\mathbf{m}}^{\ast }\mathbf{\times A}_{%
\mathbf{m}}\rightarrow \mathbf{A}_{\mathbf{m}}$) under which $\mathbf{A}_{%
\mathbf{s}}^{\ast }$ ($\mathbf{A}_{\mathbf{m}}^{\ast }$)\ is a ring and

$\mathbf{A}_{\mathbf{s}}$ ($\mathbf{A}_{\mathbf{m}}$) a left $\mathbf{A}_{%
\mathbf{s}}^{\ast }$-module ($\mathbf{A}_{\mathbf{m}}^{\ast }$-module).

For convenience,we denote\ the generalized pseudo-ring by

$\Re _{\mathbf{s}}=(A_{\mathbf{s}}^{\ast },A_{\mathbf{s}}).$\ 

\textbf{Pseudo-field }is an algebraic system consisting of\textbf{\ }an%
\textbf{\ }abelian

semigroup $\mathbf{A}_{\mathbf{s}}$\textbf{, }a distinguished subgroups $%
\mathbf{A}_{\mathbf{s}}^{\mathbf{\ast }}\subsetneqq \mathbf{A}_{\mathbf{s}}^{%
\mathbf{\#}}$ of $\mathbf{A}_{\mathbf{s}}$ and a

multiplication operations $\mathbf{A}_{\mathbf{s}}^{\ast }\mathbf{\times A}_{%
\mathbf{s}}\rightarrow \mathbf{A}_{\mathbf{s}}$and $\mathbf{A}_{\mathbf{s}}%
\mathbf{\times A}_{\mathbf{s}}^{\ast }\rightarrow \mathbf{A}_{\mathbf{s}}$
under

which $\mathbf{A}_{\mathbf{s}}^{\ast }$ is a ring,$\mathbf{A}_{\mathbf{s}}^{%
\mathbf{\#}}$ is a field and $\mathbf{A}_{\mathbf{s}}$ is a vector spase over

field $\mathbf{A}_{\mathbf{s}}^{\#}.$

\textbf{Definition 1.3.1.2}$^{\prime }$\textbf{.I.(Strong and Weak Dedekind
cut tipe}

\textbf{I in} $^{\ast }%
\mathbb{R}
$\textbf{).}

(\textbf{1})\textbf{\ }A strong Dedekind cut tipe \textbf{I} $\alpha _{%
\mathbf{s}}=\alpha _{\mathbf{s}}^{\mathbf{I}}$ in $^{\ast }%
\mathbb{R}
$ is a subset

$\alpha _{\mathbf{s}}^{\mathbf{I}}\subsetneqq $ $^{\ast }%
\mathbb{R}
$ of the hyperreals $^{\ast }%
\mathbb{R}
$ that satisfies these properties:

\textbf{1.} $\alpha _{\mathbf{s}}^{\mathbf{I}}$ is not empty.

\textbf{2.} $\beta _{\mathbf{s}}^{\mathbf{I}}=$ $^{\ast }%
\mathbb{R}
\backslash \alpha _{\mathbf{s}}^{\mathbf{I}}$ is not empty.

\textbf{3.} $\alpha _{\mathbf{s}}^{\mathbf{I}}$ contains no greatest element

\bigskip \textbf{4.} For every $x,y\in $ $^{\ast }%
\mathbb{R}
,$ if $x\in \alpha _{\mathbf{s}}^{\mathbf{I}}$ and $y<x,$ then $y\in \alpha
_{\mathbf{s}}^{\mathbf{I}}$ as well.

(\textbf{2}) A weak Dedekind cut tipe \textbf{I} $\alpha _{w}$ in $^{\ast }%
\mathbb{R}
$ is a subset

$\alpha _{w}^{\mathbf{I}}\subsetneqq $ $^{\ast }%
\mathbb{R}
$ of the hyperreals $^{\ast }%
\mathbb{R}
$ that satisfies these properties:\ 

\textbf{1.} $\alpha _{w}^{\mathbf{I}}$ is not empty.

\textbf{2. }$\beta _{w}^{\mathbf{I}}=$ $^{\ast }%
\mathbb{Q}
\backslash \alpha _{w}^{\mathbf{I}}$ is not empty.

\textbf{3. }$\alpha _{\mathbf{s}}^{\mathbf{I}}$ contains no greatest element.

\textbf{4.} For every $x,y\in $ $^{\ast }%
\mathbb{R}
,$ if $x\in \alpha _{w}^{\mathbf{I}}$ and $y<x,$ then there is exists

$z\in \alpha _{w}^{\mathbf{II}}$ such that $z<y$ as well.

\textbf{Remark. }Note that for every weak Dedekind cut $\alpha _{w}^{\mathbf{%
I}}$ in $^{\ast }%
\mathbb{R}
$ there is

exists unique strong Dedekind cut $\alpha _{\mathbf{s}}^{\mathbf{I}}\left(
\alpha _{w}^{\mathbf{I}}\right) $ in $^{\ast }\mathbf{%
\mathbb{R}
}$ such that:

$\ \ \ \ \ \ \ \ \ \ \ \ \ \ \ \ \ \ \ \ \ \ 
\begin{array}{cc}
\begin{array}{c}
\\ 
\alpha _{\mathbf{s}}^{\mathbf{I}}\left( \alpha _{w}^{\mathbf{I}}\right) =%
\left[ \alpha _{w}^{\mathbf{I}}\right] _{^{\ast }%
\mathbb{R}
}=\bigcap \left\{ \alpha _{\mathbf{s}}^{\mathbf{I}}|\alpha _{w}^{\mathbf{I}%
}\subset \text{ }\alpha _{\mathbf{s}}^{\mathbf{I}}\right\} . \\ 
\end{array}
& \left( 1.3.2\right)%
\end{array}%
$

\textbf{Example. }(\textbf{1})\textbf{\ }Let $\Delta $ denotes the set $%
\left\{ x|x\in \text{ }^{\ast }%
\mathbb{R}
\backslash ^{\ast }%
\mathbb{R}
_{+\infty }\right\} .$ It is

easy to see that $\Delta $ is a strong Dedekind cut in $^{\ast }%
\mathbb{R}
.$

(\textbf{2}) Let $\Delta \upharpoonright $ $^{\ast }%
\mathbb{Q}
$ denotes the set $\Delta \cap $ $^{\ast }%
\mathbb{Q}
=\left\{ q|\left( q\in \text{ }^{\ast }%
\mathbb{Q}
\backslash ^{\ast }%
\mathbb{Q}
_{+\infty }\right) \right\} .$

It is easy to see that $\Delta \upharpoonright $ $^{\ast }%
\mathbb{Q}
$ is a weak Dedekind cut in $^{\ast }%
\mathbb{R}
$ and

$\left[ \Omega _{k}\right] _{^{\ast }%
\mathbb{R}
}=\Delta .$

(\textbf{3}) Let $1_{^{\ast }%
\mathbb{R}
_{\mathbf{d}}}$ denotes the set $\left\{ x\in \text{ }^{\ast }%
\mathbb{R}
|x<1_{^{\ast }%
\mathbb{R}
}\right\} ,1_{^{\ast }%
\mathbb{R}
}\triangleq $ $^{\ast }1$ and $0_{^{\ast }%
\mathbb{R}
_{\mathbf{d}}}$

denotes the set $\left\{ x\in \text{ }^{\ast }%
\mathbb{R}
|x<0\right\} ,0_{^{\ast }%
\mathbb{R}
}\triangleq $ $^{\ast }0.$It is easy to see that $1_{^{\ast }%
\mathbb{R}
_{\mathbf{d}}}$

and $0_{^{\ast }%
\mathbb{R}
_{\mathbf{d}}}$ is a strong Dedekind cuts in $^{\ast }%
\mathbb{R}
.$

(\textbf{3}) Let $1_{^{\ast }%
\mathbb{R}
_{\mathbf{d}}}\upharpoonright $ $^{\ast }%
\mathbb{Q}
=\left\{ x\in \text{ }^{\ast }%
\mathbb{R}
|x<1_{^{\ast }%
\mathbb{R}
}\right\} \cap ^{\ast }%
\mathbb{Q}
$ and

$0_{^{\ast }%
\mathbb{R}
_{\mathbf{d}}}\upharpoonright $ $^{\ast }%
\mathbb{Q}
=\left\{ x\in \text{ }^{\ast }%
\mathbb{R}
|x<0_{^{\ast }%
\mathbb{R}
}\right\} \cap ^{\ast }%
\mathbb{Q}
.$It is easy to see that $1_{^{\ast }%
\mathbb{R}
_{\mathbf{d}}}\upharpoonright $ $^{\ast }%
\mathbb{Q}
$

and $0_{^{\ast }%
\mathbb{R}
_{\mathbf{d}}}\upharpoonright $ $^{\ast }%
\mathbb{Q}
$ is a strong Dedekind cuts in $^{\ast }%
\mathbb{R}
$ and

\bigskip $\left[ 1_{^{\ast }%
\mathbb{R}
_{\mathbf{d}}}\upharpoonright ^{\ast }%
\mathbb{Q}
\right] _{^{\ast }%
\mathbb{R}
_{\mathbf{d}}}=1_{^{\ast }%
\mathbb{R}
_{\mathbf{d}}},\left[ 0_{^{\ast }%
\mathbb{R}
_{\mathbf{d}}}\upharpoonright ^{\ast }%
\mathbb{Q}
\right] _{^{\ast }%
\mathbb{R}
_{\mathbf{d}}}=0_{^{\ast }%
\mathbb{R}
_{\mathbf{d}}}.$

\textbf{Definition 1.3.1.2}$^{\prime }$\textbf{.II.(Strong and Weak Dedekind
cut }

\textbf{tipe} \textbf{II in} $^{\ast }%
\mathbb{R}
$\textbf{).}

(\textbf{1})\textbf{\ }A strong Dedekind cut tipe \textbf{II} $\alpha _{%
\mathbf{s}}^{\mathbf{II}}$ in $^{\ast }%
\mathbb{R}
$ is a subset $\alpha _{\mathbf{s}}^{\mathbf{II}}\subsetneqq $ $^{\ast }%
\mathbb{R}
$

of the hyperrational numbers $^{\ast }%
\mathbb{R}
$ that satisfies these

properties:

\textbf{1.} $\alpha _{\mathbf{s}}^{\mathbf{II}}$ is not empty.

\textbf{2.} $\beta _{\mathbf{s}}^{\mathbf{II}}=$ $^{\ast }%
\mathbb{R}
\backslash \alpha _{\mathbf{s}}^{\mathbf{II}}$ is not empty.

\textbf{3.} $\alpha _{\mathbf{s}}^{\mathbf{II}}$ contains a greatest element
or $^{\ast }%
\mathbb{R}
\backslash \alpha _{\mathbf{s}}^{\mathbf{II}}$ contains no least

element.

\textbf{4.} For every $x,y\in $ $^{\ast }%
\mathbb{R}
,$ if $x\in \alpha _{\mathbf{s}}^{\mathbf{II}}$ and $y<x,$ then $y\in \alpha
_{\mathbf{s}}^{\mathbf{II}}$ as well.

\bigskip

(\textbf{2}) A weak Dedekind cut tipe \textbf{II} $\alpha _{w}^{\mathbf{II}}$
in $^{\ast }%
\mathbb{R}
$ is a subset $\alpha _{w}^{\mathbf{II}}\subsetneqq $ $^{\ast }%
\mathbb{R}
$

of the hyperrational numbers $^{\ast }%
\mathbb{R}
$ that satisfies these properties:\ 

\textbf{1.} $\alpha _{w}^{\mathbf{II}}$ is not empty.

\textbf{2. } $\beta _{w}^{\mathbf{II}}=$ $^{\ast }%
\mathbb{R}
\backslash \alpha _{w}^{\mathbf{II}}$ is not empty.

\textbf{3.} $\alpha _{\mathbf{s}}^{\mathbf{II}}$ contains a greatest element
or $^{\ast }%
\mathbb{R}
\backslash \alpha _{\mathbf{s}}^{\mathbf{II}}$ contains no least

element.

\textbf{4.} For every $x,y\in $ $^{\ast }%
\mathbb{R}
,$ if $x\in \alpha _{w}^{\mathbf{II}}$ and $y<x,$ then there is exists

$z\in \alpha _{w}^{\mathbf{II}}$ such that $z<y$ as well.

\textbf{Remark. }Note that for every weak Dedekind cut $\alpha _{w}^{\mathbf{%
II}}$ in $^{\ast }%
\mathbb{R}
$

there is exists unique strong Dedekind cut $\alpha _{\mathbf{s}}^{\mathbf{II}%
}\left( \alpha _{w}^{\mathbf{II}}\right) $ in $^{\ast }%
\mathbb{R}
$ such

that: $\alpha _{\mathbf{s}}^{\mathbf{II}}\left( \alpha _{w}^{\mathbf{II}%
}\right) =\left[ \alpha _{w}^{\mathbf{II}}\right] _{^{\ast }%
\mathbb{Q}
}=$ $\bigcap \left\{ \alpha _{\mathbf{s}}^{\mathbf{II}}|\alpha _{w}^{\mathbf{%
II}}\subset \text{ }\alpha _{\mathbf{s}}^{\mathbf{II}}\right\} .$

\bigskip

\bigskip\ $\ \ \ \ \ \ \ \ \ \ \ \ \ \ \ \ \ 
\begin{array}{cc}
\begin{array}{c}
\\ 
\alpha _{\mathbf{s}}^{\mathbf{II}}\left( \alpha _{w}^{\mathbf{II}}\right) =%
\left[ \alpha _{w}^{\mathbf{II}}\right] _{^{\ast }%
\mathbb{Q}
}=\bigcap \left\{ \alpha _{\mathbf{s}}^{\mathbf{II}}|\alpha _{w}^{\mathbf{II}%
}\subset \text{ }\alpha _{\mathbf{s}}^{\mathbf{II}}\right\} . \\ 
\end{array}
& \left( 1.3.3\right)%
\end{array}%
$

\bigskip

\bigskip

\textbf{Definition 1.3.1.2}$^{\prime }$\textbf{.a.I. (Strong and Weak
Dedekind cut }

\textbf{tipe} \textbf{I in} $^{\ast }%
\mathbb{Q}
$\textbf{).}

(\textbf{1})\textbf{\ }A strong Dedekind cut tipe \textbf{I} $\alpha _{%
\mathbf{s}}=\alpha _{\mathbf{s}}^{\mathbf{I}}$ in $^{\ast }%
\mathbb{Q}
$ is a subset

$\alpha _{\mathbf{s}}^{\mathbf{I}}\subsetneqq $ $^{\ast }%
\mathbb{Q}
$ of the hyperrational numbers $^{\ast }%
\mathbb{Q}
$ that satisfies

these properties:

\textbf{1.} $\alpha _{\mathbf{s}}^{\mathbf{I}}$ is not empty.

\textbf{2.} $\beta _{\mathbf{s}}^{\mathbf{I}}=$ $^{\ast }%
\mathbb{Q}
\backslash \alpha _{\mathbf{s}}^{\mathbf{I}}$ is not empty.

\textbf{3.} $\alpha _{\mathbf{s}}^{\mathbf{I}}$ contains no greatest element

\textbf{4.} For every $x,y\in $ $^{\ast }%
\mathbb{Q}
,$ if $x\in \alpha _{\mathbf{s}}^{\mathbf{I}}$ and $y<x,$ then $y\in \alpha
_{\mathbf{s}}^{\mathbf{I}}$

as well.

(\textbf{2}) A weak Dedekind cut tipe \textbf{I} $\alpha _{w}$ in $^{\ast }%
\mathbb{Q}
$ is a subset

$\alpha _{w}^{\mathbf{I}}\subsetneqq $ $^{\ast }%
\mathbb{Q}
$

of the hyperrational numbers $^{\ast }%
\mathbb{Q}
$ that satisfies these

properties:\ 

\textbf{1.} $\alpha _{w}^{\mathbf{I}}$ is not empty.

\textbf{2. }$\beta _{w}^{\mathbf{I}}=$ $^{\ast }%
\mathbb{Q}
\backslash \alpha _{w}^{\mathbf{I}}$ is not empty.

\textbf{3. }$\alpha _{\mathbf{s}}^{\mathbf{I}}$ contains no greatest element.

\textbf{4.} For every $x,y\in $ $^{\ast }%
\mathbb{Q}
,$ if $x\in \alpha _{w}^{\mathbf{I}}$ and $y<x,$ then there

is exists $z\in \alpha _{w}^{\mathbf{II}}$ such that $z<y$ as well.

\textbf{Remark. }Note that for every weak Dedekind cut $\alpha _{w}^{\mathbf{%
I}}$ in

$^{\ast }%
\mathbb{Q}
$ there is exists unique strong Dedekind cut $\alpha _{\mathbf{s}}^{\mathbf{I%
}}\left( \alpha _{w}^{\mathbf{I}}\right) $

in $^{\ast }%
\mathbb{Q}
$ such that:

\bigskip

$\ \ \ \ \ \ 
\begin{array}{cc}
\begin{array}{c}
\\ 
\alpha _{\mathbf{s}}^{\mathbf{I}}\left( \alpha _{w}^{\mathbf{I}}\right) =%
\left[ \alpha _{w}^{\mathbf{I}}\right] _{^{\ast }%
\mathbb{Q}
}=\bigcap \left\{ \alpha _{\mathbf{s}}^{\mathbf{I}}|\alpha _{w}^{\mathbf{I}%
}\subset \text{ }\alpha _{\mathbf{s}}^{\mathbf{I}}\right\} . \\ 
\end{array}
& \left( 1.3.4\right)%
\end{array}%
$

\bigskip

\textbf{Definition 1.3.1.2}$^{\prime }$\textbf{.a.II. (Strong and Weak
Dedekind cut }

\textbf{tipe} \textbf{II in} $^{\ast }%
\mathbb{Q}
$\textbf{).}

(\textbf{1})\textbf{\ }A strong Dedekind cut tipe \textbf{II} $\alpha _{%
\mathbf{s}}^{\mathbf{II}}$ in $^{\ast }%
\mathbb{Q}
$ is a subset $\alpha _{\mathbf{s}}^{\mathbf{II}}\subsetneqq $ $^{\ast }%
\mathbb{Q}
$

of the hyperrational numbers $^{\ast }%
\mathbb{Q}
$ that satisfies these

properties:

\textbf{1.} $\alpha _{\mathbf{s}}^{\mathbf{II}}$ is not empty.

\textbf{2.} $\beta _{\mathbf{s}}^{\mathbf{II}}=$ $^{\ast }%
\mathbb{Q}
\backslash \alpha _{\mathbf{s}}^{\mathbf{II}}$ is not empty.

\textbf{3.} $\alpha _{\mathbf{s}}^{\mathbf{II}}$ contains a greatest element
or $^{\ast }%
\mathbb{Q}
\backslash \alpha _{\mathbf{s}}^{\mathbf{II}}$ contains no least

element.

\textbf{4.} For every $x,y\in $ $^{\ast }%
\mathbb{Q}
,$ if $x\in \alpha _{\mathbf{s}}^{\mathbf{II}}$ and $y<x,$ then $y\in \alpha
_{\mathbf{s}}^{\mathbf{II}}$ as well.

\bigskip

(\textbf{2}) A weak Dedekind cut tipe \textbf{II} $\alpha _{w}^{\mathbf{II}}$
in $^{\ast }%
\mathbb{Q}
$ is a subset $\alpha _{w}^{\mathbf{II}}\subsetneqq $ $^{\ast }%
\mathbb{Q}
$

of the hyperrational numbers $^{\ast }%
\mathbb{Q}
$ that satisfies these

properties:\ 

\textbf{1.} $\alpha _{w}^{\mathbf{II}}$ is not empty.

\textbf{2. } $\beta _{w}^{\mathbf{II}}=$ $^{\ast }%
\mathbb{Q}
\backslash \alpha _{w}^{\mathbf{II}}$ is not empty.

\textbf{3.} $\alpha _{\mathbf{s}}^{\mathbf{II}}$ contains a greatest element
or $^{\ast }%
\mathbb{Q}
\backslash \alpha _{\mathbf{s}}^{\mathbf{II}}$ contains no least

element.

\textbf{4.} For every $x,y\in $ $^{\ast }%
\mathbb{Q}
,$ if $x\in \alpha _{w}^{\mathbf{II}}$ and $y<x,$ then there is exists

$z\in \alpha _{w}^{\mathbf{II}}$ such that $z<y$ as well.

\textbf{Remark. }Note that for every weak Dedekind cut $\alpha _{w}^{\mathbf{%
II}}$ in $^{\ast }%
\mathbb{Q}
$

there is exists unique strong Dedekind cut $\alpha _{\mathbf{s}}^{\mathbf{II}%
}\left( \alpha _{w}^{\mathbf{II}}\right) $ in $^{\ast }%
\mathbb{Q}
$ such

that: $\alpha _{\mathbf{s}}^{\mathbf{II}}\left( \alpha _{w}^{\mathbf{II}%
}\right) =\left[ \alpha _{w}^{\mathbf{II}}\right] _{^{\ast }%
\mathbb{Q}
}=$ $\bigcap \left\{ \alpha _{\mathbf{s}}^{\mathbf{II}}|\alpha _{w}^{\mathbf{%
II}}\subset \text{ }\alpha _{\mathbf{s}}^{\mathbf{II}}\right\} .$

\bigskip\ $\ \ 
\begin{array}{cc}
\begin{array}{c}
\\ 
\alpha _{\mathbf{s}}^{\mathbf{II}}\left( \alpha _{w}^{\mathbf{II}}\right) =%
\left[ \alpha _{w}^{\mathbf{II}}\right] _{^{\ast }%
\mathbb{Q}
}=\bigcap \left\{ \alpha _{\mathbf{s}}^{\mathbf{II}}|\alpha _{w}^{\mathbf{II}%
}\subset \text{ }\alpha _{\mathbf{s}}^{\mathbf{II}}\right\} . \\ 
\end{array}
& \left( 1.3.5\right)%
\end{array}%
$

\textbf{Definition 1.3.1.2}$^{\prime }$\textbf{.b. (Strong and Weak Dedekind
cut in} $^{\ast }%
\mathbb{Z}
$\textbf{).}\bigskip

\bigskip (\textbf{1})\textbf{\ }A \textbf{strong Dedekind cut} $\alpha _{%
\mathbf{s}}=\alpha _{^{\ast }%
\mathbb{Z}
}^{\mathbf{s}}$ in $^{\ast }%
\mathbb{Z}
$ is a subset

$\alpha _{^{\ast }%
\mathbb{Z}
}^{\mathbf{s}}\subsetneqq $ $^{\ast }%
\mathbb{Z}
$ of the hyperintegers $^{\ast }%
\mathbb{Z}
$ that satisfies these

properties:\ 

\textbf{1.} $\alpha _{^{\ast }%
\mathbb{Z}
}$ is not empty.

\textbf{2.}$\beta _{^{\ast }%
\mathbb{Z}
}^{\mathbf{s}}=$ $^{\ast }%
\mathbb{Z}
\backslash \alpha _{^{\ast }%
\mathbb{Z}
}^{\mathbf{s}}$ is not empty.

\textbf{3.} $^{\ast }%
\mathbb{Z}
\backslash \alpha _{^{\ast }%
\mathbb{Z}
}^{\mathbf{s}}$ contains no least element.

\textbf{4.} For every $x,y\in $ $^{\ast }%
\mathbb{Z}
,$ if $x\in \alpha _{^{\ast }%
\mathbb{Z}
}^{\mathbf{s}}$ and $y<x,y\notin \alpha _{^{\ast }%
\mathbb{Z}
}^{\mathbf{s}}$ then

there is exists $z\in \alpha _{^{\ast }%
\mathbb{Z}
}^{\mathbf{s}}$ such that $z<y$ as well.

(\textbf{2}) A \textbf{weak Dedekind cut} $\alpha _{w}=\alpha _{^{\ast }%
\mathbb{Z}
}^{w}$ in $^{\ast }%
\mathbb{Z}
$ is a subset

$\alpha _{^{\ast }%
\mathbb{Z}
}^{w}\subsetneqq $ $^{\ast }%
\mathbb{Z}
$ of the hyperintegers $^{\ast }%
\mathbb{Z}
$ that satisfies these

properties:\ 

\textbf{1.} $\alpha _{^{\ast }%
\mathbb{Z}
}^{w}$ is not empty.

\textbf{2.}$\beta _{^{\ast }%
\mathbb{Z}
}^{w}=$ $^{\ast }%
\mathbb{Z}
\backslash \alpha _{^{\ast }%
\mathbb{Z}
}^{w}$ is not empty.

\textbf{3.} $^{\ast }%
\mathbb{Z}
\backslash \alpha _{^{\ast }%
\mathbb{Z}
}^{w}$ contains no least element.

\textbf{4.} For every $x,y\in $ $^{\ast }%
\mathbb{Z}
,$ if $x\in \alpha _{^{\ast }%
\mathbb{Z}
}^{w}$ and $y<x,y\notin \alpha _{^{\ast }%
\mathbb{Z}
}^{w}$ then

there is exists $z\in \alpha _{^{\ast }%
\mathbb{Z}
}^{w}$ such that $z<y$ as well.

\textbf{Remark. }Note that for every weak Dedekind cut $\alpha _{w}$ in $%
^{\ast }%
\mathbb{Z}
$

there is exists unique strong Dedekind cut in $^{\ast }%
\mathbb{Z}
$ such that:

$\alpha _{\mathbf{s}}\left( \alpha _{w}\right) =\left[ \alpha _{w}\right]
_{^{\ast }%
\mathbb{Z}
}=$ $\bigcap \left\{ \alpha _{\mathbf{s}}|\alpha _{w}\subset \text{ }\alpha
_{\mathbf{s}}\right\} .$

\textbf{Example. }(\textbf{1})\textbf{\ }Let $\Omega $ denotes the set $%
\left\{ n|n\in \text{ }^{\ast }%
\mathbb{Z}
\backslash ^{\ast }%
\mathbb{Z}
_{+\infty }\right\} .$

It is easy to see that $\Omega $ is a strong Dedekind cut in $^{\ast }%
\mathbb{Z}
.$

(\textbf{2}) Let $\Omega _{k},k\in 
\mathbb{N}
,k\neq 0$ denotes the set

$\left\{ n|\left( n\in \text{ }^{\ast }%
\mathbb{Z}
\backslash ^{\ast }%
\mathbb{Z}
_{+\infty }\right) \wedge \left( n|k\right) \right\} .$

It is easy to see that $\Omega $ is a weak Dedekind cut in $^{\ast }%
\mathbb{Z}
$

and $\left[ \Omega _{k}\right] _{^{\ast }%
\mathbb{Z}
}=\Omega .$

\textbf{Definition 1.3.1.3.a.}(\textbf{Wattenberg-Dedekind hyperreal }

\textbf{numbers})

(\textbf{1}) A Wattenberg-Dedekind hyperreal number

$\alpha \in $ $^{\ast }%
\mathbb{R}
_{\mathbf{d}}$ is a strong Dedekind cut $\alpha =\alpha _{\mathbf{s}}$ in $%
^{\ast }%
\mathbb{Q}
.$

Equivalntly:

(\textbf{2}) A Wattenberg-Dedekind hyperreal number

$\alpha \in $ $^{\ast }%
\mathbb{R}
_{\mathbf{d}}$ is a weak Dedekind cut $\alpha =\alpha _{w}$ in $^{\ast }%
\mathbb{Q}
.$

(\textbf{3}) We denote the set of all Wattenberg-Dedekind

hyperreal numbers by $^{\ast }%
\mathbb{R}
_{\mathbf{d}}$ and we order them by

set-theoretic inclusion, that is to say, for any

$\alpha ,\beta \in $ $^{\ast }%
\mathbb{R}
_{\mathbf{d}},$ $\alpha <\beta $ \ if and only if $\alpha \subsetneqq \beta $
where the

inclusion is strict.

We further define $\alpha =\beta $ as real numbers

if and are equal as sets. As usual, we write

$\alpha \leqslant \beta $ if $\alpha <\beta $ or $\alpha =\beta $.

\textbf{Definition 1.3.1.3.b. }(\textbf{Wattenberg-Dedekind hyperrationals }$%
^{\ast }%
\mathbb{Q}
_{\mathbf{d}}$)

(\textbf{1}) A Wattenberg-Dedekind hyperrational is a weak Dedekind cut

$\alpha =\alpha _{^{\ast }%
\mathbb{Q}
}$ in $^{\ast }%
\mathbb{Q}
.$

\textbf{Definition 1.3.1.3.c. }(\textbf{Wattenberg-Dedekind hyperintegers }$%
^{\ast }%
\mathbb{Z}
_{\mathbf{d}}$)

(\textbf{1}) A Wattenberg-Dedekind hyperinteger is a weak Dedekind cut

$\alpha =\alpha _{^{\ast }%
\mathbb{Z}
}$ in $^{\ast }%
\mathbb{Z}
.$

(\textbf{2}) We denote the set of all Wattenberg-Dedekind

hyperintegers by $^{\ast }%
\mathbb{Z}
_{\mathbf{d}}$ and we order them by suitable

set-theoretic inclusion, that is to say, for any $\alpha _{^{\ast }%
\mathbb{Z}
},\beta _{^{\ast }%
\mathbb{Z}
}\in $ $^{\ast }%
\mathbb{Z}
_{\mathbf{d}},$

$\alpha _{^{\ast }%
\mathbb{Z}
}<\beta _{^{\ast }%
\mathbb{Z}
}$ if and only if $\left[ \alpha _{^{\ast }%
\mathbb{Z}
}\right] _{^{\ast }%
\mathbb{Q}
}\subsetneqq \left[ \beta _{^{\ast }%
\mathbb{Z}
}\right] _{^{\ast }%
\mathbb{Q}
}$ where the inclusion is

strict. We further define:

(\textbf{3}) \textbf{weak equality:}

\bigskip\ $\ \ \ \ \ \ \ \ \ \ \ \ \ \ \ \ \ \ \ \ \ \ \ \ \ \ \ \ \ \ \ \ \
\ \ \ \ \ \ \ \ \ \ \ \ \ \ 
\begin{array}{cc}
\begin{array}{c}
\\ 
\alpha _{^{\ast }%
\mathbb{Z}
}=_{w}\beta _{^{\ast }%
\mathbb{Z}
} \\ 
\end{array}
& \left( 1.3.6\right)%
\end{array}%
$

as Wattenberg-Dedekind hyperintegers iff Dedekind cut $\left[ \alpha \right]
_{^{\ast }%
\mathbb{Q}
}$

and $\left[ \beta \right] _{^{\ast }%
\mathbb{Q}
}$ are \ equal as sets,i.e.

\ \ \ \ $%
\begin{array}{cc}
\begin{array}{c}
\\ 
\forall x\left\{ x\in \left[ \alpha \right] _{^{\ast }%
\mathbb{Q}
}\iff x\in \left[ \beta \right] _{^{\ast }%
\mathbb{Q}
}\right\} \\ 
\end{array}
& \left( 1.3.7\right)%
\end{array}%
$

\bigskip

As usual, we write $\alpha \leqslant _{w}\beta $ if $\alpha <\beta $ or $%
\alpha =_{w}\beta $.

\bigskip

(\textbf{4}) \textbf{strong equality:}

\bigskip

\bigskip\ \ \ \ \ \ $%
\begin{array}{cc}
\begin{array}{c}
\\ 
\alpha _{^{\ast }%
\mathbb{Z}
}=_{s}\beta _{^{\ast }%
\mathbb{Z}
} \\ 
\end{array}
& \left( 1.3.8\right)%
\end{array}%
$

as Wattenberg-Dedekind hyperintegers iff Dedekind cut $\alpha _{^{\ast }%
\mathbb{Z}
}$

and $\beta _{^{\ast }%
\mathbb{Z}
}$ are \ equal as sets,i.e.

\bigskip

\bigskip\ \ \ \ \ \ $%
\begin{array}{cc}
\begin{array}{c}
\\ 
\forall x\left\{ x\in \alpha _{^{\ast }%
\mathbb{Z}
}\iff x\in \beta _{^{\ast }%
\mathbb{Z}
}\right\} \\ 
\end{array}
& \left( 1.3.9\right)%
\end{array}%
$

As usual, we write $\alpha _{^{\ast }%
\mathbb{Z}
}\leqslant _{s}\beta _{^{\ast }%
\mathbb{Z}
}$ if $\alpha _{^{\ast }%
\mathbb{Z}
}<\beta _{^{\ast }%
\mathbb{Z}
}$ or $\alpha _{^{\ast }%
\mathbb{Z}
}=_{s}\beta _{^{\ast }%
\mathbb{Z}
}$.

\textbf{Remark. }Note that we often write formula $\alpha _{^{\ast }%
\mathbb{Z}
}=_{s}\beta _{^{\ast }%
\mathbb{Z}
}$ as

$\alpha _{^{\ast }%
\mathbb{Z}
}=\beta _{^{\ast }%
\mathbb{Z}
}.$

\textbf{Definition 1.3.1.4.} Dedekind hyperreal $\alpha $ is said to be

\textit{Dedekind} \textit{hyperirrational }if $^{\ast }%
\mathbb{Q}
\backslash \alpha $ contains no least element.

\textbf{Theorem 1.3.1.1.} Every nonempty subset $A\subsetneqq $ $^{\ast }%
\mathbb{R}
_{\mathbf{d}}$ of Dedekind \ \ \ \ \ \ \ \ \ \ \ \ \ \ \ \ \ \ \ \ \ \ \ \ \
\ \ \ \ \ \ \ \ 

hyperreal numbers that is bounded (hyperbounded) above has a least \ \ \ \ \
\ \ \ \ \ \ \ \ \ \ \ \ \ \ \ \ \ \ \ \ \ \ 

upper bound.

\textbf{Proof.} Let $A$ be a nonempty set of hyperreal numbers, such that
for every $\alpha \in A$ we have that $\alpha \leqslant \gamma $ for some
real number $\gamma \in $ $^{\ast }%
\mathbb{R}
_{\mathbf{d}}.$Now define the set $\sup A=\dbigcup\limits_{\alpha \in
A}\alpha .$ We must show that this set is a Wattenberg-Dedekind hyperreal
number. This amounts to checking the four conditions of a Dedekind cut. $%
\sup A$ is clearly not empty, for it is the nonempty union of nonempty sets.
Because $\gamma $ is a Wattenberg-Dedekind hyperreal number, there is some
hyperrational $x\in $ $^{\ast }%
\mathbb{Q}
$ that is not in $\gamma .$ Since every $\alpha \in A$ is a subset of $%
\gamma ,x$ is not \ \ \ \ \ \ \ \ \ \ \ \ \ \ \ \ \ \ in any $\alpha ,$ so $%
x\notin \sup A$ either. Thus, $^{\ast }%
\mathbb{Q}
\backslash \sup A$ is nonempty. If $\sup A$ had a greatest element $g\in $ $%
^{\ast }%
\mathbb{Q}
,$ then $g\in \alpha $ for some $\alpha \in A.$ Then $g$ would be a greatest
element of $\alpha ,$ but $\alpha $ is a Wattenberg-Dedekind hyperreal
number, so by contrapositive law,$\sup A$ has no greatest element. Lastly,
if $x\in $ $^{\ast }%
\mathbb{Q}
$ and $x\in \sup A,$ then $x\in \alpha $ for some $\alpha ,$ so given any $%
y\in $ $^{\ast }%
\mathbb{Q}
,$ $y<x$ because $\alpha $ is a Dedekind hyperreal number $y\in \alpha $
whence $y\in \sup A.$\bigskip Thus $\sup A,$ is a Wattenberg-Dedekind
hyperreal number.Trivially,$\sup A\ $is an upper bound of $A,$ for every $%
\alpha \in A,$ $\alpha \subseteqq \sup A.$ It now suffices to prove that $%
\sup A\leqslant \gamma ,$because was an arbitrary upper bound. But this is
easy, because every $x\in \sup A,x\in $ $^{\ast }%
\mathbb{Q}
$ is an element of $\alpha $ for some $\alpha \in A,$ so because $\alpha
\subseteq \gamma ,$ $x\in \gamma .$ Thus, $\sup A$ is the least upper bound
of $A$.

\textbf{Definition 1.3.1.5.a.} Given two Wattenberg-Dedekind hyperreal
numbers $\alpha $ and $\beta $ we define:

\textbf{1.}The additive identity (zero cut) $0_{^{\ast }%
\mathbb{R}
_{\mathbf{d}}},$ denoted $0,$is

$0_{^{\ast }%
\mathbb{R}
_{\mathbf{d}}}\triangleq \left\{ x\in \text{ }^{\ast }%
\mathbb{Q}
|\text{ }x<0\right\} .$

\textbf{2.}The multiplicative identity $1_{^{\ast }%
\mathbb{R}
_{\mathbf{d}}},$ denoted $1,$is

$1_{^{\ast }%
\mathbb{R}
_{\mathbf{d}}}\triangleq \left\{ x\in \text{ }^{\ast }%
\mathbb{Q}
|\text{ }x<_{^{\ast }%
\mathbb{Q}
}1_{^{\ast }%
\mathbb{Q}
}\right\} .$

\textbf{3.} Addition $\alpha +_{^{\ast }%
\mathbb{R}
_{\mathbf{d}}}\beta $ of $\alpha $ and $\beta $ denoted $\alpha +\beta $ is

$\alpha +\beta \triangleq \left\{ x+y|\text{ }x\in \alpha ,y\in \beta
\right\} .$

It is easy to see that $\alpha +_{^{\ast }%
\mathbb{R}
_{\mathbf{d}}}0_{^{\ast }%
\mathbb{R}
_{\mathbf{d}}}=\alpha $ for all $\alpha \in $ $^{\ast }%
\mathbb{R}
_{\mathbf{d}}.$

It is easy to see that $\alpha +_{^{\ast }%
\mathbb{R}
_{\mathbf{d}}}\beta $ is a cut in $^{\ast }%
\mathbb{Q}
$ and $\alpha +_{^{\ast }%
\mathbb{R}
_{\mathbf{d}}}\beta =\beta +_{^{\ast }%
\mathbb{R}
_{\mathbf{d}}}\alpha .$

Another fundamental property of cut addition is associativity:

$\left( \alpha +_{^{\ast }%
\mathbb{R}
_{\mathbf{d}}}\beta \right) +_{^{\ast }%
\mathbb{R}
_{\mathbf{d}}}\gamma =\alpha +_{^{\ast }%
\mathbb{R}
_{\mathbf{d}}}\left( \beta +_{^{\ast }%
\mathbb{R}
_{\mathbf{d}}}\gamma \right) .$

This follows from the corresponding property of $^{\ast }%
\mathbb{Q}
.$

\textbf{4.}The opposite $-_{^{\ast }%
\mathbb{R}
_{\mathbf{d}}}\alpha $ of $\alpha ,$ denoted $-\alpha ,$ is

$-\alpha \triangleq \left\{ x\in \text{ }^{\ast }%
\mathbb{Q}
|\text{ }-x\notin \alpha ,-x\text{ is not the least element of }^{\ast }%
\mathbb{Q}
\backslash \alpha \right\} .$

\textbf{5.Remark.} We also say that the opposite $-\alpha $ of $\alpha $ is
the \textbf{additive }

\textbf{inverse} of $\alpha $ denoted $\div \alpha $ iff the next equality
is satisfied:

$\alpha +\left( \div \alpha \right) =0.$

\textbf{6.Remark. }It is easy to see that for all internal cut $\alpha ^{%
\mathbf{Int}}$ the opposite

$-\alpha ^{\mathbf{Int}}$ is the additive inverse of $\alpha ^{\mathbf{Int}%
}, $i.e. $\alpha ^{\mathbf{Int}}+\left( \div \alpha ^{\mathbf{Int}}\right)
=0.$

\textbf{7.Example}. (External cut $X$ without additive inverse $\div X$) For
any

$x,y\in $ $^{\ast }%
\mathbb{R}
$ we denote:

$x\ll \infty \triangleq \exists r\left[ \left( r\in 
\mathbb{R}
\right) \wedge \left( x<\text{ }^{\ast }r\right) \right] ,$ $y\approx
-\infty \triangleq \forall r\left[ r\in 
\mathbb{R}
\implies x<\text{ }^{\ast }r\right] .$

Let us consider two Dedekind hyperreal numbers $X$ and $Y$ defined as:

$X=\left\{ x|\left( x\in \text{ }^{\ast }%
\mathbb{R}
\right) \wedge \left( x\ll \infty \right) \right\} ,Z=\left\{ x|z<0\right\}
. $

It is easy to see that is no exist cut $Y$ such that: $X+Y=Z.$

\textbf{Proof.} Suppose that cut $Y$ such that: $X+Y=Z$ exist. It is easy to

check that $\forall y\left[ y\in Y\implies y\approx -\infty \right] .$%
Suppose that $y\in Y,$then

$\forall x\left[ x\in X\implies x+y\in Z\right] ,$i.e.$\forall \left( x\ll
\infty \right) \left[ x+y<0\right] .$Hence

$\forall \left( x\ll \infty \right) \left[ y<-x\right] ,$i.e. $y\approx
-\infty .$It is easy to check that $Z\nsubseteqq X+Y.$

If $x\in X$ and $y\in Y$ then $x\ll \infty $ and $y\approx -\infty ,$hence $%
x+y\neq -1,$i.e.

$-1\notin X+Y.$Thus $Z\nsubseteqq X+Y.$This is a contradiction.

\textbf{8.}We say that the cut $\alpha $ is positive if $0<\alpha $ or
negative if $\alpha <0.$

The absolute value of $\alpha ,$denoted $\left\vert \alpha \right\vert ,$is $%
\left\vert \alpha \right\vert \triangleq \alpha ,$if $\alpha \geq 0$ and $%
\left\vert \alpha \right\vert \triangleq -\alpha ,$

if $\alpha \leq 0$

\textbf{9.}If $\alpha ,\beta >0$ then multiplication $\alpha \times _{^{\ast
}%
\mathbb{R}
_{\mathbf{d}}}\beta $ of $\alpha $ and $\beta $ denoted $\alpha \times \beta 
$ is

$\alpha \times \beta \triangleq \left\{ z\in \text{ }^{\ast }%
\mathbb{Q}
|\text{ }z=x\times y\text{ for some }x\in \alpha ,y\in \beta \text{ with }%
x,y>0\right\} .$

In general, $\alpha \times \beta =0$ if $\alpha =\mathbf{0}$ or $\beta =0%
\mathbf{,}$

$\alpha \times \beta \triangleq \left\vert \alpha \right\vert \times
\left\vert \beta \right\vert $ if $\alpha >0,\beta >0$ or $\alpha <0,\beta <0%
\mathbf{,}$

$\alpha \times \beta \triangleq -\left( \left\vert \alpha \right\vert \cdot
\left\vert \beta \right\vert \right) $ if $\alpha >0,\beta <0\mathbf{,}$or $%
\alpha <0,\beta >0\mathbf{.}$

\textbf{10. }The cut order enjois on $^{\ast }%
\mathbb{R}
_{\mathbf{d}}$ the standard additional properties of:

(\textbf{i}) \ \ \textbf{transitivity: }$\alpha \leq \beta \leq \gamma
\implies \alpha \leq \gamma .$\textbf{\ }

(\textbf{ii}) \ \textbf{trichotomy: }eizer $\alpha <\beta ,\beta <\alpha $
or $\alpha =\beta $ but only one

of the three things is true.

(\textbf{iii}) \textbf{translation: }$\alpha \leq \beta \implies \alpha
+_{^{\ast }%
\mathbb{R}
_{\mathbf{d}}}\gamma \leq \beta +_{^{\ast }%
\mathbb{R}
_{\mathbf{d}}}\gamma .$

\textbf{11.}By definition above, this is what we mean when we say that $%
^{\ast }%
\mathbb{R}
_{\mathbf{d}}$ is an

\textbf{ordered pseudo-ring }or \textbf{ordered pseudo-field.}

\textbf{Remark.} We embed $^{\ast }%
\mathbb{R}
$ in $^{\ast }%
\mathbb{R}
_{\mathbf{d}}$ in the standard way [24].If $\alpha \in $ $^{\ast }%
\mathbb{R}
$ the

corresponding element, $\alpha ^{\#},$ of $^{\ast }%
\mathbb{R}
_{\mathbf{d}}$ is\bigskip

$\ 
\begin{array}{cc}
\begin{array}{c}
\\ 
\alpha ^{\#}\triangleq \left\{ x\in \text{ }^{\ast }%
\mathbb{R}
|x<_{^{\ast }%
\mathbb{R}
}\alpha \right\} \\ 
\end{array}
& \left( 1.3.10\right)%
\end{array}%
$

\bigskip

\textbf{Definition 1.3.1.5.b.} Given two Wattenberg-Dedekind

hyperintegers $\alpha _{^{\ast }%
\mathbb{Z}
}\in $ $^{\ast }%
\mathbb{Z}
_{\mathbf{d}}$ and $\beta _{^{\ast }%
\mathbb{Z}
}\in $ $^{\ast }%
\mathbb{Z}
_{\mathbf{d}}$ we define:

\textbf{1.}The additive identity (zero cut) denoted $0_{^{\ast }%
\mathbb{Z}
_{\mathbf{d}}}$ or $0,$is

$0_{^{\ast }%
\mathbb{Z}
_{\mathbf{d}}}\triangleq \left\{ x\in \text{ }^{\ast }%
\mathbb{Z}
|\text{ }x\leq 0\right\} .$

\textbf{2.}The multiplicative identity denoted $1_{^{\ast }%
\mathbb{Z}
_{\mathbf{d}}}$ or $1,$is

$1_{^{\ast }%
\mathbb{R}
_{\mathbf{d}}}\triangleq \left\{ x\in \text{ }^{\ast }%
\mathbb{Z}
|\text{ }x\leq _{^{\ast }%
\mathbb{Z}
}1_{^{\ast }%
\mathbb{Z}
}\right\} .$

\textbf{3.} Addition of $\alpha =\alpha _{^{\ast }%
\mathbb{Z}
}$ and $\beta =\beta _{^{\ast }%
\mathbb{Z}
}$ denoted $\alpha +_{^{\ast }%
\mathbb{Z}
_{\mathbf{d}}}\beta $

or $\alpha +\beta $ is $\alpha +_{^{\ast }%
\mathbb{Z}
_{\mathbf{d}}}\beta \triangleq \left\{ x+_{^{\ast }%
\mathbb{Z}
}y|\text{ }x\in \alpha ,y\in \beta \right\} .$

It is easy to see that $\alpha +_{^{\ast }%
\mathbb{Z}
_{\mathbf{d}}}0_{^{\ast }%
\mathbb{Z}
_{\mathbf{d}}}=\alpha $ for all $\alpha \in $ $^{\ast }%
\mathbb{Z}
_{\mathbf{d}}.$

It is easy to see that $\alpha +_{^{\ast }%
\mathbb{Z}
_{\mathbf{d}}}\beta $ is a cut in $^{\ast }%
\mathbb{Z}
$ and

$\alpha +_{^{\ast }%
\mathbb{Z}
_{\mathbf{d}}}\beta =\beta +_{^{\ast }%
\mathbb{Z}
_{\mathbf{d}}}\alpha .$

Another fundamental property of cut addition is

associativity:

$\left( \alpha +_{^{\ast }%
\mathbb{Z}
_{\mathbf{d}}}\beta \right) +_{^{\ast }%
\mathbb{Z}
_{\mathbf{d}}}\gamma =\alpha +_{^{\ast }%
\mathbb{Z}
_{\mathbf{d}}}\left( \beta +_{^{\ast }%
\mathbb{Z}
_{\mathbf{d}}}\gamma \right) .$

This follows from the corresponding property of $^{\ast }%
\mathbb{Z}
.$

\textbf{4.}The opposite of $\alpha =\alpha _{^{\ast }%
\mathbb{Z}
},$ denoted $-_{^{\ast }%
\mathbb{Z}
_{\mathbf{d}}}\alpha $ or $-\alpha ,$ is

$-_{^{\ast }%
\mathbb{Z}
_{\mathbf{d}}}\alpha \triangleq \left\{ x\in \text{ }^{\ast }%
\mathbb{Z}
|\text{ }-x\notin \alpha \right\} .$

\textbf{5.Remark.} We also say that the opposite $-_{^{\ast }%
\mathbb{Z}
_{\mathbf{d}}}\alpha _{^{\ast }%
\mathbb{Z}
}$ of $\alpha _{^{\ast }%
\mathbb{Z}
}$

is the \textbf{additive inverse} of $\alpha _{^{\ast }%
\mathbb{Z}
}$ denoted $\div _{^{\ast }%
\mathbb{Z}
_{\mathbf{d}}}\alpha $ or $\div \alpha $ iff the

next equality is satisfied: $\alpha +_{^{\ast }%
\mathbb{Z}
_{\mathbf{d}}}\left( \div _{^{\ast }%
\mathbb{Z}
_{\mathbf{d}}}\alpha \right) =0_{^{\ast }%
\mathbb{Z}
_{\mathbf{d}}}.$

\textbf{6.Remark. }It is easy to see that for all internal cut $\alpha
_{^{\ast }%
\mathbb{Z}
_{\mathbf{d}}}^{\mathbf{Int}}$ the

opposite $-_{^{\ast }%
\mathbb{Z}
_{\mathbf{d}}}\alpha _{^{\ast }%
\mathbb{Z}
}^{\mathbf{Int}}$ is the additive inverse of $\alpha _{^{\ast }%
\mathbb{Z}
}^{\mathbf{Int}},$i.e.

$\alpha _{^{\ast }%
\mathbb{Z}
}^{\mathbf{Int}}+_{^{\ast }%
\mathbb{Z}
_{\mathbf{d}}}\left( \div _{^{\ast }%
\mathbb{Z}
_{\mathbf{d}}}\alpha _{^{\ast }%
\mathbb{Z}
}^{\mathbf{Int}}\right) =0_{^{\ast }%
\mathbb{Z}
_{\mathbf{d}}}.$

\textbf{8.}We say that the cut $\alpha _{^{\ast }%
\mathbb{Z}
}$ is positive if $0_{^{\ast }%
\mathbb{Z}
_{\mathbf{d}}}<\alpha _{^{\ast }%
\mathbb{Z}
}$ or

negative if $\alpha _{^{\ast }%
\mathbb{Z}
}<0_{^{\ast }%
\mathbb{Z}
_{\mathbf{d}}}.$

The absolute value of $\alpha _{^{\ast }%
\mathbb{Z}
},$denoted $\left\vert \alpha _{^{\ast }%
\mathbb{Z}
}\right\vert ,$is $\left\vert \alpha _{^{\ast }%
\mathbb{Z}
}\right\vert \triangleq \alpha _{^{\ast }%
\mathbb{Z}
},$

if $\alpha _{^{\ast }%
\mathbb{Z}
}\geq 0_{^{\ast }%
\mathbb{Z}
_{\mathbf{d}}}$ and $\left\vert \alpha _{^{\ast }%
\mathbb{Z}
}\right\vert \triangleq -_{^{\ast }%
\mathbb{Z}
_{\mathbf{d}}}\alpha _{^{\ast }%
\mathbb{Z}
},$if $\alpha _{^{\ast }%
\mathbb{Z}
}\leq 0_{^{\ast }%
\mathbb{Z}
_{\mathbf{d}}}.$

\textbf{9.}If $\alpha =\alpha _{^{\ast }%
\mathbb{Z}
},\beta =\beta _{^{\ast }%
\mathbb{Z}
}>0_{^{\ast }%
\mathbb{Z}
_{\mathbf{d}}}$ then multiplication of $\alpha $ and $\beta $

denoted $\alpha \times _{^{\ast }%
\mathbb{Z}
_{\mathbf{d}}}\beta $ or $\alpha \times \beta $ is

$\alpha \times _{^{\ast }%
\mathbb{Z}
_{\mathbf{d}}}\beta \triangleq \left\{ z\in \text{ }^{\ast }%
\mathbb{Z}
|\text{ }z=x\times _{^{\ast }%
\mathbb{Z}
}y\text{ for some }x\in \alpha ,y\in \beta \text{ with }x,y>0_{^{\ast }%
\mathbb{Z}
}\right\} .$

In general, $\alpha \times _{^{\ast }%
\mathbb{Z}
_{\mathbf{d}}}\beta =0_{^{\ast }%
\mathbb{Z}
_{\mathbf{d}}}$ if $\alpha =0_{^{\ast }%
\mathbb{Z}
_{\mathbf{d}}}$ or $\beta =0_{^{\ast }%
\mathbb{Z}
_{\mathbf{d}}}\mathbf{,}$

$\alpha \times _{^{\ast }%
\mathbb{Z}
_{\mathbf{d}}}\beta \triangleq \left\vert \alpha \right\vert \times _{^{\ast
}%
\mathbb{Z}
_{\mathbf{d}}}\left\vert \beta \right\vert $ if $\alpha >0_{^{\ast }%
\mathbb{Z}
_{\mathbf{d}}},\beta >0_{^{\ast }%
\mathbb{Z}
_{\mathbf{d}}}$ or $\alpha <0_{^{\ast }%
\mathbb{Z}
_{\mathbf{d}}},\beta <0_{^{\ast }%
\mathbb{Z}
_{\mathbf{d}}}\mathbf{,}$

$\alpha \times _{^{\ast }%
\mathbb{Z}
_{\mathbf{d}}}\beta \triangleq -\left( \left\vert \alpha \right\vert \times
_{^{\ast }%
\mathbb{Z}
_{\mathbf{d}}}\left\vert \beta \right\vert \right) $ if $\alpha >0_{^{\ast }%
\mathbb{Z}
_{\mathbf{d}}},\beta <0_{^{\ast }%
\mathbb{Z}
_{\mathbf{d}}}\mathbf{,}$or $\alpha <0_{^{\ast }%
\mathbb{Z}
_{\mathbf{d}}},\beta >0_{^{\ast }%
\mathbb{Z}
_{\mathbf{d}}}\mathbf{.}$

\textbf{10. }The cut order enjois on $^{\ast }%
\mathbb{Z}
_{\mathbf{d}}$ the standard additional properties of:

(\textbf{i}) \ \ \textbf{transitivity: }$\alpha _{^{\ast }%
\mathbb{Z}
}\leq \beta _{^{\ast }%
\mathbb{Z}
}\leq \gamma _{^{\ast }%
\mathbb{Z}
}\implies \alpha _{^{\ast }%
\mathbb{Z}
}\leq \gamma _{^{\ast }%
\mathbb{Z}
}.$\textbf{\ }

(\textbf{ii}) \ \textbf{trichotomy: }eizer $\alpha _{^{\ast }%
\mathbb{Z}
}<\beta _{^{\ast }%
\mathbb{Z}
},\beta _{^{\ast }%
\mathbb{Z}
}<\alpha _{^{\ast }%
\mathbb{Z}
}$ or $\alpha _{^{\ast }%
\mathbb{Z}
}=\beta _{^{\ast }%
\mathbb{Z}
}$ but only one

of the three things is true.

(\textbf{iii}) \textbf{translation: }$\alpha _{^{\ast }%
\mathbb{Z}
}\leq \beta _{^{\ast }%
\mathbb{Z}
}\implies \alpha _{^{\ast }%
\mathbb{Z}
}+_{^{\ast }%
\mathbb{R}
_{\mathbf{d}}}\gamma _{^{\ast }%
\mathbb{Z}
}\leq \beta _{^{\ast }%
\mathbb{Z}
}+_{^{\ast }%
\mathbb{R}
_{\mathbf{d}}}\gamma _{^{\ast }%
\mathbb{Z}
}.$

\bigskip

\textbf{Lemma 1.3.1.1.}[24].

(\textbf{i}) Addition\textbf{\ }$\left( \circ +_{^{\ast }%
\mathbb{R}
_{\mathbf{d}}}\circ \right) $ is commutative and associative in$^{\ast }%
\mathbb{R}
_{\mathbf{d}}.$

(\textbf{ii}) $\forall \alpha \in $ $^{\ast }%
\mathbb{R}
:\alpha +_{^{\ast }%
\mathbb{R}
_{\mathbf{d}}}0_{^{\ast }%
\mathbb{R}
_{\mathbf{d}}}=\alpha .$

(\textbf{iii}) $\forall \alpha ,\beta \in $ $^{\ast }%
\mathbb{R}
:\alpha ^{\#}+_{^{\ast }%
\mathbb{R}
_{\mathbf{d}}}\beta ^{\#}=\left( \alpha +_{^{\ast }%
\mathbb{R}
_{\mathbf{d}}}\beta \right) ^{\#}.$

\textbf{Proof. }(\textbf{i}) Is clear from definitions.

(\textbf{ii}) Suppose $a\in \alpha .$ Since $a$ has no greatest element $%
\exists b\left[ \left( b>a\right) \wedge \left( b\in \alpha \right) \right]
. $

Thus $a-b\in 0_{^{\ast }%
\mathbb{R}
_{\mathbf{d}}}$ and $a=(a$ $-$ $b)+b\in 0_{^{\ast }%
\mathbb{R}
_{\mathbf{d}}}+\alpha .$

(\textbf{iii}) (a) $\alpha ^{\#}+_{^{\ast }%
\mathbb{R}
_{\mathbf{d}}}\beta ^{\#}\subseteqq \left( \alpha +_{^{\ast }%
\mathbb{R}
_{\mathbf{d}}}\beta \right) ^{\#}$ is clear since:

$\left( x<\alpha \right) \wedge \left( y<\beta \right) \implies $ $%
x+y<\alpha +\beta .$

(b) Suppose $x<\alpha +\beta .$ Thus $\alpha -\dfrac{\left( \alpha +\beta
\right) -x}{2}<\alpha $ and

$\beta -\dfrac{\left( \alpha +\beta \right) -x}{2}<\beta .$So one obtain

$x=\left[ \left( \alpha -\dfrac{\left( \alpha +\beta \right) -x}{2}\right)
+\left( \beta -\dfrac{\left( \alpha +\beta \right) -x}{2}\right) \right] \in
\alpha ^{\#}+_{^{\ast }%
\mathbb{R}
_{\mathbf{d}}}\beta ^{\#},$

$\left( \alpha +_{^{\ast }%
\mathbb{R}
_{\mathbf{d}}}\beta \right) ^{\#}\subseteqq \alpha ^{\#}+_{^{\ast }%
\mathbb{R}
_{\mathbf{d}}}\beta ^{\#}.$

Notice, here again something is lost going from $^{\ast }%
\mathbb{R}
$ to $^{\ast }%
\mathbb{R}
_{\mathbf{d}}$ since $a<\beta $ does

not imply $\alpha +\alpha <\beta +\alpha $ since $0<\varepsilon _{\mathbf{d}%
} $ but $0+\varepsilon _{\mathbf{d}}=\varepsilon _{\mathbf{d}}+\varepsilon _{%
\mathbf{d}}=\varepsilon _{\mathbf{d}}.$

\textbf{Lemma 1.3.1.2.}[24].

\textbf{(i) }$\leq _{^{\ast }%
\mathbb{R}
_{\mathbf{d}}}$a linear ordering on $^{\ast }%
\mathbb{R}
_{\mathbf{d}},$which extends the usual ordering on $^{\ast }%
\mathbb{R}
$.

\textbf{(ii) }$\left( \alpha \leq \alpha ^{\prime }\right) \wedge \left(
\beta \leq \beta ^{\prime }\right) \implies \alpha +\beta \leq \alpha
^{\prime }+\beta ^{\prime }.$

\textbf{(iii)} $\left( \alpha <\alpha ^{\prime }\right) \wedge \left( \beta
<\beta ^{\prime }\right) \implies \alpha +\beta <\alpha ^{\prime }+\beta
^{\prime }.$

\textbf{(iv) }$^{\ast }%
\mathbb{R}
$ \textbf{is dense in} $^{\ast }%
\mathbb{R}
_{\mathbf{d}}.$That is if $\alpha <\beta $ in $^{\ast }%
\mathbb{R}
_{\mathbf{d}}$ there is an $a\in $ $^{\ast }%
\mathbb{R}
$ then

$\ \ \ \ \ \ \ \alpha <a^{\#}<\beta .$

\bigskip

\textbf{Lemma 1.3.1.3.}[24].

\textbf{(i) }If $\alpha \in $ $^{\ast }%
\mathbb{R}
$ then $-_{^{\ast }%
\mathbb{R}
_{\mathbf{d}}}\left( \alpha ^{\#}\right) =\left( -_{^{\ast }%
\mathbb{R}
}\alpha \right) ^{\#}.$

\textbf{(ii) }$\mathbf{-}_{^{\ast }%
\mathbb{R}
_{\mathbf{d}}}\left( -_{^{\ast }%
\mathbb{R}
_{\mathbf{d}}}\alpha \right) =\alpha .$

\textbf{(iii) }$\alpha \leq _{^{\ast }%
\mathbb{R}
_{\mathbf{d}}}\beta \iff -_{^{\ast }%
\mathbb{R}
_{\mathbf{d}}}\beta \leq _{^{\ast }%
\mathbb{R}
_{\mathbf{d}}}-_{^{\ast }%
\mathbb{R}
_{\mathbf{d}}}\alpha .$

\textbf{(iv) }$\left( \mathbf{-}_{^{\ast }%
\mathbb{R}
_{\mathbf{d}}}\alpha \right) +_{_{^{\ast }%
\mathbb{R}
_{\mathbf{d}}}}\left( -_{_{^{\ast }%
\mathbb{R}
_{\mathbf{d}}}}\beta \right) \leq _{_{^{\ast }%
\mathbb{R}
_{\mathbf{d}}}}\mathbf{-}_{^{\ast }%
\mathbb{R}
_{\mathbf{d}}}\left( \alpha +_{_{^{\ast }%
\mathbb{R}
_{\mathbf{d}}}}\beta \right) .$

\textbf{(v) }$\forall a\in $ $^{\ast }%
\mathbb{R}
:\left( \mathbf{-}_{^{\ast }%
\mathbb{R}
}a\right) ^{\#}+_{_{^{\ast }%
\mathbb{R}
_{\mathbf{d}}}}\left( -_{_{^{\ast }%
\mathbb{R}
_{\mathbf{d}}}}\beta \right) =-_{^{\ast }%
\mathbb{R}
_{\mathbf{d}}}\left( a^{\#}+_{^{\ast }%
\mathbb{R}
_{\mathbf{d}}}\beta \right) .$

\textbf{(vi) }$\alpha +_{^{\ast }%
\mathbb{R}
_{\mathbf{d}}}\left( -_{^{\ast }%
\mathbb{R}
_{\mathbf{d}}}\alpha \right) \leq _{^{\ast }%
\mathbb{R}
_{\mathbf{d}}}0_{^{\ast }%
\mathbb{R}
_{\mathbf{d}}}.$

\bigskip

\bigskip \textbf{Lemma 1.3.1.4.}[24].

\textbf{(i) }$\forall a,b\in $ $^{\ast }%
\mathbb{R}
:\left( a\times _{^{\ast }%
\mathbb{R}
}b\right) ^{\#}=a^{\#}\times _{^{\ast }%
\mathbb{R}
_{\mathbf{d}}}b^{\#}.$

\textbf{(ii) }Multiplication $\left( \cdot \times _{^{\ast }%
\mathbb{R}
}\cdot \right) $ is associative and commutative:

\textbf{\ \ \ \ \ \ \ \ }$\alpha \times _{^{\ast }%
\mathbb{R}
_{\mathbf{d}}}\beta =\beta \times _{^{\ast }%
\mathbb{R}
_{\mathbf{d}}}\alpha ,\left( \alpha \times _{^{\ast }%
\mathbb{R}
_{\mathbf{d}}}\beta \right) \times _{^{\ast }%
\mathbb{R}
_{\mathbf{d}}}\gamma =\alpha \times _{^{\ast }%
\mathbb{R}
_{\mathbf{d}}}\left( \beta \times _{^{\ast }%
\mathbb{R}
_{\mathbf{d}}}\gamma \right) .$

\textbf{(iii) \ }$1_{^{\ast }%
\mathbb{R}
_{\mathbf{d}}}\times _{^{\ast }%
\mathbb{R}
_{\mathbf{d}}}\alpha =\alpha ;$\textbf{\ \ }$-1_{^{\ast }%
\mathbb{R}
_{\mathbf{d}}}\times _{^{\ast }%
\mathbb{R}
_{\mathbf{d}}}\alpha =-_{^{\ast }%
\mathbb{R}
_{\mathbf{d}}}\alpha ,$ where $1_{^{\ast }%
\mathbb{R}
_{\mathbf{d}}}=\left( 1_{^{\ast }%
\mathbb{R}
}\right) ^{\#}.$

\textbf{(iv) \ }$\left\vert \alpha \right\vert \times _{^{\ast }%
\mathbb{R}
_{\mathbf{d}}}\left\vert \beta \right\vert =\left\vert \beta \right\vert
\times _{^{\ast }%
\mathbb{R}
_{\mathbf{d}}}\left\vert \alpha \right\vert .$

\textbf{(v) \ \ }$\left[ \left( \alpha \geq 0\right) \wedge \left( \beta
\geq 0\right) \wedge \left( \gamma \geq 0\right) \right] \implies $

$\ \ \ \ \ \ \ \implies \alpha \times _{^{\ast }%
\mathbb{R}
_{\mathbf{d}}}\left( \beta +_{^{\ast }%
\mathbb{R}
_{\mathbf{d}}}\gamma \right) =\alpha \times _{^{\ast }%
\mathbb{R}
_{\mathbf{d}}}\beta +_{^{\ast }%
\mathbb{R}
_{\mathbf{d}}}\alpha \times _{^{\ast }%
\mathbb{R}
_{\mathbf{d}}}\gamma .$

\textbf{(vi) \ }$0_{^{\ast }%
\mathbb{R}
_{\mathbf{d}}}<_{^{\ast }%
\mathbb{R}
_{\mathbf{d}}}\alpha <_{^{\ast }%
\mathbb{R}
_{\mathbf{d}}}\alpha ^{\prime },0_{^{\ast }%
\mathbb{R}
_{\mathbf{d}}}<_{^{\ast }%
\mathbb{R}
_{\mathbf{d}}}\beta <_{^{\ast }%
\mathbb{R}
_{\mathbf{d}}}\beta ^{\prime }\implies $

$\ \ \ \ \ \ \ \implies \alpha \times _{^{\ast }%
\mathbb{R}
_{\mathbf{d}}}\beta <_{^{\ast }%
\mathbb{R}
_{\mathbf{d}}}\alpha ^{\prime }\times _{^{\ast }%
\mathbb{R}
_{\mathbf{d}}}\beta ^{\prime }.$

\textbf{Proof.(v) }Clearly $\alpha \times _{^{\ast }%
\mathbb{R}
_{\mathbf{d}}}\left( \beta +_{^{\ast }%
\mathbb{R}
_{\mathbf{d}}}\gamma \right) \leq \alpha \times _{^{\ast }%
\mathbb{R}
_{\mathbf{d}}}\beta +_{^{\ast }%
\mathbb{R}
_{\mathbf{d}}}\alpha \times _{^{\ast }%
\mathbb{R}
_{\mathbf{d}}}\gamma .$

Suppose $d\in \alpha \times _{^{\ast }%
\mathbb{R}
_{\mathbf{d}}}\beta +_{^{\ast }%
\mathbb{R}
_{\mathbf{d}}}\alpha \times _{^{\ast }%
\mathbb{R}
_{\mathbf{d}}}\gamma .$Hence:

$d=ab+a^{\prime }c,$where $a,a^{\prime }\in \alpha ,b\in \beta ,c\in \gamma
. $

Without loss of generality we may assume $a\leq a^{\prime }.$Hence:

\bigskip $d=ab+a^{\prime }c\leq a^{\prime }b+a^{\prime }c=a^{\prime }\left(
b+c\right) \in \alpha \times _{^{\ast }%
\mathbb{R}
_{\mathbf{d}}}\left( \beta +_{^{\ast }%
\mathbb{R}
_{\mathbf{d}}}\gamma \right) .$

\textbf{Definition 1.3.1.6.} Suppose $\alpha \in $ $^{\ast }%
\mathbb{R}
_{\mathbf{d}},0<_{^{\ast }%
\mathbb{R}
_{\mathbf{d}}}\alpha $ then $\alpha ^{-1_{^{\ast }%
\mathbb{R}
_{\mathbf{d}}}}$ is defined

as follows:

\textbf{(i) \ }$0_{^{\ast }%
\mathbb{R}
_{\mathbf{d}}}<_{^{\ast }%
\mathbb{R}
_{\mathbf{d}}}\alpha :\alpha ^{-1_{^{\ast }%
\mathbb{R}
_{\mathbf{d}}}}\triangleq \inf \left\{ a^{-1_{^{\ast }%
\mathbb{R}
}}|a\in \alpha \right\} ,$

\textbf{(ii) }$_{^{\ast }%
\mathbb{R}
_{\mathbf{d}}}\alpha <_{^{\ast }%
\mathbb{R}
_{\mathbf{d}}}0:\alpha ^{-1_{^{\ast }%
\mathbb{R}
_{\mathbf{d}}}}\triangleq -_{^{\ast }%
\mathbb{R}
_{\mathbf{d}}}\left( -_{^{\ast }%
\mathbb{R}
_{\mathbf{d}}}\alpha \right) ^{-1_{^{\ast }%
\mathbb{R}
_{\mathbf{d}}}}.$\bigskip

\textbf{Lemma 1.3.1.5.}[24].\textbf{\ }

\textbf{(i) \ }$\forall a\in $ $^{\ast }%
\mathbb{R}
:\left( a^{\#}\right) ^{-1_{^{\ast }%
\mathbb{R}
_{\mathbf{d}}}}=\left( a^{-1_{^{\ast }%
\mathbb{R}
}}\right) ^{\#}.$

\textbf{(ii) \ }$\left( \alpha ^{-1_{^{\ast }%
\mathbb{R}
}}\right) ^{^{-1_{^{\ast }%
\mathbb{R}
}}}=\alpha .$

\textbf{(iii) }$0_{^{\ast }%
\mathbb{R}
_{\mathbf{d}}}<_{^{\ast }%
\mathbb{R}
_{\mathbf{d}}}\alpha \leq _{^{\ast }%
\mathbb{R}
_{\mathbf{d}}}\beta \implies \beta ^{-1_{^{\ast }%
\mathbb{R}
_{\mathbf{d}}}}\leq _{^{\ast }%
\mathbb{R}
_{\mathbf{d}}}\alpha ^{-1_{^{\ast }%
\mathbb{R}
_{\mathbf{d}}}}.$

\textbf{(iv) }$\left[ \left( 0_{^{\ast }%
\mathbb{R}
_{\mathbf{d}}}<_{^{\ast }%
\mathbb{R}
_{\mathbf{d}}}\alpha \right) \wedge \left( 0_{^{\ast }%
\mathbb{R}
_{\mathbf{d}}}<_{^{\ast }%
\mathbb{R}
_{\mathbf{d}}}\beta \right) \right] \implies $

$\ \ \ \ \ \implies \left( \alpha ^{-1_{^{\ast }%
\mathbb{R}
}}\right) \times _{^{\ast }%
\mathbb{R}
_{\mathbf{d}}}\left( \beta ^{-1_{^{\ast }%
\mathbb{R}
}}\right) \leq _{^{\ast }%
\mathbb{R}
_{\mathbf{d}}}\left( \alpha \times _{^{\ast }%
\mathbb{R}
_{\mathbf{d}}}\beta \right) ^{-1_{^{\ast }%
\mathbb{R}
}}$

\textbf{(v)} \ $\forall a\in $ $^{\ast }%
\mathbb{R}
:a\neq _{^{\ast }%
\mathbb{R}
}0_{^{\ast }%
\mathbb{R}
}\implies \left( \alpha ^{\#}\right) ^{-1_{^{\ast }%
\mathbb{R}
_{\mathbf{d}}}}\times _{^{\ast }%
\mathbb{R}
_{\mathbf{d}}}\left( \beta ^{-1_{^{\ast }%
\mathbb{R}
_{\mathbf{d}}}}\right) =\left( \alpha ^{\#}\times _{^{\ast }%
\mathbb{R}
_{\mathbf{d}}}\beta \right) ^{-1_{^{\ast }%
\mathbb{R}
_{\mathbf{d}}}}.$

\textbf{(vi) }$\alpha \times _{^{\ast }%
\mathbb{R}
_{\mathbf{d}}}\alpha ^{-1_{^{\ast }%
\mathbb{R}
_{\mathbf{d}}}}\leq _{^{\ast }%
\mathbb{R}
_{\mathbf{d}}}1_{^{\ast }%
\mathbb{R}
_{\mathbf{d}}}.$

\bigskip $\mathbf{Lemma}$ $\mathbf{1.3.1.5}^{\ast }$\textbf{.}Suppose that $%
a\in $ $^{\ast }%
\mathbb{R}
,a>0,\beta ,\gamma \in $ $^{\ast }%
\mathbb{R}
_{\mathbf{d}}.$ Then

$a^{\#}\times _{^{\ast }%
\mathbb{R}
_{\mathbf{d}}}\left( \beta +_{^{\ast }%
\mathbb{R}
_{\mathbf{d}}}\gamma \right) =a^{\#}\times _{^{\ast }%
\mathbb{R}
_{\mathbf{d}}}\beta +_{^{\ast }%
\mathbb{R}
_{\mathbf{d}}}a^{\#}\times _{^{\ast }%
\mathbb{R}
_{\mathbf{d}}}\gamma .$

\textbf{Proof. }Clearly $a^{\#}\times _{^{\ast }%
\mathbb{R}
_{\mathbf{d}}}\left( \beta +_{^{\ast }%
\mathbb{R}
_{\mathbf{d}}}\gamma \right) \leq a^{\#}\times _{^{\ast }%
\mathbb{R}
_{\mathbf{d}}}\beta +_{^{\ast }%
\mathbb{R}
_{\mathbf{d}}}a^{\#}\times _{^{\ast }%
\mathbb{R}
_{\mathbf{d}}}\gamma .$

$\left( a^{\#}\right) ^{^{-1_{^{\ast }%
\mathbb{R}
}}}\left( a^{\#}\times _{^{\ast }%
\mathbb{R}
_{\mathbf{d}}}\beta +_{^{\ast }%
\mathbb{R}
_{\mathbf{d}}}a^{\#}\times _{^{\ast }%
\mathbb{R}
_{\mathbf{d}}}\gamma \right) \leq $

$\leq \left( a^{\#}\right) ^{^{-1_{^{\ast }%
\mathbb{R}
}}}\left( a^{\#}\times _{^{\ast }%
\mathbb{R}
_{\mathbf{d}}}\beta \right) +_{^{\ast }%
\mathbb{R}
_{\mathbf{d}}}\left( a^{\#}\right) ^{^{-1_{^{\ast }%
\mathbb{R}
}}}\left( a^{\#}\times _{^{\ast }%
\mathbb{R}
_{\mathbf{d}}}\gamma \right) =$

$=\beta +_{^{\ast }%
\mathbb{R}
_{\mathbf{d}}}\gamma .$Thus $\left( a^{\#}\right) ^{^{-1_{^{\ast }%
\mathbb{R}
}}}\left( a^{\#}\times _{^{\ast }%
\mathbb{R}
_{\mathbf{d}}}\beta +_{^{\ast }%
\mathbb{R}
_{\mathbf{d}}}a^{\#}\times _{^{\ast }%
\mathbb{R}
_{\mathbf{d}}}\gamma \right) \leq \beta +_{^{\ast }%
\mathbb{R}
_{\mathbf{d}}}\gamma $ and

one obtain $a^{\#}\times _{^{\ast }%
\mathbb{R}
_{\mathbf{d}}}\beta +_{^{\ast }%
\mathbb{R}
_{\mathbf{d}}}a^{\#}\times _{^{\ast }%
\mathbb{R}
_{\mathbf{d}}}\gamma \leq a^{\#}\times _{^{\ast }%
\mathbb{R}
_{\mathbf{d}}}\left( \beta +_{^{\ast }%
\mathbb{R}
_{\mathbf{d}}}\gamma \right) .$

$\mathbf{Lemma}$ $\mathbf{1.3.1.6.}$\textbf{\ }($\mathbf{General}$\textbf{\ }%
$\mathbf{Strong}$\textbf{\ }$\mathbf{Approximation}$\textbf{\ }$\mathbf{%
Property}$).

If $A$ is a nonempty subset of $^{\ast }%
\mathbb{R}
_{\mathbf{d}}$ which is bounded from above, then

$\sup (A)$ is the unique number such that:

(\textbf{i}) $\sup (A)$ is an upper bound for $A$ and

(\textbf{ii}) for any $\alpha \in \sup (A)$ there exists $x\in A$ such that $%
\alpha <x\leq \sup (A).$

\textbf{Proof. }If not, then $\alpha $ is an upper bound of $A$ less than
the least upper

bound $\sup (A)$, which is a contradiction.

\textbf{Lemma 1.3.1.7.}Let $\mathbf{A}$ and $\mathbf{B}$ be nonempty subsets
of $^{\ast }%
\mathbb{R}
\subset $ $^{\ast }%
\mathbb{R}
_{\mathbf{d}}$ and

$\mathbf{C}=$ $\left\{ a+b:a\in \mathbf{A},b\in \mathbf{B}\right\} $.If $%
\mathbf{A}$ and $\mathbf{B}$ are bounded or hyperbounded

from above,hence $\sup \left( \mathbf{A}\right) $ and $\sup \left( \mathbf{B}%
\right) $ exist, then $\mathbf{s}$-$\sup \left( \mathbf{C}\right) $ exist and

\bigskip

$\ \ \ \ 
\begin{array}{cc}
\begin{array}{c}
\\ 
\sup \left( \mathbf{C}\right) =\sup \left( \mathbf{A}\right) +\sup \left( 
\mathbf{B}\right) . \\ 
\end{array}
& \text{ }\left( 1.3.11\right) \text{\ }%
\end{array}%
$

\bigskip

\textbf{Proof.}Suppose $c<\sup \left( \mathbf{A}\right) +\sup \left( \mathbf{%
B}\right) .$From \textbf{Lemma 1.3.1.2.}(\textbf{iv}) $^{\ast }%
\mathbb{R}
$ is

dense in $^{\ast }%
\mathbb{R}
_{\mathbf{d}}.$So there is exists $x\in $ $^{\ast }%
\mathbb{R}
$ such that $c<x^{\#}<\sup \left( \mathbf{A}\right) +\sup \left( \mathbf{B}%
\right) .$

Suppose that $\alpha ,\beta \in $ $^{\ast }%
\mathbb{R}
$ and $\alpha ^{\#}<\sup \left( \mathbf{A}\right) ,\beta ^{\#}<\sup \left( 
\mathbf{B}\right) .$

From Lemma 1.3.1.6\textbf{\ }(\textbf{General Strong Approximation Property}%
) one obtain

there is exists $a\in \mathbf{A},b\in \mathbf{B}$ such that $\alpha
^{\#}<a<\sup \left( \mathbf{A}\right) ,\beta ^{\#}<b<\sup \left( \mathbf{B}%
\right) .$

Suppose $x^{\#}<\alpha ^{\#}+\beta ^{\#}.$Thus one obtain:

\bigskip

$\ \ 
\begin{array}{cc}
\begin{array}{c}
\\ 
\alpha ^{\#}-\dfrac{\left( \alpha ^{\#}+\beta ^{\#}\right) -x^{\#}}{2}%
<\alpha ^{\#}<a<\sup \left( \mathbf{A}\right) \\ 
\end{array}
& \text{\ }\left( 1.3.12\right) \text{\ }%
\end{array}%
$

\bigskip

\bigskip and

$\ \ \ \ 
\begin{array}{cc}
\begin{array}{c}
\\ 
\beta ^{\#}-\dfrac{\left( \alpha ^{\#}+\beta ^{\#}\right) -x^{\#}}{2}<\beta
^{\#}<b<\sup \left( \mathbf{B}\right) .\bigskip \\ 
\end{array}
& \text{ \ \ }\left( 1.3.13\right) \text{\ }%
\end{array}%
$

So one obtain

$\bigskip $

$\ 
\begin{array}{cc}
\begin{array}{c}
\\ 
x^{\#}=\left[ \left( \alpha ^{\#}-\dfrac{\left( \alpha ^{\#}+\beta
^{\#}\right) -x^{\#}}{2}\right) +\left( \beta ^{\#}-\dfrac{\left( \alpha
^{\#}+\beta ^{\#}\right) -x^{\#}}{2}\right) \right] \\ 
\\ 
<\alpha ^{\#}+\beta ^{\#}<a+b<\sup \left( \mathbf{A}\right) +\sup \left( 
\mathbf{B}\right) . \\ 
\end{array}
& \text{ \ }\left( 1.3.14\right) \text{\ \ \ \ \ \ \ \ }%
\end{array}%
$

\bigskip

But $a+b\in \mathbf{C,}$hence by using \textbf{Lemma 1.3.1.4 }one obtain that

$\sup \left( \mathbf{C}\right) =\sup \left( \mathbf{A}\right) +\sup \left( 
\mathbf{B}\right) .$

\textbf{Theorem 1.3.1.2. }Let $\mathbf{A}$ and $\mathbf{B}$ be nonempty
subsets of $^{\ast }%
\mathbb{R}
_{\mathbf{d}}$ and

$\mathbf{C}=$ $\left\{ a+b:a\in \mathbf{A},b\in \mathbf{B}\right\} $.If $%
\mathbf{A}$ and $\mathbf{B}$ are bounded or hyperbounded

from above,hence $\sup \left( \mathbf{A}\right) $ and $\sup \left( \mathbf{B}%
\right) $ exist, then $\mathbf{s}$-$\sup \left( \mathbf{C}\right) $ exist
and $\ \ \ \ \ \ \ \ \ \ \ \ \ \ \ \ \ \ \ \ \ \ \ \ \ \ \ \ \ \ \ \ 
\begin{array}{cc}
\begin{array}{c}
\\ 
\sup \left( \mathbf{C}\right) =\sup \left( \mathbf{A}\right) +\sup \left( 
\mathbf{B}\right) . \\ 
\end{array}
& \text{ \ \ \ \ \ }\left( 1.3.3.1\right) \text{\ \ \ \ }%
\end{array}%
$

\textbf{Proof.}Suppose $c<\sup \left( \mathbf{A}\right) +\sup \left( \mathbf{%
B}\right) .$From \textbf{Lemma 1.3.1.2.}(\textbf{iv}) $^{\ast }%
\mathbb{R}
$ is

dense in $^{\ast }%
\mathbb{R}
_{\mathbf{d}}.$So there is exists $x\in $ $^{\ast }%
\mathbb{R}
$ such that $c<x^{\#}<\sup \left( \mathbf{A}\right) +\sup \left( \mathbf{B}%
\right) .$

Suppose that $\alpha ,\beta \in $ $^{\ast }%
\mathbb{R}
$ and $\alpha ^{\#}<\sup \left( \mathbf{A}\right) ,\beta ^{\#}<\sup \left( 
\mathbf{B}\right) .$From \textbf{Lemma 1.3.1.4}

(\textbf{General Strong Approximation Property})one obtain there is exists

$a\in \mathbf{A},b\in \mathbf{B}$ such that $\alpha ^{\#}<a<\sup \left( 
\mathbf{A}\right) ,\beta ^{\#}<b<\sup \left( \mathbf{B}\right) .$ Suppose

$x^{\#}<\alpha ^{\#}+\beta ^{\#}.$Thus one obtain:

$\ \ \ \ \ \ \ \ \ \ \ \ \ \ \ \ \ \ \ \ \ \ \ \ \ \ \ \ \ \ \ \ 
\begin{array}{cc}
\begin{array}{c}
\\ 
\alpha ^{\#}-\dfrac{\left( \alpha ^{\#}+\beta ^{\#}\right) -x^{\#}}{2}%
<\alpha ^{\#}<a<\sup \left( \mathbf{A}\right) \\ 
\end{array}
& \text{ \ \ \ \ \ \ \ }%
\end{array}%
$

and

$\ \ \ \ \ \ \ \ \ \ \ \ \ \ \ \ \ \ \ \ \ \ \ \ \ \ \ \ \ \ \ 
\begin{array}{cc}
\begin{array}{c}
\\ 
\beta ^{\#}-\dfrac{\left( \alpha ^{\#}+\beta ^{\#}\right) -x^{\#}}{2}<\beta
^{\#}<b<\sup \left( \mathbf{B}\right) .\bigskip \\ 
\end{array}
& \text{ \ \ \ \ \ \ \ }%
\end{array}%
$

So one obtain

$\ \ \ \ \ \ \ \ \ \ \ \ \ 
\begin{array}{cc}
\begin{array}{c}
\\ 
x^{\#}=\left[ \left( \alpha ^{\#}-\dfrac{\left( \alpha ^{\#}+\beta
^{\#}\right) -x^{\#}}{2}\right) +\left( \beta ^{\#}-\dfrac{\left( \alpha
^{\#}+\beta ^{\#}\right) -x^{\#}}{2}\right) \right] < \\ 
\\ 
<\alpha ^{\#}+\beta ^{\#}< \\ 
\\ 
<a+b<\sup \left( \mathbf{A}\right) +\sup \left( \mathbf{B}\right) . \\ 
\end{array}
& \text{ \ \ \ \ \ \ \ }%
\end{array}%
$

But $a+b\in \mathbf{C,}$hence by using \textbf{Lemma 1.3.1.4 }one obtain that

$\sup \left( \mathbf{C}\right) =\sup \left( \mathbf{A}\right) +\sup \left( 
\mathbf{B}\right) .$

\textbf{Theorem 1.3.1.3.}Suppose that $\mathbf{S}$ is a non-empty subset of $%
^{\ast }%
\mathbb{R}
_{\mathbf{d}}$ which is

bounded or hyperbounded from above,i.e. $\sup \left( \mathbf{S}\right) $
exist and suppose that

$\xi \in $ $^{\ast }%
\mathbb{R}
,\xi >0.$

Then

$\ \ \ \ \ \ \ \ \ \ \ \ \ \ \ \ \ \ \ \ \ \ \ \ 
\begin{array}{cc}
\begin{array}{c}
\\ 
\underset{x\in \mathbf{S}}{\sup }\left\{ \xi ^{\#}\times x\right\} =\xi
^{\#}\times \left( \underset{x\in \mathbf{S}}{\sup }\left\{ x\right\}
\right) =\xi ^{\#}\times \left( \sup \mathbf{S}\right) \mathbf{.} \\ 
\end{array}
& \text{ \ \ \ \ \ \ \ \ \ \ \ \ }\left( 1.3.3.2\right) \text{\ \ \ \ \ \ \ }%
\end{array}%
$

\textbf{Proof.}Let $B=\mathbf{s}$-$\sup \mathbf{S.}$Then $B$ is the smallest
number such that, for \ \ \ \ \ \ \ \ \ \ \ \ \ \ \ \ \ \ \ \ \ \ \ \ 

any $x\in \mathbf{S,}x$ $\mathbf{\leq B.}$Let $\mathbf{T}=\left\{ \xi
^{\#}\times x|x\in \mathbf{S}\right\} .$Since $\xi ^{\#}>0,\xi ^{\#}\times
x\leq \xi ^{\#}\times B$ for any

$x\in \mathbf{S.}$Hence $\mathbf{T}$ is bounded or hyperbounded above by $%
\xi ^{\#}\times B.$Hence \ \ \ \ \ \ \ \ \ \ \ \ \ \ \ \ \ \ \ \ \ \ \ \ 

$\mathbf{T}$ has a supremum $C_{\mathbf{T}}=\mathbf{s}$-$\sup \mathbf{T.}$
Now\ we have to pruve that $C_{\mathbf{T}}=\xi ^{\#}\times B=$

$=\xi ^{\#}\times \left( \sup \mathbf{S}\right) .$Since $\xi ^{\#}\times
B=\xi ^{\#}\times \left( \sup \mathbf{S}\right) $ is an apper bound for $%
\mathbf{T}$and $C$ is the

smollest apper bound for $\mathbf{T,}C_{\mathbf{T}}\leq \xi ^{\#}\times B.$%
Now we repeat the argument above

with the roles of $\mathbf{S}$ and $\mathbf{T}$ reversed. We know that $C_{%
\mathbf{T}}$ is the smallest number

such that, for any $y\in \mathbf{T,}y\leq C_{\mathbf{T}}.$Since $\xi >0$ it
follows that

$\left( \xi ^{\#}\right) ^{-1}\times y\leq \left( \xi ^{\#}\right)
^{-1}\times C_{\mathbf{T}}$ for any $y\in \mathbf{T.}$But $\mathbf{S=}%
\left\{ \left( \xi ^{\#}\right) ^{-1}\times y|y\in \mathbf{T}\right\} .$Hence

$\left( \xi ^{\#}\right) ^{-1}\times C_{\mathbf{T}}$ is an apper bound for $%
\mathbf{S.}$But $B$ is a supremum for $\mathbf{S.}$Hence

$B\leq \left( \xi ^{\#}\right) ^{-1}\times C_{\mathbf{T}}$ and $\xi
^{\#}\times B\leq C_{\mathbf{T}}.$We have shown that $C_{\mathbf{T}}\leq \xi
^{\#}\times B$ and also

that $\xi ^{\#}\times B\leq C_{\mathbf{T}}.$Thus $\xi ^{\#}\times B=C_{%
\mathbf{T}}.$\ \ 

\bigskip

\bigskip \textbf{Theorem 1.3.1.4. }Suppose that $\alpha \in $ $^{\ast }%
\mathbb{R}
$\textbf{\ }$\alpha >0,\beta \in $ $^{\ast }%
\mathbb{R}
_{\mathbf{d}},\gamma \in $ $^{\ast }%
\mathbb{R}
_{\mathbf{d}}.$Then

$\ \ \ \ \ \ \ \ \ \ \ \ \ \ \ \ \ \ \ \ \ \ \ \ \ \ 
\begin{array}{cc}
\begin{array}{c}
\\ 
\ \ \alpha ^{\#}\times _{^{\ast }%
\mathbb{R}
_{\mathbf{d}}}\left( \beta +_{^{\ast }%
\mathbb{R}
_{\mathbf{d}}}\gamma \right) =\alpha ^{\#}\times _{^{\ast }%
\mathbb{R}
_{\mathbf{d}}}\beta +_{^{\ast }%
\mathbb{R}
_{\mathbf{d}}}\alpha ^{\#}\times _{^{\ast }%
\mathbb{R}
_{\mathbf{d}}}\gamma . \\ 
\end{array}
& \text{ \ \ \ \ \ \ \ \ \ \ \ \ \ \ \ \ \ \ }\left( 1.3.3.3\right)%
\end{array}%
$

$\ \ \ \ \ \ \ \ $

\textbf{Proof.}Let us consider any two sets $S_{\beta }\subset $ $^{\ast }%
\mathbb{R}
$ and $S_{\gamma }\subset $ $^{\ast }%
\mathbb{R}
$ such that:

$\beta =\sup \left( S_{\beta }\right) ,\gamma =\sup \left( S_{\gamma
}\right) .$Thus by using \textbf{Theorem 1.3.1.3 }and$\ $

\textbf{Theorem 1.3.1.2 }one obtain:

$\bigskip $ \ \ \ \ \ \ \ \ \ \ \ \ \ \ \ \ \ $%
\begin{array}{cc}
\begin{array}{c}
\\ 
\alpha ^{\#}\times _{^{\ast }%
\mathbb{R}
_{\mathbf{d}}}\left( \beta +_{^{\ast }%
\mathbb{R}
_{\mathbf{d}}}\gamma \right) =\alpha ^{\#}\times _{^{\ast }%
\mathbb{R}
_{\mathbf{d}}}\sup \left( S_{\beta }+S_{\gamma }\right) = \\ 
\\ 
=\sup \left[ \alpha ^{\#}\times _{^{\ast }%
\mathbb{R}
_{\mathbf{d}}}\left( S_{\beta }+S_{\gamma }\right) \right] =\sup \left[
\alpha ^{\#}\times _{^{\ast }%
\mathbb{R}
_{\mathbf{d}}}S_{\beta }+\alpha ^{\#}\times _{^{\ast }%
\mathbb{R}
_{\mathbf{d}}}S_{\gamma }\right] = \\ 
\\ 
=\sup \left( \alpha ^{\#}\times _{^{\ast }%
\mathbb{R}
_{\mathbf{d}}}S_{\beta }\right) +\sup \left( \alpha ^{\#}\times _{^{\ast }%
\mathbb{R}
_{\mathbf{d}}}S_{\gamma }\right) = \\ 
\\ 
\alpha ^{\#}\times _{^{\ast }%
\mathbb{R}
_{\mathbf{d}}}\sup \left( S_{\beta }\right) +\alpha ^{\#}\times _{^{\ast }%
\mathbb{R}
_{\mathbf{d}}}\sup \left( S_{\gamma }\right) . \\ 
\end{array}
& 
\end{array}%
$

$\bigskip $\textbf{Theorem 1.3.1.5. }Suppose that $\alpha \in $ $^{\ast }%
\mathbb{R}
,\alpha <0,\beta \in $ $^{\ast }%
\mathbb{R}
,\gamma \in $ $^{\ast }%
\mathbb{R}
_{\mathbf{d}}.$Then $\ \ \ \ \ \ \ \ \ \ \ \ \ \ \ \ \ \ \ \ $

$\ \ \ \ \ \ \ \ \ \ \ \ \ 
\begin{array}{cc}
\begin{array}{c}
\\ 
\ \ \alpha ^{\#}\times _{^{\ast }%
\mathbb{R}
_{\mathbf{d}}}\left( \beta ^{\#}+_{^{\ast }%
\mathbb{R}
_{\mathbf{d}}}\gamma \right) =\left( -1_{^{\ast }%
\mathbb{R}
_{\mathbf{d}}}\right) \times _{^{\ast }%
\mathbb{R}
_{\mathbf{d}}}\left[ \left\vert \alpha ^{\#}\right\vert \times _{^{\ast }%
\mathbb{R}
_{\mathbf{d}}}\beta ^{\#}+_{^{\ast }%
\mathbb{R}
_{\mathbf{d}}}\left\vert \alpha ^{\#}\right\vert \times _{^{\ast }%
\mathbb{R}
_{\mathbf{d}}}\gamma \right] . \\ 
\end{array}
& \text{ \ }\left( 1.3.3.4\right)%
\end{array}%
$

\textbf{Proof.}Let us consider any set $S_{\gamma }\subset $ $^{\ast }%
\mathbb{R}
$ such that:$\gamma =\sup \left( S_{\gamma }\right) .$Thus by

using \textbf{Theorem 1.3.1.3, Theorem 1.3.1.2 }and$\ $\textbf{Lemma 1.3.1.3
(v) }one obtain:\bigskip

\bigskip\ \ \ \ \ \ \ \ \ \ \ \ \ \ \ \ \ \ $%
\begin{array}{cc}
\begin{array}{c}
\\ 
\alpha ^{\#}\times _{^{\ast }%
\mathbb{R}
_{\mathbf{d}}}\left( \beta ^{\#}+_{^{\ast }%
\mathbb{R}
_{\mathbf{d}}}\gamma \right) =\left\vert \alpha ^{\#}\right\vert \times
_{^{\ast }%
\mathbb{R}
_{\mathbf{d}}}\left( -1_{^{\ast }%
\mathbb{R}
_{\mathbf{d}}}\right) \times _{^{\ast }%
\mathbb{R}
_{\mathbf{d}}}\left( \beta ^{\#}+_{^{\ast }%
\mathbb{R}
_{\mathbf{d}}}\gamma \right) = \\ 
\\ 
=\left\vert \alpha ^{\#}\right\vert \times _{^{\ast }%
\mathbb{R}
_{\mathbf{d}}}\left[ \left( -_{^{\ast }%
\mathbb{R}
_{\mathbf{d}}}\beta ^{\#}\right) +_{^{\ast }%
\mathbb{R}
_{\mathbf{d}}}\left( -_{^{\ast }%
\mathbb{R}
_{\mathbf{d}}}\gamma \right) \right] = \\ 
\\ 
=\left\vert \alpha ^{\#}\right\vert \times _{^{\ast }%
\mathbb{R}
_{\mathbf{d}}}\left( -_{^{\ast }%
\mathbb{R}
_{\mathbf{d}}}\beta ^{\#}\right) +_{^{\ast }%
\mathbb{R}
_{\mathbf{d}}}\left\vert \alpha ^{\#}\right\vert \times _{^{\ast }%
\mathbb{R}
_{\mathbf{d}}}\left( -_{^{\ast }%
\mathbb{R}
_{\mathbf{d}}}\gamma \right) = \\ 
\\ 
=\left\vert \alpha ^{\#}\right\vert \times _{^{\ast }%
\mathbb{R}
_{\mathbf{d}}}\left( -1_{^{\ast }%
\mathbb{R}
_{\mathbf{d}}}\right) \times _{^{\ast }%
\mathbb{R}
_{\mathbf{d}}}\beta ^{\#}+_{^{\ast }%
\mathbb{R}
_{\mathbf{d}}}\left\vert \alpha ^{\#}\right\vert \times _{^{\ast }%
\mathbb{R}
_{\mathbf{d}}}\left( -1_{^{\ast }%
\mathbb{R}
_{\mathbf{d}}}\right) \times _{^{\ast }%
\mathbb{R}
_{\mathbf{d}}}\gamma = \\ 
\\ 
=\alpha ^{\#}\times _{^{\ast }%
\mathbb{R}
_{\mathbf{d}}}\beta ^{\#}+_{^{\ast }%
\mathbb{R}
_{\mathbf{d}}}\alpha ^{\#}\times _{^{\ast }%
\mathbb{R}
_{\mathbf{d}}}\gamma . \\ 
\end{array}
& 
\end{array}%
$\ \ \ \ \ \ \ \ \ \ \ \ \ \ \ \ \ \ \ \ \ \ \ \ \ \ \ \ \ \ \ \ \ \ \ \ \ \
\ 

\section{I.3.2.The topology of $^{\ast }%
\mathbb{R}
_{\mathbf{d}}.$Wattenberg standard part.}

Fortunately topologically, $^{\ast }%
\mathbb{R}
_{\mathbf{d}}$ has many properties strongly reminiscent

of $%
\mathbb{R}
$ itself. We proceed as follows [24].

\textbf{Definition 1.3.2.1. }

\textbf{(i) \ }$\left( \alpha ,\beta \right) _{^{\ast }%
\mathbb{R}
_{\mathbf{d}}}\triangleq \left\{ u|\alpha <_{^{\ast }%
\mathbb{R}
_{\mathbf{d}}}u<_{^{\ast }%
\mathbb{R}
_{\mathbf{d}}}\beta \right\} ,$

\textbf{(ii) }$\left[ \alpha ,\beta \right] _{^{\ast }%
\mathbb{R}
_{\mathbf{d}}}\triangleq \left\{ u|\alpha \leq _{^{\ast }%
\mathbb{R}
_{\mathbf{d}}}u\leq _{^{\ast }%
\mathbb{R}
_{\mathbf{d}}}\beta \right\} .$

\textbf{Definition 1.3.2.2.}[24].Suppose $U\subseteqq $ $^{\ast }%
\mathbb{R}
_{\mathbf{d}}$. Then $U$ is open if and only

if for every $u\in U,$ $\exists \alpha _{\alpha \in ^{\ast }%
\mathbb{R}
_{\mathbf{d}}}\exists \beta _{\beta \in ^{\ast }%
\mathbb{R}
_{\mathbf{d}}}\left[ \alpha <_{^{\ast }%
\mathbb{R}
_{\mathbf{d}}}u<_{^{\ast }%
\mathbb{R}
_{\mathbf{d}}}\beta \right] $ such that

$u\in \left( \alpha ,\beta \right) _{^{\ast }%
\mathbb{R}
_{\mathbf{d}}}\subseteqq U.$

\textbf{Remark.1.3.2.1.}[24]. Notice this is not equivalent to:

$\forall u_{u\in U}\exists \varepsilon _{\varepsilon >0}\left[ \left(
u-\varepsilon ,u+\varepsilon \right) _{^{\ast }%
\mathbb{R}
_{\mathbf{d}}}\subseteqq U\right] .$

\textbf{Lemma} \textbf{1.3.2.1.}[24].

\textbf{(i) \ \ }$^{\ast }%
\mathbb{R}
$ is dense in $^{\ast }%
\mathbb{R}
_{\mathbf{d}}.$

\textbf{(ii)}\ \ $^{\ast }%
\mathbb{R}
_{\mathbf{d}}\backslash ^{\ast }%
\mathbb{R}
$ is dense in $^{\ast }%
\mathbb{R}
_{\mathbf{d}}.$

\bigskip

\textbf{Lemma} \textbf{1.3.2.2.}[24].\textbf{\ }Suppose $A\subseteqq $ $%
^{\ast }%
\mathbb{R}
_{\mathbf{d}}.$Then $A$ is closed if and

only if:

\textbf{(i) \ }$\forall E\left( E\subseteqq A\right) $\textbf{\ }$E$ bounded
above implies $\sup \left( E\right) \in A,$ and

\bigskip \textbf{(ii) }$\forall E\left( E\subseteqq A\right) $\textbf{\ }$E$
bounded below implies $\inf \left( E\right) \in A.$

\textbf{Proposition} \textbf{1.3.2.1.}[24].

\textbf{(i) \ \ }$^{\ast }%
\mathbb{R}
_{\mathbf{d}}$ is connected.

\textbf{(ii) \ }For $\alpha <_{^{\ast }%
\mathbb{R}
_{\mathbf{d}}}\beta $ in $^{\ast }%
\mathbb{R}
_{\mathbf{d}}$ set $\left[ \alpha ,\beta \right] _{^{\ast }%
\mathbb{R}
_{\mathbf{d}}}$ is compact.

\textbf{(iii) }Suppose $A\subseteqq $ $^{\ast }%
\mathbb{R}
_{\mathbf{d}}.$Then $A$ is compact if and

\ \ \ \ \ only if $A$ is closed and bounded.

\textbf{(iv)} $^{\ast }%
\mathbb{R}
_{\mathbf{d}}$\ is normal.

\textbf{(v) \ \ }The map $\alpha \longmapsto -_{^{\ast }%
\mathbb{R}
_{\mathbf{d}}}\alpha $ is continuous.

\textbf{(vi) \ \ }The map $\alpha \longmapsto \alpha ^{^{-1_{^{\ast }%
\mathbb{R}
_{\mathbf{d}}}}}$ is continuous.

\textbf{(vii)} \ The maps $\left( \alpha ,\beta \right) \longmapsto \left(
\alpha +_{^{\ast }%
\mathbb{R}
_{\mathbf{d}}}\beta \right) $ and $\left( \alpha ,\beta \right) \longmapsto
\left( \alpha \times _{^{\ast }%
\mathbb{R}
_{\mathbf{d}}}\beta \right) $

\ \ \ \ \ \ \ \ are not continuous.

\bigskip \textbf{Definition 1.3.2.3.}[24].(\textbf{Wattenberg Standard Part})

\textbf{(i) \ }Suppose $\alpha \in \left( -\Delta _{\mathbf{d}},\Delta _{%
\mathbf{d}}\right) _{^{\ast }%
\mathbb{R}
_{\mathbf{d}}}.$Then there is a unique

\ \ \ \ \ standard $x\in 
\mathbb{R}
$ called $WST\left( \alpha \right) ,$ such that $x\in \left[ \alpha
-\varepsilon _{\mathbf{d}},\alpha +\varepsilon _{\mathbf{d}}\right] _{^{\ast
}%
\mathbb{R}
_{\mathbf{d}}},$

\textbf{(ii) \ \ }$\alpha \leq _{^{\ast }%
\mathbb{R}
_{\mathbf{d}}}\beta $ implies $WST\left( \alpha \right) \leq WST\left( \beta
\right) ,$

\textbf{(iii) \ }the map $WST\left( \cdot \right) :$ $^{\ast }%
\mathbb{R}
_{\mathbf{d}}\rightarrow 
\mathbb{R}
$ is continuous,

\textbf{(iv) \ }$WST\left( \alpha +_{^{\ast }%
\mathbb{R}
_{\mathbf{d}}}\beta \right) =WST\left( \alpha \right) +WST\left( \beta
\right) ,$

\textbf{(v) \ \ }$WST\left( \alpha \times _{^{\ast }%
\mathbb{R}
_{\mathbf{d}}}\beta \right) =WST\left( \alpha \right) \times WST\left( \beta
\right) ,$

\textbf{(vi) \ }$WST\left( -_{^{\ast }%
\mathbb{R}
_{\mathbf{d}}}\alpha \right) =-WST\left( \alpha \right) ,$

\textbf{(vii) }$WST\left( \alpha ^{-1_{^{\ast }%
\mathbb{R}
_{\mathbf{d}}}}\right) =\left[ WST\left( \alpha \right) \right] ^{-1}$ if $%
\alpha \notin \left[ -\varepsilon _{\mathbf{d}},\varepsilon _{\mathbf{d}}%
\right] _{^{\ast }%
\mathbb{R}
_{\mathbf{d}}}.$

\bigskip

\textbf{Proposition} \textbf{1.3.2.2.}[24].Suppose $f:$ $\left[ a,b\right]
\rightarrow $ $A\subseteqq $ $^{\ast }%
\mathbb{R}
$ is internal,

$\ast $-continuous, and monotonic. Then

(\textbf{1}) $f$ has a unique continuous extension $f^{\#}$ $\left[ a,b%
\right] _{^{\ast }%
\mathbb{R}
_{\mathbf{d}}}$ $\overline{A}$ $\rightarrow $ $^{\ast }%
\mathbb{R}
_{\mathbf{d}}$, where

$\overline{A}$ denotes the closure of $A$ in $^{\ast }%
\mathbb{R}
_{\mathbf{d}}.$

(\textbf{2})The conclusion (\textbf{1}) above holds iff is piecewise
monotonic

(i.e., the domain can be decomposed into a finite (not $\ast $-finite) number

of intervals on each of which $f$ is monotonic).

\textbf{Proposition} \textbf{1.3.2.3.}[24].Suppose $f,g$ are $\ast $%
-continuous, piecewise

monotonic functions then

\bigskip \textbf{(i) }$f\circ g$ is also and

\textbf{(ii) }$\left( f\circ g\right) ^{\#}=\left( f^{\#}\right) \circ
\left( g^{\#}\right) .$

\bigskip

\ \ \ \ \ \ \ \ \ \ \ \ \ \ \ \ \ \ \ \ \ \ \ \ \ \ \ \ \ \ \ \ \ \ \ \ \ \
\ \ \ 

\section{I.3.3.Absorption numbers in $^{\ast }%
\mathbb{R}
_{\mathbf{d}}$ and \ \ \ \ \ \ \ \ \ \ \ \ \ \ \ \ \ \ \ \ \ \ \ \ \ \ \ \
idempotents.}

\bigskip\ \ \ \ \ \ \ \ \ \ \ \ \ \ \ \ \ \ \ \ \ \ \ \ \ \ \ \ \ \ \ \ \ \
\ \ \ \ \ \ \ \ \ \ \ \ 

\section{I.3.3.1.Absorption function and numbers in $^{\ast }%
\mathbb{R}
_{\mathbf{d}}$.}

One of standard ways of defining the completion of $^{\ast }%
\mathbb{R}
$ involves restricting oneself to subsets a which have the following
property $\forall \varepsilon _{\varepsilon >0}\exists x_{x\in \alpha }$ $%
\exists y_{y\in \alpha }\left[ y\text{ }-\text{ }x<\varepsilon \right] $. It
is well known that in this case we obtain a field. In fact the proof is
essentially the same as the one used in the case of ordinary Dedekind cuts
in the development of the standard real numbers, $\varepsilon _{\mathbf{d}},$%
of course, does not have the above property because no infinitesimal
works.This suggests the introduction of the concept of absorption part $%
\mathbf{ab.p.}\left( \alpha \right) $ of a number $\alpha $ for an element $%
\alpha $ of $^{\ast }%
\mathbb{R}
_{\mathbf{d}}$ which, roughly speaking, measures how much a departs from
having the above property [23]. We also introduce similar concept of an
absorption number $\alpha \left( \mathbf{ab.n.}\right) \beta \triangleq 
\mathbf{ab.n.}\left( \alpha ,\beta \right) $ (cut) for given element $\beta $
of $^{\ast }%
\mathbb{R}
_{\mathbf{d}}.$

\bigskip \textbf{Definition 1.3.3.1.1.}[23].$\mathbf{ab.p.}\left( \alpha
\right) \triangleq \left\{ d\geq 0|\forall x_{x\in \alpha }\left[ x+d\in
\alpha \right] \right\} .$

\textbf{Example 1.3.3.1.}(\textbf{i})\textbf{\ }$\forall \alpha \in $ $%
^{\ast }%
\mathbb{R}
:\mathbf{ab.p.}\left( \alpha \right) =0,$

(\textbf{ii}) $\mathbf{ab.p.}\left( \varepsilon _{\mathbf{d}}\right)
=\varepsilon _{\mathbf{d}},$ (\textbf{iii}) $\mathbf{ab.p.}\left(
-\varepsilon _{\mathbf{d}}\right) =\varepsilon _{\mathbf{d}},$

(\textbf{iv}) $\forall \alpha \in $ $^{\ast }%
\mathbb{R}
:\mathbf{ab.p.}\left( \alpha +\varepsilon _{\mathbf{d}}\right) =\varepsilon
_{\mathbf{d}},$

(\textbf{v}) $\ \ \forall \alpha \in $ $^{\ast }%
\mathbb{R}
:\mathbf{ab.p.}\left( \alpha -\varepsilon _{\mathbf{d}}\right) =\varepsilon
_{\mathbf{d}}.$

\textbf{Definition 1.3.3.2. }$\mathbf{ab.n.}\left( \alpha ,\beta \right)
\iff \alpha +\beta =\alpha .$

\textbf{Example 1.3.3.2.}(\textbf{i})$\ \forall \beta \approx 0:$ $\mathbf{%
ab.n.}\left( \varepsilon _{\mathbf{d}},\beta \right) ,$

(\textbf{ii}) $\mathbf{ab.n.}\left( \varepsilon _{\mathbf{d}},\varepsilon _{%
\mathbf{d}}\right) ,\mathbf{ab.n.}\left( -\varepsilon _{\mathbf{d}%
},\varepsilon _{\mathbf{d}}\right) ,\mathbf{ab.n.}\left( -\varepsilon _{%
\mathbf{d}},-\varepsilon _{\mathbf{d}}\right) ,$

(\textbf{iii}) $\forall \alpha \in $ $%
\mathbb{R}
:\mathbf{ab.n.}\left( \alpha +\varepsilon _{\mathbf{d}},\varepsilon _{%
\mathbf{d}}\right) ,\mathbf{ab.n.}\left( \alpha -\varepsilon _{\mathbf{d}%
},\varepsilon _{\mathbf{d}}\right) ,\mathbf{ab.n.}\left( \alpha -\varepsilon
_{\mathbf{d}},-\varepsilon _{\mathbf{d}}\right) ,$

(\textbf{iv}) $\forall \alpha \in $ $%
\mathbb{R}
:\mathbf{ab.n.}\left( \Delta _{\mathbf{d}},\beta \right) ,$

\bigskip (\textbf{v}) $\mathbf{ab.n.}\left( \Delta _{\mathbf{d}},\Delta _{%
\mathbf{d}}\right) ,\mathbf{ab.n.}\left( -\Delta _{\mathbf{d}},\Delta _{%
\mathbf{d}}\right) ,\mathbf{ab.n.}\left( -\Delta _{\mathbf{d}},-\Delta _{%
\mathbf{d}}\right) .$

\textbf{Lemma 1.3.3.1.}[23].(\textbf{i}) $c<\mathbf{ab.p.}\left( \alpha
\right) $ and $0\leq d<c\implies d\in \mathbf{ab.p.}\left( \alpha \right) $

(\textbf{ii}) $c\in \mathbf{ab.p.}\left( \alpha \right) $ and $d\in \mathbf{%
ab.p.}\left( \alpha \right) \implies c+d\in \mathbf{ab.p.}\left( \alpha
\right) .$

\textbf{Remark 1.3.3.1.} By \textbf{Lemma} \textbf{1.3.2.1} $\mathbf{ab.p.}%
\left( \alpha \right) $ may be regarded as an \ \ \ \ \ \ \ \ \ \ \ \ \ \ \
\ \ \ \ \ \ \ \ \ \ \ 

element of $^{\ast }%
\mathbb{R}
_{\mathbf{d}}$ by adding on all negative elements of $^{\ast }%
\mathbb{R}
_{\mathbf{d}}$ to $\mathbf{ab.p.}\left( \alpha \right) .$ \ \ \ \ \ \ \ \ \
\ \ \ \ \ \ \ \ \ \ \ \ \ \ \ \ \ \ \ \ \ 

Of course if the condition $d\geq 0$ in the definition of $\mathbf{ab.p.}%
\left( \alpha \right) $ is deleted we

automatically get all the negative elements to be in $\mathbf{ab.p.}\left(
\alpha \right) $ since

$x<y\in \alpha \implies x\in \alpha .$The reason for our definition is that
the real interest

lies in the non-negative numbers. A technicality occurs if $\mathbf{ab.p.}%
\left( \alpha \right) =\left\{ 0\right\} $.

We then identify $\mathbf{ab.p.}\left( \alpha \right) $ with $0.$ [$\mathbf{%
ab.p.}\left( \alpha \right) $ becomes $\{x|x<0\}$ which by

our early convention is not in $^{\ast }%
\mathbb{R}
_{\mathbf{d}}$].

\bigskip \textbf{Remark 1.3.3.1.2.}By \textbf{Lemma} \textbf{1.3.2.1}(%
\textbf{ii}), $\mathbf{ab.p.}\left( \alpha \right) $ is idempotent.

\textbf{Lemma} \textbf{1.3.3.1.2.}[23].

(\textbf{i}) $\mathbf{ab.p.}(\alpha )$ is the maximum element $\beta \in $ $%
^{\ast }%
\mathbb{R}
_{\mathbf{d}}$ such that $\alpha +\beta =\alpha .$

(\textbf{ii}) $\mathbf{ab.p.}(\alpha )\leq \alpha $ for $\alpha >0.$

(\textbf{iii}) If $\alpha $ is positive and idempotent then $\mathbf{ab.p.}%
(\alpha )=\alpha .$

\textbf{Lemma} \textbf{1.3.3.1.3.}[23]. Let $\alpha \in $ $^{\ast }%
\mathbb{R}
_{\mathbf{d}}$ satsify $\alpha >0.$ Then the following are

equivalent. In what follows assume $a,b>0.$

(\textbf{i}) $\ \ \alpha $ is idempotent,

(\textbf{ii})\ $\ a,b\in \alpha \implies a+b\in \alpha ,$

(\textbf{iii}) $a\in \alpha \implies 2a\in \alpha ,$

(\textbf{iv}) $\forall n_{n\in 
\mathbb{N}
}\left[ a\in \alpha \implies n\cdot a\in \alpha \right] ,$

(\textbf{v})\ \ $a\in \alpha \implies r\cdot a\in \alpha ,$ for all finite $%
r\in $ $^{\ast }%
\mathbb{R}
.$\ \ \ \ \ \ \ 

\bigskip\ \ \ \ \ \ \ \ \ \ \ \ \ \ \ \ \ \ \ \ \ \ \ \ \ \ \ \ \ \ \ \ \ \
\ \ \ \ \ \ \ \ \ \ \ \ 

\section{Connection with the value group.}

\textbf{Definition 1.3.3.1.2. }We define an equivalence relation on the
positive

elements of $^{\ast }%
\mathbb{R}
$ as follows: $a$ $\symbol{126}$ $b\iff \dfrac{a}{b}$ and $\dfrac{b}{a}$ are
finite.Then the

equivalence classes from a linear ordered set. We denote the order

relation by $\ll .$

\bigskip The classes may be regarded as orders of infinity.

The subring of $^{\ast }%
\mathbb{R}
$ consisting of the finite elements is a valuation ring,

and the equivalence classes may also be regarded as elements of the

value group. Condition (\textbf{v}) in Lemma 1.3.3.1.3 essentially says that 
$a$ $\in \alpha $

and $b\symbol{126}a$ $\implies $ $b\in \alpha ,$i.e. a may be regarded as a
Dedekind cut in the value

group.

\bigskip\ \ \ \ \ \ \ \ \ \ \ \ \ \ \ \ \ \ \ \ \ \ \ \ \ \ \ \ \ \ \ \ \ \
\ \ \ \ \ \ \ \ \ \ \ \ 

\section{Properties of the Absorption Function.}

\bigskip

\bigskip \textbf{Theorem} \textbf{1.3.3.1.1.}[23].\textbf{\ }$\left( -\alpha
\right) +\alpha =-\left[ \mathbf{ab.p.}(\alpha )\right] .$

\bigskip \textbf{Theorem} \textbf{1.3.3.1.2.}[23].$\mathbf{ab.p.}(\alpha
+\beta )\geq \mathbf{ab.p.}(\alpha ).$

\textbf{Theorem} \textbf{1.3.3.1.3.}[23].

(\textbf{i}) $\alpha +\beta \leq \alpha +\gamma \implies -\left[ \mathbf{%
ab.p.}(\alpha )\right] +\beta \leq \gamma .$

\bigskip (\textbf{ii}) $\alpha +\beta =\alpha +\gamma \implies -\left[ 
\mathbf{ab.p.}(\alpha )\right] +\beta =\gamma .$

\textbf{Theorem} \textbf{1.3.3.1.4.}[23].

(\textbf{i}) $\ \mathbf{ab.p.}(-\alpha )=\mathbf{ab.p.}(\alpha ),$

(\textbf{ii}) $\mathbf{ab.p.}(\alpha +\beta )=\max \left\{ \mathbf{ab.p.}%
(\alpha ),\mathbf{ab.p.}(\beta )\right\} .$

\bigskip

We now classify the elements $\beta $ such that $\alpha +\beta =\alpha $.
For positive $\beta $ we

know by Lemma \textbf{1.3.3.1.2.}(\textbf{i}) that $\alpha +\beta =\alpha $
iff $\beta \leq \mathbf{ab.p.}(\alpha ).$

\bigskip \textbf{Theorem} \textbf{1.3.3.1.5.}[23]. Assume $\beta $ $>0.$ If $%
\alpha $ absorbs $-\beta $ then a abosrbs $\beta $.

\textbf{Theorem} \textbf{1.3.3.1.6.}[23]. Let $0<\alpha \in $ $^{\ast }%
\mathbb{R}
_{\mathbf{d}}.$ Then the following are

equivalent

(\textbf{i}) $\ \alpha $ is an idempotent,

(\textbf{ii}) $\left( -\alpha \right) +\left( -\alpha \right) =-\alpha ,$

(\textbf{iii}) $\left( -\alpha \right) +\alpha =-\alpha .$

\bigskip

\bigskip\ \ \ \ \ \ \ \ \ \ \ \ \ \ \ \ \ \ \ \ \ \ \ \ \ \ \ \ \ \ \ \ \ \
\ \ \ \ \ \ \ \ \ \ \ \ 

\section{Special Equivalence Relations on $^{\ast }%
\mathbb{R}
_{\mathbf{d}}.$}

Let $\Delta $ be a positive idempotent. We define three equivalence relations

$\left( \circ \text{ }\mathbf{R\circ }\right) \mathbf{,}\left( \circ \text{ }%
\mathbf{S\circ }\right) $ and $\left( \circ \text{ }\mathbf{T\circ }\right) $
on $^{\ast }%
\mathbb{R}
_{\mathbf{d}}.$

\bigskip \textbf{Definition 1.3.3.1.3.}[23].

\bigskip (\textbf{i}) $\alpha \mathbf{R}\beta \left( \func{mod}\Delta
\right) \iff \alpha +\Delta =\beta +\Delta ,$

(\textbf{ii}) $\alpha \mathbf{S}\beta \left( \func{mod}\Delta \right) \iff
\alpha +\left( -\Delta \right) =\beta +\left( -\Delta \right) ,$

\bigskip (\textbf{iii}) $\alpha \mathbf{T}\beta \left( \func{mod}\Delta
\right) \iff \exists d\left( d\in \Delta \right) \left[ \left( \alpha
\subset \beta +d\right) \wedge \left( \beta \subset \alpha +d\right) \right]
.$

\textbf{Remark 1.3.3.1.3.}To simplify the notation $\func{mod}\Delta $ is
omitted when we are

dealing with only one $\Delta $. $\mathbf{R}$ and $\mathbf{S}$ are obviously
equivalence relations.

$\mathbf{T}$ is an equivalence relation since $\Delta $ is idempotent.

\bigskip

\textbf{Remark 1.3.3.1.4.}It is immediate that $\mathbf{R,S}$ and $\mathbf{T}
$ are congruence relations

with respect to addition. Also, if $\symbol{126}$ stands for either $\mathbf{%
R,S}$ or $\mathbf{T}$ then $\alpha <\beta <\gamma $

and $\alpha $ $\symbol{126}$ $\gamma \implies \alpha $ $\symbol{126}$ $\beta
.$ To see this it is convenient to have the following

lemma.

\textbf{Lemma} \textbf{1.3.3.1.4.}[23]. Suppose $\alpha <\beta $. Then

(\textbf{i}) $\ \ \alpha \mathbf{R}\beta \left( \func{mod}\Delta \right)
\iff \beta \leq \alpha +\Delta ,$

\bigskip (\textbf{ii}) \ $\alpha \mathbf{S}\beta \left( \func{mod}\Delta
\right) \iff \beta +\left( -\Delta \right) \leq \alpha .$

\textbf{Lemma} \textbf{1.3.3.1.5.}[23]. Let $\Delta $ be a positive
idempotent. Then

$-\left[ \alpha +\left( -\Delta \right) \right] +\left( -\Delta \right) \leq
-\alpha .$

\textbf{Remark 1.3.3.1.5.}This is not immediate since the inequality

$\left( -\alpha \right) +\left( -\beta \right) $ $\leq -\left( \alpha +\beta
\right) $ goes the wrong way. In fact, this seems

surprising at first since the first addend may be bigger than one

intuitively expects, e.g. if $\alpha =\Delta =\varepsilon _{\mathbf{d}}$
then $-\left[ \alpha +\left( -\Delta \right) \right] =$

$-\left[ \varepsilon _{\mathbf{d}}+\left( -\varepsilon _{\mathbf{d}}\right) %
\right] =\varepsilon _{\mathbf{d}}>0.$ However,$\varepsilon _{\mathbf{d}%
}+\left( -\varepsilon _{\mathbf{d}}\right) =-\varepsilon _{\mathbf{d}},$ so
the

inequality is valid after all.

\bigskip

\textbf{Theorem} \textbf{1.3.3.1.7.}[23].

(\textbf{i}) $\ \ \mathbf{S}$ is a congruence relation with respect to
negation.

(\textbf{ii}) $\ \mathbf{T}$ is a congruence relation with respect to
negation.

\bigskip (\textbf{iii}) $\mathbf{R}$ is not a congruence relation with
respect to negation.

\textbf{Theorem} \textbf{1.3.3.1.8.}[23]. $\alpha +\Delta $ is the maximum
element $\beta $ satisfying

$\beta \mathbf{R\alpha .}$

\textbf{Theorem} \textbf{1.3.3.1.9.}[23].\textbf{\ }$\alpha +\left( -\Delta
\right) $\textbf{\ }is the minimum element $\beta $ satisfying

$\beta \mathbf{S\alpha .}$

\bigskip \textbf{Theorem} \textbf{1.3.3.1.9.}[23].\textbf{\ }$\mathbf{T}%
\subsetneqq \mathbf{R}\subsetneqq \mathbf{S.}$ Both inclusions are proper.

\textbf{Theorem} \textbf{1.3.3.1.10.}[23].

(\textbf{i}) Let $\Delta _{1}$ and $\Delta _{2}$ be two positive idempotents
such that $\Delta _{2}>\Delta _{1}.$ Then:

$\Delta _{2}+\left( -\Delta _{1}\right) =\Delta _{2},$

(\textbf{ii}) Let $\Delta _{1}$ and $\Delta _{2}$ be two positive
idempotents such that $\Delta _{2}>\Delta _{1}.$ Then:

$\alpha \mathbf{S}\beta \left( \func{mod}\Delta _{1}\right) \implies \alpha 
\mathbf{R}\beta \left( \func{mod}\Delta _{2}\right) .$

\textbf{Theorem} \textbf{1.3.3.1.11.}[23].Let $\Delta _{1}$ and $\Delta _{2}$
be two positive idempotents such

that $\Delta _{2}>\Delta _{1}.$Then $\alpha \mathbf{S}\beta \left( \func{mod}%
\Delta _{1}\right) \implies \alpha \mathbf{T}\beta \left( \func{mod}\Delta
_{2}\right) $ but not conversely.

\U{f24c}

\textbf{Theorem} \textbf{1.3.3.1.12.}[23].$\mathbf{S}$ is the smallest
congruence relation with respect

to addition and negation containing $\mathbf{R.}$

\textbf{Theorem} \textbf{1.3.3.1.13.}[23].Any convex congruence relation $%
\left( \circ \text{ }\symbol{126}\circ \right) $ containing

$\mathbf{T}$ properly must contain $\mathbf{S.}$

\ \ \ \ \ \ \ \ \ \ \ \ \ \ \ \ \ \ \ \ \ \ \ \ \ \ \ \ \ \ \ \ \ \ \ 

\section{I.3.3.2.Special Kinds of Idempotents in $^{\ast }%
\mathbb{R}
_{\mathbf{d}}$.}

\bigskip

Let $a\in $ $^{\ast }%
\mathbb{R}
$ such that $a>0.$ Then $a$ gives rise to two idempotents in a

natural way.

\textbf{Definition 1.3.3.2.1.}[23].

(\textbf{i}) $\mathbf{A}_{a}$ $\triangleq $ $\left\{ x|\exists n_{n\in 
\mathbb{N}
}\left[ x<n\cdot a\right] \right\} .$

(\textbf{ii}) $\mathbf{B}_{a}$ $\triangleq $ $\left\{ x|\forall r_{r\in 
\mathbb{R}
_{+}}\left[ x<r\cdot a\right] \right\} .$

Then it is immediate that $\mathbf{A}_{a}$ and $\mathbf{B}_{a}$ are
idempotents.The usual "$\epsilon /2$

argument" shows this for $\mathbf{B}_{a}.$It is also clear that $\mathbf{A}%
_{a}$ is the smallest

idempotent containing $a$ and $\mathbf{B}_{a}$ is the largest idempotent not

containing $a.$It follows that $\mathbf{B}_{a}$ and $\mathbf{A}_{a}$ are
consecutive idempotents.

\textbf{Remark 1.3.3.2.1.}Note that $\mathbf{B}_{1}=\varepsilon _{\mathbf{d}%
}=\inf \left( 
\mathbb{R}
_{+}\right) $ (which is the set of all

infinite small positive numbers plus all negative numbers) which we \ \ \ \
\ \ \ \ \ \ \ \ \ \ \ \ \ \ \ \ \ \ \ \ \ \ \ \ \ \ \ \ \ \ \ \ \ \ \ \ \ \
\ 

have already considered above. \bigskip $\mathbf{A}_{1}=\Delta _{\mathbf{d}%
}\triangleq \sup \left( 
\mathbb{R}
_{+}\right) $ (which is the set

of all finite numbers plus all negative numbers) which we have also

already considered above.

\bigskip

\bigskip \textbf{Definition 1.3.3.2.2. }Let $a\in $ $^{\ast }%
\mathbb{N}
.$

(\textbf{i}) $\mathbf{\omega }_{\mathbf{d}}\left[ a\right] $ $\triangleq $ $%
\left\{ x|\exists n_{n\in 
\mathbb{N}
}\left[ x<n\cdot a\right] \right\} .$

(\textbf{ii}) $\mathbf{\Omega }_{\mathbf{d}}\left[ a\right] $ = $\left\{
x|\forall r_{r\in 
\mathbb{R}
_{+}}\left[ x<r\cdot a\right] \right\} ,a\in $ $^{\ast }%
\mathbb{N}
_{\infty }$

\bigskip

\textbf{Remark 1.3.3.2.2.}Then it is immediate that $\mathbf{\omega }_{%
\mathbf{d}}\left[ a\right] $ and $\mathbf{\Omega }_{\mathbf{d}}\left[ a%
\right] $ are

idempotents. It is also clear that $\mathbf{\omega }_{\mathbf{d}}\left[ a%
\right] $ is the smallest idempotent

containing hypernatural $a$ and $\mathbf{\omega }_{\mathbf{d}}\left[ a\right]
=a\cdot \mathbf{\omega }_{\mathbf{d}}.$ $\mathbf{\Omega }_{\mathbf{d}}\left[
a\right] =a\cdot \varepsilon _{\mathbf{d}}$ is the

largest idempotent not containing $a.$

It follows that $\mathbf{\Omega }_{\mathbf{d}}\left[ a\right] $ and $\mathbf{%
\omega }_{\mathbf{d}}\left[ a\right] $ are consecutive idempotents.

\textbf{Remark 1.3.3.2.3. }Note that $\mathbf{\omega }_{\mathbf{d}}\left[ 1%
\right] =\mathbf{\omega }_{\mathbf{d}}$ (which is the set of all finite

natural numbers $%
\mathbb{N}
$ plus all negative numbers) which we have also

already considered above.

\bigskip

\textbf{Theorem} \textbf{1.3.3.2.1.}[23].

(\textbf{i}) No idempotent of the form $\mathbf{A}_{a}$ has an immediate
successor.

(\textbf{ii}) All consecutive pairs of idempotents have the form $\mathbf{A}%
_{a}$ and $\mathbf{B}_{a}$

for some $a\in $ $^{\ast }%
\mathbb{R}
.$

\bigskip\ 

\ \ \ \ \ \ \ \ \ \ \ \ \ \ \ \ \ \ \ \ \ \ \ \ \ \ \ \ \ \ \ \ \ \ \ 

\section{I.3.3.3. Types of $\protect\alpha $ with a given $\mathbf{ab.p.}(%
\protect\alpha ).$}

\bigskip

Among elements of $\alpha \in $ $^{\ast }%
\mathbb{R}
_{\mathbf{d}}$ such that $\mathbf{ab.p.}(\alpha )=\Delta $ we can distinguish

two types.

\bigskip

\textbf{Definition 1.3.3.3.1.}[23].\textbf{\ }Assume $\Delta >0.$

(\textbf{i}) $\ \alpha \in $ $^{\ast }%
\mathbb{R}
_{\mathbf{d}}$ has type $1$ if $\exists x\left( x\in \alpha \right) \forall y%
\left[ x+y\in \alpha \implies y\in \Delta \right] ,$

(\textbf{ii}) $\alpha \in $ $^{\ast }%
\mathbb{R}
_{\mathbf{d}}$ has type $2$ if $\forall x\left( x\in \alpha \right) \exists
y\left( y\notin \Delta \right) \left[ x+y\in \alpha \right] ,$i.e.

$\alpha \in $ $^{\ast }%
\mathbb{R}
_{\mathbf{d}}$ has type $2$ iff $\alpha $ does not have type $1.$

A similar classification exists from above.

\bigskip \textbf{Definition 1.3.3.3.2.}[23]. Assume $\Delta >0.$

\bigskip (\textbf{i}) $\ \alpha \in $ $^{\ast }%
\mathbb{R}
_{\mathbf{d}}$ has type $1\mathbf{A}$ if $\exists x\left( x\notin \alpha
\right) \forall y\left[ x-y\notin \alpha \implies y\in \Delta \right] ,$

(\textbf{ii}) $\alpha \in $ $^{\ast }%
\mathbb{R}
_{\mathbf{d}}$ has type $2\mathbf{A}$ if $\forall x\left( x\notin \alpha
\right) \exists y\left( y\notin \alpha \right) \left[ x-y\notin \alpha %
\right] .$

\bigskip

\bigskip \textbf{Theorem} \textbf{1.3.3.3.3.}[23].

(\textbf{i}) $\ \alpha \in $ $^{\ast }%
\mathbb{R}
_{\mathbf{d}}$ has type $1$ iff $-\alpha $ has type $1\mathbf{A},$

\bigskip (\textbf{ii}) $\alpha \in $ $^{\ast }%
\mathbb{R}
_{\mathbf{d}}$ cannot have type $1$ and type $1\mathbf{A}$ simultaneously.

\textbf{Theorem} \textbf{1.3.3.3.4.}[23].Suppose $\mathbf{ab.p.}(\alpha
)=\Delta >0.$ Then $\alpha $ has type $1$

iff $\alpha $\ has the form $a+$ $\Delta $ for some $a\in $ $^{\ast }%
\mathbb{R}
.$

\textbf{Theorem} \textbf{1.3.3.3.5.}[23].$\alpha \in $ $^{\ast }%
\mathbb{R}
_{\mathbf{d}}$ has type $1\mathbf{A}$ iff $\alpha $ has the form $a+$ $%
\left( -\Delta \right) $

for some $a\in $ $^{\ast }%
\mathbb{R}
.$

\bigskip \textbf{Theorem} \textbf{1.3.3.3.6.}[23].

(\textbf{i}) \ If $\mathbf{ab.p.}(\alpha )>\mathbf{ab.p.}(\beta )$ then $%
\alpha +\beta $ has type $1$ iff $\alpha $ has type $1.$

(\textbf{ii}) If $\mathbf{ab.p.}(\alpha )=\mathbf{ab.p.}(\beta )$ then $%
\alpha +\beta $ has type $2$ iff either $\alpha $ or $\beta $

has type $2.$

\textbf{Theorem} \textbf{1.3.3.3.7.}[23].\ If $\mathbf{ab.p.}(\alpha )$ has
the form $\mathbf{B}_{a}$ then $\alpha $ has

type $1$ or type $1\mathbf{A.}$

\bigskip

\bigskip\ 

\ \ \ \ \ \ \ \ \ \ \ \ \ \ \ \ \ \ \ \ \ \ \ \ \ \ \ \ \ \ \ \ \ \ \ 

\section{I.3.3.4. $\protect\varepsilon $-Part of $\protect\alpha $ with $%
\mathbf{ab.p.}(\protect\alpha )\neq 0.$}

\bigskip

$\mathbf{Theorem}$ $\mathbf{1.3.3.4.8.}$\textbf{\ }(\textbf{i}) Suppose:

1) $-\Delta _{\mathbf{d}}<\alpha <\Delta _{\mathbf{d}},$

2) $\mathbf{ab.p.}(\alpha )=\varepsilon _{\mathbf{d}}$ i.e. $\alpha $ has
type $1.$

Then there is exist \textit{unique} $a\in 
\mathbb{R}
$ such that\bigskip\ 

\ $\ \ 
\begin{array}{cc}
\begin{array}{c}
\\ 
\alpha =\left( ^{\ast }a\right) ^{\#}+\varepsilon _{\mathbf{d}}, \\ 
\\ 
a=WST\left( \alpha \right) . \\ 
\end{array}
& \text{ \ \ }\left( 1.3.3.5\right)%
\end{array}%
$

\bigskip

(\textbf{ii}) Suppose:

1) $-\Delta _{\mathbf{d}}<\alpha _{1}<\Delta _{\mathbf{d}},-\Delta _{\mathbf{%
d}}<\alpha _{2}<\Delta _{\mathbf{d}},$

2) $\mathbf{ab.p.}(\alpha _{1})=\varepsilon _{\mathbf{d}},\mathbf{ab.p.}%
(\alpha _{2})=\varepsilon _{\mathbf{d}}$ i.e. $\alpha _{1}$ and $\alpha _{2}$
has type $1.$

Then:

$\ \ \ \ \ \ \ \ \ \ \ \ \ \ \ \ \ \ \ \ \ \ \ \ \ \ \ \ \ \ \ \ \ \ \ \ \ \ 
$

$\ \ 
\begin{array}{cc}
\begin{array}{c}
\\ 
\alpha _{1}+\alpha _{2}=WST\left( \alpha _{1}\right) +WST\left( \alpha
_{2}\right) +\varepsilon _{\mathbf{d}}. \\ 
\end{array}
& \text{ }\left( 1.3.3.6\right)%
\end{array}%
$

(\textbf{iii}) Suppose:

1) $-\Delta _{\mathbf{d}}<\alpha <\Delta _{\mathbf{d}},$

2) $\mathbf{ab.p.}(\alpha )=\varepsilon _{\mathbf{d}}$ i.e. $\alpha $ has
type $1.$

Then $\forall b\left( b\in \text{ }^{\ast }%
\mathbb{R}
\right) $:

$\ \ \ \ \ \ \ \ \ \ \ \ \ \ \ \ \ \ \ \ \ $

$\ \ 
\begin{array}{cc}
\begin{array}{c}
\\ 
b^{\#}\times \alpha =b^{\#}\times \left( ^{\ast }WST\left( \alpha \right)
\right) ^{\#}+b^{\#}\times \varepsilon _{\mathbf{d}}. \\ 
\end{array}
& \text{ \ }\left( 1.3.3.7\right)%
\end{array}%
$

$.$

(\textbf{iv}) Suppose:

1) $-\Delta _{\mathbf{d}}<\alpha _{1}<\Delta _{\mathbf{d}},-\Delta _{\mathbf{%
d}}<\alpha _{2}<\Delta _{\mathbf{d}},$

2) $\mathbf{ab.p.}(\alpha _{1})=\varepsilon _{\mathbf{d}},\mathbf{ab.p.}%
(\alpha _{2})=\varepsilon _{\mathbf{d}}$ i.e. $\alpha _{1}$ and $\alpha _{2}$
has type $1.$

Then $\forall b\left( b\in \text{ }^{\ast }%
\mathbb{R}
\right) $:

\bigskip

$\ \ 
\begin{array}{cc}
\begin{array}{c}
\\ 
b^{\#}\times \left( \alpha _{1}+\alpha _{2}\right) =b^{\#}\times \left(
^{\ast }WST\left( \alpha \right) \right) ^{\#}+ \\ 
\\ 
b^{\#}\times \left( ^{\ast }WST\left( \alpha _{2}\right) \right)
^{\#}+b^{\#}\times \varepsilon _{\mathbf{d}}. \\ 
\end{array}
& \text{ \ \ \ \ }\left( 1.3.3.8\right)%
\end{array}%
$

$\bigskip $

$\mathbf{Theorem}$ $\mathbf{1.3.3.4.9.}$\textbf{\ }(\textbf{i}) Suppose:

1) $-\Delta _{\mathbf{d}}<\alpha <\Delta _{\mathbf{d}},$

2) $\mathbf{ab.p.}(\alpha )=-\varepsilon _{\mathbf{d}}$ i.e. $\alpha $ has
type $1A.$

Then there is exist \textit{unique} $a\in 
\mathbb{R}
$ such that\bigskip\ \ 

$\ \ \ \ \ \ \ \ \ \ \ \ \ \ \ \ \ \ \ \ $

$\ 
\begin{array}{cc}
\begin{array}{c}
\\ 
\alpha =\left( ^{\ast }a\right) ^{\#}-\varepsilon _{\mathbf{d}}. \\ 
\\ 
a=WST\left( \alpha \right) . \\ 
\end{array}
& \text{ \ }\left( 1.3.3.9\right)%
\end{array}%
$

\bigskip

(\textbf{ii}) Suppose:

1) $-\Delta _{\mathbf{d}}<\alpha _{1}<\Delta _{\mathbf{d}},-\Delta _{\mathbf{%
d}}<\alpha _{2}<\Delta _{\mathbf{d}},$

2) $\mathbf{ab.p.}(\alpha _{1})=-\varepsilon _{\mathbf{d}},\mathbf{ab.p.}%
(\alpha _{2})=-\varepsilon _{\mathbf{d}}$ i.e. $\alpha _{1}$ and $\alpha
_{2} $ has type $1A$ or

3) $\mathbf{ab.p.}(\alpha _{1})=\varepsilon _{\mathbf{d}},\mathbf{ab.p.}%
(\alpha _{2})=-\varepsilon _{\mathbf{d}}$ i.e. $\alpha _{1}$ has type $1$
and $\alpha _{2}$ has

type $1A.$Then:

\bigskip

$%
\begin{array}{cc}
\begin{array}{c}
\\ 
\alpha _{1}+\alpha _{2}=WST\left( \alpha _{1}\right) +WST\left( \alpha
_{2}\right) -\varepsilon _{\mathbf{d}}. \\ 
\end{array}
& \text{ \ \ }\left( 1.3.3.10\right)%
\end{array}%
$

\bigskip

(\textbf{iii}) Suppose:

1) $-\Delta _{\mathbf{d}}<\alpha <\Delta _{\mathbf{d}},$

2) $\mathbf{ab.p.}(\alpha )=-\varepsilon _{\mathbf{d}}$ i.e. $\alpha $ has
type $1A.$

Then $\forall b\left( b\in \text{ }^{\ast }%
\mathbb{R}
\right) $:

$\ 
\begin{array}{cc}
\begin{array}{c}
\\ 
b^{\#}\times \alpha =b^{\#}\times \left( ^{\ast }WST\left( \alpha \right)
\right) ^{\#}-b^{\#}\times \varepsilon _{\mathbf{d}}. \\ 
\end{array}
& \text{ \ \ }\left( 1.3.3.11\right)%
\end{array}%
$

$.$

\bigskip

(\textbf{iv}) Suppose:

1) $-\Delta _{\mathbf{d}}<\alpha _{1}<\Delta _{\mathbf{d}},-\Delta _{\mathbf{%
d}}<\alpha _{2}<\Delta _{\mathbf{d}},$

2) $\mathbf{ab.p.}(\alpha _{1})=-\varepsilon _{\mathbf{d}},\mathbf{ab.p.}%
(\alpha _{2})=-\varepsilon _{\mathbf{d}}$ i.e. $\alpha _{1}$ and $\alpha
_{2} $ has type $1A$ or

3) $\mathbf{ab.p.}(\alpha _{1})=\varepsilon _{\mathbf{d}},\mathbf{ab.p.}%
(\alpha _{2})=-\varepsilon _{\mathbf{d}}$ i.e. $\alpha _{1}$ has type $1$
and $\alpha _{2}$ has

type $1A.$Then $\forall b\left( b\in \text{ }^{\ast }%
\mathbb{R}
\right) $:

\bigskip

$\ 
\begin{array}{cc}
\begin{array}{c}
\\ 
b^{\#}\times \left( \alpha _{1}+\alpha _{2}\right) =b^{\#}\times \left(
^{\ast }WST\left( \alpha \right) \right) ^{\#}+ \\ 
\\ 
b^{\#}\times \left( ^{\ast }WST\left( \alpha _{2}\right) \right)
^{\#}-b^{\#}\times \varepsilon _{\mathbf{d}}. \\ 
\end{array}
& \text{ \ \ \ }\left( 1.3.3.12\right)%
\end{array}%
$

\bigskip $\bigskip $

$\mathbf{Definition}$ $\mathbf{1.3.3.4.3.}$Suppose $\alpha >0$ then

$\alpha ^{+}\triangleq \left[ \alpha \right] ^{+}\triangleq \left[ x|\left(
x\in \alpha \right) \wedge \left( x\geq 0\right) \right] .$

Suppose (1) $\alpha >0$ and (2) $\mathbf{ab.p.}(\alpha )=\varepsilon _{%
\mathbf{d}}$ i.e.

$\alpha $ has type $1,$ i.e. $\alpha =\left( ^{\ast }a\right)
^{\#}+\varepsilon _{\mathbf{d}},a=WST\left( \alpha \right) ,a\in $ $%
\mathbb{R}
.$

Then \bigskip $\beta \triangleq \left[ \alpha \right] _{\varepsilon },$($%
\varepsilon \approx 0,\varepsilon \in $ $^{\ast }%
\mathbb{R}
_{+}$) is an $\varepsilon $-part of $\alpha $ iff:

$\ 
\begin{array}{cc}
\begin{array}{c}
\\ 
\ \ \forall y\left( y\geq 0\right) \left[ \left( \left[ \left( ^{\ast
}a\right) ^{\#}\right] ^{+}+y\in \alpha ^{+}\right) \wedge \left( \left[
\left( ^{\ast }a\right) ^{\#}\right] ^{+}+y\in \beta \right) \right. \\ 
\\ 
\left. \iff y\in \varepsilon ^{\#}\times \varepsilon _{\mathbf{d}}^{+}\right]
. \\ 
\end{array}
& \text{ }\left( 1.3.3.13\right)%
\end{array}%
$

$\bigskip $

$\mathbf{Theorem}$ $\mathbf{1.3.3.4.10.}$Suppose $0<\alpha <\Delta _{\mathbf{%
d}},\mathbf{ab.p.}(\alpha )=\varepsilon _{\mathbf{d}}$

i.e. $\alpha $ has type $1,$i.e. $\alpha =\left( ^{\ast }a\right)
^{\#}+\varepsilon _{\mathbf{d}},a=WST\left( \alpha \right) ,a\in $ $%
\mathbb{R}
.$

$0<b<\Delta _{\mathbf{d}},b\in $ $^{\ast }%
\mathbb{R}
,c\in 
\mathbb{R}
.$

Then:

\textbf{(i) \ \ }$\left[ \alpha \right] _{\varepsilon }=\left[ \left( ^{\ast
}a\right) ^{\#}\right] ^{+}+\varepsilon ^{\#}\times \varepsilon _{\mathbf{d}%
}^{+}.$

\textbf{(ii) \ }$\left[ b^{\#}+\alpha \right] _{\varepsilon }=\left[ \left(
^{\ast }\left( \mathbf{st}\left( b\right) \right) \right) ^{\#}\right] ^{+}+%
\left[ \left( ^{\ast }a\right) ^{\#}\right] ^{+}+\varepsilon ^{\#}\times
\varepsilon _{\mathbf{d}}^{+}.$

\textbf{(iii) }$\left[ \left( ^{\ast }c\right) ^{\#}+\alpha \right]
_{\varepsilon }=\left[ \left( ^{\ast }\left( c\right) \right) ^{\#}\right]
^{+}+\left[ \left( ^{\ast }a\right) ^{\#}\right] ^{+}+\varepsilon
^{\#}\times \varepsilon _{\mathbf{d}}^{+}.$

$\mathbf{Theorem}$ $\mathbf{1.3.3.4.11.}$\textbf{(i) }Suppose $0<\alpha
_{1}<\Delta _{\mathbf{d}},0<\alpha _{2}<\Delta _{\mathbf{d}},$

$\mathbf{ab.p.}(\alpha _{1})=\mathbf{ab.p.}(\alpha _{2})=\varepsilon _{%
\mathbf{d}},WST\left( \alpha _{1}\right) =a\in 
\mathbb{R}
,WST\left( \alpha _{2}\right) =b\in 
\mathbb{R}
.$

Then: \textbf{\ }

\textbf{(i) \ \ }$\left[ \alpha _{1}+\alpha _{2}\right] _{\varepsilon }=%
\left[ \left( ^{\ast }a\right) ^{\#}\right] ^{+}+\left[ \left( ^{\ast
}b\right) ^{\#}\right] ^{+}+\varepsilon ^{\#}\times \varepsilon _{\mathbf{d}%
}^{+}.$

\textbf{(ii) \ }$\left[ \alpha _{1}-\alpha _{2}\right] _{\varepsilon }=\left[
\left( ^{\ast }a\right) ^{\#}\right] ^{+}-\left[ \left( ^{\ast }b\right)
^{\#}\right] ^{+}-\varepsilon ^{\#}\times \varepsilon _{\mathbf{d}}^{+}.$

$\mathbf{Theorem}$ $\mathbf{1.3.3.4.11.\forall }\varepsilon \left(
\varepsilon \approx 0\right) \left[ \alpha ^{+}=\varepsilon _{\mathbf{d}%
}^{+}\iff \left[ \alpha \right] _{\varepsilon }=\varepsilon ^{\#}\times
\varepsilon _{\mathbf{d}}^{+}\right] .$

$\mathbf{Definition}$ $\mathbf{1.3.3.4.4.}$Suppose $\alpha \geq 0$ $\mathbf{%
ab.p.}(\alpha )=\Delta \geq \varepsilon _{\mathbf{d}},\alpha \neq \Delta $

and $\alpha $ has type $1,$i.e. $\alpha $ has representation $\alpha
=a^{\#}+\Delta $ for some

$a\in $ $^{\ast }%
\mathbb{R}
_{+},a^{\#}\notin \Delta .$

Then \bigskip $\beta \triangleq \left[ \alpha |a^{\#}\right] _{\varepsilon
}, $($\varepsilon \approx 0,\varepsilon \in $ $^{\ast }%
\mathbb{R}
_{+}$) is an $\varepsilon $-part of $\alpha $ \textit{for a} \textit{given}

$a\in $ $^{\ast }%
\mathbb{R}
_{+}$ iff:

\bigskip

$\ 
\begin{array}{cc}
\begin{array}{c}
\\ 
\ \ \forall y\left( y\geq 0\right) \left[ \left( \left[ a^{\#}\right]
^{+}+y\in \alpha ^{+}\right) \wedge \left( \left[ a^{\#}\right] ^{+}+y\in
\beta \right) \iff y\in \varepsilon ^{\#}\times \Delta ^{+}\right] . \\ 
\end{array}
& \text{ \ }\left( 1.3.3.14\right)%
\end{array}%
$ \ \ \ \ \ 

\bigskip\ \ \ 

\textbf{Remark. }Suppose $\mathbf{ab.p.}(\alpha )=\Delta \geq \varepsilon _{%
\mathbf{d}},\alpha =\Delta .$Then \bigskip $\beta \triangleq \left[ \alpha
|\Delta \right] _{\varepsilon },$

($\varepsilon \approx 0,\varepsilon \in $ $^{\ast }%
\mathbb{R}
_{+}$) is an $\varepsilon $-part of $\alpha $ \textit{for a} \textit{given} $%
a\in $ $^{\ast }%
\mathbb{R}
$ iff:

\bigskip

$\ 
\begin{array}{cc}
\begin{array}{c}
\\ 
\ \ \forall y\left( y\geq 0\right) \left[ \left( y\in \alpha ^{+}\right)
\wedge \left( y\in \beta \right) \iff y\in \varepsilon ^{\#}\times \Delta
^{+}\right] . \\ 
\end{array}
& \text{ }\left( 1.3.3.15\right)%
\end{array}%
$

\bigskip

Note if $\mathbf{ab.p.}(\alpha )=\varepsilon _{\mathbf{d}}$ and $\alpha
=\left( ^{\ast }a\right) ^{\#}+\varepsilon _{\mathbf{d}},a\in $ $%
\mathbb{R}
$

then \bigskip $\beta \triangleq \left[ \alpha |\left( ^{\ast }a\right) ^{\#}%
\right] _{\varepsilon }=\left[ \alpha \right] _{\varepsilon },$ $\left[
\alpha \right] _{\varepsilon }$is an $\varepsilon $-part of $\alpha .$

\textbf{Definition 1.3.3.4.4. }Suppose $\alpha >0,$ $\mathbf{ab.p.}(\alpha
)=\Delta \leq -\varepsilon _{\mathbf{d}}$ and $\alpha $

has type $1A,$i.e. $\alpha $ has representation $\alpha =a^{\#}-\Delta $ for
some

$a\in $ $^{\ast }%
\mathbb{R}
,a^{\#}\notin \Delta .$Then \bigskip $\beta \triangleq \left[ \alpha |a^{\#}%
\right] _{\varepsilon }$ is an $\varepsilon $-part of $\alpha $ \textit{for a%
} \textit{given}

$a\in $ $^{\ast }%
\mathbb{R}
$ iff:

\bigskip $%
\begin{array}{cc}
\begin{array}{c}
\\ 
\ \ \forall y\left( y\geq 0\right) \left[ \left( \left[ a^{\#}\right]
^{+}-y\in \alpha ^{+}\right) \wedge \left( \left[ a^{\#}\right] ^{+}-y\in
\beta \right) \right. \\ 
\\ 
\left. \iff y\in -\varepsilon ^{\#}\times \left( -\Delta \right) ^{+}\right]
. \\ 
\end{array}
& \text{ \ }\left( 1.3.3.16\right)%
\end{array}%
$

\bigskip

Note if $\mathbf{ab.p.}(\alpha )=-\varepsilon _{\mathbf{d}}$ i.e. $\alpha
=\left( ^{\ast }a\right) ^{\#}-\varepsilon _{\mathbf{d}},a\in $ $%
\mathbb{R}
$

then \bigskip $\beta \triangleq \left[ \alpha |\left( ^{\ast }a\right) ^{\#}%
\right] _{\varepsilon }=\left[ \alpha \right] _{\varepsilon },\left[ \alpha %
\right] _{\varepsilon }$ is an $\varepsilon $-part of $-\alpha $.$\ \ $

$\ \ \ $

\textbf{Theorem} \textbf{1.3.3.4.10. }

(\textbf{1}) Suppose $\alpha >0,$ $\mathbf{ab.p.}(\alpha )=\Delta \geq
\varepsilon _{\mathbf{d}}$ and $\alpha $ has type $1,$

i.e. $\alpha =a^{\#}+\Delta $ for some $a\in $ $^{\ast }%
\mathbb{R}
.$

Then \bigskip $\left[ \alpha |a^{\#}\right] _{\varepsilon }$ has the form

$%
\begin{array}{cc}
\begin{array}{c}
\\ 
\left[ \alpha |a^{\#}\right] _{\varepsilon }=a^{\#}+\varepsilon ^{\#}\times
\Delta ^{+} \\ 
\end{array}
& \text{ \ }\left( 1.3.3.17\right)%
\end{array}%
$\ \ \ \ \ \ \ \ \ 

\bigskip\ \ \ \ \ \ \ \ \ \ \ \ \ \ \ \ \ 

for a given $a\in $ $^{\ast }%
\mathbb{R}
.$

(\textbf{2}) Suppose $\alpha >0,$ $\mathbf{ab.p.}(\alpha )=\Delta \leq
-\varepsilon _{\mathbf{d}}$ and $\alpha $ has type $1A,$

i.e. $\alpha =a^{\#}-\Delta $ for some $a\in $ $^{\ast }%
\mathbb{R}
_{+}.$

Then \bigskip $\left[ \alpha |a^{\#}\right] _{\varepsilon }$ has the form

\ \ \ \ \ \ \ \ \ 

$\ \ \ \ 
\begin{array}{cc}
\begin{array}{c}
\\ 
\left[ \alpha |a^{\#}\right] _{\varepsilon }=\left[ a^{\#}\right]
^{+}-\varepsilon ^{\#}\times \left( -\Delta \right) ^{+} \\ 
\end{array}
& \text{\ }\left( 1.3.3.18\right)%
\end{array}%
$\ \ \ \ \ \ \ \ \ \ \ \ \ \ \ \ \ \ \ \ \ \ \ \ \ \ 

\bigskip

for a given $a\in $ $^{\ast }%
\mathbb{R}
_{+}.$

\bigskip

$\mathbf{Theorem}$ $\mathbf{1.3.3.4.11.}$\textbf{\ }

(\textbf{1}) Suppose $\alpha >0,$ $\mathbf{ab.p.}(\alpha )=\varepsilon _{%
\mathbf{d}}$ i.e. $\alpha $ has type $1$ and

$\alpha $ has representation $\alpha =\left( ^{\ast }a\right)
^{\#}+\varepsilon _{\mathbf{d}},$ for some \textit{unique} $a\in $ $%
\mathbb{R}
.$

Then \bigskip $\left[ \alpha \right] _{\varepsilon }$ has the unique form:\ 

\ \ \ \ \ \ 

$\ \ 
\begin{array}{cc}
\begin{array}{c}
\\ 
\left[ \alpha \right] _{\varepsilon }=\left[ \left( ^{\ast }a\right) ^{\#}%
\right] ^{+}+\varepsilon ^{\#}\times \varepsilon _{\mathbf{d}}^{+}. \\ 
\end{array}
& \text{ }\left( 1.3.3.19\right)%
\end{array}%
$

\bigskip

(\textbf{2}) Suppose $\alpha >0,$ $\mathbf{ab.p.}(\alpha )=-\varepsilon _{%
\mathbf{d}}$ i.e. $\alpha $ has type $1A$ and

$\alpha $ has representation $\alpha =\left( ^{\ast }a\right)
^{\#}-\varepsilon _{\mathbf{d}},$ for some \textit{unique} $a\in $ $%
\mathbb{R}
_{+}.$

Then \bigskip $\left[ \alpha \right] _{\varepsilon }$ has the unique form:\ 

$\ \ \ 
\begin{array}{cc}
\begin{array}{c}
\\ 
\left[ \alpha \right] _{\varepsilon }=\left[ \left( ^{\ast }a\right) ^{\#}%
\right] ^{+}-\varepsilon ^{\#}\times \varepsilon _{\mathbf{d}}^{+}. \\ 
\end{array}
& \text{ }\left( 1.3.3.20\right)%
\end{array}%
$

\bigskip

$\mathbf{Theorem}$ $\mathbf{1.3.3.4.12.}$\textbf{\ }(\textbf{1}) Suppose $%
\mathbf{ab.p.}(\alpha )=\varepsilon _{\mathbf{d}},WST\left( \alpha \right)
\geq 0$

i.e. $\alpha $ has type $1$ and $\alpha $ has representation $\alpha =\left(
^{\ast }a\right) ^{\#}+\varepsilon _{\mathbf{d}},$

for some \textit{unique} $a\in $ $%
\mathbb{R}
_{+}.$Then \bigskip for every $M\in $ $^{\ast }%
\mathbb{R}
_{+}$\ \ 

\ \ 

$\ \ \ \ \ \ 
\begin{array}{cc}
\begin{array}{c}
\\ 
M\times \left[ \alpha \right] _{\varepsilon }=\left[ M\times \alpha
\left\vert M\times \left( ^{\ast }a\right) ^{\#}\right. \right]
_{\varepsilon }= \\ 
\\ 
=M\times \left[ \left( ^{\ast }a\right) ^{\#}\right] ^{+}+\left( \varepsilon
^{\#}\times M\right) \times \varepsilon _{\mathbf{d}}^{+}. \\ 
\end{array}
& \text{ \ }\left( 1.3.3.21\right)%
\end{array}%
$

\bigskip

(\textbf{2}) Suppose $\mathbf{ab.p.}(\alpha )=-\varepsilon _{\mathbf{d}%
},WST\left( \alpha \right) \geq 0$

and $\alpha $ has type $1A$ i.e. $\alpha =\left( ^{\ast }a\right)
^{\#}-\varepsilon _{\mathbf{d}},$ for some \textit{unique} $a\in $ $%
\mathbb{R}
_{+}.$

Then \bigskip for every $M\in $ $^{\ast }%
\mathbb{R}
_{+}$

\ \ \ \ \ \ 

$\ \ 
\begin{array}{cc}
\begin{array}{c}
\\ 
M\times \left[ \alpha \right] _{\varepsilon }=\left[ M\times \alpha
\left\vert M\times \left( ^{\ast }a\right) ^{\#}\right. \right]
_{\varepsilon }= \\ 
\\ 
=M\times \left[ \left( ^{\ast }a\right) ^{\#}\right] ^{+}-\left( \varepsilon
^{\#}\times M\right) \times \varepsilon _{\mathbf{d}}^{+}. \\ 
\end{array}
& \text{ \ }\left( 1.3.3.22\right)%
\end{array}%
$

\bigskip

\textbf{Theorem} \textbf{1.3.3.4.13. (i) }Suppose $\mathbf{ab.p.}(\alpha
)=\varepsilon _{\mathbf{d}}$ i.e. $\alpha $ has type $1.$

Then

$\alpha =\varepsilon _{\mathbf{d}}\iff \forall y\left( y\geq 0\right)
\forall \varepsilon \left( \varepsilon \approx 0\right) \left[ \left( y\in
\alpha \right) \wedge \left( y\in \left[ \alpha \right] _{\varepsilon
}\right) \iff y\in \varepsilon ^{\#}\times \varepsilon _{\mathbf{d}}^{+}%
\right] .$

\textbf{(ii) }Suppose $\mathbf{ab.p.}(\alpha )=-\varepsilon _{\mathbf{d}}$
i.e. $\alpha $ has type $1A.$

Then

$\alpha =-\varepsilon _{\mathbf{d}}\iff \forall y\left( y\geq 0\right)
\forall \varepsilon \left( \varepsilon \approx 0\right) \left[ \left( y\in
\alpha \right) \wedge \left( y\in \left[ \alpha \right] _{\varepsilon
}\right) \iff y\in -\varepsilon ^{\#}\times \varepsilon _{\mathbf{d}}^{+}%
\right] .$

\bigskip

\ \ \ \ \ \ \ \ \ \ \ \ \ \ \ \ \ \ \ \ \ \ \ \ \ \ \ \ \ \ \ \ \ \ \ \ \ \
\ \ \ \ \ \ \ \ 

\section{I.3.3.5.Multiplicative idempotents.}

\bigskip

\textbf{Definition 1.3.3.5.1.}[23]. Let $\left[ S\right] _{\mathbf{d}%
}=\left\{ x|\exists y\left( y\in S\right) \left[ x\leq y\right] \right\} $%
.Then $\left[ S\right] _{\mathbf{d}}$

satisfies the usual axioms for a closure operation. \bigskip

Let $f$ be a continuous strictly increasing function in each variable from a

subset of $%
\mathbb{R}
^{n}$ into $%
\mathbb{R}
$. Specifically, we want the domain to be the cartesian

product $\prod_{i=1}^{n}A_{i},$ where $A_{i}$ = $\left\{ x|x>a_{i}\right\} $
for some $a_{i}\in 
\mathbb{R}
.$By transfer $f$

extends to a function $^{\ast }f$ from the corresponding subset of $^{\ast }%
\mathbb{R}
^{n}$ into $^{\ast }%
\mathbb{R}
$

which is also strictly increasing in each variable and continuous in the $Q$

topology (i.e. $\varepsilon $ and $\delta $ range over arbitrary positive
elements in $^{\ast }%
\mathbb{R}
$).

\bigskip

\bigskip \textbf{Definition 1.3.3.5.2.}[23]. Let $\alpha _{i}\in $ $^{\ast }%
\mathbb{R}
_{\mathbf{d}},b_{i}\in $ $^{\ast }%
\mathbb{R}
,$ then

\bigskip\ $\ 
\begin{array}{cc}
\begin{array}{c}
\\ 
\left[ f\right] _{\mathbf{d}}\left( \alpha _{1},\alpha _{2},...,\alpha
_{n}\right) = \\ 
\\ 
\left[ \left\{ ^{\ast }f\left( b_{1},b_{2},...,b_{n}\right) |\text{ }%
b_{i}\in \alpha _{i}\right\} \right] _{\mathbf{d}} \\ 
\end{array}
& 1.3.3.23%
\end{array}%
$

\textbf{Theorem} \textbf{1.3.3.5.1.}[23].\ If $f$ and $g$ are functions of
one variable

then $\left[ f\cdot g\right] _{\mathbf{d}}\left( \alpha \right) =\left( %
\left[ f\right] _{\mathbf{d}}\left( \alpha \right) \right) \cdot \left( %
\left[ g\right] _{\mathbf{d}}\left( \alpha \right) \right) .$

\textbf{Theorem} \textbf{1.3.3.5.2.}[23].Let $f$ and $g$ be any two terms
obtained by

compositions of strictly increasing continuous functions possibly

containing parameters in $^{\ast }%
\mathbb{R}
$. Then any relation $^{\ast }f=$ $^{\ast }g$ or $^{\ast }f<$ $^{\ast }g$

valid in $^{\ast }%
\mathbb{R}
$ extends to $^{\ast }%
\mathbb{R}
_{\mathbf{d}},$i.e. $\left[ f\right] _{\mathbf{d}}\left( \alpha \right) =$ $%
\left[ g\right] _{\mathbf{d}}\left( \alpha \right) $ or $\left[ f\right] _{%
\mathbf{d}}\left( \alpha \right) <$ $\left[ g\right] _{\mathbf{d}}\left(
\alpha \right) .$\bigskip

\textbf{Theorem} \textbf{1.3.3.5.3.}[23].The map $\alpha \longmapsto $ $%
\left[ \exp \right] _{\mathbf{d}}\left( \alpha \right) $ maps the set of

additive idempotents onto the set of all multiplicative idempotents

other than $0.$

\bigskip

Similarly, multiplicative absorption can be defined and reduced to

the study of additive absorption. Incidentally the map $\alpha \longmapsto $ 
$\left[ \exp \right] _{\mathbf{d}}\left( \alpha \right) $

is essentially the same as the map in [34, Theorem 6] which is the

map from the set of ideals onto the set of all prime ideals of the

valuation ring consisting of the finite elements of $^{\ast }%
\mathbb{R}
$.

\bigskip \bigskip

\ \ \ \ \ \ \ \ \ \ \ \ \ \ \ \ \ \ \ \ \ \ \ \ \ \ \ \ \ \ \ \ \ \ \ \ \ \
\ \ \ \ \ \ \ \ 

\section{I.3.3.6. Additive monoid of Dedekind hyperreal \ \ \ \ \ \ \ \ \ \
\ \ \ \ \ \ \ \ \ \ integers $^{\ast }\breve{%
\mathbb{Z}%
}_{\mathbf{d}}.$}

\bigskip

Well-order relation $\left( \cdot \preceq _{\mathbf{s}}\cdot \right) $ (or
strong well-ordering) on a set $S$ is a total order

on $S$ with the property: that every non-empty subset $S^{\prime }$ of $S$
has a least element in this ordering.The set $S$ together with the
well-order relation $\preceq _{\mathbf{s}}$ is then called a (strong) well-
ordered set. The natural numbers of $^{\ast }%
\mathbb{N}
$ with the well-order relation $\left( \cdot \leq _{^{\ast }%
\mathbb{N}
}\cdot \right) $ are not strong well-ordered set,for there is no smallest
infinite one.

\textbf{Definition 1.3.3.6.1. }Weak well-order relation $\left( \cdot
\preceq _{w}\cdot \right) $ (or weak well-ordering) on

a set $S$ is a total order on $S$ with the property:every non-empty subset $%
S^{\prime }\subseteqq $ $S$ has

a least element in this ordering or $S^{\prime }$ has a greatest lower bound
($\inf \left( S^{\prime }\right) $) in this

ordering.The set $S$ together with the weak well-order relation $\preceq
_{w} $is then called a

weak well-ordered set.

The natural numbers of $^{\ast }%
\mathbb{N}
$ with the well-order relation $\left( \cdot \leq _{^{\ast }%
\mathbb{N}
}\cdot \right) $ are not iven

weak well-ordered set,for there is no $\inf \left( S^{\prime }\right) $ in $%
S.$

Let us considered completion of the ring $^{\ast }%
\mathbb{Z}
$.Possible standard completion \ \ \ \ \ \ \ \ \ \ \ \ \ \ \ \ \ \ \ \ \ \ \
\ \ \ \ 

of the ring $^{\ast }%
\mathbb{Z}
$ can be constructed by Dedekind sections. Making a semantic \ \ \ \ \ \ \ \
\ \ \ \ \ \ \ \ \ \ \ \ 

leap, we now answer the question:"what is a \textit{Dedekind hyperintegers }%
?"

\textbf{Definition 1.3.3.6.2. }A \textbf{Dedekind hyperinteger }is a cut in $%
^{\ast }%
\mathbb{Z}
.$

$^{\ast }\breve{%
\mathbb{Z}%
}_{\mathbf{d}}=\left( ^{\ast }%
\mathbb{Z}
_{\mathbf{d}},+\right) $ is the class of all Dedekind \textbf{hyperinteger}s 
$x=A|B,A\subsetneqq $ $^{\ast }%
\mathbb{Z}
,B\subsetneqq $ $^{\ast }%
\mathbb{Z}
$ \ \ \ \ \ \ \ \ \ \ \ \ \ \ \ 

$(x=A,A\subsetneqq ^{\ast }%
\mathbb{Z}
).$

We will show that in a natural way $^{\ast }%
\mathbb{Z}
_{\mathbf{d}}$ is a complete ordered additive\textbf{\ }

monoid\textbf{\ }containing $^{\ast }%
\mathbb{Z}
.$\textbf{\ \ \ \ \ \ \ \ \ \ \ \ \ \ \ \ \ \ \ \ \ \ \ \ \ \ \ \ \ }

Before spelling out what this means, here are some examples of cuts in $%
^{\ast }%
\mathbb{Z}
.$

(\textbf{i}) \ \ 

$A|B=\left. \left\{ n\in \text{ }^{\ast }%
\mathbb{Z}
\text{ }|\text{ }n<1\right\} \right\vert \left\{ n\in \text{ }^{\ast }%
\mathbb{Z}
\text{ }|\text{ }n\geq 1\right\} .$

(\textbf{ii})

$A|B=\left. \left\{ n\in \text{ }^{\ast }%
\mathbb{Z}
\text{ }|\text{ }n<\omega \right\} \right\vert \left\{ n\in \text{ }^{\ast }%
\mathbb{Z}
\text{ }|\text{ }n\geq \omega \right\} ,$where $\omega \in $ $^{\ast }%
\mathbb{N}
_{\infty }.$\ 

(\textbf{iii})

\bigskip $A|B=\left. \left\{ n\in \text{ }^{\ast }%
\mathbb{Z}
\text{ }|\left( n\leq 0\right) \vee \left( \text{ }n\in 
\mathbb{Z}
_{+}\right) \right\} \right\vert \left\{ n\in \text{ }^{\ast }%
\mathbb{Z}
\text{ }|\left( n\in \text{ }^{\ast }%
\mathbb{N}
_{\infty }\right) \right\} .$

(i\textbf{v})

$A|B=$

$\left. \left\{ n\in \text{ }^{\ast }%
\mathbb{Z}
\text{ }|\left( n\leq 0\right) \vee \left[ \left( n\in \text{ }^{\ast }%
\mathbb{Z}
_{+}\right) \wedge \left( \underset{i\in 
\mathbb{N}
}{\bigwedge }\left( n\leq \omega +i\right) \right) \right] \right\}
\right\vert $

$\left\{ n\in \text{ }^{\ast }%
\mathbb{Z}
\text{ }|\left( n\in \text{ }^{\ast }%
\mathbb{Z}
_{+}\right) \wedge \text{ }\left( \underset{i\in 
\mathbb{N}
}{\bigwedge }\left( n>\omega +i\right) \right) \right\} ,$

where $\omega \in $ $^{\ast }%
\mathbb{N}
_{\infty }.$

\textbf{Remark. 1.3.3.6.1.} It is convenient to say that $A|B\in $ $^{\ast }%
\mathbb{Z}
_{\mathbf{d}}$ is an \textbf{integer}

(hyperinteger)\textbf{\ cut }in $^{\ast }%
\mathbb{Z}
$\textit{\ }if it is like the cut in examples (\textbf{i}),(\textbf{ii}):
fore some

fixed integer (hyperinteger) number $c\in $ $^{\ast }%
\mathbb{Z}
,A$ is the set of all integer $n\in $ $^{\ast }%
\mathbb{Z}
$

such that $n<c$ while $B$ is the rest of $^{\ast }%
\mathbb{Z}
.$ \ \ \ \ \ \ \ \ \ \ \ \ \ \ \ \ \ \ \ \ \ \ \ \ \ \ \ \ 

The $B$-set of an integer (hyperinteger) cut contains a smollest $c,$ and \
\ \ \ \ \ \ \ \ \ \ \ \ \ \ \ \ \ 

conversaly if $A|B$ is a cut in $^{\ast }%
\mathbb{Z}
$ and $B$ contains a smollest element $c$ then

$A|B$ is an integer (hyperinteger)\ cut at $c.$We write $\breve{c}$ for the
integer and

hyperinteger cut at $c.$This lets us think of $^{\ast }%
\mathbb{Z}
\subset $ $^{\ast }%
\mathbb{Z}
_{\mathbf{d}}$ by identifying $c$ with $\breve{c}.$

\textbf{Remark.1.3.3.6.2. }It is convenient to say that:

(\textbf{1}) $A|B\in $ $^{\ast }%
\mathbb{Z}
_{\mathbf{d}}$ is an \textbf{standard cut} in $^{\ast }%
\mathbb{Q}
$\textit{\ }if it is like the cut in example (\textbf{i}):

fore some cut $A^{\prime }|B^{\prime }\in 
\mathbb{Z}
$ the next equality is satisfied:$A|B=$ $^{\ast }\left( A^{\prime }\right)
|^{\ast }\left( B^{\prime }\right) ,$

i.e. $A$-set of a cut is an standard set.

(\textbf{2}) $A|B\in $ $^{\ast }\breve{%
\mathbb{Z}%
}_{\mathbf{d}}$ is an \textbf{internal} \textbf{cut }or \textbf{nonstandard
cut} in $^{\ast }%
\mathbb{Z}
$\textit{\ }if it is like the cut

in example (\textbf{ii}), i.e. $A$-set of a cut is an \textit{internal
nonstandard set.}

(\textbf{3}) $A|B\in $ $^{\ast }\breve{%
\mathbb{Z}%
}_{\mathbf{d}}$ is an \textbf{external} \textbf{cut }in $^{\ast }%
\mathbb{Z}
$\textit{\ }if it is like the cut in examples (\textbf{iii})-(i\textbf{v}),

i.e. $A$-set of a cut is an \textit{external set. }

\textbf{Definition 1.3.3.6.3. }A Dedekind cut $\alpha $ in $^{\ast }%
\mathbb{Z}
$ is a subset $\alpha \subset $ $^{\ast }%
\mathbb{Z}
$ of the \ \ \ \ \ \ \ \ \ \ \ \ \ \ \ \ \ \ \ \ \ \ 

hyperinteger numbers $^{\ast }%
\mathbb{Z}
$ that satisfies these properties:

\textbf{1.} $\alpha $ is not empty.

\textbf{2.} $\beta =$ $^{\ast }%
\mathbb{Z}
\backslash \alpha $ is not empty.

\textbf{3.} $\alpha $ contains no greatest element.

\textbf{4.} For $x,y\in $ $^{\ast }%
\mathbb{Z}
,$ if $x\in \alpha $ and $y<x,$ then $y\in \alpha $ as well.

\textbf{Definition 1.3.3.6.4.} A Dedekind hyperinteger $\alpha \in $ $^{\ast
}\breve{%
\mathbb{Z}%
}_{\mathbf{d}}$ is a

Dedekind cut $\alpha $ in $^{\ast }%
\mathbb{Z}
.$ We denote the set of all Dedekind hyperinteger

by $^{\ast }\breve{%
\mathbb{Z}%
}_{\mathbf{d}}$ and we order them by set-theoretic inclusion, that is to

say, for any $\alpha ,\beta \in $ $^{\ast }\breve{%
\mathbb{Z}%
}_{\mathbf{d}},$ $\alpha <_{^{\ast }\breve{%
\mathbb{Z}%
}_{\mathbf{d}}}\beta $ (or $\alpha <\beta $) if and only if $\alpha
\subsetneqq \beta $ where the

inclusion is strict. We further define $\alpha =_{^{\ast }%
\mathbb{Z}
_{\mathbf{d}}}\beta $ (or $\alpha =\beta $) as hyperinteger

if and are equal as sets. As usual, we write $\alpha \leqslant _{^{\ast }%
\mathbb{Z}
_{\mathbf{d}}}\beta $ if $\alpha <_{^{\ast }%
\mathbb{Z}
_{\mathbf{d}}}\beta $

or $\alpha =_{^{\ast }%
\mathbb{Z}
_{\mathbf{d}}}\beta $.

\textbf{Definition 1.3.3.6.5. }$M\in $ $^{\ast }%
\mathbb{Z}
_{\mathbf{d}}$ is an \textbf{upper bound }for a set $S\subset $ $^{\ast }%
\mathbb{R}
_{\mathbf{d}}$ if each

$s\in S$ satisfies $s\leq _{^{\ast }%
\mathbb{Z}
_{\mathbf{d}}}M.$ We also say that the set $S$ is \textbf{bounded above }by

$M$ iff $M\in $ $\mathbf{L}\left( ^{\ast }%
\mathbb{R}
\right) .$We also say that the set $S$ is \textbf{hyperbounded above }iff

$M\notin $ $\mathbf{L}\left( ^{\ast }%
\mathbb{R}
\right) ,$i.e. $\left\vert M\right\vert \in $ $^{\ast }%
\mathbb{R}
_{+}\backslash 
\mathbb{R}
_{+}.$

\textbf{Definition 1.3.3.6.6. }An\textbf{\ }upper bound for $S$ that is less
than all other upper \ \ \ \ \ \ \ \ \ \ \ \ \ \ \ \ \ \ \ \ \ \ 

bound for $S$ is a \textbf{least upper bound} for $S.$

\textbf{Theorem 1.3.3.6.1} Every nonempty subset $A\subsetneqq $ $^{\ast }%
\mathbb{Z}
_{\mathbf{d}}$ of Dedekind \ \ \ \ \ \ \ \ \ \ \ \ \ \ \ \ \ \ \ \ \ \ \ \ \
\ \ \ \ \ \ \ \ 

hyperinteger that is bounded (hyperbounded) above has a least \ \ \ \ \ \ \
\ \ \ \ \ \ \ \ \ \ \ \ \ \ \ \ \ \ \ \ 

upper bound.

\textbf{Definition 1.3.3.6.7.} Given two Dedekind hyperinteger $\alpha $ and 
$\beta $ we

define:

\textbf{1.}The additive identity (zero cut) denoted $0_{^{\ast }%
\mathbb{Z}
_{\mathbf{d}}}$ or $\mathbf{0},$is

$0_{^{\ast }%
\mathbb{Z}
_{\mathbf{d}}}\triangleq \left\{ x\in \text{ }^{\ast }%
\mathbb{Z}
|\text{ }x<0\right\} .$

\textbf{2.}The multiplicative identity denoted $1_{^{\ast }%
\mathbb{Z}
_{\mathbf{d}}}$ or $1,$is

$1_{^{\ast }\breve{%
\mathbb{Z}%
}_{\mathbf{d}}}\triangleq \left\{ x\in \text{ }^{\ast }%
\mathbb{Z}
|\text{ }x<1\right\} .$

\textbf{3.} Addition $\alpha +_{^{\ast }\breve{%
\mathbb{Z}%
}_{\mathbf{d}}}\beta $ of $\alpha $ and $\beta $ also denoted $\alpha +\beta 
$ is

$\alpha +_{^{\ast }\breve{%
\mathbb{Z}%
}_{\mathbf{d}}}\beta \triangleq \left\{ x+y|\text{ }x\in \alpha ,y\in \beta
\right\} .$

It is easy to see that $\alpha +_{^{\ast }\breve{%
\mathbb{Z}%
}_{\mathbf{d}}}0_{^{\ast }\breve{%
\mathbb{Z}%
}_{\mathbf{d}}}=0_{^{\ast }\breve{%
\mathbb{Z}%
}_{\mathbf{d}}}$ for all $\alpha \in $ $^{\ast }\breve{%
\mathbb{Z}%
}_{\mathbf{d}}.$

It is easy to see that $\alpha +_{^{\ast }\breve{%
\mathbb{Z}%
}_{\mathbf{d}}}\beta $ is a cut in $^{\ast }%
\mathbb{Z}
$ and $\alpha +_{^{\ast }\breve{%
\mathbb{Z}%
}_{\mathbf{d}}}\beta =\beta +_{^{\ast }\breve{%
\mathbb{Z}%
}_{\mathbf{d}}}\alpha .$

Another fundamental property of cut addition is associativity:

$\left( \alpha +_{^{\ast }\breve{%
\mathbb{Z}%
}_{\mathbf{d}}}\beta \right) +_{^{\ast }\breve{%
\mathbb{Z}%
}_{\mathbf{d}}}\gamma =\alpha +_{^{\ast }\breve{%
\mathbb{Z}%
}_{\mathbf{d}}}\left( \beta +_{^{\ast }\breve{%
\mathbb{Z}%
}_{\mathbf{d}}}\gamma \right) .$

This follows from the corresponding property of $^{\ast }%
\mathbb{Z}
.$

\textbf{4.}The opposite $-_{^{\ast }\breve{%
\mathbb{Z}%
}_{\mathbf{d}}}\alpha $ of $\alpha ,$ also denoted $-\alpha ,$is

$-\alpha \triangleq \left\{ x\in \text{ }^{\ast }%
\mathbb{Z}
|\text{ }-x\notin \alpha ,-x\text{ is not the least element of }^{\ast }%
\mathbb{Z}
\backslash \alpha \right\} .$

\textbf{5.Remark 1.3.3.6.3.} We also say that the opposite $-\alpha $ of $%
\alpha $ is the \textbf{additive }

\textbf{inverse} of $\alpha $ denoted $\div \alpha $ iff the next equality
is satisfied: $\alpha +\left( \div \alpha \right) =0.$

\textbf{6.Remark 1.3.3.6.4. }It is easy to see that for all standard and
internal cut $\alpha ^{\mathbf{Int}}$ the

opposite $-\alpha ^{\mathbf{Int}}$ is the additive inverse of $\alpha ^{%
\mathbf{Int}},$i.e. $\alpha ^{\mathbf{Int}}+\left( \div \alpha ^{\mathbf{Int}%
}\right) =0.$

\textbf{7.}We say that the cut $\alpha $ is positive if $0<\alpha $ or
negative if $\alpha <0.$

The absolute value of $\alpha ,$denoted $\left\vert \alpha \right\vert ,$is $%
\left\vert \alpha \right\vert \triangleq \alpha ,$if $\alpha \geq 0$ and $%
\left\vert \alpha \right\vert \triangleq -\alpha ,$if $\alpha \leq 0.$

\textbf{8. }The cut order enjois on $^{\ast }%
\mathbb{Z}
_{\mathbf{d}}$ the standard additional properties of:

(\textbf{i}) \ \ \textbf{transitivity: }$\alpha \leq _{^{\ast }%
\mathbb{Z}
_{\mathbf{d}}}\beta \leq _{^{\ast }\breve{%
\mathbb{Z}%
}_{\mathbf{d}}}\gamma \implies \alpha \leq _{^{\ast }\breve{%
\mathbb{Z}%
}_{\mathbf{d}}}\gamma .$\textbf{\ }

(\textbf{ii}) \ \textbf{trichotomy: }eizer $\alpha <_{^{\ast }\breve{%
\mathbb{Z}%
}_{\mathbf{d}}}\beta ,\beta <_{^{\ast }\breve{%
\mathbb{Z}%
}_{\mathbf{d}}}\alpha $ or $\alpha =_{^{\ast }\breve{%
\mathbb{Z}%
}_{\mathbf{d}}}\beta $ but only one of the

\ \ \ \ \ three things is true.

(\textbf{iii}) \textbf{translation: }$\alpha \leq _{^{\ast }\breve{%
\mathbb{Z}%
}_{\mathbf{d}}}\beta \implies \alpha +_{^{\ast }\breve{%
\mathbb{Z}%
}_{\mathbf{d}}}\gamma \leq _{^{\ast }\breve{%
\mathbb{Z}%
}_{\mathbf{d}}}\beta +_{^{\ast }\breve{%
\mathbb{Z}%
}_{\mathbf{d}}}\gamma .$

\textbf{9.}By definition above, this is what we mean when we say that

$^{\ast }\breve{%
\mathbb{Z}%
}_{\mathbf{d}}$ is an \textbf{complete} \textbf{ordered additive monoid.}

\textbf{Remark 1.3.3.6.5.} Let us considered Dedekind integer cut $c\in $ $%
^{\ast }\breve{%
\mathbb{Z}%
}_{\mathbf{d}}$ as subset of

$c\subset $ $^{\ast }%
\mathbb{R}
_{\mathbf{d}}.$We write $\widetilde{c}=$ $^{\ast }%
\mathbb{R}
_{\mathbf{d}}$-$\sup \left( c\right) =\underset{x}{\sup }\left\{ x|x\in
c\subset \text{ }^{\ast }%
\mathbb{R}
_{\mathbf{d}}\right\} $ for the cut $c\in $ $^{\ast }\breve{%
\mathbb{Z}%
}_{\mathbf{d}}.$

This lets us think of canonical imbeding $^{\ast }\breve{%
\mathbb{Z}%
}_{\mathbf{d}}$ $\underset{\mathbf{j}_{\mathbf{d}}}{\longmapsto }$ $^{\ast }%
\mathbb{R}
_{\mathbf{d}}$ monoid

$^{\ast }\breve{%
\mathbb{Z}%
}_{\mathbf{d}}$ into generalized pseudo-field $^{\ast }%
\mathbb{R}
_{\mathbf{d}}$

$\ \ \ \ \ \ \ \ \ \ \ \ \ \ \ \ \ \ \ \ \ \ \ \ \ \ \ \ \ \ \ \ \ \ \ \ \ \
\ \ \ \ \ \ \ \ \ \ \ \ \ \ \ \ 
\begin{array}{cc}
\begin{array}{c}
\\ 
^{\ast }\breve{%
\mathbb{Z}%
}_{\mathbf{d}}\subset _{\longrightarrow }\text{ }^{\ast }%
\mathbb{R}
_{\mathbf{d}} \\ 
\end{array}
& 
\end{array}%
$

by identifying $c$ with it image $\widetilde{c}=j_{\mathbf{d}}\left(
c\right) .$

\textbf{Remark 1.3.3.6.6.} It is convenient to identify monoid $^{\ast }%
\breve{%
\mathbb{Z}%
}_{\mathbf{d}}$ with it

image $j_{\mathbf{d}}\left( ^{\ast }\breve{%
\mathbb{Z}%
}_{\mathbf{d}}\right) \subset $ $^{\ast }%
\mathbb{R}
_{\mathbf{d}}.$

\ \ \ \ \ \ \ \ \ \ \ \ \ \ \ \ \ \ \ \ \ \ \ \ \ \ \ \ \ \ \ \ \ \ \ \ \ \
\ \ \ \ \ \ \ 

\section{I.3.5.Pseudo-ring of Wattenberg hyperreal \ \ \ \ \ \ \ \ \ \ \ \ \
\ \ \ \ \ \ \ \ \ \ \ \ \ \ \ \ integers $^{\ast }%
\mathbb{Z}
_{\mathbf{d}}.$}

\bigskip

The set $^{\ast }%
\mathbb{R}
_{\mathbf{d}}$ has within it a set $^{\ast }%
\mathbb{Z}
_{\mathbf{d}}\subsetneqq $ $^{\ast }\breve{%
\mathbb{Z}%
}_{\mathbf{d}}$ of Wattenberg hyperreal integers

which behave very much like hyperreals $^{\ast }%
\mathbb{Z}
$ inside $^{\ast }%
\mathbb{R}
.$In particular the

greatest integer function $^{\ast }\left[ \cdot \right] :$ $^{\ast }%
\mathbb{R}
\rightarrow $ $^{\ast }%
\mathbb{Z}
$ extends in a natural way

to $\left[ \alpha \right] _{^{\ast }%
\mathbb{R}
_{\mathbf{d}}}:$ $^{\ast }%
\mathbb{R}
_{\mathbf{d}}\rightarrow $ $^{\ast }%
\mathbb{Z}
_{\mathbf{d}}.$

\textbf{Lemma 1.3.5.1.}[24].Suppose $\alpha \in $ $^{\ast }%
\mathbb{R}
_{\mathbf{d}}$. Then the following two

conditions on $\alpha $ are equivalent:

(\textbf{i}) $\alpha =\sup \left\{ \nu ^{\#}|\left( \nu \in \text{ }^{\ast }%
\mathbb{Z}
\right) \wedge \left( \nu \leq _{^{\ast }%
\mathbb{R}
_{\mathbf{d}}}\alpha \right) \right\} ,$

(\textbf{ii}) $\alpha =\inf \left\{ \nu ^{\#}|\left( \nu \in \text{ }^{\ast }%
\mathbb{Z}
\right) \wedge \left( \alpha \leq _{^{\ast }%
\mathbb{R}
_{\mathbf{d}}}\nu \right) \right\} .$

\textbf{Definition 1.3.5.7.}If $\alpha $ satisfies conditions (\textbf{i})
or (\textbf{ii}) from lemma 1.3.5.1

$\alpha $ is said to be a $^{\ast }%
\mathbb{Z}
_{\mathbf{d}}$-integer or Wattenberg (hyperreal) integer.

\textbf{Lemma 1.3.5.2.}[24].(\textbf{i}) $^{\ast }%
\mathbb{Z}
_{\mathbf{d}}$ is the closure in $^{\ast }%
\mathbb{R}
_{\mathbf{d}}$ of $^{\ast }%
\mathbb{Z}
$,

(\textbf{ii}) $\ ^{\ast }%
\mathbb{N}
_{\mathbf{d}}$ is the closure in $^{\ast }%
\mathbb{R}
_{\mathbf{d}}$ of $^{\ast }%
\mathbb{N}
,$

(\textbf{iii}) both $^{\ast }%
\mathbb{Z}
_{\mathbf{d}}$ and $^{\ast }%
\mathbb{N}
_{\mathbf{d}}$ are closed with respect to taking $\sup $

\ \ \ \ \ \ and $\inf $.

(\textbf{iv}) $\ ^{\ast }%
\mathbb{N}
_{\mathbf{d}}$ is a weak well-ordered set.

\textbf{Lemma 1.3.5.3.}[24].Suppose that $\lambda ,\nu \in $ $^{\ast }%
\mathbb{Z}
_{\mathbf{d}}.$ Then,

(\textbf{i}) $\ \ \ \ \lambda +_{^{\ast }%
\mathbb{R}
_{\mathbf{d}}}\nu \in $ $^{\ast }%
\mathbb{Z}
_{\mathbf{d}}.$

(\textbf{ii}) \ \ \ \ $-_{^{\ast }%
\mathbb{R}
_{\mathbf{d}}}\lambda \in $ $^{\ast }%
\mathbb{Z}
_{\mathbf{d}}.$

(\textbf{iii}) $\ \ \lambda \times _{^{\ast }%
\mathbb{R}
_{\mathbf{d}}}\nu \in $ $^{\ast }%
\mathbb{Z}
_{\mathbf{d}}.$\bigskip

\textbf{Definition 1.3.5.8. }Suppose $\alpha \in $ $^{\ast }%
\mathbb{R}
_{\mathbf{d}}.$ Then, we define

$\left[ \cdot \right] _{^{\ast }%
\mathbb{R}
_{\mathbf{d}}}:$ $^{\ast }%
\mathbb{R}
_{\mathbf{d}}\rightarrow $ $^{\ast }%
\mathbb{Z}
_{\mathbf{d}}$ by: $\left[ \alpha \right] _{^{\ast }%
\mathbb{R}
_{\mathbf{d}}}\triangleq \sup \left\{ \nu |\left( \nu \in \text{ }^{\ast }%
\mathbb{Z}
\right) \wedge \left( \nu \leq _{^{\ast }%
\mathbb{R}
_{\mathbf{d}}}\alpha \right) \right\} .$

\textbf{Remark 1.3.5.7.}There are two possibilities:

(\textbf{i}) collection $\left\{ \nu |\left( \nu \in \text{ }^{\ast }%
\mathbb{Z}
\right) \wedge \left( \nu \leq _{^{\ast }%
\mathbb{R}
_{\mathbf{d}}}\alpha \right) \right\} $ has no greatest element. In this

case $\left[ \alpha \right] _{^{\ast }%
\mathbb{R}
_{\mathbf{d}}}=\alpha $ since $\left[ \alpha \right] _{^{\ast }%
\mathbb{R}
_{\mathbf{d}}}<_{^{\ast }%
\mathbb{R}
_{\mathbf{d}}}\alpha $ implies $\exists a\left( a\in \text{ }^{\ast }%
\mathbb{R}
\right) \left[ \left[ \alpha \right] _{^{\ast }%
\mathbb{R}
_{\mathbf{d}}}<_{^{\ast }%
\mathbb{R}
_{\mathbf{d}}}a<_{^{\ast }%
\mathbb{R}
_{\mathbf{d}}}\alpha \right] .$

But then $[a]_{^{\ast }%
\mathbb{R}
_{\mathbf{d}}}<_{^{\ast }%
\mathbb{R}
_{\mathbf{d}}}\alpha $ which implies $[a]_{^{\ast }%
\mathbb{R}
_{\mathbf{d}}}+_{^{\ast }%
\mathbb{R}
_{\mathbf{d}}}1_{^{\ast }%
\mathbb{R}
_{\mathbf{d}}}<_{^{\ast }%
\mathbb{R}
_{\mathbf{d}}}\alpha $ contradicting

with $\left[ \alpha \right] _{^{\ast }%
\mathbb{R}
_{\mathbf{d}}}<_{^{\ast }%
\mathbb{R}
_{\mathbf{d}}}a\leq _{^{\ast }%
\mathbb{R}
_{\mathbf{d}}}[a]_{^{\ast }%
\mathbb{R}
_{\mathbf{d}}}+_{^{\ast }%
\mathbb{R}
_{\mathbf{d}}}1_{^{\ast }%
\mathbb{R}
_{\mathbf{d}}}$

(\textbf{ii}) collection $\left\{ \nu |\left( \nu \in \text{ }^{\ast }%
\mathbb{Z}
\right) \wedge \left( \nu \leq _{^{\ast }%
\mathbb{R}
_{\mathbf{d}}}\alpha \right) \right\} $ has a greatest element,$\nu .$In this

case $\left[ \alpha \right] _{^{\ast }%
\mathbb{R}
_{\mathbf{d}}}=\nu \in $ $^{\ast }%
\mathbb{N}
.$

\textbf{Definition 1.3.5.9. }$^{\ast }%
\mathbb{Z}
_{\mathbf{d}}$-integer $\sup \left( 
\mathbb{N}
\right) =\inf \left( ^{\ast }%
\mathbb{N}
_{\infty }\right) $ we denote $\omega _{\mathbf{d}}.$

\textbf{Definition 1.3.5.10. }Suppose $\nu \in $ $^{\ast }%
\mathbb{N}
_{\infty }.$ Then, we define $\nu $-\textbf{block} $\mathbf{bk}\left[ \nu %
\right] $

as a set of hyper integers such that $\mathbf{bk}\left[ \nu \right] =\left\{
\nu \pm n|n\in 
\mathbb{N}
\right\} .$

For $\nu ,\lambda \in $ $^{\ast }%
\mathbb{N}
_{\infty }$ there are two possibilities:

(\textbf{i}) $\nu -\lambda \in 
\mathbb{Z}
.$ In this case $\mathbf{bk}\left[ \nu \right] =\mathbf{bk}\left[ \lambda %
\right] $ and we write $\mathbf{bk}\left[ \nu \right] =\mathbf{bk}\left[ 
\widetilde{\nu }\right] $ where

$\widetilde{\nu }\in $ $^{\ast }%
\mathbb{Z}
_{+}/%
\mathbb{Z}
.$

(\textbf{ii}) $\left\vert \nu -\lambda \right\vert \in $ $^{\ast }%
\mathbb{N}
_{\infty }.$In this case $\mathbf{bk}\left[ \nu \right] \neq \mathbf{bk}%
\left[ \lambda \right] $ and $\mathbf{bk}\left[ \widetilde{\nu }\right] \neq 
\mathbf{bk}\left[ \widetilde{\lambda }\right] .$

\textbf{Lemma 1.3.5.11. }$^{\ast }%
\mathbb{N}
=%
\mathbb{N}
\cup \left( \bigcup_{\widetilde{\nu }}\mathbf{bk}\left[ \widetilde{\nu }%
\right] \right) .$

\textbf{Proof. }Clear\textbf{\ }by\textbf{\ }using\textbf{\ }%
[25,Chapt.1,section 9].\textbf{\ \ \ \ \ \ \ \ \ \ \ \ \ \ \ \ \ \ \ }

\bigskip\ \ \ \ \ \ \ \ \ \ \ \ \ \ \ \ \ \ \ \ \ \ \ \ \ \ \ \ \ \ \ \ \ \
\ \ \ 

\section{I.3.6.External summation of countable and hyperfinite sequences in $%
^{\ast }%
\mathbb{R}
_{\mathbf{d}}.$}

\textbf{Definition 1.3.6.1. }Let $\mathbf{S}_{X}$ denote the group of
permutations of the set $X$

and $\mathbf{H}_{X}$ denote ultrafilter on the set $X.$ Permutation $\sigma
\in \mathbf{S}_{X}$ is \textit{admissible}

iff $\sigma $ preserv $\mathbf{H}_{X},$i.e. for any $A\in \mathbf{H}_{X}$
the next condition is satisfied:

$\sigma \left( A\right) \in \mathbf{H}_{X}.$ \ \ \ \ \ \ \ \ \ \ \ \ \ \ \ \
\ \ \ \ \ \ \ \ \ \ \ 

Below we denote by $\widehat{\mathbf{S}}_{X,\mathbf{H}_{X}}$ the subgroup $%
\widehat{\mathbf{S}}_{X,\mathbf{H}_{X}}\subsetneqq \mathbf{S}_{X}$ of the all

admissible permutations.

\textbf{Definition 1.3.6.2. }Let us consider countable sequence $\mathbf{s}%
_{n}:%
\mathbb{N}
\rightarrow 
\mathbb{R}
,$

such that:

(\textbf{a}) $\forall n\left( \mathbf{s}_{n}\geq 0\right) $ or

(\textbf{b}) $\forall n\left( \mathbf{s}_{n}<0\right) $ and hyperreal number
denoted $\left[ \mathbf{s}_{n}\right] $ which \ \ \ \ \ \ \ \ \ \ \ \ \ 

formed from sequence $\ \left\{ \mathbf{s}_{n}\right\} _{n\in 
\mathbb{N}
}$ by the law

\bigskip 

\bigskip $\ 
\begin{array}{cc}
\begin{array}{c}
\\ 
\left[ \mathbf{s}_{n}\right] = \\ 
\\ 
\left( \mathbf{s}_{0},\mathbf{s}_{0}+\mathbf{s}_{1},\mathbf{s}_{0}+\mathbf{s}%
_{1}+\mathbf{s}_{2},...,\sum_{0}^{i}\mathbf{s}_{i},...\right) \in \text{ }%
^{\ast }%
\mathbb{R}
. \\ 
\end{array}
& \text{\ }\left( 1.3.6.1\right) 
\end{array}%
$

\bigskip

\bigskip Then external sum of the countable sequence $\mathbf{s}_{n}$ denoted

\bigskip 

$\ 
\begin{array}{cc}
\begin{array}{c}
\\ 
Ext\text{-}\dsum\limits_{n\in 
\mathbb{N}
}\mathbf{s}_{n} \\ 
\end{array}
& \text{ \ }\left( 1.3.6.2\right) 
\end{array}%
$

\bigskip is

$\ 
\begin{array}{cc}
\begin{array}{c}
\\ 
\left( a\right) :Ext-\dsum\limits_{n\in 
\mathbb{N}
}\mathbf{s}_{n}\triangleq \inf \left\{ \left[ \mathbf{s}_{\sigma \left(
n\right) }\right] |\text{ }\sigma \in \widehat{\mathbf{S}}_{%
\mathbb{N}
,\mathbf{H}_{%
\mathbb{N}
}}\right\} , \\ 
\\ 
\left( b\right) :Ext-\dsum\limits_{n\in 
\mathbb{N}
}\mathbf{s}_{n}\triangleq \sup \left\{ \left[ \mathbf{s}_{\sigma \left(
n\right) }\right] |\text{ }\sigma \in \widehat{\mathbf{S}}_{%
\mathbb{N}
,\mathbf{H}_{%
\mathbb{N}
}}\right\} \\ 
\end{array}
& \text{\ }\left( 1.3.6.3\right)%
\end{array}%
$

\bigskip accordingly.

\textbf{Example 1.3.6.1. }Let us consider countable sequence $\left\{ 
\mathbf{1}_{n}\right\} _{n\in 
\mathbb{N}
}$ such \ \ \ \ \ \ \ \ \ \ \ \ \ \ \ \ \ \ \ \ \ \ \ \ \ \ \ \ \ \ 

that: $\forall n\left( \mathbf{1}_{n}=1\right) .$Hence $\left[ \mathbf{1}_{n}%
\right] =\left( 1,2,3,....,i,...\right) =\varpi \in $ $^{\ast }%
\mathbb{R}
$ and using \ \ \ \ \ \ \ \ \ \ \ \ \ \ \ \ \ \ \ \ \ \ \ \ \ 

Eq.(1.3.3) one obtain $\ \ \ \ \ \ \ \ \ \ \ \ \ \ \ \ \ \ \ \ \ \ \ \ \ \ \
\ \ \ $

$\bigskip $

$%
\begin{array}{cc}
\begin{array}{c}
\\ 
Ext\text{ -}\dsum\limits_{n\in 
\mathbb{N}
}\mathbf{1}_{n}=\varpi \in \text{ }^{\ast }%
\mathbb{R}
. \\ 
\end{array}
& \text{ \ \ }\left( 1.3.6.4\right)%
\end{array}%
$

\ \ \ \ \ \ \ \ \ \ \ \ \ \ \ \ \ \ \ \ \ \ \ \ \ \ \ \ \ \ \ \ \ \ \ \ \ \
\ \ \ \ \ \ \ \ \ \ \ \ \ \ \ \ 

\textbf{Example 1.3.6.2. }Let us consider countable sequence $\left\{ 
\mathbf{1}_{n}^{\blacktriangledown }\right\} _{n\in 
\mathbb{N}
}$ such that: $\ \ \ \ \ \ \ $

$\left\{ n|\mathbf{1}_{n}^{\blacktriangledown }=1\right\} \in \mathbf{H}_{%
\mathbb{N}
}.$Hence $\left[ \mathbf{1}_{n}^{\blacktriangledown }\right] =\left(
1,2,3,....,i,...\right) \left( \text{mod}\mathbf{H}_{%
\mathbb{N}
}\right) $ $=\varpi \in $ $^{\ast }%
\mathbb{R}
$ and \ \ \ \ \ \ \ \ \ \ \ \ \ \ \ \ \ \ \ \ \ \ 

using Eq.(1.3.3) one obtain

\bigskip

$\ \ 
\begin{array}{cc}
\begin{array}{c}
\\ 
Ext-\dsum\limits_{n\in 
\mathbb{N}
}\mathbf{1}_{n}^{\blacktriangledown }=\varpi \in \text{ }^{\ast }%
\mathbb{R}
. \\ 
\end{array}
& \text{ \ }\left( 1.3.6.5\right)%
\end{array}%
$

\bigskip

\textbf{Example 1.3.6.3. }(\textbf{Euler's infinite number }$E^{\#}$). Let
us consider

countable sequence $\mathbf{h}_{n}=n^{-1}$. Hence

$\ \left[ \mathbf{h}_{n}\right] =\left( 1,1+\dfrac{1}{2},1+\dfrac{1}{2}+%
\dfrac{1}{3},....,1+\dfrac{1}{2}+\dfrac{1}{3}+...+\dfrac{1}{i},...\right)
\in $ $^{\ast }%
\mathbb{R}
$

and using Eq.(1.3.3) one obtain

\bigskip

\ $\ 
\begin{array}{cc}
\begin{array}{c}
\\ 
Ext-\dsum\limits_{n=1}^{\infty }\mathbf{h}_{n}=E^{\#}\in \text{ }^{\ast }%
\mathbb{R}
_{\mathbf{d}}. \\ 
\end{array}
& \text{ \ \ }\left( 1.3.6.6\right)%
\end{array}%
$

\bigskip

\textbf{Definition 1.3.6.8. }Let us consider countable sequence $\mathbf{s}%
_{n}:%
\mathbb{N}
\rightarrow 
\mathbb{R}
$ and

two subsequences denoted $\mathbf{s}_{n}^{+}:%
\mathbb{N}
\rightarrow 
\mathbb{R}
,\mathbf{s}_{n}^{-}:%
\mathbb{N}
\rightarrow 
\mathbb{R}
$ which formed from \ \ \ \ \ \ \ \ \ \ \ \ \ \ \ \ \ 

sequence $\left\{ \mathbf{s}_{n}\right\} _{n\in 
\mathbb{N}
}$ by the law

\bigskip

$\ 
\begin{array}{cc}
\begin{array}{c}
\\ 
\mathbf{s}_{n}^{+}=\mathbf{s}_{n}\iff \mathbf{s}_{n}\geq 0, \\ 
\\ 
\mathbf{s}_{n}^{+}=0\iff \mathbf{s}_{n}<0 \\ 
\end{array}
& \text{ }\left( 1.3.6.7\right)%
\end{array}%
$

\bigskip

and accordingly by the law \bigskip

\bigskip $\ 
\begin{array}{cc}
\begin{array}{c}
\\ 
\mathbf{s}_{n}^{-}=\mathbf{s}_{n}\iff \mathbf{s}_{n}<0, \\ 
\\ 
\mathbf{s}_{n}^{-}=0\iff \mathbf{s}_{n}\geq 0 \\ 
\end{array}
& \text{ \ }\left( 1.3.6.8\right)%
\end{array}%
$

\bigskip Hence $\left\{ \mathbf{s}_{n}\right\} _{n\in 
\mathbb{N}
}=\left\{ \mathbf{s}_{n}^{+}+\mathbf{s}_{n}^{-}\right\} _{n\in 
\mathbb{N}
}.$

\textbf{Example 1.3.6.4. }Let us consider countable sequence$\ $

$\bigskip \ \ \ \ \ \ \ \ \ \ \ \ \ \ \ \ \ \ \ \ \ \ \ \ \ \ \ $

\bigskip $\ 
\begin{array}{cc}
\begin{array}{c}
\\ 
\ \left\{ \mathbf{1}_{n}^{\pm }\right\} _{n\in 
\mathbb{N}
}=\left\{ 1,-1,1,-1,...,1,-1,...\right\} . \\ 
\end{array}
& \text{ \ \ }\left( 1.3.6.9\right)%
\end{array}%
$

Hence $\left\{ \mathbf{1}_{n}^{\pm }\right\} _{n\in 
\mathbb{N}
}=$ $\left\{ \mathbf{1}_{n}^{+}+\mathbf{1}_{n}^{-}\right\} _{n\in 
\mathbb{N}
}\ \ $where$\ \ \ \ \ \ \ \ \ \ \ \ \ \ \ \ \ \ \ \ \ \ $

\bigskip\ \ \ \ \ \ \ \ \ \ \ \ \ \ \ \ \ \ \ \ \ \ \ \ \ \ \ \ \ \ \ \ \ \
\ \ \ \ \ \ \ \ \ \ \ \ \ \ \ \ \ $\ \ 
\begin{array}{cc}
\begin{array}{c}
\\ 
\left\{ \mathbf{1}_{n}^{+}\right\} _{n\in 
\mathbb{N}
}=\left\{ 1,0,1,0,...,1,0,...\right\}  \\ 
\\ 
\left\{ \mathbf{1}_{n}^{-}\right\} _{n\in 
\mathbb{N}
}=\left\{ 0,-1,0,-1,...,0,-1,...\right\} . \\ 
\end{array}
& \text{ \ \ }\left( 1.3.6.10\right) 
\end{array}%
$\ \ 

\textbf{Definition 1.3.6.9.}The\textbf{\ }external sum of the arbitrary
countable \ \ \ \ \ \ \ \ \ \ \ \ \ \ \ \ \ \ \ \ \ \ \ \ \ \ \ \ \ \ \ \ \
\ \ \ 

sequence $\left\{ \mathbf{s}_{n}\right\} _{n\in 
\mathbb{N}
}$ denoted

\bigskip

$\ \ \ 
\begin{array}{cc}
\begin{array}{c}
\\ 
Ext\text{ -}\dsum\limits_{n\in 
\mathbb{N}
}\mathbf{s}_{n} \\ 
\end{array}
& \text{ \ }\left( 1.3.6.12\right)%
\end{array}%
$

\bigskip

\bigskip is $\ \ \ \ \ \ \ \ \ \ \ \ \ \ \ \ \ \ \ \ \ \ \ \ \ \ \ \ \ \ \ \
\ \ \ \ $

$\ \ \ \ 
\begin{array}{cc}
\begin{array}{c}
\\ 
Ext-\dsum\limits_{n\in 
\mathbb{N}
}\mathbf{s}_{n}\triangleq \\ 
\\ 
\left( Ext-\dsum\limits_{n\in 
\mathbb{N}
}\mathbf{s}_{n}^{+}\right) +\left( Ext-\dsum\limits_{n\in 
\mathbb{N}
}\mathbf{s}_{n}^{-}\right) . \\ 
\end{array}
& \text{ }\left( 1.3.6.13\right)%
\end{array}%
$

\bigskip

\textbf{Example 1.3.6.5. }Let us consider countable sequence (1.3.9)$\ \ $%
Using \ \ \ \ \ \ \ \ \ \ \ \ \ \ \ \ \ \ \ \ \ \ 

Eq.(1.3.3),Eq.(1.3.13) and Eq.(1.3.5) one obtain$\ $

$\bigskip $

$\ \ \ \ \ \ \ \ \ \ \ \ \ \ \ \ \ \ $

$\ \ 
\begin{array}{cc}
\begin{array}{c}
\\ 
Ext\text{ -}\dsum\limits_{n\in 
\mathbb{N}
}\ \mathbf{1}_{nn\in 
\mathbb{N}
}^{\pm }= \\ 
\\ 
\left( Ext\text{ -}\dsum\limits_{n\in 
\mathbb{N}
}\mathbf{1}_{n}^{+}\right) +\left( Ext\text{ -}\dsum\limits_{n\in 
\mathbb{N}
}\mathbf{1}_{n}^{-}\right) = \\ 
\\ 
=\varpi -\varpi =0. \\ 
\end{array}
& \text{ \ \ \ \ \ \ \ \ }\left( 1.3.6.14\right)%
\end{array}%
$

\bigskip

\textbf{Definition 1.3.6.10. }Let us consider countable sequence $\mathbf{s}%
_{n}^{\#}:%
\mathbb{N}
\rightarrow $ $^{\ast }%
\mathbb{R}
_{\mathbf{d}},$ \ \ \ \ \ \ \ \ \ \ \ \ \ \ \ \ \ \ \ \ \ \ \ \ \ \ \ \ \ 

such that: (\textbf{a}) $\forall n\left( \mathbf{s}_{n}^{\#}\geq 0\right) $
or

(\textbf{b}) $\forall n\left( \mathbf{s}_{n}^{\#}<0\right) .$

\bigskip Then external sum of the countable sequence $\mathbf{s}_{n}^{\#}$
denoted

\bigskip

$\ 
\begin{array}{cc}
\begin{array}{c}
\\ 
\#Ext\text{-}\dsum\limits_{n\in 
\mathbb{N}
}\mathbf{s}_{n}^{\#} \\ 
\end{array}
& \text{ }\left( 1.3.6.15\right) 
\end{array}%
$

\bigskip is

$\ $

$\ \ 
\begin{array}{cc}
\begin{array}{c}
\\ 
\left( a\right) :\#Ext\text{ -}\dsum\limits_{n\in 
\mathbb{N}
}\mathbf{s}_{n}^{\#}\triangleq \\ 
\\ 
\text{ }\underset{k\in 
\mathbb{N}
}{\sup }\left\{ \left. \dsum\limits_{n\leq k}\mathbf{s}_{\sigma \left(
n\right) }^{\#}\right\vert \text{ }\sigma \in \mathbf{S}_{%
\mathbb{N}
}\right\} , \\ 
\\ 
\left( b\right) :\#Ext\text{ -}\dsum\limits_{n\in 
\mathbb{N}
}\mathbf{s}_{n}^{\#}\triangleq \\ 
\\ 
\text{ }\underset{k\in 
\mathbb{N}
}{\inf }\left\{ \left. \dsum\limits_{n\leq k}\mathbf{s}_{\sigma \left(
n\right) }^{\#}|\text{ }\right\vert \sigma \in \mathbf{S}_{%
\mathbb{N}
}\right\} . \\ 
\end{array}
& \text{ \ }\left( 1.3.6.16\right)%
\end{array}%
$

\textbf{Definition 1.3.6.11. }Let us consider countable sequence $\mathbf{s}%
_{n}^{\#}:%
\mathbb{N}
\rightarrow $ $^{\ast }%
\mathbb{R}
_{\mathbf{c}}$ \ \ \ \ \ \ \ \ \ \ \ \ \ \ \ \ \ \ \ \ \ \ \ \ \ \ 

and two subsequences denoted $^{\#}\mathbf{s}_{n}^{+}:%
\mathbb{N}
\rightarrow $ $^{\ast }%
\mathbb{R}
_{\mathbf{c}},^{\#}\mathbf{s}_{n}^{-}:%
\mathbb{N}
\rightarrow $ $^{\ast }%
\mathbb{R}
_{\mathbf{c}}$ which \ \ \ \ \ \ \ \ \ \ \ \ \ \ \ \ \ \ \ \ \ \ \ \ \ 

formed from sequence $\left\{ \mathbf{s}_{n}^{\#}\right\} _{n\in 
\mathbb{N}
}$ by the law

\bigskip

$\ 
\begin{array}{cc}
\begin{array}{c}
\\ 
^{\#}\mathbf{s}_{n}^{+}=\mathbf{s}_{n}\iff \mathbf{s}_{n}^{\#}\geq 0, \\ 
\\ 
^{\#}\mathbf{s}_{n}^{+}=0\iff \mathbf{s}_{n}^{\#}<0 \\ 
\end{array}
& \text{ \ \ }\left( 1.3.6.17\right)%
\end{array}%
$

\bigskip

and accordingly by the law

\bigskip

$\ 
\begin{array}{cc}
\begin{array}{c}
\\ 
^{\#}\mathbf{s}_{n}^{-}=\mathbf{s}_{n}^{\#}\iff \mathbf{s}_{n}^{\#}<0, \\ 
\\ 
^{\#}\mathbf{s}_{n}^{-}=0\iff \mathbf{s}_{n}^{\#}\geq 0 \\ 
\end{array}
& \text{ }\left( 1.3.6.18\right)%
\end{array}%
$

\bigskip

Hence $\left\{ \mathbf{s}_{n}^{\#}\right\} _{n\in 
\mathbb{N}
}=\left\{ ^{\#}\mathbf{s}_{n}^{+}+\text{ }^{\#}\mathbf{s}_{n}^{-}\right\}
_{n\in 
\mathbb{N}
}.$

\textbf{Definition 1.3.6.12.}The\textbf{\ }external sum of the arbitrary
countable

sequence $\left\{ \mathbf{s}_{n}\right\} _{n\in 
\mathbb{N}
}$ denoted $\ \ \ \ $

$\bigskip $

$\ \ 
\begin{array}{cc}
\begin{array}{c}
\\ 
\#Ext-\dsum\limits_{n\in 
\mathbb{N}
}\mathbf{s}_{n}^{\#} \\ 
\end{array}
& \text{ \ \ }\left( 1.3.6.19\right)%
\end{array}%
$

\bigskip

\bigskip is

$\ \ \ \ \ \ \ \ \ \ \ \ \ \ \ \ \ \ \ \ \ \ \ \ \ \ \ \ \ \ \ \ \ \ \ \ $

$\ 
\begin{array}{cc}
\begin{array}{c}
\\ 
\#Ext\text{ -}\dsum\limits_{n\in 
\mathbb{N}
}\mathbf{s}_{n}^{\#}\triangleq \\ 
\\ 
\left( \#Ext\text{ -}\dsum\limits_{n\in 
\mathbb{N}
}\left( ^{\#}\mathbf{s}_{n}^{+}\right) \right) +\left( \#Ext\text{ -}%
\dsum\limits_{n\in 
\mathbb{N}
}\left( ^{\#}\mathbf{s}_{n}^{-}\right) \right) . \\ 
\end{array}
& \text{ \ \ }\left( 1.3.6.20\right)%
\end{array}%
$

\bigskip

\textbf{Definition 1.3.6.13. }Let us consider an nonempty subset $\mathbf{A}%
\subsetneqq $ $^{\ast }%
\mathbb{R}
_{\mathbf{d}}$\ which \ \ \ \ \ \ \ \ \ \ \ \ \ \ \ \ \ \ \ \ \ \ \ \ \ \ \
\ \ \ 

is bounded or hyperbounded from above and such that:

$\sup \left( \mathbf{A}\right) \pm \varepsilon \neq \sup \left( \mathbf{A}%
\right) $ for any $\varepsilon \approx 0.$We call this least upper bound $%
\sup \left( \mathbf{A}\right) $ \ \ \ \ \ \ \ \ \ \ \ \ \ \ \ \ \ \ \ \ \ \
\ \ \ \ \ 

the \textit{strong least upper bound} or \textit{strong supremum}, written
as $\mathbf{s}$-$\sup \left( \mathbf{A}\right) $.

\textbf{Proposition 1.3.6.1.} If $\mathbf{A}$ is a nonempty subset of $%
^{\ast }%
\mathbb{R}
_{\mathbf{d}}$ which is bounded \ \ \ \ \ \ \ \ \ \ \ \ \ \ \ \ \ \ \ \ \ \
\ \ \ \ 

from above and strong supremum $\mathbf{s}$-$\sup \left( \mathbf{A}\right) $
exist, then:

(\textbf{1}) $\mathbf{s}$-$\sup \left( \mathbf{A}\right) $ is the unique
number such that $\mathbf{s}$-$\sup \left( \mathbf{A}\right) $ is an upper
bound \ \ \ \ \ \ \ \ \ \ \ \ \ \ \ \ \ \ \ \ \ \ \ \ \ \ 

for $\mathbf{A}$ and $\mathbf{s}$-$\sup \left( \mathbf{A}\right)
-\varepsilon $ is not a upper bound for $\mathbf{A}$ for any $\varepsilon
\approx 0,\varepsilon >0;$

(\textbf{2}) (\textbf{The Strong Approximation Property}) let $\varepsilon
\approx 0,\varepsilon >0$ there exist $x\in \mathbf{A}$ \ \ \ \ \ \ \ \ \ \
\ \ \ \ \ \ \ \ \ \ 

such that $\mathbf{s}$-$\sup \left( \mathbf{A}\right) -\varepsilon <x\leq 
\mathbf{s}$-$\sup \left( \mathbf{A}\right) .$

\textbf{Proof.}(\textbf{2})\textbf{\ }If not, then $\mathbf{s}$-$\sup \left( 
\mathbf{A}\right) -\varepsilon $ is an upper bound of $\mathbf{A}$ less than
the least \ \ \ \ \ \ \ \ \ \ \ \ \ \ \ \ \ 

upper bound, which is a contradiction.

\textbf{Corollary} \textbf{1.3.6.1. }Let $\mathbf{A}$ be bounded or
hyperbounded from above and \ \ \ \ \ \ \ \ \ \ \ \ \ \ \ \ \ \ \ \ \ \ 

non-empty set such that $\mathbf{s}$-$\sup \left( \mathbf{A}\right) $ exist.
There is a function

$\alpha \left( \circ \right) :$ $^{\ast }\mathbf{%
\mathbb{N}
}_{\infty }\rightarrow $ $^{\ast }%
\mathbb{R}
_{\mathbf{d}}$ such that for all $\mathbf{n\in }$ $^{\ast }\mathbf{%
\mathbb{N}
}_{\infty }$ we have

$\mathbf{s}$-$\sup \left( \mathbf{A}\right) -\mathbf{n}^{-1}<\alpha \left( 
\mathbf{n}\right) \leq \mathbf{s}$-$\sup \left( \mathbf{A}\right) $

\textbf{Theorem 1.3.6.2. }Let $\mathbf{A}$ be a non-empty set which is
bounded or \ \ \ \ \ \ \ \ \ \ \ \ \ \ \ \ \ \ \ \ \ \ \ \ \ \ \ \ 

hyperbounded from below. Then the set of lower bounds of $\mathbf{A}$ has a
\ \ \ \ \ \ \ \ \ \ \ \ \ \ \ \ \ \ \ \ \ \ \ \ \ \ 

greatest element.

\textbf{Proof. }Let $-\mathbf{A\triangleq }\left\{ -x|x\in \mathbf{A}%
\right\} .$We know that (\textbf{i}) $\forall x_{x\in \text{ }^{\ast }%
\mathbb{R}
_{\mathbf{d}}}\forall y_{y\in \text{ }^{\ast }%
\mathbb{R}
_{\mathbf{d}}}\left( x\leq y\iff -y\leq -x\right) .$

Let $l_{\mathbf{A}}$ be a lower bound of $\mathbf{A.}$Then $l_{\mathbf{A}%
}\leq x$ for all $x\in \mathbf{A.}$So $-x\leq -l_{\mathbf{A}}$ for all

$x\in \mathbf{A,}$that is $y\leq l_{\mathbf{A}}$ for all $y\in -\mathbf{A.}$
So $-\mathbf{A}$ is bounded above, and non-empty,

so by the Theorem 1.3.1 $\sup \left( -\mathbf{A}\right) $ exists.

We shall prove now that: (\textbf{ii}) $-\sup \left( -\mathbf{A}\right) $ is
a lower bound of $\mathbf{A,}$(\textbf{iii}) if $l_{\mathbf{A}}$ is a lower
\ \ \ \ \ \ \ \ \ \ \ \ 

bound of $\mathbf{A}$ then $l_{\mathbf{A}}\leq -\sup \left( -\mathbf{A}%
\right) .$ (\textbf{ii}) if $x\in \mathbf{A}$ then $-x\in -\mathbf{A}$ and
so $-x\leq \sup \left( -\mathbf{A}\right) $

Hence by statement (\textbf{i}) $x\geq -\sup \left( -\mathbf{A}\right) $ and
we see that $-\sup \left( -\mathbf{A}\right) $ is a lower

bound of $\mathbf{A.}$(\textbf{iii}) If $l_{\mathbf{A}}\leq x$ for all $x\in 
\mathbf{A}$ then $-l_{\mathbf{A}}\geq y$ for all $y\in -\mathbf{A.}$Hence

$-l_{\mathbf{A}}\geq \sup \left( -\mathbf{A}\right) $ by virtue of $\sup
\left( -\mathbf{A}\right) $ being the least upper bound of $-\mathbf{A.}$

Finally we obtain: $l_{\mathbf{A}}\leq -\sup \left( -\mathbf{A}\right) .$

\textbf{Definition 1.3.6.14. }We call this greatest element a \textit{%
greatest lower bound} $\ \ \ \ \ \ \ \ \ \ \ \ \ \ \ \ \ \ \ \ \ $

or \textit{infinum of }$\mathbf{A}$\textit{,}written is\textit{\ }$\inf
\left( \mathbf{A}\right) .$

\textbf{Definition 1.3.6.15.}Let us consider an nonempty subset $\mathbf{A}%
\subsetneqq $ $^{\ast }%
\mathbb{R}
_{\mathbf{d}}$\ which \ \ \ \ \ \ \ \ \ \ \ \ \ \ \ \ \ \ \ \ \ \ \ \ \ \ \
\ \ \ 

is bounded or hyperbounded from below and such that:

$\inf \left( \mathbf{A}\right) \pm \varepsilon \neq \inf \left( \mathbf{A}%
\right) $ for any $\varepsilon >0,$ $\varepsilon \approx 0.$

We call this greatest lower bound $\inf \left( \mathbf{A}\right) $ a \textit{%
strong greatest lower bound} or

\textit{strong infinum}, written is $\mathbf{s}$-$\inf \left( \mathbf{A}%
\right) $.

\textbf{Definition 1.3.6.16.}Let us consider an nonempty subset $\mathbf{A}%
\subsetneqq $ $^{\ast }%
\mathbb{R}
_{\mathbf{d}}$\ which \ \ \ \ \ \ \ \ \ \ \ \ \ \ \ \ \ \ \ \ \ \ \ \ \ \ \
\ \ \ 

is bounded or hyperbounded from below and such that:

(\textbf{1}) there exist $\varepsilon _{0}\approx 0$ such that $\inf \left( 
\mathbf{A}\right) \pm \varepsilon _{0}=\inf \left( \mathbf{A}\right) ,$

(\textbf{2}) $\inf \left( \mathbf{A}\right) \pm \varepsilon \neq \inf \left( 
\mathbf{A}\right) $ for any $\varepsilon >0$ such that $\varepsilon \geq
\varepsilon _{0}\approx 0.$

We call this greatest lower bound $\inf \left( \mathbf{A}\right) $ \textit{%
almost} \textit{strong greatest lower bound}

or \textit{almost} \textit{strong infinum}, written is $o\mathbf{s}$-$\inf
\left( \mathbf{A}\right) $.

\textbf{Definition 1.3.6.17. }Let us consider an nonempty subset $\mathbf{A}%
\subsetneqq $ $^{\ast }%
\mathbb{R}
_{\mathbf{d}}$\ which \ \ \ \ \ \ \ \ \ \ \ \ \ \ \ \ \ \ \ \ \ \ \ \ \ \ \
\ \ \ 

is bounded or hyperbounded from below and such that:

(\textbf{1}) $\inf \left( \mathbf{A}\right) \pm \varepsilon =\inf \left( 
\mathbf{A}\right) $ for any $\varepsilon >0,$ $\varepsilon \approx 0,$

(\textbf{2}) $\inf \left( \mathbf{A}\right) \pm \varepsilon \neq \inf \left( 
\mathbf{A}\right) $ for any $\varepsilon >0$ such that $\varepsilon
\not\approx 0.$

We call this greatest lower bound $\inf \left( \mathbf{A}\right) $ \textit{%
weak greatest lower bound}

or \textit{weak infinum}, written is $w$-$\inf \left( \mathbf{A}\right) $.

\textbf{Definition 1.3.6.18. }Let us consider an nonempty subset $\mathbf{A}%
\subsetneqq $ $^{\ast }%
\mathbb{R}
_{\mathbf{d}}$\ which \ \ \ \ \ \ \ \ \ \ \ \ \ \ \ \ \ \ \ \ \ \ \ \ \ \ \
\ \ \ 

is hyperbounded from below and such that:

(\textbf{1}) $\inf \left( \mathbf{A}\right) \pm \alpha =\inf \left( \mathbf{A%
}\right) $ for any $\alpha >0,$ $\alpha \in 
\mathbb{R}
,$

(\textbf{2}) $\inf \left( \mathbf{A}\right) \pm \Gamma \neq \inf \left( 
\mathbf{A}\right) $ for any $\Gamma >0$ such that $\Gamma \in $ $^{\ast }%
\mathbb{N}
_{\infty }.$

We call this greatest lower bound $\inf \left( \mathbf{A}\right) $ \textit{%
ultra} \textit{weak greatest lower bound}

or \textit{ultra} \textit{weak infinum}, written is $uw$-$\inf \left( 
\mathbf{A}\right) $.

\textbf{Proposition 1.3.6.2. }(\textbf{1}) If $\mathbf{A}$ is a nonempty
subset of $^{\ast }%
\mathbb{R}
_{\mathbf{d}}$ which is bounded \ \ \ \ \ \ \ \ \ \ \ \ \ \ \ \ \ \ \ \ \ \
\ \ \ \ 

from below and strong infinum $\mathbf{s}$-$\inf \left( \mathbf{A}\right) $
exist, then:

$\mathbf{s}$-$\inf \left( \mathbf{A}\right) $ is the unique number such that 
$\mathbf{s}$-$\inf \left( \mathbf{A}\right) $ is an upper bound \ \ \ \ \ \
\ \ \ \ \ \ \ \ \ \ \ \ \ \ \ \ \ \ \ \ 

for $\mathbf{A}$ and $\mathbf{s}$-$\inf \left( \mathbf{A}\right)
+\varepsilon $ is not a lower bound for $\mathbf{A}$ for any $\varepsilon
\approx 0,\varepsilon >0;$

(\textbf{2}) If $\mathbf{A}$ is a nonempty subset of $^{\ast }%
\mathbb{R}
_{\mathbf{d}}$ which is bounded from above, then:

$\mathbf{s}$-$\sup \left( \mathbf{A}\right) $ is the unique number such that 
$\mathbf{s}$-$\sup \left( \mathbf{A}\right) $ is an upper bound for $\mathbf{%
A}$

and $\mathbf{s}$-$\sup \left( \mathbf{A}\right) -\varepsilon $ is not an
upper bound for $\mathbf{A}$ for any $\varepsilon \approx 0,\varepsilon >0.$

\textbf{Proposition 1.3.6.3.(a). }(\textbf{Strong Approximation Property.})

(\textbf{1}) If $\mathbf{A}$ is a nonempty subset of $^{\ast }%
\mathbb{R}
_{\mathbf{d}}$ which is bounded (hyperbounded)

from above and such that strong supremum $\mathbf{s}$-$\sup \left( \mathbf{A}%
\right) $ exist, and let

$\varepsilon \approx 0,\varepsilon >0$ there exist $x\in \mathbf{A}$ such
that $\mathbf{s}$-$\sup \left( \mathbf{A}\right) -\varepsilon <x\leq \mathbf{%
s}$-$\sup \left( \mathbf{A}\right) .$

(\textbf{2}) If $\mathbf{A}$ is a nonempty subset of $^{\ast }%
\mathbb{R}
_{\mathbf{d}}$ which is bounded (hyperbounded)

from below and such that strong infinum $\mathbf{s}$-$\inf \left( \mathbf{A}%
\right) $ exist, and let

$\varepsilon \approx 0,\varepsilon >0$ there exist $x\in \mathbf{A}$ such
that $\mathbf{s}$-$\inf \left( \mathbf{A}\right) \leq x<\mathbf{s}$-$\inf
\left( \mathbf{A}\right) +\varepsilon .$

\textbf{Proof.} (\textbf{1}) If not, then $\mathbf{s}$-$\sup \left( \mathbf{A%
}\right) -\varepsilon $ is an upper bound of $\mathbf{A}$ less than the

strong\textit{\ }upper bound $\mathbf{s}$-$\sup \left( \mathbf{A}\right) $,
which is a contradiction.

(\textbf{2}) If not, then $\mathbf{s}$-$\inf \left( \mathbf{A}\right)
+\varepsilon $ is an lower bound of $\mathbf{A}$ bigger than the

strong lower bound $\mathbf{s}$-$\inf \left( \mathbf{A}\right) $, which is a
contradiction.

\textbf{Proposition 1.3.6.3.(b) (The Almost Strong Approximation Property.) }

(\textbf{1}) If $\mathbf{A}$ is a nonempty subset of $^{\ast }%
\mathbb{R}
_{\mathbf{d}}$ which is bounded (hyperbounded)

from above and such that almost strong supremum $o\mathbf{s}$-$\sup \left( 
\mathbf{A}\right) $ exist, and

let $\varepsilon \approx 0,\varepsilon >0,o\mathbf{s}$-$\sup \left( \mathbf{A%
}\right) \pm \varepsilon \neq o\mathbf{s}$-$\sup \left( \mathbf{A}\right) $
there exist $x\in \mathbf{A}$ such that

$o\mathbf{s}$-$\sup \left( \mathbf{A}\right) -\varepsilon <x\leq o\mathbf{s}$%
-$\sup \left( \mathbf{A}\right) .$

(\textbf{2}) If $\mathbf{A}$ is a nonempty subset of $^{\ast }%
\mathbb{R}
_{\mathbf{d}}$ which is bounded (hyperbounded)

from below and such that almost strong infinum $o\mathbf{s}$-$\inf \left( 
\mathbf{A}\right) $ exist, and let

$\varepsilon \approx 0,\varepsilon >0,o\mathbf{s}$-$\inf \left( \mathbf{A}%
\right) \pm \varepsilon \neq o\mathbf{s}$-$\inf \left( \mathbf{A}\right) $
there exist $x\in \mathbf{A}$ such that

$o\mathbf{s}$-$\inf \left( \mathbf{A}\right) \leq x<o\mathbf{s}$-$\inf
\left( \mathbf{A}\right) +\varepsilon .$

\textbf{Proof.} (\textbf{1}) If not, then $o\mathbf{s}$-$\sup \left( \mathbf{%
A}\right) -\varepsilon $ is an upper bound of $\mathbf{A}$ less than the

almost strong upper bound $o\mathbf{s}$-$\sup \left( \mathbf{A}\right) $,
which is a contradiction.

(\textbf{2}) If not, then $o\mathbf{s}$-$\inf \left( \mathbf{A}\right)
+\varepsilon $ is an lower bound of $\mathbf{A}$ bigger than the

almost strong \textit{\ }lower bound $o\mathbf{s}$-$\inf \left( \mathbf{A}%
\right) $, which is a contradiction.

\textbf{Proposition 1.3.6.3.(c) (The Weak Approximation Property.) }

(\textbf{1}) If $\mathbf{A}$ is a nonempty subset of $^{\ast }%
\mathbb{R}
_{\mathbf{d}}$ which is bounded (hyperbounded)

from above and such that weak supremum $w$-$\sup \left( \mathbf{A}\right) $
exist, and let

$\varepsilon \in 
\mathbb{R}
,\varepsilon >0,w$-$\sup \left( \mathbf{A}\right) \pm \varepsilon \neq w$-$%
\sup \left( \mathbf{A}\right) $ there exist $x\in \mathbf{A}$ such that

$w$-$\sup \left( \mathbf{A}\right) -\varepsilon <x\leq w$-$\sup \left( 
\mathbf{A}\right) .$

(\textbf{2}) If $\mathbf{A}$ is a nonempty subset of $^{\ast }%
\mathbb{R}
_{\mathbf{d}}$ which is bounded (hyperbounded)

from below and such that weak infinum $w$-$\inf \left( \mathbf{A}\right) $
exist, and let

$\varepsilon \in 
\mathbb{R}
,\varepsilon >0,w$-$\inf \left( \mathbf{A}\right) \pm \varepsilon \neq w$-$%
\inf \left( \mathbf{A}\right) $ there exist $x\in \mathbf{A}$ such that

$w$-$\inf \left( \mathbf{A}\right) \leq x<w$-$\inf \left( \mathbf{A}\right)
+\varepsilon .$

\textbf{Proof.} (\textbf{1}) If not, then $w$-$\sup \left( \mathbf{A}\right)
-\varepsilon $ is an upper bound of $\mathbf{A}$ less than the

weak upper bound $w$-$\sup \left( \mathbf{A}\right) $, which is a
contradiction.

(\textbf{2}) If not, then

\bigskip

\textbf{Corollary} \textbf{1.3.6.2. }(\textbf{1}) Let $\mathbf{A}$ be
bounded or hyperbounded from below and \ \ \ \ \ \ \ \ \ \ \ \ \ \ \ \ \ \ \
\ \ \ 

non-empty set such that $\mathbf{s}$-$\inf \left( \mathbf{A}\right) $ exist.
There is a function \textbf{\ }

$\beta \left( \circ \right) :$ $^{\ast }\mathbf{%
\mathbb{N}
}_{\infty }\rightarrow $ $^{\ast }%
\mathbb{R}
_{\mathbf{d}}$ such that for all $\mathbf{n\in }$ $^{\ast }\mathbf{%
\mathbb{N}
}_{\infty }$ we have

$\mathbf{s}$-$\inf \left( \mathbf{A}\right) \leq \beta \left( \mathbf{n}%
\right) <\mathbf{s}$-$\inf \left( \mathbf{A}\right) +\mathbf{n}^{-1}.$

(\textbf{2}) Let $\mathbf{A}$ be bounded or hyperbounded from above and \ \
\ \ \ \ \ \ \ \ \ \ \ \ \ \ \ \ \ \ \ \ 

non-empty set such that $\mathbf{s}$-$\sup \left( \mathbf{A}\right) $ exist.
There is a function \textbf{\ }

$\alpha \left( \circ \right) :$ $^{\ast }\mathbf{%
\mathbb{N}
}_{\infty }\rightarrow $ $^{\ast }%
\mathbb{R}
_{\mathbf{d}}$ such that for all $\mathbf{n\in }$ $^{\ast }\mathbf{%
\mathbb{N}
}_{\infty }$ we have

$\mathbf{s}$-$\sup \left( \mathbf{A}\right) -\mathbf{n}^{-1}<\alpha \left( 
\mathbf{n}\right) \leq \mathbf{s}$-$\sup \left( \mathbf{A}\right) .$

\textbf{Example 1.3.6.6. }(\textbf{a})\textbf{\ }The subset $\left\{ \mathbf{%
n}^{-1}\right\} _{\mathbf{n\in }^{\ast }\mathbf{%
\mathbb{N}
}_{\infty }}\triangleq \left\{ \left. \dfrac{1}{\mathbf{n}}\right\vert 
\mathbf{n\in }\text{ }^{\ast }\mathbf{%
\mathbb{N}
}_{\infty }\right\} \subsetneqq $ $^{\ast }%
\mathbb{R}
_{\mathbf{d}}$

has a \textit{strong greatest lower bound} $\mathbf{s}$-$\inf \left( \left\{ 
\mathbf{n}^{-1}\right\} _{\mathbf{n\in }^{\ast }\mathbf{%
\mathbb{N}
}_{\infty }}\right) =0$ in $^{\ast }%
\mathbb{R}
_{\mathbf{d}}.$

(\textbf{b}) The subset $\left\{ n^{-1}\right\} _{\mathbf{n\in 
\mathbb{N}
}}\triangleq \left\{ \left. \dfrac{1}{n}\right\vert n\mathbf{\in }\text{ }%
\mathbf{%
\mathbb{N}
}\right\} \subsetneqq $ $^{\ast }%
\mathbb{R}
_{\mathbf{d}}$ has a \textit{greatest} \textit{lower bound }

$\inf \left( \left\{ n^{-1}\right\} _{n\mathbf{\in 
\mathbb{N}
}}\right) $ in $^{\ast }%
\mathbb{R}
_{\mathbf{d}}$ but has not a strong greatest lower bound in $^{\ast }%
\mathbb{R}
_{\mathbf{d}}.$

\textbf{Example 1.3.6.7. }The subset $%
\mathbb{R}
\subsetneqq $ $^{\ast }%
\mathbb{R}
_{\mathbf{d}}$ has the \textit{least upper bound} $\sup \left( 
\mathbb{R}
\right) $

in $^{\ast }%
\mathbb{R}
_{\mathbf{d}}$ but has not \textit{strong least} \textit{upper bound} in $%
^{\ast }%
\mathbb{R}
_{\mathbf{d}}.$

\textbf{Example 1.3.6.8. }The subset $\mathbf{I}_{\ast }$ the all
infinitesimal members of the $^{\ast }%
\mathbb{R}
,$

$\mathbf{I}_{\ast }\subsetneqq $ $^{\ast }%
\mathbb{R}
_{\mathbf{d}}$ has least upper bound $\sup \left( \mathbf{I}_{\ast }\right) $
in $^{\ast }%
\mathbb{R}
_{\mathbf{d}}$ but has not \textit{strong least}

\textit{upper bound} \ in $^{\ast }%
\mathbb{R}
_{\mathbf{d}}.$

\textbf{Example 1.3.6.9. }The subset $%
\mathbb{R}
_{+}\subsetneqq $ $^{\ast }%
\mathbb{R}
_{\mathbf{d}}$ has lower bound $\inf \left( 
\mathbb{R}
_{+}\right) $ \ \ \ \ \ \ \ \ \ \ \ \ \ \ \ \ \ \ \ \ \ \ 

in $^{\ast }%
\mathbb{R}
_{\mathbf{d}}$ but has not \textit{strong} \textit{lower bound} in $^{\ast }%
\mathbb{R}
_{\mathbf{d}}.$

\textbf{Example 1.3.6.10. }The subset $^{\ast }%
\mathbb{R}
_{+}\subsetneqq $ $^{\ast }%
\mathbb{R}
_{\mathbf{d}}$ has lower bound $\inf \left( ^{\ast }%
\mathbb{R}
_{+}\right) $

in $^{\ast }%
\mathbb{R}
_{\mathbf{d}}$ but \ has not \textit{strong} \textit{lower bound} in $^{\ast
}%
\mathbb{R}
_{\mathbf{d}}.$

\textbf{Proposition 1.3.6.4.}Let $\mathbf{A}$ and $\mathbf{B}$ be nonempty
subsets of $^{\ast }%
\mathbb{R}
_{\mathbf{d}}$

\textbf{Theorem 1.3.6.3.A. }Let $\mathbf{A}$ and $\mathbf{B}$ be nonempty
subsets of $^{\ast }%
\mathbb{R}
_{\mathbf{d}}$ and

$\mathbf{C}=$ $\left\{ a+b:a\in \mathbf{A},b\in \mathbf{B}\right\} $.

(\textbf{1.a}) If $\mathbf{A}$ and $\mathbf{B}$ are bounded or hyperbounded
from above

(hence $\mathbf{s}$-$\sup \left( \mathbf{A}\right) $ and $\mathbf{s}$-$\sup
\left( \mathbf{B}\right) $ exist) then $\mathbf{s}$-$\sup \left( \mathbf{C}%
\right) $ exist

and

\bigskip

$\ 
\begin{array}{cc}
\begin{array}{c}
\\ 
\mathbf{s}\text{-}\sup \left( \mathbf{C}\right) =\mathbf{s}\text{-}\sup
\left( \mathbf{A}\right) +\mathbf{s}\text{-}\sup \left( \mathbf{B}\right) \\ 
\end{array}
& \text{ \ }\left( 1.3.6.21.\mathbf{a}\right)%
\end{array}%
$

\bigskip

(\textbf{2.a}) If $\mathbf{A}$ and $\mathbf{B}$ are bounded or hyperbounded
from below

(hence $\mathbf{s}$-$\inf \left( \mathbf{A}\right) $ and $\mathbf{s}$-$\inf
\left( \mathbf{B}\right) $ exist) then $\mathbf{s}$-$\inf \left( \mathbf{C}%
\right) $ exist and

\bigskip

\bigskip $\ 
\begin{array}{cc}
\begin{array}{c}
\\ 
\mathbf{s}\text{-}\inf \left( \mathbf{C}\right) =\mathbf{s}\text{-}\inf
\left( \mathbf{A}\right) +\mathbf{s}\text{-}\inf \left( \mathbf{B}\right) \\ 
\end{array}
& \text{ \ }\left( 1.3.6.22.\mathbf{a}\right)%
\end{array}%
$

\bigskip

(\textbf{1.b}) If $\mathbf{A}$ and $\mathbf{B}$ are bounded or hyperbounded
from above,

(hence $o\mathbf{s}$-$\sup \left( \mathbf{A}\right) $ and $o\mathbf{s}$-$%
\sup \left( \mathbf{B}\right) $ exist) then $o\mathbf{s}$-$\sup \left( 
\mathbf{C}\right) $ exist

and $\ \ \ \ \ \ \ \ \ \ \ \ \ \ \ \ \ \ \ \ \ \ \ \ \ \ \ \ \ \ \ \ \ \ \ \
\ \ \ \ \ $

$\bigskip $

$\ \ 
\begin{array}{cc}
\begin{array}{c}
\\ 
o\mathbf{s}\text{-}\sup \left( \mathbf{C}\right) =o\mathbf{s}\text{-}\sup
\left( \mathbf{A}\right) +o\mathbf{s}\text{-}\sup \left( \mathbf{B}\right)
\\ 
\end{array}
& \text{ \ }\left( 1.3.6.21.\mathbf{b}\right)%
\end{array}%
$

\bigskip

(\textbf{2.b}) If $\mathbf{A}$ and $\mathbf{B}$ are bounded or hyperbounded
from below

(hence $o\mathbf{s}$-$\inf \left( \mathbf{A}\right) $ and $o\mathbf{s}$-$%
\inf \left( \mathbf{B}\right) $ exist) then $o\mathbf{s}$-$\inf \left( 
\mathbf{C}\right) $ exist

and

\bigskip

$%
\begin{array}{cc}
\begin{array}{c}
\\ 
o\mathbf{s}\text{-}\inf \left( \mathbf{C}\right) =o\mathbf{s}\text{-}\inf
\left( \mathbf{A}\right) +o\mathbf{s}\text{-}\inf \left( \mathbf{B}\right)
\\ 
\end{array}
& \text{ \ }\left( 1.3.6.22.\mathbf{b}\right)%
\end{array}%
$\bigskip

(\textbf{1.c}) If $\mathbf{A}$ and $\mathbf{B}$ are bounded or hyperbounded
from above,

hence and $w$-$\sup \left( \mathbf{A}\right) $ and $w$-$\sup \left( \mathbf{B%
}\right) $ exist, then $w$-$\sup \left( \mathbf{C}\right) $ exist

and

$\bigskip \ 
\begin{array}{cc}
\begin{array}{c}
\\ 
w\text{-}\sup \left( \mathbf{C}\right) =w\text{-}\sup \left( \mathbf{A}%
\right) +w\text{-}\sup \left( \mathbf{B}\right) \\ 
\end{array}
& \text{ \ }\left( 1.3.6.21.\mathbf{c}\right)%
\end{array}%
$

(\textbf{2.c}) If $\mathbf{A}$ and $\mathbf{B}$ are bounded or hyperbounded
from below

(hence $w$-$\inf \left( \mathbf{A}\right) $ and $w$-$\inf \left( \mathbf{B}%
\right) $ exist) then $w$-$\inf \left( \mathbf{C}\right) $ exist and

\bigskip

$\ \ 
\begin{array}{cc}
\begin{array}{c}
\\ 
w\text{-}\inf \left( \mathbf{C}\right) =w\text{-}\inf \left( \mathbf{A}%
\right) +w\text{-}\inf \left( \mathbf{B}\right) \\ 
\end{array}
& \text{ \ }\left( 1.3.6.22.\mathbf{c}\right)%
\end{array}%
$

\bigskip

(\textbf{1.d}) If $\mathbf{A}$ and $\mathbf{B}$ are bounded or hyperbounded
from above

(then $uw$-$\sup \left( \mathbf{A}\right) $ and $uw$-$\sup \left( \mathbf{B}%
\right) $ exist) then $uw$-$\sup \left( \mathbf{C}\right) $ exist

and $\ \ \ \ \ \ \ \ \ \ \ \ \ \ \ \ \ \ \ \ \ \ \ \ \ \ \ \ \ \ \ \ \ \ \ \
\ \ \ \ $

$\bigskip $

$\ \ 
\begin{array}{cc}
\begin{array}{c}
\\ 
uw\text{-}\sup \left( \mathbf{C}\right) =uw\text{-}\sup \left( \mathbf{A}%
\right) +uw\text{-}\sup \left( \mathbf{B}\right) \\ 
\end{array}
& \text{\ }\left( 1.3.6.21.\mathbf{d}\right)%
\end{array}%
$

\bigskip

(\textbf{2.d}) If $\mathbf{A}$ and $\mathbf{B}$ are bounded or hyperbounded
from below

(hence $uw$-$\inf \left( \mathbf{A}\right) $ and $uw$-$\inf \left( \mathbf{B}%
\right) $ exist) then $uw$-$\inf \left( \mathbf{C}\right) $ exist

and $\ \ \ \ \ \ \ \ \ \ \ \ \ \ \ \ \ \ \ \ \ \ \ \ \ \ \ \ \ \ \ \ \ \ \ \
\ \ \ \ \ \ $

$\bigskip $

$\ \ 
\begin{array}{cc}
\begin{array}{c}
\\ 
uw\text{-}\inf \left( \mathbf{C}\right) =uw\text{-}\inf \left( \mathbf{A}%
\right) +uw\text{-}\inf \left( \mathbf{B}\right) \\ 
\end{array}
& \text{ \ \ }\left( 1.3.6.22.\mathbf{d}\right)%
\end{array}%
$

\bigskip

\textbf{Proof. }(\textbf{1.a}) Suppose that $\mathbf{A}$ and $\mathbf{B}$
are bounded or hyperbounded from \ \ \ \ \ \ \ \ \ \ \ \ \ \ \ \ \ \ \ \ \ \
\ \ \ 

above,hence $\mathbf{s}$-$\sup (\mathbf{A})$ and $\mathbf{s}$-$\sup (\mathbf{%
B})$ exist. Let $c$ $\in $ $\mathbf{C}$. Then $c=a+b$ for

numbers $a\in $ $\mathbf{A}$ and $b$ $\in $ $\mathbf{B}$. Since $a\leq 
\mathbf{s}$-$\sup (\mathbf{A})$ and $b\leq \mathbf{s}$-$\sup (\mathbf{B})$,

$c=a+b\leq \mathbf{s}$-$\sup (\mathbf{A})$ $+$ $\mathbf{s}$-$\sup (\mathbf{B}%
)$. This shows that $\mathbf{s}$-$\sup (\mathbf{A})+\mathbf{s}$-$\sup (%
\mathbf{B})$ is an

upper bound for $\mathbf{C}$, in particular, $\mathbf{C}$ is bounded or
hyperbounded from above.

Given $\varepsilon >0,$ $\mathbf{s}$-$\sup (\mathbf{A})-\varepsilon /2$ is
not an a strong upper bound for $\mathbf{A}$ hence there \ \ \ \ \ \ 

exists $a^{\prime }\in \mathbf{A}$ such that $\mathbf{s}$-$\sup (\mathbf{A}%
)-\varepsilon /2<a^{\prime }.$ Similarly, $\mathbf{s}$-$\sup (\mathbf{B})$ $%
-\varepsilon /2$ is not an \ \ \ \ \ \ \ \ \ \ \ \ \ \ \ \ 

upper bound for $\mathbf{B}$ and there exists $b^{\prime }$ $\in $ $\mathbf{B%
}$ such that $\mathbf{s}$-$\sup (\mathbf{B})-\varepsilon /2<b^{\prime }.$ So

for $c^{\prime }=a^{\prime }+b^{\prime }$ $\in $ $\mathbf{C}$ we have $%
\mathbf{s}$-$\sup (\mathbf{A})$ $+$ $\mathbf{s}$-$\sup (\mathbf{B}%
)-\varepsilon <c^{\prime }.$ This shows that

$\mathbf{s}$-$\sup (\mathbf{A})$ $+$ $\mathbf{s}$-$\sup (\mathbf{B}%
)-\varepsilon $ is not an a strong upper bound for $\mathbf{C}$ for any $%
\varepsilon >0.$

Hence by the Proposition 1.3.1 one obtain: $\mathbf{s}$-$\sup (\mathbf{C})=$ 
$\mathbf{s}$-$\sup (\mathbf{A})$ $+$ $\mathbf{s}$-$\sup (\mathbf{B})$.

By using \textbf{Theorem 1.3.6.3.A }one obtain:

\textbf{Theorem 1.3.6.3.B. }Let $\mathbf{A}$ and $\mathbf{B}$ be nonempty
subsets of $^{\ast }%
\mathbb{R}
_{\mathbf{d}}$ and

$\mathbf{C}=$ $\left\{ a+b:a\in \mathbf{A},b\in \mathbf{B}\right\} $.

(\textbf{1}) If $\mathbf{A}$ and $\mathbf{B}$ are bounded or hyperbounded
from above \ \ \ \ \ \ \ \ \ \ \ \ \ \ \ \ \ \ \ \ \ \ \ \ \ \ \ \ \ \ \ \ \
\ \ \ \ \ \ \ \ \ \ \ \ \ \ \ \ \ 

(hence $\sup \left( \mathbf{A}\right) $ and $\sup \left( \mathbf{B}\right) $
exist) then $\sup \left( \mathbf{C}\right) $ exist and

\bigskip

$\ \ 
\begin{array}{cc}
\begin{array}{c}
\\ 
\sup \left( \mathbf{C}\right) =\sup \left( \mathbf{A}\right) +\sup \left( 
\mathbf{B}\right) \\ 
\end{array}
& \text{ \ \ }\left( 1.3.6.23\right) \text{\ }%
\end{array}%
$

\bigskip

(\textbf{2}) If $\mathbf{A}$ and $\mathbf{B}$ are bounded or hyperbounded
from below \ \ \ \ \ \ \ \ \ \ \ \ \ \ \ \ \ \ \ \ \ \ \ \ \ \ \ \ \ \ \ \ \
\ \ \ \ \ \ \ \ \ 

(hence $\inf \left( \mathbf{A}\right) $ and $\inf \left( \mathbf{B}\right) $
exist) then $\inf \left( \mathbf{C}\right) $ exist and $\ \ \ \ \ \ \ \ \ \
\ \ \ \ \ \ \ \ \ \ \ \ \ \ \ \ \ \ \ \ \ \ \ \ \ \ \ \ \ \ \ \ \ \ \ \ \ $

$\bigskip $

$\ 
\begin{array}{cc}
\begin{array}{c}
\\ 
\inf \left( \mathbf{C}\right) =\inf \left( \mathbf{A}\right) +\inf \left( 
\mathbf{B}\right) . \\ 
\end{array}
& \text{ \ \ }\left( 1.3.6.24\right) \text{\ \ }%
\end{array}%
$

\bigskip

\textbf{Theorem 1.3.6.3.C. }Setting (\textbf{1}). Suppose that $\mathbf{S}$
is a non-empty subset \ \ \ \ \ \ \ \ \ \ \ \ \ \ \ \ \ \ \ \ \ \ \ \ \ \ \
\ \ \ \ \ \ \ \ \ 

of $^{\ast }%
\mathbb{R}
_{\mathbf{d}}$ which is bounded or hyperbounded from above and $\mathbf{s}$-$%
\sup \mathbf{S}$ exist

and $\xi \in $ $^{\ast }%
\mathbb{R}
,\xi >0.$Then

\bigskip $\ \ 
\begin{array}{cc}
\begin{array}{c}
\\ 
\mathbf{s}\text{-}\underset{x\in \mathbf{S}}{\sup }\left\{ \xi \times
x\right\} = \\ 
\\ 
\xi \times \left( \mathbf{s}\text{-}\underset{x\in \mathbf{S}}{\sup }\left\{
x\right\} \right) =\xi \times \left( \mathbf{s}\text{-}\sup \mathbf{S}%
\right) \mathbf{.} \\ 
\end{array}
& \text{ \ }\left( 1.3.6.25\right) \text{\ }%
\end{array}%
$

Setting (\textbf{2}).Suppose that $\mathbf{S}$ is a non-empty subset of $%
^{\ast }%
\mathbb{R}
_{\mathbf{d}}$ which is

bounded or hyperbounded from above and $o\mathbf{s}$-$\sup \mathbf{S}$ exist
and $\xi \in $ $^{\ast }%
\mathbb{R}
,$

$\xi >0.$Then

\bigskip

$\ \ \ \ 
\begin{array}{cc}
\begin{array}{c}
\\ 
o\mathbf{s}\text{-}\underset{x\in \mathbf{S}}{\sup }\left\{ \xi \times
x\right\} = \\ 
\\ 
\xi \times \left( o\mathbf{s}\text{-}\underset{x\in \mathbf{S}}{\sup }%
\left\{ x\right\} \right) =\xi \times \left( o\mathbf{s}\text{-}\sup \mathbf{%
S}\right) \mathbf{.} \\ 
\end{array}
& \text{ \ }\left( 1.3.6.26\right) \text{\ }%
\end{array}%
$

\bigskip

\textbf{Proof.} (\textbf{1}) Let $B=\mathbf{s}$-$\sup \mathbf{S.}$Then $B$
is the smallest number such that, for \ \ \ \ \ \ \ \ \ \ \ \ \ \ \ \ \ \ \
\ \ \ \ \ 

any $x\in \mathbf{S,}x$ $\mathbf{\leq }B\mathbf{.}$Let $\mathbf{T}=\left\{
\xi \times x|x\in \mathbf{S}\right\} .$Since $\xi >0,\xi \times x\leq \xi
\times B$ for any

$x\in \mathbf{S.}$Hence $\mathbf{T}$ is bounded or hyperbounded above by $%
\xi \times B.$By the \ \ \ \ \ \ \ \ \ \ \ \ \ \ \ \ \ \ \ \ \ \ \ \ 

Theorem 1 and setting (1), $\mathbf{T}$ has a strong supremum $C,C=\mathbf{s}
$-$\sup \mathbf{T.}$ \ \ \ \ \ \ \ \ \ \ \ \ \ \ \ \ \ \ \ \ \ \ \ \ \ \ \ \ 

Now\ we have to pruve that $C=\xi \times B.$Since $\xi \times B$ is an apper
bound for $\mathbf{T}$

and $C$ is the smollest apper bound for $\mathbf{T,}C\leq \xi \times B.$Now
we repeat the \ \ \ \ \ \ \ \ \ \ \ \ \ \ \ \ \ \ \ \ 

argument above with the roles of $\mathbf{S}$ and $\mathbf{T}$ reversed. We
know that $C$ is the

smallest number such that, for any $y\in \mathbf{T,}y\leq C.$Since $\xi >0$
it follows that

$\xi ^{-1}\times y\leq \xi ^{-1}\times C$ for any $y\in \mathbf{T.}$But $%
\mathbf{S=}\left\{ \xi ^{-1}\times y|y\in \mathbf{T}\right\} .$Hence $\xi
^{-1}\times C$ is an

apper bound for $\mathbf{S.}$But $B$ is a strong supremum for $\mathbf{S.}$%
Hence $B\leq \xi ^{-1}\times C$

and $\xi \times B\leq C.$We have shown that $C\leq \xi \times B$ and also
that $\xi \times B\leq C.$Thus

$\xi \times B=C.$

\bigskip

\textbf{Theorem 1.3.6.3.D. }Let $\mathbf{A}$ and $\mathbf{B}$ be nonempty
subsets of $^{\ast }%
\mathbb{R}
_{\mathbf{d}}$ such that \ \ \ \ \ \ \ \ \ \ \ \ \ \ \ \ \ \ \ \ \ \ \ \ \ \ 

$0\leq \mathbf{A,}0\leq \mathbf{B}$ and $\mathbf{C}=$ $\left\{ a\times
b:a\in \mathbf{A},b\in \mathbf{B}\right\} $.

(\textbf{1.a}) If $\mathbf{A}$ and $\mathbf{B}$ are bounded or hyperbounded
from above,hence $\mathbf{s}$-$\sup \left( \mathbf{A}\right) $ \ \ \ \ \ \ \
\ \ \ \ \ \ \ \ \ \ \ \ \ 

and $\mathbf{s}$-$\sup \left( \mathbf{B}\right) $ exist, then $\mathbf{s}$-$%
\sup \left( \mathbf{C}\right) $ exist and

\bigskip $\ 
\begin{array}{cc}
\begin{array}{c}
\\ 
\mathbf{s}\text{-}\sup \left( \mathbf{C}\right) =\left[ \mathbf{s}\text{-}%
\sup \left( \mathbf{A}\right) \right] \times \left[ \mathbf{s}\text{-}\sup
\left( \mathbf{B}\right) \right]  \\ 
\end{array}
& \text{ \ \ }\left( 1.3.6.21^{\prime }.\mathbf{a}\right) 
\end{array}%
$

\textbf{Proposition 1.3.6.5. }Let $\mathbf{A}$ and $\mathbf{B}$ be nonempty
subsets of $^{\ast }%
\mathbb{R}
_{\mathbf{d}}$.

(\textbf{i}) If for every $a\in \mathbf{A}$ there exists $b\in \mathbf{B}$
with $a\leq b$ and $\mathbf{B}$ is bounded \ \ \ \ \ \ \ \ \ \ \ \ \ \ \ \ \
\ \ \ \ \ \ \ \ \ \ \ \ \ \ \ 

from above,then so is $\mathbf{A}$ and $\sup \left( \mathbf{A}\right) \leq
\sup \left( \mathbf{B}\right) .$

(\textbf{ii}) If for every $b\in \mathbf{B}$ there exists $a\in \mathbf{A}$
with $a\leq b$ and $\mathbf{A}$ is bounded \ \ \ \ \ \ \ \ \ \ \ \ \ \ \ \ \
\ \ \ \ \ \ \ \ \ \ \ \ \ \ \ 

from below,then so is $\mathbf{B}$ and $\inf \left( \mathbf{A}\right) \leq
\inf \left( \mathbf{B}\right) .$

\textbf{Proof. }(\textbf{ii})\textbf{\ }Suppose that for every $b\in \mathbf{%
B}$ there exists $a\in \mathbf{A}$ with $a\leq b$ \ \ \ \ \ \ \ \ \ \ \ \ \
\ \ \ \ \ \ \ \ \ \ \ \ \ \ \ 

and $\mathbf{A}$ is bounded from below. Then $\inf (\mathbf{A})$ exists. For
every $b\in \mathbf{B}$ there \ \ \ \ \ \ \ \ \ \ \ \ \ \ \ \ \ \ \ \ \ \ \
\ \ \ \ \ \ 

is a $\in $ $\mathbf{A}$ such that $a\leq b.$ So $\inf \left( \mathbf{A}%
\right) $ $\leq a\leq b.$Therefore $\inf \left( \mathbf{A}\right) $ is a
lower \ \ \ \ \ \ \ \ \ \ \ \ \ \ \ \ \ \ \ \ \ \ \ \ \ \ \ \ 

bound for $\mathbf{B.}$Hence $\mathbf{B}$ is bounded from below and $\inf
\left( \mathbf{B}\right) $ exists. By \ \ \ \ \ \ \ \ \ \ \ \ \ \ \ \ \ \ \
\ \ \ \ \ \ \ \ \ \ 

definition of the infimum (greatest lower bound) one obtain:

$\inf \left( \mathbf{B}\right) \geq \inf \left( \mathbf{A}\right) .$

\textbf{Lemma} \textbf{1.3.6.1. }(\textbf{a})\textbf{\ }$\mathbf{s}$\textbf{-%
}$\inf \left( ^{\ast }%
\mathbb{Z}
\right) $ and $\mathbf{s}$\textbf{-}$\sup \left( ^{\ast }%
\mathbb{Z}
\right) $ is not exists in $^{\ast }%
\mathbb{R}
_{\mathbf{d}}$.

(\textbf{b}) $\mathbf{s}$\textbf{-}$\sup \left( ^{\ast }%
\mathbb{N}
\right) $ is not exists in $^{\ast }%
\mathbb{R}
_{\mathbf{d}}$.

(\textbf{c}) $^{\ast }%
\mathbb{Z}
$ is bounded neither from below nor from above in $^{\ast }%
\mathbb{R}
_{\mathbf{d}}$.

(\textbf{d}) $^{\ast }%
\mathbb{N}
$ is not bounded from above in $^{\ast }%
\mathbb{R}
_{\mathbf{d}}$.

\textbf{Proof. }(\textbf{a})\textbf{\ }Assume that $\mathbf{s}$\textbf{-}$%
\sup \left( ^{\ast }%
\mathbb{Z}
\right) $ exists in $^{\ast }%
\mathbb{R}
_{\mathbf{d}}$.Then $\mathbf{s}$\textbf{-}$\sup \left( ^{\ast }%
\mathbb{Z}
\right) -1$

is not an upper bound and hence there exists $n\in $ $^{\ast }%
\mathbb{Z}
$ such that

$\mathbf{s}$\textbf{-}$\sup \left( ^{\ast }%
\mathbb{Z}
\right) -1<n$ hence $\mathbf{s}$\textbf{-}$\sup \left( ^{\ast }%
\mathbb{Z}
\right) <n+1.$But since $n+1\in $ $^{\ast }%
\mathbb{Z}
$

this is a contradiction. Therefore $\mathbf{s}$\textbf{-}$\sup \left( ^{\ast
}%
\mathbb{Z}
\right) $ is not exists in $^{\ast }%
\mathbb{R}
_{\mathbf{d}}$.

(\textbf{d}) Assume that $^{\ast }%
\mathbb{N}
$ has an upper bound, call it $\mathbf{\Theta .}$Hence $\mathbf{\Theta }%
^{-1} $

is a lower bound for the set $\left\{ \mathbf{n}^{-1}\right\} _{\mathbf{n\in 
}^{\ast }\mathbf{%
\mathbb{N}
}_{\infty }}$ and consequently

$\inf \left( \left\{ \mathbf{n}^{-1}\right\} _{\mathbf{n\in }^{\ast }\mathbf{%
\mathbb{N}
}_{\infty }}\right) \geq \mathbf{\Theta }^{-1}\neq 0.$But we know that $\ast 
$-$\lim_{\mathbf{n}\rightarrow \text{ }^{\ast }\infty }\mathbf{n}^{-1}=0$

which is a contradiction.

\textbf{Theorem 1.3.6.4. (Generalyzed Archimedean Property of }$^{\ast }%
\mathbb{R}
_{\mathbf{d}}$\textbf{). }

For any $\varepsilon \in $ $^{\ast }%
\mathbb{R}
_{\mathbf{d}},\varepsilon \approx 0,\varepsilon >0$ there exists $\mathbf{%
n\in }^{\ast }%
\mathbb{N}
$ such that $\mathbf{n}^{-1}<\varepsilon .$

\textbf{Proof. }Since $^{\ast }%
\mathbb{N}
$ is not hyperbounded from above $\varepsilon ^{-1}$ is not an upper \ \ \ \
\ \ \ \ \ \ \ \ \ \ \ \ \ \ \ \ \ \ \ \ \ \ \ \ \ \ \ 

bound for $^{\ast }%
\mathbb{N}
.$ Hence there exists $\mathbf{n}\in $ $^{\ast }%
\mathbb{N}
_{\infty }$ such that $\mathbf{n>}$ $\varepsilon ^{-1}$ and \ \ \ \ \ \ \ \
\ \ \ \ \ \ \ \ \ \ \ \ \ \ \ \ \ \ \ 

consequently $\mathbf{n}^{-1}<\varepsilon .$

\textbf{Theorem 1.3.6.5. }For every $x\in $ $^{\ast }%
\mathbb{R}
_{\mathbf{d}}$ such that for the set

$\left\{ n\in \text{ }^{\ast }%
\mathbb{Z}
|n\leq x\right\} $one of the next conditions is satisfied:

(\textbf{i}) strong supremum $\mathbf{s}$-$\sup \left( \left\{ n\in \text{ }%
^{\ast }%
\mathbb{Z}
|n\leq x\right\} \right) $ exists in $^{\ast }%
\mathbb{R}
_{\mathbf{d}}$ or

(\textbf{ii}) almost strong supremum $\mathbf{os}$-$\sup \left( \left\{ n\in 
\text{ }^{\ast }%
\mathbb{Z}
|n\leq x\right\} \right) $ exists in $^{\ast }%
\mathbb{R}
_{\mathbf{d}}$ or

(\textbf{iii}) weak supremum $w$-$\sup \left( \left\{ n\in \text{ }^{\ast }%
\mathbb{Z}
|n\leq x\right\} \right) $ exists in $^{\ast }%
\mathbb{R}
_{\mathbf{d}}$ \ \ \ \ \ \ \ \ \ \ \ \ \ \ \ \ \ \ \ \ \ \ \ \ \ \ \ \ \ \ \
\ \ \ \ \ \ \ \ \ \ \ \ 

there exists a unique $m\in $ $^{\ast }%
\mathbb{Z}
$ such that $m\leq x<m+1.$

\textbf{Proof. }Let $x\in $ $^{\ast }%
\mathbb{R}
_{\mathbf{d}}.$

\textit{Existence:}\textbf{\ }Since $^{\ast }%
\mathbb{Z}
$ is not hyperbounded from below $x$ is not a lower

bound for $^{\ast }%
\mathbb{Z}
$ hence the set $\mathbf{A=}\left\{ n\in \text{ }^{\ast }%
\mathbb{Z}
|n\leq x\right\} $is not empty. Moreover,

$x$ is an upper bound for $\mathbf{A}$ by definition of $\mathbf{A.}$Hence,
as a subset of $^{\ast }%
\mathbb{R}
_{\mathbf{d}},$

$\mathbf{A}$ has a supremum $\sup \left( \left\{ n\in \text{ }^{\ast }%
\mathbb{Z}
|n\leq x\right\} \right) $,call it $\Delta \left( x\right) .$ $\Delta \left(
x\right) -1$ is not an

upper bound for $\mathbf{A}$ hence there exists \ $m\in \mathbf{A\subset }$ $%
^{\ast }%
\mathbb{Z}
$ such that $\Delta \left( x\right) -1<m$

and consequently $\Delta \left( x\right) <m+1.$So $m+1\notin \mathbf{A,}$so $%
m+1.$Therefore

$m\leq x<m+1.$

\textit{Uniqueness:} Suppose that $m^{\prime }\leq x<m^{\prime }+1$ for $%
m^{\prime }\in $ $^{\ast }%
\mathbb{Z}
.$If $m^{\prime }<m,$ then

$m^{\prime }+1\leq m$ implying $m^{\prime }\leq x<m^{\prime }+1\leq m\leq x,$
a contradiction. $m<m^{\prime }$

leads to a similar contradiction. So $m=m^{\prime }.$

Let $E=\left\{ x\in \text{ }^{\ast }%
\mathbb{R}
_{\mathbf{d}}|x^{2}=x\times _{^{\ast }%
\mathbb{R}
_{\mathbf{d}}}x<_{\text{ }^{\ast }%
\mathbb{R}
_{\mathbf{d}}}2\right\} .$Note that $1^{2}=1<2,$ so that

$1\in E$ and in particular $E$ is non-empty. Further if $x>2$ then:

$x^{2}=x\times _{^{\ast }%
\mathbb{R}
_{\mathbf{d}}}x>2x>4>2.$

Hence $2$ is an upper bound for $E$ and so we may define $\zeta \triangleq
\sup E.$

\textbf{Theorem} \textbf{1.3.6.6. }Suppose that $\zeta =\mathbf{s}$-$\sup E$%
\ exist.There exists a unique

positive number $\zeta \triangleq Ext$-$\sqrt{2}\triangleq \#$-$\sqrt{2}\in $
$^{\ast }%
\mathbb{R}
_{\mathbf{d}}$ such that $\zeta ^{2}=\zeta \times _{^{\ast }%
\mathbb{R}
_{\mathbf{d}}}\zeta =_{\text{ }^{\ast }%
\mathbb{R}
_{\mathbf{d}}}2.$

\textbf{Proof. }Note further that $\zeta $ $>1>0$ is positive. We split the
remainder of

the proof into showing that $\zeta ^{2}<2$ and $\zeta ^{2}>2$ both lead to
contradictions.

Suppose for a contradiction that $\zeta ^{2}<2.$Let $h=\dfrac{1}{2}\min
\left( \zeta ,\dfrac{2-\zeta ^{2}}{3\zeta }\right) >0.\ $

Then $(\zeta +h)^{2}=\zeta ^{2}+2h\times \zeta +h^{2}<\zeta ^{2}+3h\times
\zeta <\zeta ^{2}+(2-\zeta ^{2})=2.$Since $h<\zeta $

and $h<\dfrac{2-\zeta ^{2}}{3\zeta }.$Hence $\zeta +h\in E$ and since $\zeta
=\mathbf{s}$-$\sup E$ we get $\zeta +h<\zeta ,$ a

contradiction.

Suppose instead that $\zeta ^{2}>2.$ Let $h=\dfrac{1}{2}\left( \dfrac{\zeta
^{2}-2}{2\zeta }\right) >0.\ $As $\zeta -h<\zeta $ there

exists $\epsilon \in E$ with $\zeta -h<\epsilon $ by the Strong
Approximation Property; then

$(\zeta -h)^{2}<\epsilon ^{2}<2\Rightarrow \zeta ^{2}-2h\times \zeta
+h^{2}<2.\ \ \ \ \ \ \ \ \ \ \ \ \ \ \ \ \ \ \ \ \ \ \ \ \ \ \ \ \ \ \ \ \ \
\ \ $

As $h^{2}>0$ this gives $\zeta ^{2}-2h\times \zeta <2,$ and so, since $\zeta
>0,$ we have

$h>\left( \zeta ^{2}-2\right) /2\zeta $ which contradicts our choice of $h.$%
Finally, by trichotomy,

$\zeta ^{2}=2$ follows as the only remaining possibility.

Let $E_{<}=\left\{ x\in \text{ }^{\ast }%
\mathbb{Q}
|x^{2}=x\times _{^{\ast }%
\mathbb{R}
}x<_{\text{ }^{\ast }%
\mathbb{R}
}2\right\} $ and

$E_{>}=\left\{ x\in \text{ }^{\ast }%
\mathbb{Q}
|x^{2}=x\times _{^{\ast }%
\mathbb{R}
}x>_{\text{ }^{\ast }%
\mathbb{R}
}2\right\} .$

Hence a \textit{Dedekind hyperreal }$\#$-$\sqrt{2}\in $ $^{\ast }%
\mathbb{R}
_{\mathbf{d}}$ is a pair $(U,V)\in \mathbf{P}\left( ^{\ast }%
\mathbb{Q}
\right) \times $ $\mathbf{P}\left( ^{\ast }%
\mathbb{Q}
\right) $

where $U=E_{<},V=E_{>}.$

\textbf{Theorem} \textbf{1.3.6.7.}Let $a\in $ $^{\ast }%
\mathbb{R}
$ be any positive hyperreal number. Then for any

$n\in $ $^{\ast }%
\mathbb{N}
$ there exists a unique Dedekind hyperreal number $\alpha \in $ $^{\ast }%
\mathbb{R}
_{\mathbf{d}}$ \ \ \ \ \ \ \ \ \ \ \ \ \ \ \ \ \ \ \ \ \ \ \ \ \ \ \ \ \ \ 

(denoted by$\left( \sqrt[n]{a}\right) _{\mathbf{d}}$) such that $\alpha
^{n}=a.$

\textbf{Theorem} \textbf{1.3.6.8.}($\ast $\textbf{-Density of} $^{\ast }%
\mathbb{Q}
$ \textbf{in} $^{\ast }%
\mathbb{R}
_{\mathbf{d}}$). Let $x\in $ $^{\ast }%
\mathbb{R}
_{\mathbf{d}}$ be a Dedekind \ \ \ \ \ \ \ \ \ \ \ \ \ \ \ \ \ \ \ \ \ \ \ 

hyperreal number such that $x\pm \varepsilon \neq x$ for any $\varepsilon
>0, $ $\varepsilon \approx 0.$. For every $\epsilon >0$ \ \ \ \ \ \ \ \ \ \
\ \ \ \ \ \ \ \ \ 

there exists a hyperrational number $r\in $ $^{\ast }%
\mathbb{Q}
$ such that $x-\epsilon <r<x+\epsilon .$ \ \ \ \ \ \ \ \ \ \ \ \ \ \ \ \ \ \
\ \ \ \ \ \ \ \ \ \ \ \ \ \ \ \ \ \ \ \ \ 

\textbf{Proof.}\bigskip\ Let $\epsilon >0$ be given. By the Generalyzed
Archimedean Property of $^{\ast }%
\mathbb{R}
_{\mathbf{d}}$ \ \ \ \ \ \ \ \ \ \ \ \ \ \ \ \ \ \ \ 

we can pick $n\in $ $^{\ast }%
\mathbb{N}
$ with $n^{-1}<\epsilon .$Let $q=\left[ \left\vert nx\right\vert \right] \in 
$ $^{\ast }%
\mathbb{N}
.$Since $q\leq nx<q+1,$we \ \ \ \ \ \ \ \ \ \ \ \ \ \ \ \ \ \ \ \ \ \ 

have $\dfrac{q}{n}\leq x<\dfrac{q}{n}+\dfrac{1}{n}<\dfrac{q}{n}+\epsilon .$%
Now let $r=\dfrac{q}{n}\in $ $^{\ast }%
\mathbb{Q}
.$Then $r\leq x<r+\epsilon $ and \ \ \ \ \ \ \ \ \ \ \ \ \ \ \ \ 

hence $x-\epsilon <r<x+\epsilon .$

\ \ \ \ \ \ \ \ \ \ \ \ \ \ \ \ \ \ \ \ \ \ \ \ \ \ \ \ \ \ \ \ \ \ \ \ \ \
\ \ \ \ \ \ \ \ 

\section{Rearrangements of countable infinite series.}

\bigskip

\textbf{Definition 1.3.6.19.}(\textbf{i})\textbf{\ }Let be $\left\{ \mathbf{s%
}_{n}\right\} _{n=1}^{\infty }$ countable sequence $\mathbf{s}_{n}:%
\mathbb{N}
\rightarrow $ $%
\mathbb{R}
.$

such that: (\textbf{a}) $\forall n\left( \mathbf{s}_{n}\geq 0\right) $ or (%
\textbf{b}) $\forall n\left( \mathbf{s}_{n}<0\right) $ or

(\textbf{c}) $\left\{ \mathbf{s}_{n}\right\} _{n=1}^{\infty }=\left\{ 
\mathbf{s}_{n_{1}}\right\} _{n_{1}\in 
\mathbb{N}
_{1}}^{\infty }\cup \left\{ \mathbf{s}_{n_{2}}\right\} _{n_{2}\in 
\mathbb{N}
_{2}}^{\infty },\forall n_{1}\left( n_{1}\in \widehat{%
\mathbb{N}
}_{1}\right) \left[ \mathbf{s}_{n_{1}}\geq 0\right] ,$

$\forall n_{2}\left( n_{2}\in \widehat{%
\mathbb{N}
}_{2}\right) \left[ \mathbf{s}_{n_{2}}<0\right] ,%
\mathbb{N}
=\widehat{%
\mathbb{N}
}_{1}\cup \widehat{%
\mathbb{N}
}_{2}.$

Then external $\flat $-sum of the countable sequence $\mathbf{s}_{n}$ denoted

\bigskip

$\ \ 
\begin{array}{cc}
\begin{array}{c}
\\ 
\left( \#Ext\text{-}\dsum\limits_{n\in 
\mathbb{N}
}\left( ^{\ast }\mathbf{s}_{n}\right) \right) ^{\flat } \\ 
\end{array}
& \text{\ }\left( 1.3.6.23\right) 
\end{array}%
$

\bigskip is$\ \ \ $

$\ \ \ \ \ \ \ \ \ \ \ \ \ \ \ \ \ \ \ \ \ \ \ \ \ \ \ $

$\ \ \ \ \ \ 
\begin{array}{cc}
\begin{array}{c}
\left( \mathbf{a}\right) \text{ \ \ \ \ }\forall n\left( \mathbf{s}_{n}\geq
0\right) : \\ 
\\ 
\left( \#Ext\text{-}\dsum\limits_{n\in 
\mathbb{N}
}\left( ^{\ast }\mathbf{s}_{n}\right) \right) ^{\flat }\triangleq \\ 
\\ 
\triangleq \text{ }\underset{k\in 
\mathbb{N}
}{\sup }\left\{ \dsum\limits_{n\leq k}\left( ^{\ast }\mathbf{s}_{n}\right)
^{\#}\right\} , \\ 
\\ 
\left( \mathbf{b}\right) \text{ \ \ \ \ \ \ }\forall n\left( \mathbf{s}%
_{n}<0\right) : \\ 
\\ 
\left( \#Ext\text{-}\dsum\limits_{n\in 
\mathbb{N}
}\left( ^{\ast }\mathbf{s}_{n}\right) \right) ^{\flat }\triangleq \\ 
\\ 
\triangleq \text{ }\underset{k\in 
\mathbb{N}
}{\inf }\left\{ \dsum\limits_{n\leq k}\left( ^{\ast }\mathbf{s}_{n}\right)
^{\#}\right\} . \\ 
\\ 
\left( \mathbf{c}\right) \text{ \ \ }\forall n_{1}\left( n_{1}\in \widehat{%
\mathbb{N}
}_{1}\right) \left[ \mathbf{s}_{n_{1}}\geq 0\right] , \\ 
\\ 
\forall n_{2}\left( n_{2}\in \widehat{%
\mathbb{N}
}_{2}\right) \left[ \mathbf{s}_{n_{2}}<0\right] ,%
\mathbb{N}
=\widehat{%
\mathbb{N}
}_{1}\cup \widehat{%
\mathbb{N}
}_{2}: \\ 
\\ 
\left( \#Ext\text{-}\dsum\limits_{n\in 
\mathbb{N}
}\left( ^{\ast }\mathbf{s}_{n}\right) \right) ^{\flat }\triangleq \\ 
\\ 
\triangleq \left( \#Ext\text{-}\dsum\limits_{n_{1}\in \widehat{%
\mathbb{N}
}_{1}}\left( ^{\ast }\mathbf{s}_{n_{1}}\right) \right) ^{\flat }+\left( \#Ext%
\text{-}\dsum\limits_{n_{2}\in \widehat{%
\mathbb{N}
}_{2}}\left( ^{\ast }\mathbf{s}_{n_{2}}\right) \right) ^{\flat }. \\ 
\end{array}
& \text{ \ }\left( 1.3.6.24\right)%
\end{array}%
$

\bigskip

\textbf{Definition 1.3.6.20.}(\textbf{i})\textbf{\ }Let be $\left\{ \mathbf{s%
}_{n}\right\} _{n=1}^{\infty }$ countable sequence $\mathbf{s}_{n}:%
\mathbb{N}
\rightarrow $ $^{\ast }%
\mathbb{R}
$

such that: (\textbf{a}) $\forall n\left( \mathbf{s}_{n}\geq 0\right) $ or (%
\textbf{b}) $\forall n\left( \mathbf{s}_{n}<0\right) $ or

(\textbf{c}) $\left\{ \mathbf{s}_{n}\right\} _{n=1}^{\infty }=\left\{ 
\mathbf{s}_{n_{1}}\right\} _{n_{1}\in 
\mathbb{N}
_{1}}^{\infty }\cup \left\{ \mathbf{s}_{n_{2}}\right\} _{n_{2}\in 
\mathbb{N}
_{2}}^{\infty },\forall n_{1}\left( n_{1}\in \widehat{%
\mathbb{N}
}_{1}\right) \left[ \mathbf{s}_{n_{1}}\geq 0\right] ,$

$\forall n_{2}\left( n_{2}\in \widehat{%
\mathbb{N}
}_{2}\right) \left[ \mathbf{s}_{n_{2}}<0\right] ,%
\mathbb{N}
=\widehat{%
\mathbb{N}
}_{1}\cup \widehat{%
\mathbb{N}
}_{2}.$

Then external $\flat $-sum of the countable sequence $\mathbf{s}_{n}$ denoted

\bigskip

$\ \ \ \ \ \ \ \ \ \ \ 
\begin{array}{cc}
\begin{array}{c}
\\ 
\left( \#Ext\text{-}\dsum\limits_{n\in 
\mathbb{N}
}\mathbf{s}_{n}^{\#}\right) ^{\flat } \\ 
\end{array}
& \text{ \ }\left( 1.3.6.23^{\prime }\right) 
\end{array}%
$

\bigskip is$\ \ $

$\ \ \ \ \ \ \ \ \ \ 
\begin{array}{cc}
\begin{array}{c}
\left( \mathbf{a}\right) \text{ \ \ \ \ }\forall n\left( \mathbf{s}_{n}\geq
0\right) : \\ 
\\ 
\left( \#Ext\text{-}\dsum\limits_{n\in 
\mathbb{N}
}\mathbf{s}_{n}^{\#}\right) ^{\flat }\triangleq \\ 
\\ 
\triangleq \text{ }\underset{k\in 
\mathbb{N}
}{\sup }\left\{ \dsum\limits_{n\leq k}s_{n}^{\#}\right\} , \\ 
\\ 
\left( \mathbf{b}\right) \text{ \ \ \ \ \ \ }\forall n\left( \mathbf{s}%
_{n}<0\right) : \\ 
\\ 
\left( \#Ext\text{-}\dsum\limits_{n\in 
\mathbb{N}
}\mathbf{s}_{n}^{\#}\right) ^{\flat }\triangleq \\ 
\\ 
\triangleq \text{ }\underset{k\in 
\mathbb{N}
}{\inf }\left\{ \dsum\limits_{n\leq k}s_{n}^{\#}\right\} . \\ 
\\ 
\left( \mathbf{c}\right) \text{ \ \ }\forall n_{1}\left( n_{1}\in 
\mathbb{N}
_{1}\right) \left[ \mathbf{s}_{n_{1}}\geq 0\right] , \\ 
\\ 
\forall n_{2}\left( n_{2}\in 
\mathbb{N}
_{2}\right) \left[ \mathbf{s}_{n_{2}}<0\right] ,%
\mathbb{N}
=%
\mathbb{N}
_{1}\cup 
\mathbb{N}
_{2}: \\ 
\\ 
\left( \#Ext\text{-}\dsum\limits_{n\in 
\mathbb{N}
}\mathbf{s}_{n}^{\#}\right) ^{\flat }\triangleq \\ 
\\ 
\triangleq \left( \#Ext\text{-}\dsum\limits_{n_{1}\in \widehat{%
\mathbb{N}
}_{1}}\mathbf{s}_{n_{1}}^{\#}\right) ^{\flat }+\left( \#Ext\text{-}%
\dsum\limits_{n_{2}\in \widehat{%
\mathbb{N}
}_{2}}\mathbf{s}_{n_{2}}^{\#}\right) ^{\flat }. \\ 
\end{array}
& \text{ \ }\left( 1.3.6.24^{\prime }\right)%
\end{array}%
$

(\textbf{ii}) Let be $\left\{ \mathbf{s}_{n}\right\} _{n=1}^{\infty }$
countable sequence $\mathbf{s}_{n}:%
\mathbb{N}
\rightarrow $ $^{\ast }%
\mathbb{R}
_{\mathbf{d}},$ \ \ \ \ \ \ \ \ \ \ \ \ \ \ \ \ \ \ \ \ \ \ \ \ \ \ \ \ \ 

such that (\textbf{a}) $\forall n\left( \mathbf{s}_{n}\geq 0\right) $ or (%
\textbf{b}) $\forall n\left( \mathbf{s}_{n}<0\right) $ or

(\textbf{c}) $\left\{ \mathbf{s}_{n}\right\} _{n=1}^{\infty }=\left\{ 
\mathbf{s}_{n_{1}}\right\} _{n_{1}\in 
\mathbb{N}
_{1}}^{\infty }\cup \left\{ \mathbf{s}_{n_{2}}\right\} _{n_{2}\in 
\mathbb{N}
_{2}}^{\infty },\forall n_{1}\left( n_{1}\in \widehat{%
\mathbb{N}
}_{1}\right) \left[ \mathbf{s}_{n_{1}}\geq 0\right] ,$

$\forall n_{2}\left( n_{2}\in \widehat{%
\mathbb{N}
}_{2}\right) \left[ \mathbf{s}_{n_{2}}<0\right] ,%
\mathbb{N}
=\widehat{%
\mathbb{N}
}_{1}\cup \widehat{%
\mathbb{N}
}_{2}.$

\bigskip Then external $\flat $-sum of the countable sequence $\mathbf{s}_{n}
$ denoted

$\ \ 
\begin{array}{cc}
\begin{array}{c}
\\ 
\left( \#Ext\text{-}\dsum\limits_{n\in 
\mathbb{N}
}\mathbf{s}_{n}\right) ^{\flat } \\ 
\end{array}
& \text{ \ \ \ }\left( 1.3.6.23^{\prime \prime }\right) 
\end{array}%
$ \ 

\bigskip is$\ $

$\ \ \ \ \ \ \ \ \ \ \ \ \ \ \ \ \ \ \ \ \ \ \ \ \ \ \ \ \ \ $

$\ \ \ 
\begin{array}{cc}
\begin{array}{c}
\left( \mathbf{a}\right) \text{ \ \ \ \ }\forall n\left( \mathbf{s}_{n}\geq
0\right) : \\ 
\\ 
\left( \#Ext\text{-}\dsum\limits_{n\in 
\mathbb{N}
}\mathbf{s}_{n}\right) ^{\flat }\triangleq \\ 
\\ 
\triangleq \text{ }\underset{k\in 
\mathbb{N}
}{\sup }\left\{ \dsum\limits_{n\leq k}\mathbf{s}_{n}\right\} , \\ 
\\ 
\left( \mathbf{b}\right) \text{ \ \ \ \ \ \ }\forall n\left( \mathbf{s}%
_{n}<0\right) : \\ 
\\ 
\left( \#Ext\text{-}\dsum\limits_{n\in 
\mathbb{N}
}\mathbf{s}_{n}\right) ^{\flat }\triangleq \\ 
\\ 
\triangleq \text{ }\underset{k\in 
\mathbb{N}
}{\inf }\left\{ \dsum\limits_{n\leq k}\mathbf{s}_{n}\right\} . \\ 
\\ 
\left( \mathbf{c}\right) \text{ \ }\forall n_{1}\left( n_{1}\in \widehat{%
\mathbb{N}
}_{1}\right) \left[ \mathbf{s}_{n_{1}}\geq 0\right] , \\ 
\\ 
\forall n_{2}\left( n_{2}\in \widehat{%
\mathbb{N}
}_{2}\right) \left[ \mathbf{s}_{n_{2}}<0\right] ,%
\mathbb{N}
=\widehat{%
\mathbb{N}
}_{1}\cup \widehat{%
\mathbb{N}
}_{2}: \\ 
\\ 
\left( \#Ext\text{-}\dsum\limits_{n\in 
\mathbb{N}
}\mathbf{s}_{n}\right) ^{\flat }\triangleq \\ 
\\ 
\triangleq \left( \#Ext\text{-}\dsum\limits_{n_{1}\in \widehat{%
\mathbb{N}
}_{1}}\mathbf{s}_{n_{1}}\right) ^{\flat }+\left( \#Ext\text{-}%
\dsum\limits_{n_{2}\in \widehat{%
\mathbb{N}
}_{2}}\mathbf{s}_{n_{2}}\right) ^{\flat }.%
\end{array}
& \text{ \ \ }\left( 1.3.6.24^{\prime \prime }\right)%
\end{array}%
$

\bigskip

\textbf{Theorem} \textbf{1.3.6.9.}(\textbf{i}) Let be $\left\{ \mathbf{s}%
_{n}\right\} _{n=1}^{\infty }$ countable sequence $\mathbf{s}_{n}:%
\mathbb{N}
\rightarrow $ $%
\mathbb{R}
$

such that $\forall n\left( n\in 
\mathbb{N}
\right) \left[ \mathbf{s}_{n}\geq 0\right] ,$ $\sum_{n=1}^{\infty }\mathbf{s}%
_{n}=\eta <\infty ,$ i.e. infinite series

$\sum_{n=1}^{\infty }\mathbf{s}_{n}$ converges to $\eta $ in $%
\mathbb{R}
.$Then

\ \ \ $\ \ 
\begin{array}{cc}
\begin{array}{c}
\\ 
\left( \#Ext\text{-}\dsum\limits_{n\in 
\mathbb{N}
}\left( ^{\ast }\mathbf{s}_{n}\right) \right) ^{\flat }\triangleq \text{ }%
\underset{k\in 
\mathbb{N}
}{\sup }\left\{ \dsum\limits_{n\leq k}\left( ^{\ast }\mathbf{s}_{n}\right)
^{\#}\right\} = \\ 
\\ 
=\text{ }\left( ^{\ast }\eta \right) ^{\#}-\varepsilon _{\mathbf{d}}\in 
\text{ }^{\ast }%
\mathbb{R}
_{\mathbf{d}}, \\ 
\end{array}
& \left( 1.3.6.25.a\right)%
\end{array}%
$

\bigskip

(\textbf{ii}) Let be $\left\{ \mathbf{s}_{n}\right\} _{n=1}^{\infty }$
countable sequence $\mathbf{s}_{n}:%
\mathbb{N}
\rightarrow $ $%
\mathbb{R}
$

such that $\forall n\left( n\in 
\mathbb{N}
\right) \left[ \mathbf{s}_{n}<0\right] ,$ $\sum_{n=1}^{\infty }\mathbf{s}%
_{n}=\eta <\infty ,$ i.e. infinite series

$\sum_{n=1}^{\infty }\mathbf{s}_{n}$ converges to $\eta $ in $%
\mathbb{R}
.$Then\ \ \ \ \ \ \ \ \ \ \ \ \ \ \ 

\bigskip

\ $\ 
\begin{array}{cc}
\begin{array}{c}
\\ 
\left( \#Ext\text{-}\dsum\limits_{n\in 
\mathbb{N}
}\left( ^{\ast }\mathbf{s}_{n}\right) \right) ^{\flat }\triangleq \text{ }%
\underset{k\in 
\mathbb{N}
}{\inf }\left\{ \dsum\limits_{n\leq k}\left( ^{\ast }\mathbf{s}_{n}\right)
^{\#}\right\} = \\ 
\\ 
=\text{ }\left( ^{\ast }\eta \right) ^{\#}+\varepsilon _{\mathbf{d}}\in 
\text{ }^{\ast }%
\mathbb{R}
_{\mathbf{d}}, \\ 
\end{array}
& \left( 1.3.6.25.b\right)%
\end{array}%
$

\bigskip

(\textbf{iii}) Let be $\left\{ \mathbf{s}_{n}\right\} _{n=1}^{\infty }$
countable sequence $\mathbf{s}_{n}:%
\mathbb{N}
\rightarrow $ $%
\mathbb{R}
$ such that

(1) $\left\{ \mathbf{s}_{n}\right\} _{n=1}^{\infty }=\left\{ \mathbf{s}%
_{n_{1}}\right\} _{n_{1}\in 
\mathbb{N}
_{1}}^{\infty }\cup \left\{ \mathbf{s}_{n_{2}}\right\} _{n_{2}\in 
\mathbb{N}
_{2}}^{\infty },\forall n_{1}\left( n_{1}\in \widehat{%
\mathbb{N}
}_{1}\right) \left[ \mathbf{s}_{n_{1}}\geq 0\right] ,$

$\forall n_{2}\left( n_{2}\in \widehat{%
\mathbb{N}
}_{2}\right) \left[ \mathbf{s}_{n_{2}}<0\right] ,%
\mathbb{N}
=\widehat{%
\mathbb{N}
}_{1}\cup \widehat{%
\mathbb{N}
}_{2}$ and

(2) \ $\dsum\limits_{n_{1}\in \widehat{%
\mathbb{N}
}_{1}}\mathbf{s}_{n_{1}}=\eta _{1}<\infty ,\dsum\limits_{n_{2}\in \widehat{%
\mathbb{N}
}_{2}}\mathbf{s}_{n_{2}}=\eta _{2}>-\infty .$

\bigskip 

Then \ \ \ \ \ 

\ \ \ \ \ \ \ \ \ $\ \ \ \ \ \ \ \ \ $

$\ \ \ 
\begin{array}{cc}
\begin{array}{c}
\\ 
\left( \#Ext\text{-}\dsum\limits_{n\in 
\mathbb{N}
}\left( ^{\ast }\mathbf{s}_{n}\right) \right) ^{\flat }\triangleq \\ 
\\ 
\triangleq \left( \#Ext\text{-}\dsum\limits_{n_{1}\in \widehat{%
\mathbb{N}
}_{1}}\left( ^{\ast }\mathbf{s}_{n_{1}}\right) \right) ^{\flat }\text{ }%
+\left( \#Ext\text{-}\dsum\limits_{n_{2}\in \widehat{%
\mathbb{N}
}_{2}}\left( ^{\ast }\mathbf{s}_{n_{2}}\right) \right) ^{\flat }= \\ 
\\ 
=\left( ^{\ast }\eta _{1}\right) ^{\#}-\varepsilon _{\mathbf{d}}+\text{ }%
\left( ^{\ast }\eta _{2}\right) ^{\#}+\varepsilon _{\mathbf{d}}= \\ 
\\ 
=\left( ^{\ast }\eta _{1}\right) ^{\#}+\text{ }\left( ^{\ast }\eta
_{2}\right) ^{\#}-\varepsilon _{\mathbf{d}}\in \text{ }^{\ast }%
\mathbb{R}
_{\mathbf{d}}. \\ 
\end{array}
& \left( 1.3.6.25.c\right)%
\end{array}%
$

\bigskip

\textbf{Theorem} \textbf{1.3.6.10.}Let be $\left\{ a_{n}\right\}
_{n=1}^{\infty }$ countable sequence $a_{n}:%
\mathbb{N}
\rightarrow $ $^{\ast }%
\mathbb{R}
_{\mathbf{d}},$ \ \ \ \ \ \ \ \ \ \ \ \ \ \ \ \ \ \ \ \ \ \ \ \ \ \ \ \ \ 

such that $\forall n\left( a_{n}\geq 0\right) $ and $\left( \#Ext\text{-}%
\dsum\limits_{n\in 
\mathbb{N}
}a_{n}\right) ^{\flat }$ external sum of the

countable sequence $\left\{ a_{n}\right\} _{n=1}^{\infty }$ denoted by $%
\mathbf{s}.$\ \ \ \ \ \ \ \ \ \ \ \ \ \ \ \ \ \ \ \ \ \ \ \ \ \ \ \ \ \ \ \
\ \ \ \ \ \ \ \ \ \ \ \ \ \ \ \ \ \ \ \ \ \ \ \ \ \ \ \ \ \ \ \ \ \ \ \ \ \
\ \ \ \ \ \ \ \ \ \ \ \ \ \ \ 

Let be $\left\{ b_{n}\right\} _{n=1}^{\infty }$ countable sequence where $%
b_{n}=a_{m\left( n\right) }$ any

rearrangement of \ terms of the sequence $\left\{ a_{n}\right\}
_{n=1}^{\infty }.$

Then external sum $\mathbf{\sigma }=\left( \#Ext\text{-}\dsum\limits_{n\in 
\mathbb{N}
}b_{n}\right) ^{\flat }$ of the countable

sequence $\left\{ b_{n}\right\} _{m=1}^{\infty }$ \ has the same value $%
\mathbf{s}$ as external sum of the

countable sequence $\left\{ a_{n}\right\} \mathbf{,}$i.e.

\bigskip

$\ \ 
\begin{array}{cc}
\begin{array}{c}
\\ 
\mathbf{\sigma =s.} \\ 
\end{array}
& \text{ \ }\left( 1.3.6.25\right)%
\end{array}%
$

\bigskip

\bigskip

\textbf{Proof.}\bigskip Let be $\mathbf{\sigma }_{n}=b_{1}+b_{2}+...+b_{n}$
\ the $n$-th partial sum of the \ \ \ \ \ \ \ \ \ \ \ \ \ \ \ \ \ \ \ \ \ \ 

sequence $\left\{ b_{n}\right\} _{n=1}^{\infty }$ and $\mathbf{s}%
_{m}=a_{1}+a_{2}+...+a_{m}$ the $m$-th partial

sum of the sequence $\left\{ a_{n}\right\} _{n=1}^{\infty }$.

It is easy to see that for any given $n$-th partial sum

$\mathbf{\sigma }_{n}=b_{1}+b_{2}+...+b_{n}$ there is exist $\ m$-th partial
sum

$\mathbf{s}_{m\left( n\right) }=a_{1}+a_{2}+...+a_{m\left( n\right) }$ such
that:

\bigskip

$%
\begin{array}{cc}
\begin{array}{c}
\\ 
\left\{ a_{m}\right\} _{m=1}^{m\left( n\right) }\supseteqq \left\{
b_{i}\right\} _{i=1}^{n}, \\ 
\end{array}
& \text{ }\left( 1.3.6.26\right)%
\end{array}%
$

\bigskip

and there is exist $\ N$-th partial sum

$\mathbf{\sigma }_{N\left( m\right) }=b_{1}+b_{2}+...+b_{n}+...+b_{N\left(
m\right) }$ such that:

\bigskip

$\ \ \ 
\begin{array}{cc}
\begin{array}{c}
\\ 
\left\{ b_{j}\right\} _{j=1}^{N\left( m\right) }\supseteqq \left\{
a_{i}\right\} _{i=1}^{m\left( n\right) }, \\ 
\end{array}
& \text{ \ }\left( 1.3.6.27\right)%
\end{array}%
$

\bigskip

By using setting and Eqs.(1.3.26)-(1.3.27) one

obtain inequality \bigskip\ 

$\ \ \ \ \ \ \ \ \ \ \ \ \ \ \ \ \ \ \ \ \ \ \ \ \ \ \ \ \ \ \ \ \ \ \ \ \ \
\ \ \ \ \ \ \ \ \ \ \ \ \ \ \ \ $

$\ \ \ \ 
\begin{array}{cc}
\begin{array}{c}
\\ 
\mathbf{\sigma }_{n}\leq \mathbf{s}_{m\left( n\right) }\leq \mathbf{\sigma }%
_{N\left( m\right) }. \\ 
\end{array}
& \text{ \ \ }\left( 1.3.6.28\right)%
\end{array}%
$

\bigskip

\bigskip By using \textbf{Proposition 1.3.6.5. }one obtain

\bigskip

$\ \ \ \ 
\begin{array}{cc}
\begin{array}{c}
\\ 
\underset{n\in 
\mathbb{N}
}{\sup }\left\{ \mathbf{\sigma }_{n}\right\} \leq \text{ }\underset{n\in 
\mathbb{N}
}{\sup }\left\{ \mathbf{s}_{m\left( n\right) }\right\} \leq \text{ }\underset%
{n\in 
\mathbb{N}
}{\sup }\left\{ \mathbf{\sigma }_{N\left( m\right) }\right\} . \\ 
\end{array}
& \text{ \ \ }\left( 1.3.6.29\right)%
\end{array}%
$

\bigskip

Hence $\mathbf{\sigma }\leq \mathbf{s}\leq \mathbf{\sigma }$ and finally we
obtain $\mathbf{\sigma =s.}$

\bigskip \textbf{Theorem} \textbf{1.3.21.}

\textbf{Theorem} \textbf{1.3.22.(i) }Let be $\left\{ a_{n}\right\}
_{n=1}^{\infty }$ countable sequence $a_{n}:%
\mathbb{N}
\rightarrow $ $^{\ast }%
\mathbb{R}
_{\mathbf{d}},$ \ \ \ \ \ \ \ \ \ \ \ \ \ \ \ \ \ \ \ \ \ \ \ \ \ \ \ \ \ 

such that $\forall n\left( a_{n}\geq 0\right) $ and $\#Ext$-$%
\dsum\limits_{n\in 
\mathbb{N}
}a_{n}$ external sum of the sequence $\left\{ a_{n}\right\} _{n=1}^{\infty
}. $ \ \ \ \ \ \ \ \ \ \ \ \ \ \ \ \ \ 

Then for any $c\in $ $^{\ast }%
\mathbb{R}
_{+}$ the next equality is satisfied:

\bigskip

$\ \ 
\begin{array}{cc}
\begin{array}{c}
\\ 
c^{\#}\times \left( \#Ext\text{-}\dsum\limits_{n\in 
\mathbb{N}
}a_{n}\right) =\left( \#Ext\text{-}\dsum\limits_{n\in 
\mathbb{N}
}c^{\#}\times a_{n}\right) \\ 
\end{array}
& \text{\ }\left( 1.3.6.30\right)%
\end{array}%
$

\bigskip

\textbf{(ii) }Let be $\left\{ a_{n}\right\} _{n=1}^{\infty }$ countable
sequence $a_{n}:%
\mathbb{N}
\rightarrow $ $^{\ast }%
\mathbb{R}
_{\mathbf{d}},$ \ \ \ \ \ \ \ \ \ \ \ \ \ \ \ \ \ \ \ \ \ \ \ \ \ \ \ \ \ 

such that $\forall n\left( a_{n}<0\right) $ and $\#Ext$-$\dsum\limits_{n\in 
\mathbb{N}
}a_{n}$ external sum of the

sequence $\left\{ a_{n}\right\} _{n=1}^{\infty }.$Then for any $c\in $ $%
^{\ast }%
\mathbb{R}
_{+}$ the next equality is

satisfied:

\bigskip

$\ \ 
\begin{array}{cc}
\begin{array}{c}
\\ 
c^{\#}\times \left( \#Ext\text{-}\dsum\limits_{n\in 
\mathbb{N}
}a_{n}\right) =\left( \#Ext\text{-}\dsum\limits_{n\in 
\mathbb{N}
}c^{\#}\times a_{n}\right) \\ 
\end{array}
& \text{ \ }\left( 1.3.6.30^{\prime }\right)%
\end{array}%
$

\bigskip

\textbf{(iii) }Let be $\left\{ \mathbf{s}_{n}\right\} _{n=1}^{\infty }$
countable sequence $\mathbf{s}_{n}:%
\mathbb{N}
\rightarrow $ $^{\ast }%
\mathbb{R}
_{\mathbf{d}}$ such that

$\left\{ \mathbf{s}_{n}\right\} _{n=1}^{\infty }=\left\{ \mathbf{s}%
_{n_{1}}\right\} _{n_{1}\in 
\mathbb{N}
_{1}}^{\infty }\cup \left\{ \mathbf{s}_{n_{2}}\right\} _{n_{2}\in 
\mathbb{N}
_{2}}^{\infty },\forall n_{1}\left( n_{1}\in \widehat{%
\mathbb{N}
}_{1}\right) \left[ \mathbf{s}_{n_{1}}\geq 0\right] ,$

$\forall n_{2}\left( n_{2}\in \widehat{%
\mathbb{N}
}_{2}\right) \left[ \mathbf{s}_{n_{2}}<0\right] ,%
\mathbb{N}
=\widehat{%
\mathbb{N}
}_{1}\cup \widehat{%
\mathbb{N}
}_{2}$

Then the next equality is satisfied:

\bigskip

$\ \ 
\begin{array}{cc}
\begin{array}{c}
\\ 
\ \ \ \#Ext\text{-}\dsum\limits_{n\in 
\mathbb{N}
}\mathbf{s}_{n}= \\ 
\\ 
\#Ext\text{-}\dsum\limits_{n_{1}\in \widehat{%
\mathbb{N}
}_{1}}\mathbf{s}_{n_{1}}+\#Ext\text{-}\dsum\limits_{n_{2}\in \widehat{%
\mathbb{N}
}_{2}}\mathbf{s}_{n_{2}} \\ 
\end{array}
& \text{ \ \ }\left( 1.3.6.30^{\prime \prime }\right)%
\end{array}%
$

\bigskip $\ \ \ \ \ \ \ \ \ \ \ \ \ \ \ \ \ \ \ \ \ \ \ \ \ \ \ $

\textbf{Proof.(i) }By using \textbf{Definition} \textbf{1.3.20 (ii)} and 
\textbf{Theorem 1.3.1.3 }

one obtain

\bigskip

\ \ $\ \ 
\begin{array}{cc}
\begin{array}{c}
\\ 
c\times \left( \#Ext\text{-}\dsum\limits_{n\in 
\mathbb{N}
}a_{n}\right) =c\times \sup \left\{ \left. \dsum\limits_{n\leq
m}a_{n}\right\vert m\in 
\mathbb{N}
\right\} = \\ 
\\ 
\sup \left[ c\times \left\{ \left. \dsum\limits_{n\leq m}a_{n}\right\vert
m\in 
\mathbb{N}
\right\} \right] = \\ 
\\ 
\sup \left[ \left\{ \left. c\times \dsum\limits_{n\leq m}a_{n}\right\vert
m\in 
\mathbb{N}
\right\} \right] = \\ 
\\ 
\sup \left[ \left\{ \left. \dsum\limits_{n\leq m}c\times a_{n}\right\vert
m\in 
\mathbb{N}
\right\} \right] = \\ 
\\ 
\left( \#Ext\text{-}\dsum\limits_{n\in 
\mathbb{N}
}c\times a_{n}\right) . \\ 
\end{array}
& \text{ \ }\left( 1.3.6.31\right)%
\end{array}%
$

\textbf{Theorem} \textbf{1.3.6.13.}Let be $\left\{ a_{n}\right\}
_{n=1}^{\infty }$ countable sequence $a_{n}:%
\mathbb{N}
\rightarrow $ $^{\ast }%
\mathbb{R}
_{\mathbf{d}},$ \ \ \ \ \ \ \ \ \ \ \ \ \ \ \ \ \ \ \ \ \ \ \ \ \ \ \ \ \ 

and $\#Ext$-$\dsum\limits_{n\in 
\mathbb{N}
}a_{n}$ external sum of the sequence $\left\{ a_{n}\right\} _{n=1}^{\infty
}. $

\bigskip

\textbf{Definition 1.3.6.20. }Let be $\left\{ a_{n}\right\} _{n=1}^{\infty }$
arbitrary countable Cauchy sequence

$a_{n}:%
\mathbb{N}
\rightarrow 
\mathbb{R}
.$The \textbf{upper limit in }$^{\ast }%
\mathbb{R}
_{\mathbf{d}}$ of the countable sequence $\left\{ a_{n}\right\}
_{n=1}^{\infty }$

denoted $^{\ast }%
\mathbb{R}
_{\mathbf{d}}$-$\overline{\overline{\lim a_{n}}}$ is

\bigskip

$\ 
\begin{array}{cc}
\begin{array}{c}
\\ 
^{\ast }%
\mathbb{R}
_{\mathbf{d}}\text{-}\overline{\overline{\lim a_{n}}}=\text{ }\underset{m\in 
\mathbb{N}
}{\inf }\left( \underset{n\geq m}{\sup }\left( a_{n}^{\#}\right) \right) .
\\ 
\end{array}
& \text{ \ \ }\left( 1.3.6.32\right)%
\end{array}%
$

\bigskip

The \textbf{lower limit in }$^{\ast }%
\mathbb{R}
_{\mathbf{d}}$ of the countable sequence $\left\{ a_{n}^{\#}\right\}
_{n=1}^{\infty }$

denoted $^{\ast }%
\mathbb{R}
_{\mathbf{d}}$-$\underline{\underline{\lim a_{n}}}$ is

\bigskip

$\ 
\begin{array}{cc}
\begin{array}{c}
\\ 
^{\ast }%
\mathbb{R}
_{\mathbf{d}}\text{-}\underline{\underline{\lim a_{n}}}=\text{ }\underset{%
m\in 
\mathbb{N}
}{\sup }\left( \underset{n\geq m}{\inf }\left( a_{n}^{\#}\right) \right) .
\\ 
\end{array}
& \text{ \ \ }\left( 1.3.6.33\right)%
\end{array}%
$

\bigskip

\bigskip \textbf{Theorem} \textbf{1.3.6.14. }Suppose that $%
\lim_{n\rightarrow \infty }a_{n}=\zeta \in 
\mathbb{R}
.$ Then

\bigskip

$\ 
\begin{array}{cc}
\begin{array}{c}
\\ 
^{\ast }%
\mathbb{R}
_{\mathbf{d}}\text{-}\overline{\overline{\lim a_{n}}}=\text{ }\left( ^{\ast
}\zeta \right) ^{\#}+\varepsilon _{\mathbf{d}}, \\ 
\\ 
^{\ast }%
\mathbb{R}
_{\mathbf{d}}\text{-}\underline{\underline{\lim a_{n}}}=\text{ }\left(
^{\ast }\zeta \right) ^{\#}-\varepsilon _{\mathbf{d}}. \\ 
\end{array}
& \text{ \ }\left( 1.3.6.34\right) 
\end{array}%
$

\bigskip 

\textbf{Definition 1.3.6.21. }Let be $\left\{ b_{n}\right\} _{n=1}^{\infty }$
countable sequence $b_{n}:%
\mathbb{N}
\rightarrow $ $%
\mathbb{R}
$

such that $\sum_{n=1}^{\infty }b_{n}<\infty ,$ i.e. infinite series $%
\sum_{n=1}^{\infty }b_{n}$ converges in $%
\mathbb{R}
.$

\bigskip The \textbf{upper sum in }$^{\ast }%
\mathbb{R}
_{\mathbf{d}}$ of the infinite series $\sum_{n=1}^{\infty }b_{n}$ denoted

\bigskip $\ 
\begin{array}{cc}
\begin{array}{c}
\\ 
\overline{\overline{\#Ext\text{-}\dsum\limits_{n\in 
\mathbb{N}
}b_{n}}} \\ 
\end{array}
& \text{ \ \ }\left( 1.3.6.32\right)%
\end{array}%
$

\bigskip \bigskip is$\ $\ $\ \ \ \ \ \ \ \ \ \ \ \ $

$\bigskip \ \ \ \ 
\begin{array}{cc}
\begin{array}{c}
\\ 
\overline{\overline{\#Ext\text{-}\dsum\limits_{n\in 
\mathbb{N}
}b_{n}}}\text{ }\triangleq \text{ }^{\ast }%
\mathbb{R}
_{\mathbf{d}}\text{-}\overline{\overline{\lim \left(
\dsum\limits_{i=1}^{n}b_{i}^{\#}\right) }}= \\ 
\\ 
\text{ }\underset{m\in 
\mathbb{N}
}{\inf }\left( \underset{n\geq m}{\sup }\left(
\dsum\limits_{i=1}^{n}b_{i}^{\#}\right) \right) . \\ 
\end{array}
& \text{ \ \ \ \ \ \ }\left( 1.3.6.33\right)%
\end{array}%
$

\bigskip

The \textbf{lower sum in }$^{\ast }%
\mathbb{R}
_{\mathbf{d}}$ of the infinite series $\sum_{n=1}^{\infty }b_{n}$ denoted

\bigskip\ $\ \ \ \ \ 
\begin{array}{cc}
\begin{array}{c}
\\ 
\underline{\underline{\#Ext\text{ -}\dsum\limits_{n\in 
\mathbb{N}
}b_{n}}} \\ 
\end{array}
& \text{ \ }\left( 1.3.6.34\right)%
\end{array}%
$

\bigskip is

$\ \ \ \ \ \ \ \ \ \ \ \ \ \ \ \ \ $

$\ \ \ \ \ \ \ \ \ 
\begin{array}{cc}
\begin{array}{c}
\\ 
\underline{\underline{\#Ext\text{ -}\dsum\limits_{n\in 
\mathbb{N}
}b_{n}}}\triangleq \text{ } \\ 
\\ 
^{\ast }%
\mathbb{R}
_{\mathbf{d}}\text{-}\underline{\underline{\lim \left(
\dsum\limits_{i=1}^{n}b_{i}\right) }}=\text{ }\underset{m\in 
\mathbb{N}
}{\sup }\left( \underset{n\geq m}{\inf }\left(
\dsum\limits_{i=1}^{n}b_{i}\right) \right) . \\ 
\end{array}
& \text{\ }\left( 1.3.6.35\right) 
\end{array}%
$

\bigskip \textbf{Theorem} \textbf{1.3.6.15. }Suppose that $%
\lim_{n\rightarrow \infty }\dsum\limits_{i=1}^{n}b_{i}=\zeta \in 
\mathbb{R}
.$ Then

$\ \ \ \ \ \ \ \ \ \ \ 
\begin{array}{cc}
\begin{array}{c}
\\ 
\overline{\overline{\#Ext\text{-}\dsum\limits_{n\in 
\mathbb{N}
}b_{n}}}=\left( ^{\ast }\zeta \right) ^{\#}+\varepsilon _{\mathbf{d}}, \\ 
\\ 
\underline{\underline{\#Ext-\dsum\limits_{n\in 
\mathbb{N}
}b_{n}}}=\left( ^{\ast }\zeta \right) ^{\#}-\varepsilon _{\mathbf{d}}, \\ 
\end{array}
& \text{ \ }\left( 1.3.6.36\right)%
\end{array}%
$

\bigskip

\textbf{Definition 1.3.6.22. }Let be $\left\{ a_{n}\right\} _{n=1}^{\infty }$
arbitrary countable sequence $a_{n}:%
\mathbb{N}
\rightarrow $ $^{\ast }%
\mathbb{R}
_{\mathbf{d}}.$

\bigskip The \textbf{upper sum} of the countable sequence $\left\{
a_{n}\right\} _{n=1}^{\infty }$ denoted

\bigskip

$\ \ \ 
\begin{array}{cc}
\begin{array}{c}
\\ 
\overline{\overline{\#Ext\text{-}\dsum\limits_{n\in 
\mathbb{N}
}a_{n}}} \\ 
\end{array}
& \text{ \ }\left( 1.3.6.37\right)%
\end{array}%
$

\bigskip

is$\ $\ 

$\bigskip $

$\ \ \ \ \ \ \ $

$\bigskip \ \ 
\begin{array}{cc}
\begin{array}{c}
\\ 
\overline{\overline{\#Ext\text{-}\dsum\limits_{n\in 
\mathbb{N}
}a_{n}}}\text{ }\triangleq \text{ }\underset{m\in 
\mathbb{N}
}{\inf }\left( \underset{n\geq m}{\sup }\left(
\dsum\limits_{i=1}^{n}a_{i}\right) \right) . \\ 
\end{array}
& \text{ \ }\left( 1.3.6.38\right)%
\end{array}%
$\ 

\bigskip

\bigskip The \textbf{lower sum} of the countable sequence $a_{n}$ denoted

\bigskip $\ \ 
\begin{array}{cc}
\begin{array}{c}
\\ 
\underline{\underline{\#Ext-\dsum\limits_{n\in 
\mathbb{N}
}a_{n}}} \\ 
\end{array}
& \text{ \ \ }\left( 1.3.6.39\right)%
\end{array}%
$\bigskip \bigskip

is$\ $\ $\ \ \ \ \ \ \ \ \ $

$\bigskip $

$\ \ 
\begin{array}{cc}
\begin{array}{c}
\\ 
\underline{\underline{\#Ext-\dsum\limits_{n\in 
\mathbb{N}
}a_{n}}}\triangleq \text{ }\underset{m\in 
\mathbb{N}
}{\sup }\left( \underset{n\geq m}{\inf }\left(
\dsum\limits_{i=1}^{n}b_{i}\right) \right) . \\ 
\end{array}
& \text{ \ }\left( 1.3.6.40\right)%
\end{array}%
$

\bigskip

\textbf{Theorem} \textbf{1.3.6.16. }Let be $\left\{ a_{n}\right\}
_{n=1}^{\infty }$ arbitrary countable sequence

$a_{n}:%
\mathbb{N}
\rightarrow $ $^{\ast }%
\mathbb{R}
.$Then for every $b$ such that $b\in $ $^{\ast }%
\mathbb{R}
,b>0:$

\bigskip

$\ 
\begin{array}{cc}
\begin{array}{c}
\\ 
\left( b\right) ^{\#}\times \left( \overline{\overline{\#Ext\text{-}%
\dsum\limits_{n\in 
\mathbb{N}
}\left( a_{n}\right) ^{\#}}}\right) = \\ 
\\ 
\overline{\overline{\#Ext\text{-}\dsum\limits_{n\in 
\mathbb{N}
}\left( b\right) ^{\#}\times \left( a_{n}\right) ^{\#}}}, \\ 
\\ 
\left( b\right) ^{\#}\times \left( \underline{\underline{\#Ext-\dsum%
\limits_{n\in 
\mathbb{N}
}\left( a_{n}\right) ^{\#}}}\right) = \\ 
\\ 
\underline{\underline{\#Ext-\dsum\limits_{n\in 
\mathbb{N}
}\left( b\right) ^{\#}\times \left( a_{n}\right) ^{\#}}}. \\ 
\end{array}
& \text{ }\left( 1.3.6.41\right)%
\end{array}%
$

\textbf{Theorem} \textbf{1.3.6.17. }Suppose that $\lim_{n\rightarrow \infty
}\dsum\limits_{i=1}^{n}b_{i}=\zeta \in 
\mathbb{R}
.$Then for every $b$

such that $b\in $ $^{\ast }%
\mathbb{R}
,b>0:$\bigskip\ $\ \ \ $

$\ \ \ \ 
\begin{array}{cc}
\begin{array}{c}
\\ 
\left( b\right) ^{\#}\times \left( \overline{\overline{\#Ext\text{-}%
\dsum\limits_{n\in 
\mathbb{N}
}\left( ^{\ast }a_{n}\right) ^{\#}}}\right) = \\ 
\\ 
\overline{\overline{\#Ext\text{-}\dsum\limits_{n\in 
\mathbb{N}
}\left( b\right) ^{\#}\times \left( ^{\ast }a_{n}\right) ^{\#}}}= \\ 
\\ 
=\left( b\right) ^{\#}\times \left( ^{\ast }\zeta \right) ^{\#}+\left(
b\right) ^{\#}\times \varepsilon _{\mathbf{d}}, \\ 
\\ 
\left( b\right) ^{\#}\times \left( \underline{\underline{\#Ext\text{-}%
\dsum\limits_{n\in 
\mathbb{N}
}\left( ^{\ast }a_{n}\right) ^{\#}}}\right) = \\ 
\\ 
\underline{\underline{\#Ext\text{-}\dsum\limits_{n\in 
\mathbb{N}
}\left( b\right) ^{\#}\times \left( ^{\ast }a_{n}\right) ^{\#}}}= \\ 
\\ 
=\left( b\right) ^{\#}\times \left[ \left( ^{\ast }\zeta \right)
^{\#}-\varepsilon _{\mathbf{d}}\right] . \\ 
\end{array}
& \text{ \ \ \ \ \ }\left( 1.3.6.42\right)%
\end{array}%
$

\section{I.3.7.The construction non-archimedean field $\ \ \ \ \ \ \ \ \ \ \
\ \ \ \ \ \ \ \ \ \ $ \ \ \ \ \ \ \ \ \ \ \ \ \ \ \ \ \ \ \ \ \ \ $\ ^{\ast }%
\mathbb{R}
_{\mathbf{d}}^{\protect\omega }$ as Dedekind completion of countable \ \ \ \
\ \ \ \ \ \ \ \ \ \ \ \ \ \ \ non-standard \ models of field $%
\mathbb{R}
.$\ }

\bigskip

Let $^{\ast }%
\mathbb{R}
_{\omega }$ be a countable field which is elementary equivalent, but not\ \
\ \ \ \ \ \ \ \ \ \ \ \ \ \ \ \ 

isomorphic to $%
\mathbb{R}
.$

\textbf{Remark.1.3.7.1. }The \textquotedblleft elementary
equivalence\textquotedblright\ means that an (arithmetic) \ 

expression of first order is true in field $^{\ast }%
\mathbb{R}
_{\omega }$ if and only if it is true in field $%
\mathbb{R}
.$

Note that any non-standard model of $%
\mathbb{R}
$ contains an element $\mathbf{\upsilon \in }^{\ast }%
\mathbb{R}
_{\omega }$ such \ \ \ \ \ \ \ \ \ \ \ \ \ \ \ \ \ 

that $\mathbf{\upsilon }>x$ for each $x\in 
\mathbb{R}
.$

The canonical way to construct a model for $^{\ast }%
\mathbb{R}
_{\omega }$ uses model theory [31]. We

simply take as axioms all axioms of $%
\mathbb{R}
$ and additionally the following countable

number of axioms: the existence of an element $\mathbf{\upsilon }$ with $%
\mathbf{\upsilon }>1,\mathbf{\upsilon }>2,...,\mathbf{\upsilon }>n,...%
\mathbf{.}$

Each finite subset of this axioms is satisfied by the standard $%
\mathbb{R}
$.By the

compactness theorem in first order model theory, there exists a model which

also satisfies the given infinite set of axioms. By the theorem of L\"{o}%
wenheim-

Skolem, we can choose such models of countable cardinality.

\bigskip Each non-standard model $^{\ast }%
\mathbb{R}
$ contains the (externally defined) subset

\bigskip\ $%
\begin{array}{cc}
\begin{array}{c}
\\ 
^{\ast }%
\mathbb{R}
_{\mathbf{fin}}\triangleq \left\{ x\in \text{ }^{\ast }%
\mathbb{R}
|\exists n_{n\in 
\mathbb{Q}
}\left[ -n\leq x\leq n\right] \right\} . \\ 
\end{array}
& \text{ \ \ }\left( 1.3.7.1\right)%
\end{array}%
$

\bigskip

Every element $x\in $ $^{\ast }%
\mathbb{R}
_{\mathbf{fin}}$ defines a Dedekind cut:

\bigskip

$\ \ 
\begin{array}{cc}
\begin{array}{c}
\\ 
\mathbb{R}
=\left\{ y\in 
\mathbb{R}
|\text{ }y\leq x\right\} \cup \left\{ y\in 
\mathbb{R}
|y>x\right\} . \\ 
\end{array}
& \text{ \ \ \ }\left( 1.3.7.2\right)%
\end{array}%
$

\bigskip

We therefore get a order preserving map $\mathbf{j}_{p}\mathbf{:}^{\ast }%
\mathbb{R}
_{\mathbf{fin}}\rightarrow 
\mathbb{R}
$ which restricts to

the standard inclusion of the standard irrationals and which

respects addition and multiplication. An element of $^{\ast }%
\mathbb{R}
_{\mathbf{fin}}$ is called

infinitesimal,if it is mapped to $0$ under the map $\mathbf{j}_{p}.$

\textbf{Proposition} [30].\textbf{1.3.7.1.}Choose an arbitrary subset $%
M\subset 
\mathbb{R}
.$Then

(\textbf{i}) \ there is a model $^{\ast }%
\mathbb{R}
^{M}$ such that $\mathbf{j}_{p}\left( ^{\ast }%
\mathbb{R}
_{\mathbf{fin}}^{M}\right) \supset M.$

(\textbf{ii}) the cardinality of $^{\ast }%
\mathbb{Q}
^{M}$ can be chosen to coincide with $\mathbf{card}\left( M\right) $,if $M$
is

\ \ \ \ infinite.

\textbf{Proof.} Choose $M\subset 
\mathbb{R}
$. For each $m\in M$ choose $q_{1}^{m}<q_{2}^{m}<...<...<p_{2}^{m}<p_{1}^{m}$

with $\lim_{k\rightarrow \infty }q_{k}^{m}=\lim_{k\rightarrow \infty
}p_{k}^{m}=m.$

We add to the axioms of $%
\mathbb{R}
$ the following axioms:$\forall m\in M$ $\exists e_{m}$ such that

$q_{k}^{m}<e_{m}<p_{k}^{m}$ for all $k\in 
\mathbb{N}
.$

Again, the standard $%
\mathbb{R}
$ is a model for each finite subset of these axioms,

so that the compactness theorem implies the existence of $^{\ast }%
\mathbb{R}
^{M}$ as required,

where the cardinality of $^{\ast }%
\mathbb{R}
^{M}$ can be chosen to be the cardinality of the set of

axioms, i.e. of $M,$ if $M$ is infinite. Note that by construction $\mathbf{j%
}_{p}\left( e_{m}\right) =e_{m}.$

$\mathbf{Remark.}$\textbf{1.4.2.2.3.}It follows in particular that for each
countable subset of $%
\mathbb{R}
$

we can find a countable model of $^{\ast }%
\mathbb{R}
$ such that the image of fp contains this

subset. Note, on the other hand, that the image will only be countable, so

that the different models will have very different ranges.

\bigskip

\textbf{Definition 1.4.2.2.1.}[30].\textbf{\ }A Cauchy sequence in $^{\ast }%
\mathbb{R}
_{\omega }$ is a sequence $\left( a_{k}\right) _{k\in 
\mathbb{N}
}$

such that for every $\varepsilon \in $ $^{\ast }%
\mathbb{R}
_{\omega },$ $\varepsilon >0$ there is an $n_{\varepsilon }\in 
\mathbb{N}
$ such that:

$\forall m_{m\text{ }>\text{ }n_{\varepsilon }}\forall n_{n\text{ }>\text{ }%
n_{\varepsilon }}\left[ \left\vert \text{ }a_{m}-a_{n}\right\vert
<\varepsilon \right] .$

\textbf{Definition 1.4.2.2.2. }We define Cauchy completion $^{\ast }%
\mathbb{R}
_{\mathbf{c}}^{\omega }\triangleq \left[ ^{\ast }%
\mathbb{R}
_{\omega }\right] _{\mathbf{c}}$ in the

canonical way as equivalence classes of Cauchy sequences.

\ \ \ \ \ \ \ \ \ \ \ \ \ \ 

$\mathbf{Remark.}$\textbf{1.4.2.2.4. }This is a standard construction and
works for all ordered

fields.The result is again a field, extending the original field. Note that,
in

our case,each point in $^{\ast }%
\mathbb{R}
_{\mathbf{c}}^{\omega }$ is infinitesimally close to a point in $^{\ast }%
\mathbb{R}
$.

$\mathbf{Remark.}$\textbf{1.4.2.2.5. }In many non-standard models of $%
\mathbb{R}
$, there are no countable

sequences $\left( a_{k}\right) _{k\in 
\mathbb{N}
}$ tending to zero which are not eventually zero.

\bigskip \textbf{Proposition} [30].\textbf{1.3.7.2. }Assume that $^{\ast }%
\mathbb{R}
$ is countable.

\bigskip $\ \ \ \ \ \ \ \ \ \ \ \ \ \ \ \ \ \ \ $

\bigskip \bigskip \bigskip

\section{I.4.The construction non-archimedean field $^{\ast }%
\mathbb{R}
_{\mathbf{c}}.$}

\bigskip \bigskip

\section{I.4.1.Completion of ordered group and fields in general by using
'Cauchy pregaps'.}

We cketch here the aspects of the general theory that is concerned with

completion ordered group and fields, to be constructed by using 'Cauchy

pregaps' [32].

Throughout in this section we shall only consider fields which are
associative,

linear algebra over ground field $\Bbbk ,$ where $\Bbbk =%
\mathbb{Q}
,%
\mathbb{R}
,^{\ast }%
\mathbb{Q}
,^{\ast }%
\mathbb{R}
.$

\bigskip

\section{I.4.1.1.Totally ordered group and fields}

\textbf{Definition 1.4.1.1.1. }Let $\left( K,+,\cdot \right) $ be a field
and let $\left( \circ \leq \circ \right) $ be a binary

relation on $K.$Then $\left( K,+,\cdot ,\leq \right) $ is an ordered field if

(\textbf{i}) $\ \ \left( K,\leq \right) $ is totally ordered set,

(\textbf{ii})$\ \ \left( K,+,\leq \right) $ is an ordered group and

\bigskip (\textbf{iii})$\ a,b\in K^{+}\implies a\cdot b\in K^{+}.$

Note the standard convention that the order $\leq $ on an ordered field $K$

is necessarily a total order.

Let $K$ be an ordered field. It is easy to see from Definition 1.4.1.1. that

$\left\vert a\right\vert \cdot \left\vert b\right\vert =\left\vert a\cdot
b\right\vert .$

Let $K,L$ be an ordered fields. An imbeding of $K$ in $L$ is an algebra

monomorphism from $K$ into $L$ which is isotonic. A surjective embedding

is an isomorphism. In the case where exist an isomorphism $K$ onto $L$ then

$K$ and $L$ are isomorphic, and we write $K\cong L.$

\textbf{Definition 1.4.1.1.2. }Let $A$ be an algebra. Then:

\textbf{1.} $A\left[ X\right] $ denotes the algebra of polinomials $p\left(
X\right) $ with coefficients in $A;$

\textbf{2.} $^{\ast }A\left[ X\right] $ denotes the algebra of
hyperpolinomials $P\left( X\right) $ with coefficients in

$^{\ast }A;$

\textbf{3. }in the case where $A$ is a subalgebra of an algebra $B$ and $%
b\in B,$

$\ \ \ A\left[ b\right] \triangleq \left\{ p\left( b\right) |p\in A\left[ X%
\right] \right\} ;$

\textbf{4.} in the case where $A$ is a subalgebra of an algebra $B$ and

$\ \ \ b\in $ $^{\ast }B,^{\ast }A\left[ b\right] \triangleq \left\{ P\left(
b\right) |P\in \text{ }^{\ast }A\left[ X\right] \right\} ;$

\textbf{Definition 1.4.1.1.3.}Let\textbf{\ }$A$ be a subalgebra of an
algebra $B$ and $b\in B.$

Then:

\textbf{1.} $b\in B$ is algebraic over $A$ if there exist $p\in A\left[ X%
\right] \backslash \left\{ 0\right\} $ with $p\left( b\right) =0;$

\textbf{2.} $b\in B$ is transcendental over $A$ if is not algebraic over $A;$

\textbf{3.} $b\in $ $^{\ast }B$ is hyperalgebraic over $A$ if there exist $%
P\in $ $^{\ast }A\left[ X\right] \backslash \left\{ 0\right\} $ with

$P\left( ^{\ast }b\right) =0;$

\textbf{4.} $b\in $ $^{\ast }B$ is hypertranscendental over $A$ if is not
hyperalgebraic over $A;$

\textbf{Definition 1.4.1.1.4.}Let\textbf{\ }$A$ be a subalgebra of an
algebra $B$ and $b\in B.$

Then:

\textbf{1.}$b\in B$ is $w$-$\mathbf{transcendental}$ $\mathbf{over}$ $A$ if:

(a) $b\in B$ is transcendental over $A$ and

(b) there exist $P\in $ $^{\ast }A\left[ X\right] \backslash \left\{
0\right\} $ with $P\left( ^{\ast }b\right) =0$

\ \ \ \ or with $P\left( ^{\ast }b\right) \approx 0;$

\textbf{2. }$b\in B$ is $\#$-$\mathbf{transcendental}$ $\mathbf{over}$ $A$
if:

(a) $b\in B$ is transcendental over $A$ and

\bigskip (b) there is no exist $P\in $ $^{\ast }A\left[ X\right] \backslash
\left\{ 0\right\} $ with $P\left( ^{\ast }b\right) \approx 0.$\ \ \ 

\textbf{Example 1.4.1.1.1. }Number\textbf{\ }$\pi \in 
\mathbb{R}
$ is $w$-transcendental over$\ \ 
\mathbb{Q}
.$

There exist $P_{\pi }\left( X\right) \in $ $^{\ast }%
\mathbb{Q}
\left[ X\right] \backslash \left\{ 0\right\} $ with $P_{\pi }\left( ^{\ast
}\pi \right) \approx 0$ where\ 

\bigskip

$\ \ 
\begin{array}{cc}
\begin{array}{c}
\\ 
P_{\pi }\left( X\right) =\left[ \sin \left( X\right) \right] _{N\in \text{ }%
^{\ast }%
\mathbb{N}
_{\infty }}= \\ 
\\ 
\left( ^{\ast }\sum_{m=1}^{N\in \text{ }^{\ast }%
\mathbb{N}
_{\infty }}\dfrac{\left( -1\right) ^{m-1}X^{2m-1}}{\left( 2m-1\right) !}%
\right) . \\ 
\end{array}
& \left( 1.4.1.1.1\right)%
\end{array}%
$

\bigskip

\textbf{Example 1.4.1.1.2. }Number\textbf{\ }$\ln 2\in 
\mathbb{R}
$ is $w$-transcendental over$\ \ 
\mathbb{Q}
.$

There exist $P_{\ln 2}\left( X\right) \in $ $^{\ast }%
\mathbb{Q}
\left[ X\right] \backslash \left\{ 0\right\} $ with $P_{\ln 2}\left( \ln
2\right) -2\approx 0$ where

\bigskip

\ \ $\ \ 
\begin{array}{cc}
\begin{array}{c}
\\ 
P_{\ln 2}\left( X\right) =\left[ \exp \left( X\right) \right] _{N\in \text{ }%
^{\ast }%
\mathbb{N}
_{\infty }} \\ 
\\ 
=\left( ^{\ast }\sum_{m=1}^{N\in \text{ }^{\ast }%
\mathbb{N}
_{\infty }}\dfrac{X^{m}}{m!}\right) . \\ 
\end{array}
& \text{ \ \ }\left( 1.4.1.1.2\right)%
\end{array}%
$

\bigskip

\textbf{Definition 1.4.1.1.4.}Let\textbf{\ }$A$ be a subalgebra of an
algebra $B.$Then:

(\textbf{i}) The algebra $B$ is an algebraic extension of $A$ if each

$\ b\in B$ is algebraic over $A;$otherwise $B$ is transcendental extension
of $A.$

(\textbf{ii}) The algebra $^{\ast }B$ is an hyperalgebraic extension of $A$
if each

$\ b\in $ $^{\ast }B$ is hyperalgebraic over $A;$otherwise $B$ is
hypertranscendental

extension of $A.$

(\textbf{iii}) The algebra $^{\ast }B$ is an $w$-transcendental extension of 
$A$ if each

$\ b\in $ $^{\ast }B$ is $w$-transcendental over $A.$

(\textbf{iv}) \ The algebra $^{\ast }B$ is an $\#$-transcendental extension
of $A$ if each

$\ b\in $ $^{\ast }B$ is $\#$-transcendental over $A.$

\textbf{Definition 1.4.1.1.5. }Let\textbf{\ }$K$ be a field.

(\textbf{i}) A field $K$ is algebraically closed iff there is no

\ \ \ field $L$ which is a proper algebraic extension of $K,$or,equivalently,

$\ \ \ K$ is algebraically closed iff each non constant $p\in K\left[ X%
\right] $ has a root in $K.$

(\textbf{ii}) A field $K$ is hyperalgebraically closed iff there is no

\ \ \ field $L$ which is a proper hyperalgebraic extension of $K,$%
or,equivalently,

$\ \ \ K$ is hyperalgebraically closed iff each non constant $P\in $ $^{\ast
}K\left[ X\right] $ has a root

in $K,$ i.e. $P\left( ^{\ast }b\right) =0$ fore some $b\in K.$

\textbf{Definition 1.4.1.1.6. }An ordered field $K$ is real-closed if

(a) it has no proper algebraic extension to an ordered field,

or, equivalently,if

(b) the complexification $K_{%
\mathbb{C}
}$ of $K$ is algebraically closed,

or, equivalently,if

(c) every positive element in $K$ is a square and every polinomial over $K$

of odd degree has a root in $K.$

\bigskip

Let $K$ be a real-closed ordered field and take some $c\in K^{+}\backslash
\left\{ 0\right\} $ and

$n\in 
\mathbb{N}
.$ There is a anique element $b\in K^{+}\backslash \left\{ 0\right\} $ such
that $b^{n}=c,$ and

so there is a map $\psi :%
\mathbb{Q}
_{+}\rightarrow K^{+}\backslash \left\{ 0\right\} $\bigskip\ where $\psi
:\alpha \longmapsto c^{\alpha }.$ Thus\ 

\bigskip

\ \ \ \ \ \ \ \ \ \ \ \ \ \ \ \ \ \ \ \ \ \ \ \ \ \ 

\bigskip\ $\ \ \ 
\begin{array}{cc}
\begin{array}{c}
\\ 
\psi \left( 0\right) =1,\psi \left( 1\right) =c, \\ 
\\ 
\psi \left( \alpha +\beta \right) =\psi \left( \alpha \right) \cdot \psi
\left( \beta \right) ,\left( \alpha ,\beta \in 
\mathbb{Q}
_{+}\right) . \\ 
\end{array}
& \text{ \ }\left( 1.4.1.1.3\right)%
\end{array}%
$

\bigskip

\textbf{\ }

\section{I.4.1.2.Cauchy completion of ordered group and fields.}

\textbf{Definition 1.4.1.2.1.}Let $\left\langle A,B\right\rangle $ be a
pregap in totally ordered group $G$. Then

$\left\langle A,B\right\rangle $ is a Cauchy pregap if $A$ has no maximum,$B$
has no minimum, and,for

each $\varepsilon >0,\varepsilon \in G,$there exist $a\in A$ and $b\in B$
with $b<a+\varepsilon .$

\textbf{Definition 1.4.1.2.2. }The group $G$ is Cauchy complete if, for each
Cauchy

pregap there exists $x\in G$ with $a<<x<<b.$

\textbf{Remark.1.4.1.1.}Thus\textbf{\ }totally ordered group $G$ is Cauchy
complete iff there

are no Cauchy gaps.

\textbf{Remark.1.4.1.2.}The element $x$ arising in the above definition is
necessarily

unique.

\textbf{Example} \textbf{1.4.1.1. (i) }The group $\left( 
\mathbb{R}
,+\right) $ is certainly Cauchy complete.

\textbf{(ii) }The group $\left( ^{\ast }%
\mathbb{R}
,+\right) $ is Cauchy complete.

\textbf{(iii) }The monoid $\left( ^{\ast }%
\mathbb{R}
_{\mathbf{d}},+\right) $ is certainly Cauchy complete.

\textbf{Definition 1.4.1.3. }The set of the all Cauchy pregaps in totally
ordered

group $G$ we denote by $C\left( G\right) .$

\textbf{Definition 1.4.1.4. }A\textbf{\ }totally ordered group $G$ is
discrete if the set $G^{+}\backslash \left\{ 0\right\} $

is empty or has a minimum element and is non-discrete otherwise.

\bigskip

For any $\left\langle A_{1},B_{1}\right\rangle \in C\left( G\right) $ and $%
\left\langle A_{2},B_{2}\right\rangle \in Cp\left( G\right) $ we hawe

$\left\langle A_{1}+A_{2},B_{1}+B_{2}\right\rangle \in Cp\left( G\right) .$%
Let us define sum of the classes

$\left[ \left\langle A_{1},B_{1}\right\rangle \right] \in Cl\left[ Cp\left(
G\right) \right] $ and $\left[ \left\langle A_{2},B_{2}\right\rangle \right]
\in Cl\left[ Cp\left( G\right) \right] $ by formula \bigskip

\bigskip

\ $\ 
\begin{array}{cc}
\begin{array}{c}
\\ 
\left[ \left\langle A_{1},B_{1}\right\rangle \right] +\left[ \left\langle
A_{2},B_{2}\right\rangle \right] =\left[ \left\langle
A_{1}+A_{2},B_{1}+B_{2}\right\rangle \right] . \\ 
\end{array}
& \text{ \ \ }\left( 1.4.1.2.1\right)%
\end{array}%
$

\bigskip

Then $+$ is well defined in $H_{G}=Cl\left[ Cp\left( G\right) \right] $ and $%
\left\{ Cl\left[ Cp\left( G\right) \right] ,+\right\} =$

$\left\{ H_{G},+\right\} $ is an abelian group.The map

$\ 
\begin{array}{cc}
\begin{array}{c}
\\ 
\iota :\left\{ G,+\right\} \hookrightarrow \left\{ H_{G},+\right\} \\ 
\end{array}
& \text{ \ \ \ }\left( 1.4.1.2.2\right)%
\end{array}%
$

\bigskip

is a canonical group morphism. It is easy to see that $\left\{
H_{G},+\right\} $ is a

totally ordered group. Let $\left\langle \hat{A},\hat{B}\right\rangle $ be a
Cauchy pregap in $H_{G},$and

define

\bigskip

\ \ \ \ \ \ \ \ \ \ \ \ \ \ \ \ \ \ \ \ \ \ \ \ \ \ \ \ \ \ \ \ \ \ \ \ \ \
\ \ \ \ $\ \ \ \ \ \ 
\begin{array}{cc}
\begin{array}{c}
\\ 
\mathbf{A}=\bigcup \left\{ A_{\alpha }|\left[ \left( A_{\alpha },B_{\alpha
}\right) \right] \in \hat{A}\text{ for some }B_{\alpha }\right\} \\ 
\\ 
\mathbf{B}=\bigcup \left\{ B_{\alpha }|\left[ \left( A_{\alpha },B_{\alpha
}\right) \right] \in \hat{B}\text{ for some }A_{\alpha }\right\} , \\ 
\end{array}
& \text{ \ \ }\left( 1.4.1.2.3\right)%
\end{array}%
$

\bigskip

hence $\left\langle \mathbf{A,B}\right\rangle \in Cp\left( G\right) $ and $%
\hat{A}<<\left[ \left( \mathbf{A,B}\right) \right] <<\hat{B}$ and $\left\{
H_{G},+\right\} $ is

Cauchy complete. On it we have the following result:

\textbf{Theorem} \textbf{1.4.1.1.}[32].Let $G$ be a totally ordered
non-discrete group.

Then the group $\left\{ H_{G},+\right\} $ is defined by formula (1.4.1.1) is
a Cauchy

completion of $G.$

\textbf{Definition 1.4.1.5. }An ordered field $K$ is Cauchy complete if the

totally ordered group $\left\{ K,+\right\} $ is Cauchy complete.

\textbf{Theorem} \textbf{1.4.1.2.}[32].Let $K$ be an ordered field.Then the
Cauchy

completion $\widetilde{K}$ of the group $\left\{ K,+\right\} $ can be made
into an ordered field

in such a way that $K$ is a subfield of $\widetilde{K}.$ If $K$ is
real-closed, then so is $\widetilde{K}.$

\textbf{Proof.} Suppose $a_{1},a_{2}\in \widetilde{K},a_{1}>0,a_{2}>0$ and $%
a_{1}=\left[ \left( A_{1},B_{1}\right) \right] ,$

$a_{2}=\left[ \left( A_{2},B_{2}\right) \right] .$We may suppose that $%
A_{1},A_{2}\subset K^{+}\backslash \left\{ 0\right\} .$Set

\bigskip

\ $\ \ 
\begin{array}{cc}
\begin{array}{c}
\\ 
A=\left\{ x_{1}\cdot x_{2}|x_{1}\in A_{1},x_{2}\in A_{2}\right\} , \\ 
\\ 
B=\left\{ y_{1}\cdot y_{2}|y_{1}\in B_{1},y_{2}\in B_{2}\right\} . \\ 
\end{array}
& \text{ \ \ }\left( 1.4.1.2.4\right)%
\end{array}%
$

\bigskip

\bigskip Then $\left\langle A,B\right\rangle \in Cp\left( \left\{
K,+\right\} \right) .$Define

$\bigskip $

$\ 
\begin{array}{cc}
\begin{array}{c}
\\ 
a_{1}\cdot a_{2}=\left[ \left\langle A_{1},B_{1}\right\rangle \right] \cdot %
\left[ \left\langle A_{2},B_{2}\right\rangle \right] \\ 
\end{array}
& \text{ \ \ \ }\left( 1.4.1.2.5\right)%
\end{array}%
$

\bigskip

The operation $\left( \circ \cdot \circ \right) $ is well defined in $%
\widetilde{K}$ and $a_{1}\cdot a_{2}>0.$If $a_{1}<0,$

$a_{2}>0$ set $a_{1}\cdot a_{2}=-\left( \left( -a_{1}\right) \cdot
a_{2}\right) ,$ etc.

It is simple to check that the group $\left\{ \widetilde{K},+\right\} $
together with product

$\left( \circ \cdot \circ \right) $ is an ordered field $\widetilde{K}%
\triangleq \left\{ \widetilde{K},+,\cdot \right\} $ with the required
properties.

Let $K^{\prime }$ be any ordered field containing $K$ as an order-dense

subfield. Then there is an isotonic morphism from $K^{\prime }$ into $%
\widetilde{K},$ and so

$\widetilde{K}$ is the maximum ordered field containing $K$ as an proper

order-dense subfield. On it we have the following result:

\textbf{Theorem} \textbf{1.4.1.3.}[32].Let $K$ be an ordered field.Then $K$
is Cauchy

complete iff no ordered field $L$ containing $K$ as an proper

order-dense subfield.

Standard main tool for anderstanding structure of the totally ordered

external group well be Hahn's embedding theorem. We shall associate to

a totally ordered external and internal group $\left( G,+,\leq \right) $ a
'value set' $\Gamma ^{\#},$

define a collections $\tciFourier \left( 
\mathbb{R}
,\Gamma \right) ,\tciFourier \left( ^{\ast }%
\mathbb{R}
,\Gamma ^{\#}\right) $ of 'formal power series' over $\Gamma $ and

$\Gamma ^{\#}$ and imbed $G$ into $\tciFourier \left( 
\mathbb{R}
,\Gamma \right) $ or $\tciFourier \left( ^{\ast }%
\mathbb{R}
,\Gamma ^{\#}\right) $ correspondingly.

Let us define the value set of external group $G.$

\textbf{Definition 1.4.1.6.}[32].\textbf{\ }Let $\left( G,+,\leq \right) $
be a totally ordered external group,

i.e. $\left( G,+,\leq \right) \in V^{Ext}$ and let $x,y\in G.$Set:

(\textbf{i}) \ \ $x=o\left( y\right) $ iff $\forall n_{n\in 
\mathbb{N}
}\left[ n\cdot \left\vert x\right\vert \leq \left\vert y\right\vert \right]
; $

(\textbf{ii}) \ \ $x=O\left( y\right) $ iff $\exists m_{m\in 
\mathbb{N}
}\left[ \left\vert x\right\vert \leq m\cdot \left\vert y\right\vert \right]
; $

\bigskip (\textbf{iii}) \ $x$ $\symbol{126}$ $y$ iff $\left[ x=O\left(
y\right) \right] \wedge \left[ y=O\left( x\right) \right] .$

For each $y\in G$ the sets $\left\{ x|x=o\left( y\right) \right\} $ and $%
\left\{ x|x=O\left( y\right) \right\} $ are absolutely

convex subsets of $G.$Is clear that $\left( \circ \text{ }\symbol{126}\circ
\right) $ is an equivalence relation on $G.$

Each $\symbol{126}$ -equivalence class (other than $\left\{ 0\right\} $) is
the union of an interval

contained in $G^{+}\backslash \left\{ 0\right\} $ and an interval contained
in $G^{-}\backslash \left\{ 0\right\} .$

\textbf{Definition 1.4.1.7.}Let $\left( G,+,\leq \right) $ be a totally
ordered external group.

The set $\Gamma =\Gamma _{G}=\left( G\backslash \left\{ 0\right\} \right) /%
\symbol{126}$ \ of equivalence classes in the value set

of $G$ and the elements of $\Gamma $ are the archimedian classes of $G.$

The quotient map from $G\backslash \left\{ 0\right\} $ onto $\Gamma $ is
denoted by $v:$ $G\backslash \left\{ 0\right\} \rightarrow \Gamma .$It is

the archimedian valuation on $G.$Set $v\left( x\right) \leq v\left( y\right) 
$ for $x,y\in G\backslash \left\{ 0\right\} $ iff

$y=O\left( x\right) .$

It is easy checked that $\leq $ is well defined on $\Gamma ,$hence $\left(
\Gamma ,\leq \right) $ is a totally

ordered set, such that

\bigskip

$\ \ \ 
\begin{array}{cc}
\begin{array}{c}
\\ 
v\left( x+y\right) \geq \min \left\{ v\left( x\right) ,v\left( y\right)
\right\} , \\ 
\\ 
v\left( x+y\right) =\min \left\{ v\left( x\right) ,v\left( y\right) \right\} 
\text{ if } \\ 
\\ 
\left[ v\left( x\right) \neq v\left( y\right) \right] \vee \left[ \left(
x\in G^{+}\right) \wedge \left( y\in G^{+}\right) \right] . \\ 
\end{array}
& \text{\ }\left( 1.4.1.2.6\right)%
\end{array}%
$

\bigskip

\textbf{Definition 1.4.1.8. }Let $\left( \breve{G},+,\leq \right) $ a
totally ordered internal group,i.e.

$\left( \breve{G},+,\leq \right) \in V^{Int}$ and let $x,y\in G.$In
particular $\breve{G}=$ $^{\ast }\left( G,+,\leq \right) $ for some

$\left( G,+,\leq \right) ,$i.e. in particular $\breve{G}$ is an standard
group.

Set:\bigskip

\bigskip

\section{I.4.2.1.The construction non-archimedean \ \ \ \ \ \ \ \ \ \ \ \ \
\ \ \ \ \ \ \ \ \ \ \ \ \ \ \ \ \ \ \ field $^{\ast }%
\mathbb{R}
_{\mathbf{c}}$ by using \ Cauchy hypersequence \ \ \ \ \ \ \ \ \ \ \ \ \ \ \
\ \ \ \ \ \ \ \ \ \ \ \ \ \ \ \ \ \ \ \ in ancountable field $^{\ast }%
\mathbb{Q}
.$}

\bigskip

Let $^{\ast }%
\mathbb{Q}
_{\omega _{\alpha }}\triangleq $ $^{\ast }%
\mathbb{Q}
,\omega <\omega _{\alpha }$ be a ancountable field which is elementary

equivalent, to $%
\mathbb{Q}
.$The \textquotedblleft elementary equivalence\textquotedblright\ means that
an

(arithmetic) expression of first order is true in field $^{\ast }%
\mathbb{Q}
_{\omega _{\alpha }}$ if and only

if it is true in field $%
\mathbb{Q}
.$Note that any non-standard model of $%
\mathbb{Q}
$ contains

an element $\mathbf{e\in }^{\ast }%
\mathbb{Q}
_{\omega }$ such that $\mathbf{e}>q$ for each $q\in 
\mathbb{Q}
.$

We define Cauchy completion $\left[ ^{\ast }%
\mathbb{Q}
_{\omega _{\alpha }}\right] _{\mathbf{c}}\triangleq $ $^{\ast }%
\mathbb{R}
_{\mathbf{c}}$ in the canonical

way as equivalence classes of Cauchy hypersequences.

$\mathbf{Remark.}$\textbf{1.4.2.1.1.}This is a general construction and
works for all

nonstandard ordered fields $^{\ast }\Bbbk $.The result is again a field $%
\left[ ^{\ast }\Bbbk \right] _{\mathbf{c}}$ which

is potentially different from extending the original field $\Bbbk $, and we

actually see that $^{\ast }%
\mathbb{R}
_{\mathbf{c}}$ is different from $^{\ast }%
\mathbb{Q}
.$

$\mathbf{Remark.}$\textbf{1.4.2.1.2.}In many non-standard ancountable models
of $%
\mathbb{Q}
$,

there are no countable sequences tending to zero which are not

eventually zero.Thus dealing with analysis over field $^{\ast }%
\mathbb{R}
$ we are

compelled to enter into consideration hypersequences of various

classes: $\mathbf{s}_{\mathbf{n\in }^{\ast }%
\mathbb{N}
}:$ $^{\ast }%
\mathbb{N}
\rightarrow $ $^{\ast }%
\mathbb{Q}
,\mathbf{s}_{\mathbf{n\in }^{\ast }%
\mathbb{N}
}:$ $^{\ast }%
\mathbb{N}
\rightarrow $ $^{\ast }%
\mathbb{R}
$ and $\mathbf{s}_{\mathbf{n\in }^{\ast }%
\mathbb{N}
}:$ $^{\ast }%
\mathbb{N}
\rightarrow $ $^{\ast }%
\mathbb{R}
_{\mathbf{d}}$

\textbf{Definition 1.4.2.1.1. }A hypersequence

$\mathbf{s}_{\mathbf{n}}:$ $^{\ast }%
\mathbb{N}
\rightarrow $ $^{\ast }%
\mathbb{R}
_{\mathbf{d}}\supset $ $^{\ast }%
\mathbb{R}
\supset $ $^{\ast }%
\mathbb{Q}
,\mathbf{n\in }^{\ast }%
\mathbb{N}
$ tends

to a $\ast $-limit $\alpha $\ ($\alpha \in ^{\ast }%
\mathbb{Q}
,^{\ast }%
\mathbb{R}
$ or $^{\ast }%
\mathbb{R}
_{\mathbf{d}}$) in $^{\ast }%
\mathbb{Q}
,^{\ast }%
\mathbb{R}
$ or $^{\ast }%
\mathbb{R}
_{\mathbf{d}}$ iff

\bigskip

$%
\begin{array}{cc}
\begin{array}{c}
\\ 
\exists \alpha \left( \alpha \in \text{ }^{\ast }%
\mathbb{R}
_{\mathbf{d}}\right) \left[ \forall \varepsilon _{\varepsilon >0}\left(
\varepsilon \in \text{ }^{\ast }%
\mathbb{R}
_{\mathbf{d}}\right) \exists \mathbf{n}_{0}\left( \mathbf{n}_{0}\in \text{ }%
^{\ast }%
\mathbb{N}
_{\infty }\right) \right. \\ 
\\ 
\left. \forall \mathbf{n}\left[ \mathbf{n\geqslant n}_{0}\implies \left\vert
\alpha -\mathbf{s}_{\mathbf{n}_{0}}\right\vert <\varepsilon \right] \right] .
\\ 
\end{array}
& \left( 1.4.2.1.1\right)%
\end{array}%
$

\bigskip \bigskip

We write $\ast $-$\lim_{\mathbf{n}\text{ }\longrightarrow \text{ \ }^{\ast
}\infty }\mathbf{s}_{\mathbf{n}_{\mathbf{n\in }^{\ast }%
\mathbb{N}
}}=\alpha $ $\ $or$\ \underset{\mathbf{n}\text{ }\longrightarrow \text{ \ }%
^{\ast }\infty }{\ast \text{-}\lim }\mathbf{s}_{\mathbf{n}}=\alpha \ \ $iff
condition

(1.4.1.1.1) is \ satisfied.

\textbf{Definition 1.4.2.1.2. }A hypersequence $\mathbf{s}_{\mathbf{n}}:$ $%
^{\ast }%
\mathbb{N}
\rightarrow $ $^{\ast }%
\mathbb{R}
_{\mathbf{d}}\supseteqq $ $^{\ast }%
\mathbb{Q}
$ is

divergent in $^{\ast }%
\mathbb{R}
_{\mathbf{d}},$ or tends to $^{\ast }\infty $ iff

$\bigskip $

$\ \ 
\begin{array}{cc}
\begin{array}{c}
\\ 
\forall r_{r>0}\left( r\in \text{ }^{\ast }%
\mathbb{R}
\right) \exists \mathbf{n}_{0}\left( \mathbf{n}_{0}\in \text{ }^{\ast }%
\mathbb{N}
_{\infty }\right) \forall \mathbf{n}\left[ \mathbf{n}\geqslant \mathbf{n}%
_{0}\right. \\ 
\\ 
\left. \implies \left\vert \mathbf{s}_{\mathbf{n}}\right\vert >r\right] .\ \ 
\\ 
\end{array}
& \text{\ }\left( 1.4.2.1.2\right)%
\end{array}%
$

\ \ \ \ \ 

\textbf{Lemma 1.4.2.1.1. }Suppose that $\mathbf{s}_{\mathbf{n}}:$ $^{\ast }%
\mathbb{N}
\rightarrow $ $^{\ast }%
\mathbb{R}
\supset $ $^{\ast }%
\mathbb{Q}
,\mathbf{n\in }^{\ast }%
\mathbb{N}
.$

(a) If\textbf{\ }$\underset{\mathbf{n}\text{ }\longrightarrow \text{ \ }%
^{\ast }\infty }{\ast \text{-}\lim }\mathbf{s}_{\mathbf{n}}$ exists in $%
^{\ast }%
\mathbb{R}
,$ then it is unique.

(b) If\textbf{\ }$\underset{\mathbf{n}\text{ }\longrightarrow \text{ \ }%
^{\ast }\infty }{\ast \text{-}\lim }\mathbf{s}_{\mathbf{n}}$ exists in $%
^{\ast }%
\mathbb{R}
_{\mathbf{d}},$ then it is unique.

That is if $\ast $-$\lim_{\mathbf{n}\text{ }\longrightarrow \text{ \ }^{\ast
}\infty }\mathbf{s}_{\mathbf{n}}=\alpha _{1},\ast -\lim_{\mathbf{n}\text{ }%
\longrightarrow \text{ \ }^{\ast }\infty }\mathbf{s}_{\mathbf{n}}=\alpha
_{2} $ then $\alpha _{1}=\alpha _{2}.$

\textbf{Proof. }(a) Let $\varepsilon $ be any positive number $\varepsilon
>0,\varepsilon \in $ $^{\ast }%
\mathbb{R}
\supset $ $^{\ast }%
\mathbb{Q}
$.Then, by

definition, we must be able to find a number $N_{1}$ so that $\left\vert 
\mathbf{s}_{\mathbf{n}}-\alpha _{1}\right\vert <\varepsilon $

whenever $\mathbf{n}\geq N_{1}.$

We must also be able to find a number $N_{2}$ so that $\left\vert \mathbf{s}%
_{\mathbf{n}}-\alpha _{2}\right\vert <\varepsilon $

whenever $\mathbf{n}\geq N_{2}.$ Take $\mathbf{m}$ to be the maximum of $%
N_{1}$ and $N_{2}.$ Then both

assertions $\left\vert \mathbf{s}_{\mathbf{m}}-\alpha _{1}\right\vert
<\varepsilon $ and $\left\vert \mathbf{s}_{\mathbf{m}}-\alpha
_{2}\right\vert <\varepsilon $ are true.

This by using triangle inequality allows us to conclude that

$\left\vert \alpha _{1}-\alpha _{2}\right\vert =\left\vert \left( \alpha
_{1}-\mathbf{s}_{\mathbf{m}}\right) +\left( \mathbf{s}_{\mathbf{m}}-\alpha
_{2}\right) \right\vert \leq \left\vert \alpha _{1}-\mathbf{s}_{\mathbf{m}%
}\right\vert +\left\vert \mathbf{s}_{\mathbf{m}}-\alpha _{2}\right\vert
<2\varepsilon .$ So that

$\left\vert \alpha _{1}-\alpha _{2}\right\vert <2\varepsilon .$ But $%
\varepsilon $ can be any positive infinite small number whatsoever.

This could only be true if $\alpha _{1}=\alpha _{2}$, which is what we
wished to show.

\textbf{Definition 1.4.2.1.3. }A \textit{Cauchy hypersequence in }$^{\ast }%
\mathbb{Q}
,^{\ast }%
\mathbb{R}
$ and $^{\ast }%
\mathbb{R}
_{\mathbf{d}}$ is a sequence $\ \ \ \ \ \ \ \ \ \ \ \ \ \ \ \ \ \ \ \ \ $

$\ \mathbf{s}_{\mathbf{n}}:$ $^{\ast }%
\mathbb{N}
\rightarrow $ $^{\ast }%
\mathbb{R}
_{\mathbf{d}}\supseteqq $ $^{\ast }%
\mathbb{Q}
$ with\ the following property: for every\ $\varepsilon \in $ $^{\ast }%
\mathbb{R}
_{\mathbf{d}}$ such that $\ \ \ \ \ \ \ \ \ \ \varepsilon >0,$there exists
an $\mathbf{n}_{0}\in $ $^{\ast }%
\mathbb{N}
_{\infty }$ such that $\mathbf{m,n}\geqslant n_{0}$ implies $|\mathbf{s}_{%
\mathbf{m}}-\mathbf{s}_{\mathbf{n}}|$ $<\varepsilon ,$

i.e.

\bigskip

\ 

\bigskip $\ \ 
\begin{array}{cc}
\begin{array}{c}
\\ 
\forall \varepsilon _{(\varepsilon \in ^{\ast }%
\mathbb{R}
_{\mathbf{d}})}\left( \varepsilon >0\right) \exists \mathbf{n}_{0}\left( 
\mathbf{n}_{0}\in \text{ }^{\ast }%
\mathbb{N}
_{\infty }\right) \left[ \mathbf{m,n}\geqslant n_{0}\right. \\ 
\\ 
\left. \implies |\mathbf{s}_{\mathbf{m}}-\mathbf{s}_{\mathbf{n}%
}|<\varepsilon \right] . \\ 
\end{array}
& \text{ \ \ }\left( 1.4.2.1.3\right)%
\end{array}%
$

\bigskip

\textbf{Lemma 1.4.2.1.2. }A hypersequence of numbers $^{\ast }%
\mathbb{Q}
,^{\ast }%
\mathbb{R}
$ and $^{\ast }%
\mathbb{R}
_{\mathbf{d}},$ that

converges is Cauchy hypersequence.\textbf{\ }

\textbf{Lemma 1.4.2.1.3. }A Cauchy hypersequence $\left( \mathbf{s}_{\mathbf{%
n}}\right) _{\mathbf{n\in }^{\ast }%
\mathbb{N}
}$ in $^{\ast }%
\mathbb{Q}
$ and $^{\ast }%
\mathbb{R}
$,is

bounded or hyperbounded.

\textbf{Proof.}Choose in (1.4.2.1) $\varepsilon =1.$ Since the sequence $%
\left( \mathbf{s}_{\mathbf{n}}\right) _{\mathbf{n\in }^{\ast }%
\mathbb{N}
}$ is Cauchy,

there exists a positive hyperinteger $N\in $ $^{\ast }%
\mathbb{N}
$ such that $\left\vert \mathbf{s}_{\mathbf{i}}-\mathbf{s}_{\mathbf{j}%
}\right\vert <1$whenever

$\mathbf{i,j\geq }N.$In particular,$\left\vert \mathbf{s}_{\mathbf{i}}-%
\mathbf{s}_{N}\right\vert <1$ whenever $\mathbf{i\geq }N$ By the triangle
inequality,

\bigskip $\left\vert \mathbf{s}_{\mathbf{i}}\right\vert -\left\vert \mathbf{s%
}_{N}\right\vert \leq \left\vert \mathbf{s}_{\mathbf{i}}-\mathbf{s}%
_{N}\right\vert $ and therefore,$\left\vert \mathbf{s}_{\mathbf{i}%
}\right\vert <\left\vert \mathbf{s}_{N}\right\vert +1$ for all $\mathbf{%
i\geq }$ $N.$

\textbf{Definition 1.4.2.1.4. }Cauchy hypersequences $(x_{\mathbf{n}})_{%
\mathbf{n}\in ^{\ast }%
\mathbb{N}
}$ and $(y_{\mathbf{n}})_{\mathbf{n}\in ^{\ast }%
\mathbb{N}
},$

can be added, multiplied and compared as follows:

(\textbf{a}) $(x_{\mathbf{n}})_{\mathbf{n}\in \text{ }^{\ast }%
\mathbb{N}
}+(y_{\mathbf{n}})_{\mathbf{n}\in \text{ }^{\ast }%
\mathbb{N}
}=(x_{\mathbf{n}}+y_{\mathbf{n}})_{\mathbf{n}\in \text{ }^{\ast }%
\mathbb{N}
},$

(\textbf{b}) $(x_{\mathbf{n}})_{\mathbf{n}\in \text{ }^{\ast }%
\mathbb{N}
}\times (y_{\mathbf{n}})_{\mathbf{n}\in \text{ }^{\ast }%
\mathbb{N}
}=(x_{\mathbf{n}}\times y_{\mathbf{n}})_{\mathbf{n}\in \text{ }^{\ast }%
\mathbb{N}
},$

(\textbf{c}) $\dfrac{(x_{\mathbf{n}})_{\mathbf{n}\in \text{ }^{\ast }%
\mathbb{N}
}}{(y_{\mathbf{n}})_{\mathbf{n}\in \text{ }^{\ast }%
\mathbb{N}
}}=\left( \dfrac{x_{\mathbf{n}}}{y_{\mathbf{n}}}\right) _{\mathbf{n}\in 
\text{ }^{\ast }%
\mathbb{N}
}$ iff $\forall \mathbf{n}\left( \mathbf{n}\in \text{ }^{\ast }%
\mathbb{N}
\right) $ $\left[ y_{n}\neq 0\right] ,$

(\textbf{d}) $(x_{\mathbf{n}})_{\mathbf{n}\in \text{ }^{\ast }%
\mathbb{N}
}^{-1}=(x_{\mathbf{n}}^{-1})_{\mathbf{n}\in \text{ }^{\ast }%
\mathbb{N}
}$ \ iff $\ \forall \mathbf{n}\in $ $^{\ast }%
\mathbb{N}
(y_{\mathbf{n}}\neq 0),$

(\textbf{e}) $(x_{\mathbf{n}})_{\mathbf{n}\in \text{ }^{\ast }%
\mathbb{N}
}\geq (y_{\mathbf{n}})_{\mathbf{n}\in \text{ }^{\ast }%
\mathbb{N}
}$ if and only if for every $\epsilon >0,\epsilon \in $ $^{\ast }%
\mathbb{Q}
$

there exists an integer $\mathbf{n}_{0}$ such that \bigskip $x_{\mathbf{n}%
}\geq y_{\mathbf{n}}-\epsilon $ for all $\mathbf{n>n}_{0}\mathbf{.}$

\textbf{Definition 1.4.2.1.5.}Two Cauchy hypersequences $(x_{\mathbf{n}})_{%
\mathbf{n}\in \text{ }^{\ast }%
\mathbb{N}
}$ and $(y_{\mathbf{n}})_{\mathbf{n}\in \text{ }^{\ast }%
\mathbb{N}
}$

are called equivalent: $(x_{\mathbf{n}})_{\mathbf{n}\in \text{ }^{\ast }%
\mathbb{N}
}$ $\approx _{\mathbf{c}}(y_{n})_{\mathbf{n}\in \text{ }^{\ast }%
\mathbb{N}
}$ if the hypersequence $\ \ \ \ \ \ \ \ \ \ \ \ \ \ \ \ \ \ \ \ \ \ \ \ \ \
\ \ \ \ \ \ \ \ \ \ \ \ \ \ $ $\ $

$(x_{\mathbf{n}}-y_{\mathbf{n}})_{\mathbf{n}\in \text{ }^{\ast }%
\mathbb{N}
}$has $\ast $-limit zero$,$i.e. $\ast $-$\lim_{\mathbf{n\rightarrow }^{\ast
}\infty }(x_{\mathbf{n}}-y_{\mathbf{n}})_{\mathbf{n}\in \text{ }^{\ast }%
\mathbb{N}
}=0.$

\textbf{Lemma 1.4.2.1.4. }If $(x_{\mathbf{n}})_{\mathbf{n}\in \text{ }^{\ast
}%
\mathbb{N}
}\approx _{\mathbf{c}}(x_{\mathbf{n}}^{\prime })_{\mathbf{n}\in \text{ }%
^{\ast }%
\mathbb{N}
}$ and $(y_{\mathbf{n}})_{\mathbf{n}\in \text{ }^{\ast }%
\mathbb{N}
}\approx _{\mathbf{c}}(y_{\mathbf{n}}^{\prime })_{\mathbf{n}\in \text{ }%
^{\ast }%
\mathbb{N}
},$

are two pairs of equivalent Cauchy hypersequences, then:

\textbf{(a)} hypersequence $(x_{\mathbf{n}}+y_{\mathbf{n}})_{\mathbf{n}\in 
\text{ }^{\ast }%
\mathbb{N}
}$ is Cauchy and\bigskip

$\ 
\begin{array}{cc}
\begin{array}{c}
\\ 
(x_{\mathbf{n}}+y_{\mathbf{n}})_{\mathbf{n}\in \text{ }^{\ast }%
\mathbb{N}
}\approx _{\mathbf{c}}(x_{\mathbf{n}}^{\prime }+y_{\mathbf{n}}^{\prime })_{%
\mathbf{n}\in \text{ }^{\ast }%
\mathbb{N}
}, \\ 
\end{array}
& \text{ \ }\left( 1.4.2.1.4\right)%
\end{array}%
\ \ $

\bigskip

\textbf{(b)} hypersequence $(x_{\mathbf{n}}-y_{\mathbf{n}})_{\mathbf{n}\in 
\text{ }^{\ast }%
\mathbb{N}
}$ is Cauchy and

\bigskip

$\ \ \ 
\begin{array}{cc}
\begin{array}{c}
\\ 
\ (x_{\mathbf{n}}-y_{\mathbf{n}})_{\mathbf{n}\in \text{ }^{\ast }%
\mathbb{N}
}\approx _{\mathbf{c}}(x_{\mathbf{n}}^{\prime }-y_{\mathbf{n}}^{\prime })_{%
\mathbf{n}\in \text{ }^{\ast }%
\mathbb{N}
}, \\ 
\end{array}
& \text{ \ }\left( 1.4.2.1.5\right)%
\end{array}%
$

$\ \ \ \ \ \ $

\textbf{(c)} hypersequence $(x_{\mathbf{n}}\times y_{\mathbf{n}})_{\mathbf{n}%
\in \text{ }^{\ast }%
\mathbb{N}
}$ \ is Cauchy and

\bigskip

\bigskip $\ \ 
\begin{array}{cc}
\begin{array}{c}
\\ 
\ (x_{\mathbf{n}}\times y_{\mathbf{n}})_{\mathbf{n}\in \text{ }^{\ast }%
\mathbb{N}
}\approx _{\mathbf{c}}(x_{\mathbf{n}}^{\prime }\times y_{\mathbf{n}}^{\prime
})_{\mathbf{n}\in \text{ }^{\ast }%
\mathbb{N}
}, \\ 
\end{array}
& \text{ \ \ }\left( 1.4.2.1.6\right)%
\end{array}%
\ \ \ \ \ \ $

\bigskip

\textbf{(d)} hypersequence $\left( \dfrac{x_{\mathbf{n}}}{y_{\mathbf{n}}}%
\right) _{\mathbf{n}\in \text{ }^{\ast }%
\mathbb{N}
}$ is Cauchy and

\bigskip

$\ \ \ 
\begin{array}{cc}
\begin{array}{c}
\\ 
\ \left( \dfrac{x_{\mathbf{n}}}{y_{\mathbf{n}}}\right) _{\mathbf{n}\in \text{
}^{\ast }%
\mathbb{N}
}\approx _{\mathbf{c}}\left( \dfrac{x_{\mathbf{n}}^{\prime }}{y_{\mathbf{n}%
}^{\prime }}\right) _{\mathbf{n}\in \text{ }^{\ast }%
\mathbb{N}
} \\ 
\end{array}
& \text{\ }\left( 1.4.2.1.7\right) 
\end{array}%
$

\ \ \ \ iff $\forall \mathbf{n}_{(\mathbf{n\in }^{\ast }%
\mathbb{N}
)}\left[ (y_{\mathbf{n}}\not\equiv 0)\wedge \left( y_{\mathbf{n}}^{\prime
}\not\equiv 0\right) \wedge \left( y_{\mathbf{n}}\not\approx _{\mathbf{c}%
}0\right) \right] ,$

\textbf{(e) }hypersequence $(x_{\mathbf{n}}+0_{\mathbf{n}})_{\mathbf{n}\in 
\text{ }^{\ast }%
\mathbb{N}
}$ where $\forall \mathbf{n}_{(\mathbf{n\in }^{\ast }%
\mathbb{N}
)}\left[ 0_{\mathbf{n}}=0\right] $ \ \ \ \ \ \ \ \ \ \ \ \ \ \ \ \ \ \ \ \ \
\ \ \ \ \ \ \ \ \ \ \ \ \ \ \ \ \ \ \ \ \ \ \ \ \ \ \ 

\ \ \ \ is Cauchy and

$\ \ \ \ \ 
\begin{array}{cc}
\begin{array}{c}
\\ 
(x_{\mathbf{n}})_{\mathbf{n}\in \text{ }^{\ast }%
\mathbb{N}
}+(0_{\mathbf{n}})_{\mathbf{n}\in \text{ }^{\ast }%
\mathbb{N}
}\approx _{\mathbf{c}}(x_{\mathbf{n}})_{\mathbf{n\in }^{\ast }%
\mathbb{N}
}, \\ 
\end{array}
& \text{ \ \ }\left( 1.4.2.1.8\right)%
\end{array}%
$

\bigskip

\bigskip here $(0_{\mathbf{n}})_{\mathbf{n}\in \text{ }^{\ast }%
\mathbb{N}
}$ is a \textit{null hypersequence,}

\textbf{(f) }hypersequence $(x_{\mathbf{n}}\times 1_{\mathbf{n}})_{\mathbf{%
n\in }^{\ast }%
\mathbb{N}
}$ where $\forall \mathbf{n}_{(\mathbf{n\in }^{\ast }%
\mathbb{N}
)}\left[ 1_{\mathbf{n}}=1\right] $ is

Cauchy and

\bigskip

$\ \ 
\begin{array}{cc}
\begin{array}{c}
\\ 
\left( x_{\mathbf{n}}\right) _{\mathbf{n\in }^{\ast }%
\mathbb{N}
}\times (1_{\mathbf{n}})_{\mathbf{n\in }^{\ast }%
\mathbb{N}
}\approx _{\mathbf{c}}(x_{\mathbf{n}})_{\mathbf{n\in }^{\ast }%
\mathbb{N}
}, \\ 
\end{array}
& \text{ \ \ }\left( 1.4.2.1.9\right)%
\end{array}%
$

\ \ \ \ \ \ \ \ \ \ \ \ \ \ 

here $(1_{\mathbf{n}})_{\mathbf{n\in }^{\ast }%
\mathbb{N}
}$ is a unit hypersequence.

\textbf{(g) }hypersequence $(x_{\mathbf{n}})_{\mathbf{n\in }^{\ast }%
\mathbb{N}
}\times (x_{\mathbf{n}})_{\mathbf{n\in }^{\ast }%
\mathbb{N}
}^{-1}$ \ is Cauchy and

\bigskip

\bigskip $\ \ 
\begin{array}{cc}
\begin{array}{c}
\\ 
(x_{\mathbf{n}})_{\mathbf{n\in }^{\ast }%
\mathbb{N}
}\times (x_{\mathbf{n}})_{\mathbf{n\in }^{\ast }%
\mathbb{N}
}^{-1}\approx _{\mathbf{c}}(1_{\mathbf{n}})_{\mathbf{n\in }^{\ast }%
\mathbb{N}
} \\ 
\end{array}
& \text{ \ \ }\left( 1.4.2.1.10\right)%
\end{array}%
$

\ \ \ \ iff $\forall \mathbf{n}_{(\mathbf{n\in }^{\ast }%
\mathbb{N}
)}\left[ (x_{\mathbf{n}}\not\equiv 0)\wedge \left( x_{\mathbf{n}}\not\approx
_{\mathbf{c}}\left( 0_{\mathbf{n}}\right) _{\mathbf{n\in }^{\ast }%
\mathbb{N}
}\right) \right] .$

\bigskip \textbf{Proof. (a) }From definition of the Cauchy hypersequences
one obtain:

\bigskip

$%
\begin{array}{cc}
\begin{array}{c}
\\ 
\exists \varepsilon _{1}\exists \mathbf{m}_{\left( \mathbf{m\in }^{\ast }%
\mathbb{N}
_{\infty }\right) }\mathbf{\forall k}\left( \mathbf{k\geqslant m}\right) 
\mathbf{\forall l}\left( \mathbf{l\geqslant m}\right) \left[ \left(
\left\vert x_{\mathbf{k}}-x_{\mathbf{l}}\right\vert <\varepsilon _{1}\right)
\right. \\ 
\\ 
\left. \wedge \left( \left\vert y_{\mathbf{k}}-y_{\mathbf{l}}\right\vert
<\varepsilon _{1}\right) \right] . \\ 
\end{array}
& \text{ }\left( 1.4.2.1.11\right)%
\end{array}%
$

\bigskip

Suppose $\varepsilon _{1}=\varepsilon /2,$ then from formula above we can to
choose $\mathbf{m=m}\left( \varepsilon _{1}\right) $

such that for all $\mathbf{k\geqslant m,l\geqslant m}$ valid the next
inequalities:

\bigskip

$\ 
\begin{array}{cc}
\begin{array}{c}
\\ 
\left\vert (x_{\mathbf{k}}+y_{\mathbf{k}})-(x_{\mathbf{l}}+y_{\mathbf{l}%
})\right\vert =\left\vert (x_{\mathbf{k}}-x_{\mathbf{l}})+\left( y_{\mathbf{k%
}}-y_{\mathbf{l}}\right) \right\vert \leqslant \\ 
\\ 
\leqslant \left\vert (x_{\mathbf{k}}-x_{\mathbf{l}})\right\vert +\left\vert
\left( y_{\mathbf{k}}-y_{\mathbf{l}}\right) \right\vert <\varepsilon
/2+\varepsilon /2=\varepsilon , \\ 
\\ 
\left\vert (x_{\mathbf{k}}^{\prime }+y_{\mathbf{k}}^{\prime })-(x_{\mathbf{l}%
}^{\prime }+y_{\mathbf{l}}^{\prime })\right\vert =\left\vert (x_{\mathbf{k}%
}^{\prime }-x_{\mathbf{l}}^{\prime })+\left( y_{\mathbf{k}}^{\prime }-y_{%
\mathbf{l}}^{\prime }\right) \right\vert \leqslant \\ 
\\ 
\leqslant \left\vert (x_{\mathbf{k}}^{\prime }-x_{\mathbf{l}}^{\prime
})\right\vert +\left\vert \left( y_{\mathbf{k}}^{\prime }-y_{\mathbf{l}%
}^{\prime }\right) \right\vert <\varepsilon /2+\varepsilon /2=\varepsilon .
\\ 
\end{array}
& \text{ \ }\left( 1.4.2.1.12\right)%
\end{array}%
$

\bigskip

From \textbf{Definition 1.4.2.1.5.} and inequalities (1.4.12) we have \ \ \
\ \ \ \ \ \ \ \ \ \ \ \ \ \ \ \ \ \ \ \ \ \ \ \ \ \ \ \ \ \ \ \ \ \ \ \ 

the statement (\textbf{a}).

\bigskip \textbf{(b)} Similarly proof the statement (a) we have the next
inequalities:\ $\ \ \ \ $

$\ \ \ 
\begin{array}{cc}
\begin{array}{c}
\\ 
\left\vert (x_{\mathbf{k}}-y_{\mathbf{k}})-(x_{\mathbf{l}}-y_{\mathbf{l}%
})\right\vert =\left\vert (x_{\mathbf{k}}-x_{\mathbf{l}})+\left( y_{\mathbf{k%
}}-y_{\mathbf{l}}\right) \right\vert \leqslant \\ 
\\ 
\leqslant \left\vert (x_{\mathbf{k}}-x_{\mathbf{l}})\right\vert +\left\vert
\left( y_{\mathbf{k}}-y_{\mathbf{l}}\right) \right\vert <\varepsilon
/2+\varepsilon /2=\varepsilon , \\ 
\\ 
\left\vert (x_{\mathbf{k}}^{\prime }-y_{\mathbf{k}}^{\prime })-(x_{\mathbf{l}%
}^{\prime }-y_{\mathbf{l}}^{\prime })\right\vert =\left\vert (x_{\mathbf{k}%
}^{\prime }-x_{\mathbf{l}}^{\prime })-\left( y_{\mathbf{k}}^{\prime }-y_{%
\mathbf{l}}^{\prime }\right) \right\vert \leqslant \\ 
\\ 
\leqslant \left\vert (x_{\mathbf{k}}^{\prime }-x_{\mathbf{l}}^{\prime
})\right\vert +\left\vert \left( y_{\mathbf{k}}^{\prime }-y_{\mathbf{l}%
}^{\prime }\right) \right\vert <\varepsilon /2+\varepsilon /2=\varepsilon .\
\ \ \  \\ 
\end{array}
& \text{ \ }\left( 1.4.2.1.13\right)%
\end{array}%
$

From \textbf{Definition 1.4.2.1.5} and inequalities (1.4.13) we have the \ \
\ \ \ \ \ \ \ \ \ \ \ \ \ \ \ \ \ \ \ \ \ \ \ \ \ \ \ \ \ \ \ \ \ 

statement (\textbf{b}).

\bigskip \textbf{(c)} $\mathbf{\forall k}\left( \mathbf{k\geqslant m}\right) 
$ and $\mathbf{\forall l}\left( \mathbf{l\geqslant m}\right) $ we have the
next inequalities:

$\ \ \ \ \ \ \ \ \ \ \ \ \ \ \ 
\begin{array}{cc}
\begin{array}{c}
\\ 
\left\vert x_{\mathbf{k}}\cdot y_{\mathbf{k}}-x_{\mathbf{l}}\cdot y_{\mathbf{%
l}}\right\vert =\left\vert \left( x_{\mathbf{k}}\cdot y_{\mathbf{k}}-x_{%
\mathbf{l}}\cdot y_{\mathbf{k}}\right) +\left( x_{\mathbf{l}}\cdot y_{%
\mathbf{k}}-x_{\mathbf{l}}\cdot y_{\mathbf{l}}\right) \right\vert \leqslant
\\ 
\\ 
\leqslant \left\vert x_{\mathbf{k}}-x_{\mathbf{l}}\right\vert \cdot
\left\vert y_{\mathbf{l}}\right\vert +\left\vert y_{\mathbf{k}}-y_{\mathbf{l}%
}\right\vert \cdot \left\vert x_{\mathbf{l}}\right\vert , \\ 
\\ 
\left\vert x_{\mathbf{k}}^{\prime }\cdot y_{\mathbf{k}}^{\prime }-x_{\mathbf{%
l}}^{\prime }\cdot y_{\mathbf{l}}^{\prime }\right\vert =\left\vert \left( x_{%
\mathbf{k}}^{\prime }\cdot y_{\mathbf{k}}^{\prime }-x_{\mathbf{l}}^{\prime
}\cdot y_{\mathbf{k}}^{\prime }\right) +\left( x_{\mathbf{l}}^{\prime }\cdot
y_{\mathbf{k}}^{\prime }-x_{\mathbf{l}}^{\prime }\cdot y_{\mathbf{l}%
}^{\prime }\right) \right\vert \leqslant \\ 
\\ 
\leqslant \left\vert x_{\mathbf{k}}^{\prime }-x_{\mathbf{l}}^{\prime
}\right\vert \cdot \left\vert y_{\mathbf{l}}^{\prime }\right\vert
+\left\vert y_{\mathbf{k}}^{\prime }-y_{\mathbf{l}}^{\prime }\right\vert
\cdot \left\vert x_{\mathbf{l}}^{\prime }\right\vert , \\ 
\\ 
\left\vert x_{\mathbf{k}}\cdot y_{\mathbf{k}}-x_{\mathbf{k}}^{\prime }\cdot
y_{\mathbf{k}}^{\prime }\right\vert =\left\vert \left( x_{\mathbf{k}}\cdot
y_{\mathbf{k}}-x_{\mathbf{k}}\cdot y_{\mathbf{k}}^{\prime }\right) +\left(
x_{\mathbf{k}}\cdot y_{\mathbf{k}}^{\prime }-x_{\mathbf{k}}^{\prime }\cdot
y_{\mathbf{k}}^{\prime }\right) \right\vert \leqslant \\ 
\\ 
\leqslant \left\vert x_{\mathbf{k}}-x_{\mathbf{k}}^{\prime }\right\vert
\cdot \left\vert y_{\mathbf{k}}^{\prime }\right\vert +\left\vert y_{\mathbf{k%
}}-y_{\mathbf{k}}^{\prime }\right\vert \cdot \left\vert x_{\mathbf{k}%
}\right\vert . \\ 
\end{array}
& \text{ }\left( 1.4.2.1.14\right)%
\end{array}%
$

From definition Cauchy hypersequences one obtain $\exists c\forall \mathbf{k:%
}\left\vert x_{\mathbf{k}}\right\vert \leqslant c,\left\vert y_{\mathbf{k}%
}\right\vert \leqslant c,$

$\left\vert x_{\mathbf{k}}^{\prime }\right\vert \leqslant c,\left\vert y_{%
\mathbf{k}}^{\prime }\right\vert \leqslant c.$ From \textbf{Definition
1.4.2.1.5} and inequalities (1.4.14) we

have the statement (\textbf{c}).

Let $\Re _{\mathbf{c}}^{\ast }$ denote the set of the all equivalence
classes $\left\{ \left( x_{n}\right) _{\mathbf{n\in }^{\ast }%
\mathbb{N}
}\right\} \in \Re _{\mathbf{c}}^{\ast }$

Using Lemma 1.4.1. one can define an equivalence relation $\approx _{\mathbf{%
c}},$which is \ \ \ \ \ \ \ \ \ \ \ \ \ \ \ 

compatible with the operations defined above, and the set $^{\ast }%
\mathbb{R}
_{\mathbf{c}}=\Re _{\mathbf{c}}^{\ast }/\approx _{\mathbf{c}}$ \ \ \ \ \ \ \
\ \ \ \ \ \ \ \ \ \ \ 

is satisfy of the all usual field axioms of the hyperreal numbers.

\textbf{Lemma 1.4.2.1.5. }Suppose that

$\left\{ \left( x_{\mathbf{n}}\right) _{\mathbf{n\in }^{\ast }%
\mathbb{N}
}\right\} ,\left\{ \left( y_{\mathbf{n}}\right) _{\mathbf{n\in }^{\ast }%
\mathbb{N}
}\right\} ,\left\{ \left( z_{\mathbf{n}}\right) _{\mathbf{n\in }^{\ast }%
\mathbb{N}
}\right\} \in $ $^{\ast }%
\mathbb{R}
_{\mathbf{c}},$ then:$\ \ \ \ \ \ \ \ \ \ \ \ \ \ \ \ \ \ \ $

$\bigskip \ \ \ \ \ \ \ \ \ \ \ \ \ \ 
\begin{array}{cc}
\begin{array}{c}
\\ 
\text{(\textbf{a}) } \\ 
\\ 
\text{(\textbf{b})} \\ 
\\ 
\\ 
\\ 
\text{(\textbf{c}) } \\ 
\\ 
\\ 
\\ 
\text{(\textbf{d})} \\ 
\\ 
\\ 
\\ 
\text{(\textbf{e})} \\ 
\\ 
\\ 
\text{(\textbf{f})} \\ 
\\ 
\text{(\textbf{g})} \\ 
\\ 
\text{(\textbf{i})} \\ 
\\ 
\text{(\textbf{j})} \\ 
\\ 
\text{(\textbf{k})} \\ 
\\ 
\\ 
\end{array}%
\begin{array}{c}
\\ 
\left\{ \left( x_{\mathbf{n}}\right) _{\mathbf{n\in }^{\ast }%
\mathbb{N}
}\right\} +\left\{ \left( y_{\mathbf{n}}\right) _{\mathbf{n\in }^{\ast }%
\mathbb{N}
}\right\} =\left\{ \left( y_{\mathbf{n}}\right) _{\mathbf{n\in }^{\ast }%
\mathbb{N}
}\right\} +\left\{ \left( x_{\mathbf{n}}\right) _{\mathbf{n\in }^{\ast }%
\mathbb{N}
}\right\} , \\ 
\\ 
\left[ \left\{ \left( x_{\mathbf{n}}\right) _{\mathbf{n\in }^{\ast }%
\mathbb{N}
}\right\} +\left\{ \left( y_{\mathbf{n}}\right) _{\mathbf{n\in }^{\ast }%
\mathbb{N}
}\right\} \right] +\left\{ \left( z_{\mathbf{n}}\right) _{\mathbf{n\in }%
^{\ast }%
\mathbb{N}
}\right\} = \\ 
\\ 
=\left\{ \left( x_{\mathbf{n}}\right) _{\mathbf{n\in }^{\ast }%
\mathbb{N}
}\right\} +\left[ \left\{ \left( y_{\mathbf{n}}\right) _{\mathbf{n\in }%
^{\ast }%
\mathbb{N}
}\right\} +\left\{ \left( z_{\mathbf{n}}\right) _{\mathbf{n\in }^{\ast }%
\mathbb{N}
}\right\} \right] , \\ 
\\ 
\left\{ \left( z_{\mathbf{n}}\right) _{\mathbf{n\in }^{\ast }%
\mathbb{N}
}\right\} \times \left[ \left\{ \left( x_{\mathbf{n}}\right) _{\mathbf{n\in }%
^{\ast }%
\mathbb{N}
}\right\} +\left\{ \left( y_{\mathbf{n}}\right) _{\mathbf{n\in }^{\ast }%
\mathbb{N}
}\right\} \right] = \\ 
\\ 
=\left\{ \left( z_{\mathbf{n}}\right) _{\mathbf{n\in }^{\ast }%
\mathbb{N}
}\right\} \times \left\{ \left( x_{\mathbf{n}}\right) _{\mathbf{n\in }^{\ast
}%
\mathbb{N}
}\right\} +\left\{ \left( z_{\mathbf{n}}\right) _{\mathbf{n\in }^{\ast }%
\mathbb{N}
}\right\} \times \left\{ \left( y_{\mathbf{n}}\right) _{\mathbf{n\in }^{\ast
}%
\mathbb{N}
}\right\} , \\ 
\\ 
\left[ \left\{ \left( x_{\mathbf{n}}\right) _{\mathbf{n\in }^{\ast }%
\mathbb{N}
}\right\} \times \left\{ \left( y_{\mathbf{n}}\right) _{\mathbf{n\in }^{\ast
}%
\mathbb{N}
}\right\} \right] \times \left\{ \left( z_{\mathbf{n}}\right) _{\mathbf{n\in 
}^{\ast }%
\mathbb{N}
}\right\} = \\ 
\\ 
=\left\{ \left( x_{\mathbf{n}}\right) _{\mathbf{n\in }^{\ast }%
\mathbb{N}
}\right\} \times \left[ \left\{ \left( y_{\mathbf{n}}\right) _{\mathbf{n\in }%
^{\ast }%
\mathbb{N}
}\right\} \times \left\{ \left( z_{\mathbf{n}}\right) _{\mathbf{n\in }^{\ast
}%
\mathbb{N}
}\right\} \right] , \\ 
\\ 
\left[ \left\{ \left( x_{\mathbf{n}}\right) _{\mathbf{n\in }^{\ast }%
\mathbb{N}
}\right\} \times \left\{ \left( y_{\mathbf{n}}\right) _{\mathbf{n\in }^{\ast
}%
\mathbb{N}
}\right\} \right] \times \left\{ \left( z_{\mathbf{n}}\right) _{\mathbf{n\in 
}^{\ast }%
\mathbb{N}
}\right\} = \\ 
=\left\{ \left( x_{\mathbf{n}}\right) _{\mathbf{n\in }^{\ast }%
\mathbb{N}
}\right\} \times \left[ \left\{ \left( y_{\mathbf{n}}\right) _{\mathbf{n\in }%
^{\ast }%
\mathbb{N}
}\right\} \times \left\{ \left( z_{\mathbf{n}}\right) _{\mathbf{n\in }^{\ast
}%
\mathbb{N}
}\right\} \right] , \\ 
\\ 
\left\{ \left( x_{\mathbf{n}}\right) _{\mathbf{n\in }^{\ast }%
\mathbb{N}
}\right\} +\left\{ \left( 0_{\mathbf{n}}\right) _{\mathbf{n\in }^{\ast }%
\mathbb{N}
}\right\} =\left\{ \left( x_{\mathbf{n}}\right) _{\mathbf{n\in }^{\ast }%
\mathbb{N}
}\right\} , \\ 
\\ 
\left\{ \left( x_{\mathbf{n}}\right) _{\mathbf{n\in }^{\ast }%
\mathbb{N}
}\right\} \cdot \left\{ \left( x_{\mathbf{n}}\right) _{\mathbf{n\in }^{\ast }%
\mathbb{N}
}\right\} ^{-1}=\left\{ \left( 1_{\mathbf{n}}\right) _{\mathbf{n\in }^{\ast }%
\mathbb{N}
}\right\} , \\ 
\\ 
\left\{ \left( x_{\mathbf{n}}\right) _{\mathbf{n\in }^{\ast }%
\mathbb{N}
}\right\} \times \left\{ \left( 0_{\mathbf{n}}\right) _{\mathbf{n\in }^{\ast
}%
\mathbb{N}
}\right\} =\left\{ \left( 0_{\mathbf{n}}\right) _{\mathbf{n\in }^{\ast }%
\mathbb{N}
}\right\} , \\ 
\\ 
\left\{ \left( x_{\mathbf{n}}\right) _{\mathbf{n\in }^{\ast }%
\mathbb{N}
}\right\} \times \left\{ \left( 1_{\mathbf{n}}\right) _{{}}\right\} =\left\{
\left( x_{\mathbf{n}}\right) _{{}}\right\} , \\ 
\\ 
\left\{ \left( x_{\mathbf{n}}\right) _{\mathbf{n\in }^{\ast }%
\mathbb{N}
}\right\} <\left\{ \left( y_{\mathbf{n}}\right) _{\mathbf{n\in }^{\ast }%
\mathbb{N}
}\right\} \wedge \left\{ \left( 0_{\mathbf{n}}\right) _{\mathbf{n\in }^{\ast
}%
\mathbb{N}
}\right\} <\left\{ \left( z_{\mathbf{n}}\right) _{\mathbf{n\in }^{\ast }%
\mathbb{N}
}\right\} \implies \\ 
\\ 
\implies \left\{ \left( z_{\mathbf{n}}\right) _{\mathbf{n\in }^{\ast }%
\mathbb{N}
}\right\} \times \left\{ \left( x_{\mathbf{n}}\right) _{\mathbf{n\in }^{\ast
}%
\mathbb{N}
}\right\} <\left\{ \left( z_{\mathbf{n}}\right) _{\mathbf{n\in }^{\ast }%
\mathbb{N}
}\right\} \times \left\{ \left( y_{\mathbf{n}}\right) _{\mathbf{n\in }^{\ast
}%
\mathbb{N}
}\right\} . \\ 
\end{array}
& \text{ \ \ }\left( 1.4.2.1.15\right) \text{\ }%
\end{array}%
$

\textbf{Proof. }Statements\textbf{\ (a),(b),(c),(d),(e),(f),(g),(i) (j)} and 
\textbf{(k)} is\textbf{\ }evidently from \ \ \ \ \ \ \ \ \ \ \ \ \ \ \ 

Lemma.1.4.1\textbf{\ }and definition of the equivalence relation $\approx _{%
\mathbf{c}}.$

\bigskip

\bigskip

\section{I.4.2.2.The construction non-archimedean field $^{\ast }%
\mathbb{R}
_{\mathbf{c}}^{\protect\omega }$ as \ Cauchy completion of countable
non-standard models of field $%
\mathbb{Q}
.$\ }

\bigskip

Let $^{\ast }%
\mathbb{Q}
_{\omega }$ be a countable field which is elementary equivalent, but not \ \
\ \ \ \ \ \ \ \ \ \ \ \ \ \ \ \ \ \ \ 

isomorphic to $%
\mathbb{Q}
.$

\textbf{Remark.1.4.2.2.1. }The \textquotedblleft elementary
equivalence\textquotedblright\ means that an

(arithmetic) expression of first order is true in field $^{\ast }%
\mathbb{Q}
_{\omega }$ if and only if it

is true in field $%
\mathbb{Q}
.$

Note that any non-standard model of $%
\mathbb{Q}
$ contains an element $\mathbf{e\in }^{\ast }%
\mathbb{Q}
_{\omega }$ such \ \ \ \ \ \ \ \ \ \ \ \ \ \ \ \ \ 

that $\mathbf{e}>q$ for each $q\in 
\mathbb{Q}
.$

The canonical way to construct a model for $^{\ast }%
\mathbb{Q}
_{\omega }$ uses model theory [30],[31].

We simply take as axioms all axioms of $%
\mathbb{Q}
$ and additionally the following

countable number of axioms: the existence of an element $\mathbf{e}$ with

$\mathbf{e}>1,\mathbf{e}>2,...,\mathbf{e}>n,...\mathbf{.}$Each finite subset
of this axioms is satisfied by the

standard $%
\mathbb{Q}
$.By the compactness theorem in first order model theory, there

exists a model which also satisfies the given infinite set of axioms. By the

theorem of L\"{o}wenheim-Skolem, we can choose such models of countable

cardinality.\bigskip

Each non-standard model $^{\ast }%
\mathbb{Q}
$ contains the (externally defined) subset\ 

\bigskip

$\ \ \ 
\begin{array}{cc}
\begin{array}{c}
\\ 
^{\ast }%
\mathbb{Q}
_{\mathbf{fin}}\triangleq \left\{ x\in \text{ }^{\ast }%
\mathbb{Q}
|\exists n_{n\in 
\mathbb{Q}
}\left[ -n\leq x\leq n\right] \right\} . \\ 
\end{array}
& \text{ }\left( 1.4.2.2.1\right)%
\end{array}%
$

\bigskip

\bigskip Every element $x\in $ $^{\ast }%
\mathbb{Q}
_{\mathbf{fin}}$ defines a Dedekind cut:

\bigskip

\ $\ \ 
\begin{array}{cc}
\begin{array}{c}
\\ 
\mathbb{Q}
=\left\{ q\in 
\mathbb{Q}
|\text{ }q\leq x\right\} \cup \left\{ q\in 
\mathbb{Q}
|q>x\right\} . \\ 
\end{array}
& \text{\ }\left( 1.4.2.2.2\right)%
\end{array}%
$

\bigskip

We therefore get a order preserving map $\mathbf{j}_{op}\mathbf{:}^{\ast }%
\mathbb{Q}
_{\mathbf{fin}}\rightarrow 
\mathbb{R}
$ which restricts to

the standard inclusion of the standard rationals and which

respects addition and multiplication. An element of $^{\ast }%
\mathbb{Q}
_{\mathbf{fin}}$ is called

infinitesimal,if it is mapped to $0$ under the map $\mathbf{j}_{op}.$

\textbf{Proposition} [30].\textbf{1.4.2.2.1.}Choose an arbitrary subset $%
M\subset 
\mathbb{R}
.$Then

(\textbf{i}) \ there is a model $^{\ast }%
\mathbb{Q}
^{M}$ such that $\mathbf{j}_{op}\left( ^{\ast }%
\mathbb{Q}
_{\mathbf{fin}}^{M}\right) \supset M.$

(\textbf{ii}) the cardinality of $^{\ast }%
\mathbb{Q}
^{M}$ can be chosen to coincide with $\mathbf{card}\left( M\right) $,if $M$
is

\ \ \ \ infinite.

\textbf{Proof.} Choose $M\subset 
\mathbb{R}
$. For each $m\in M$ choose $q_{1}^{m}<q_{2}^{m}<...<...$

$<p_{2}^{m}<p_{1}^{m}$ with $\lim_{k\rightarrow \infty
}q_{k}^{m}=\lim_{k\rightarrow \infty }p_{k}^{m}=m.$

We add to the axioms of $%
\mathbb{Q}
$ the following axioms:$\forall m\in M$ $\exists e_{m}$ such that

$q_{k}^{m}<e_{m}<p_{k}^{m}$ for all $k\in 
\mathbb{N}
.$

Again, the standard $%
\mathbb{R}
$ is a model for each finite subset of these axioms,

so that the compactness theorem implies the existence of $^{\ast }%
\mathbb{Q}
^{M}$ as required,

where the cardinality of $^{\ast }%
\mathbb{Q}
^{M}$ can be chosen to be the cardinality of the set

of axioms, i.e. of $M,$ if $M$ is infinite. Note that by construction $%
\mathbf{j}_{op}\left( e_{m}\right) =e_{m}.$

\textbf{Remark.1.4.2.2.2. }It follows in particular that for each countable
subset

of $%
\mathbb{R}
$ we can find a countable model $^{\ast }%
\mathbb{Q}
_{\omega }$ of $^{\ast }%
\mathbb{Q}
$ such that the image of

$\mathbf{j}_{op}\left( \circ \right) $ contains this subset.Note, on the
other hand, that the image will

only be countable, so that the different models will have very different

ranges.

\textbf{Definition 1.4.2.2.1.}[30].\textbf{\ }A Cauchy sequence in $^{\ast }%
\mathbb{Q}
_{\omega }$ is a sequence $\left( a_{k}\right) _{k\in 
\mathbb{N}
}$

such that for every $\varepsilon \in $ $^{\ast }%
\mathbb{Q}
_{\omega },$ $\varepsilon >0$ there is an $n_{\varepsilon }\in 
\mathbb{N}
$ such that:

$\forall m_{m\text{ }>\text{ }n_{\varepsilon }}\forall n_{n\text{ }>\text{ }%
n_{\varepsilon }}\left[ \left\vert \text{ }a_{m}-a_{n}\right\vert
<\varepsilon \right] .$

\textbf{Definition 1.4.2.2.2. }We define Cauchy completion $^{\ast }%
\mathbb{R}
_{\mathbf{c}}^{\omega }\triangleq \left[ ^{\ast }%
\mathbb{Q}
_{\omega }\right] _{\mathbf{c}}$ in the

canonical way as equivalence classes of Cauchy sequences.

\ \ \ \ \ \ \ \ \ \ \ \ \ \ 

$\mathbf{Remark.}$\textbf{1.4.2.2.3. }This is a standard construction and
works for all

ordered fields.The result is again a field, extending the original field.

Note that, in our case,each point in $\left[ ^{\ast }%
\mathbb{Q}
_{\omega }\right] _{\mathbf{c}}$ is infinitesimally close to a

point in $^{\ast }%
\mathbb{Q}
$.

$\mathbf{Remark.}$\textbf{1.4.2.2.4. }In many non-standard models of $%
\mathbb{Q}
$, there are no

countable sequences $\left( a_{k}\right) _{k\in 
\mathbb{N}
}$ tending to zero which are not eventually

zero.

\bigskip

$\ \ \ \ \ \ \ \ \ \ \ \ \ \ \ \ \ \ \ \ \ \ \ \ \ \ \ \ \ \ \ \ $

\section{II.Euler's proofs by using non-archimedean analysis on the
pseudo-ring $^{\ast }%
\mathbb{R}
_{\mathbf{d}}$ revisited.\ \ \ \ \ \ \ \ \ \ \ \ \ \ \ \ \ \ \ II.1.Euler's
original proof of the Goldbach-Euler Theorem revisited.}

\ \ \ \ \ \ \ \ \ \ \ \ \ \ \ \ \ \ \ \ \ \ \ \ \ \ \ \ \ \ \ \ \ \ \ \ \ \
\ \ \ \ \ \ \ \ \ \ \ \ 

\bigskip\ \ \ \ \ \ \ \ \ \ \ \ \ \ \ \ \ \ \ \ \ \ \ \ \ \ \ \ \ \ \ \ \ \
\ \ \ \ \ \ \ \ \ \ \ \ 

That's what he's infected me with,he thought. His madness.That's

why I've come here.That's what I want here. A strange and very new

feeling overwhelmed him. He was aware that the feeling was really not

new at all, that it had been hidden in him for a long time, but that he was

acknowledging it only now, and everything was falling into place.

\bigskip\ \ \ \ \ \ \ \ \ \ \ \ \ \ \ \ \ \ \ \ \ \ \ \ \ \ \ \ \ \ \ \ \ \
\ \ \ \ \ \ \ \ \ \ \ \ \ \ \ \ \ \ \ \ \ \ \ \ \ \ \ \ \ \ \ Arkady and
Boris Strugatsky

\bigskip\ \ \ \ \ \ \ \ \ \ \ \ \ \ \ \ \ \ \ \ \ \ \ \ \ \ \ \ \ \ \ \ \ \
\ \ \ \ \ \ \ \ \ \ \ \ \ \ \ \ \ \ \ \ \ \ \ \ \ \ \ \ \ \ \ \ \ \ \ \ \ \
\ \ \ \ \ "Roadside Picnic"

Euler's paper of 1737 \textquotedblleft Variae Observationes Circa Series
Infinitas,\textquotedblright\ is Euler's first paper that closely follows
the modern Theorem-Proof format. There are no definitions in the paper, or
it would probably follow the Definition-Theorem-Proof format. After an
introductory paragraph in which Euler tells part of the story of the
problem, Euler gives us a theorem and a "proof". Euler's "proof" begins with
an 18-th century step that treats \textit{infinity as a number.} Such steps
became unpopular among rigorous mathematicians about a hundred years later.
He takes $x$ to be the "sum" of the harmonic series:

\bigskip

$\ \ 
\begin{array}{cc}
\begin{array}{c}
\\ 
x=1+\dfrac{1}{2}+\dfrac{1}{3}+\dfrac{1}{4}+\dfrac{1}{5}+\dfrac{1}{6}+ \\ 
\\ 
+\dfrac{1}{7}+\dfrac{1}{8}+\dfrac{1}{9}+...+\dfrac{1}{n}+... \\ 
\end{array}
& \text{ \ \ \ \ \ \ \ }\left( 2.4\right)%
\end{array}%
$

\bigskip

The Euler's original proofs is one of those examples of \textit{completely} 
\textit{misuse} \ \ \ \ \ \ \ \ \ \ \ \ \ \ \ \ \ \ \ \ \ \ \ \ \ of
divergent series to obtain \textit{completely correct results} so frequent
during the \ seventeenth and eighteenth centuries.The acceptance of Euler's
proofs seems to lie in the fact that,at the time,Euler (and most of his
contemporaries) actually manipulated a model of real numbers which included
infinitely large and infinitely small numbers.A model that much later
Bolzano would try to build on solid grounds and that today is called
\textquotedblleft nonstandard\textquotedblright\ after A.Robinson definitely
established it in the 1960's [1],[2],[3],[4],[5]. This last approach,
though, is completely in tune with Euler's proof [7] Nevertheless using
ideas borrowed from modern nonstandard analysis the same reconstruction
rigorous by modern Robinsonian standards \textit{is not found.} In
particular "nonstandard" proof proposed in paper [7] is not completely
nonstandard becourse authors use the solution Catalan's conjecture [9]

Unfortunately completely correct proofs of the Goldbach-Euler Theorem, was
presented many authors as rational reconstruction only in terms which could
be considered rigorous by modern Weierstrassian standards.

In this last section we show how, a few simple ideas from non-archimedean
analysis on the pseudoring $^{\ast }%
\mathbb{R}
_{\mathbf{d}},$ vindicate Euler's work.

\bigskip

\textbf{Theorem 2.1.1. }(Euler [6],[8]) Consider the following series,
infinitely \ \ \ \ \ \ \ \ \ \ \ \ \ \ \ \ \ \ \ \ \ \ \ \ \ \ \ \ \ 

continued:

\ \ $\ 
\begin{array}{cc}
\begin{array}{c}
\\ 
\dfrac{1}{3}+\dfrac{1}{7}+\dfrac{1}{8}+\dfrac{1}{15}+\dfrac{1}{24}+ \\ 
\\ 
+\dfrac{1}{26}+\dfrac{1}{31}+\dfrac{1}{35}+... \\ 
\end{array}
& \text{\ }\left( 2.4.1\right) 
\end{array}%
$

whose denominators, increased by one, are all the numbers which \ \ \ \ \ \
\ \ \ \ \ \ \ \ \ \ \ \ \ \ \ \ \ \ \ \ \ \ \ \ \ \ \ \ \ \ \ \ \ \ \ \ \ \
\ \ \ \ \ 

are powers of the integers, either squares or any other higher \ \ \ \ \ \ \
\ \ \ \ \ \ \ \ \ \ \ \ \ \ \ \ \ \ \ \ \ \ \ \ \ \ \ \ \ \ 

degree.Thus each term may be expressed by the formula

\bigskip

$\ \ \ 
\begin{array}{cc}
\begin{array}{c}
\\ 
\dfrac{1}{m^{n}-1} \\ 
\end{array}
& \text{ \ }\left( 2.4.2\right)%
\end{array}%
$

\bigskip

where $m$ and $n$ are integers greater than one. The sum of this series is $%
1.$

\textbf{Proof.} Let

$\ \ \ \ \ \ \ \ \ \ \ \ \ \ \ \ \ 
\begin{array}{cc}
\begin{array}{c}
\\ 
\ \mathbf{h}=\mathbf{cl}\left( 1+\dfrac{1}{2},1+\dfrac{1}{2}+\dfrac{1}{3},1+%
\dfrac{1}{2}+\dfrac{1}{3}+\right. \\ 
\\ 
\left. \dfrac{1}{4},1+\dfrac{1}{2}+\dfrac{1}{3}+\dfrac{1}{4}+\dfrac{1}{5}%
,...\right) \\ 
\end{array}
& \text{ \ \ }\left( 2.4.3\right)%
\end{array}%
$

\bigskip

from Eq.(2.4.3), we obtain

\ $\ \ \ \ \ \ \ 
\begin{array}{cc}
\begin{array}{c}
\\ 
1=\mathbf{cl}\left( \dfrac{1}{2},\dfrac{1}{2}+\dfrac{1}{4},\dfrac{1}{2}+%
\dfrac{1}{4}+\dfrac{1}{8},\dfrac{1}{2}+\dfrac{1}{4}+\dfrac{1}{8}+\dfrac{1}{16%
},\right. \\ 
\\ 
\left. \dfrac{1}{2}+\dfrac{1}{4}+\dfrac{1}{8}+\dfrac{1}{16}+\dfrac{1}{32}%
,...,\dfrac{1}{2}+\dfrac{1}{4}+\dfrac{1}{8}+\dfrac{1}{16}+\dfrac{1}{32}+...+%
\dfrac{1}{2^{i}},...\right) -\varepsilon _{1}, \\ 
\\ 
\varepsilon _{1}\approx 0, \\ 
\\ 
\varepsilon _{1}=\mathbf{cl}\left( \dfrac{1}{2^{M}},\dfrac{1}{2^{M+1}},...,%
\dfrac{1}{2^{M+i}},...\right) \\ 
\end{array}
& \left( 2.4.4\right)%
\end{array}%
$

\bigskip

Thus we obtain

$\ \ 
\begin{array}{cc}
\begin{array}{c}
\\ 
\ \mathbf{h}-1= \\ 
\\ 
\mathbf{cl}\left( 1,1+\dfrac{1}{3},1+\dfrac{1}{3}+\dfrac{1}{5},1+\dfrac{1}{3}%
+\dfrac{1}{5}+\dfrac{1}{6},1+\dfrac{1}{3}+\dfrac{1}{5}+\dfrac{1}{6}+\dfrac{1%
}{7},\right. \\ 
\\ 
\left. 1+\dfrac{1}{3}+\dfrac{1}{5}+\dfrac{1}{6}+\dfrac{1}{7}+\dfrac{1}{9}+%
\dfrac{1}{10},...\right) -\varepsilon _{1}. \\ 
\end{array}
& \text{ }\left( 1.4.5\right)%
\end{array}%
$

\bigskip

From Eq.(2.4.5), we obtain

\bigskip 

$\ \ \ \ \ \ \ \ \ \ \ \ \ 
\begin{array}{cc}
\begin{array}{c}
\\ 
\dfrac{1}{2}=\mathbf{cl}\left( \dfrac{1}{3},\dfrac{1}{3}+\dfrac{1}{9},\dfrac{%
1}{3}+\dfrac{1}{9}+\dfrac{1}{27},...,\right. \\ 
\\ 
\left. \dfrac{1}{3}+\dfrac{1}{9}+\dfrac{1}{27}+...+\dfrac{1}{3^{i}}%
,...\right) -\varepsilon _{2}, \\ 
\\ 
\varepsilon _{2}\approx 0, \\ 
\\ 
\varepsilon _{2}=\dfrac{1}{2}\mathbf{cl}\left( \dfrac{1}{3^{M}},\dfrac{1}{%
3^{M+1}},...,\dfrac{1}{3^{M+i}},...\right) \\ 
\end{array}
& \text{ \ \ \ \ }\left( 2.4.6\right)%
\end{array}%
$

\bigskip

we obtain

$\ \ \ 
\begin{array}{cc}
\begin{array}{c}
\\ 
\ \mathbf{h}-\left( 1+\dfrac{1}{2}\right) =\mathbf{cl}\left( 1,1+\dfrac{1}{5}%
,1+\dfrac{1}{5}+\dfrac{1}{6},1+\dfrac{1}{5}+\dfrac{1}{6}+\dfrac{1}{7},\right.
\\ 
\\ 
\left. 1+\dfrac{1}{5}+\dfrac{1}{6}+\dfrac{1}{7}+\dfrac{1}{10}+\dfrac{1}{11}%
,...\right) -\left( \varepsilon _{1}+\varepsilon _{2}\right) . \\ 
\end{array}
& \text{ \ \ }\left( 2.4.7\right)%
\end{array}%
$

\bigskip\ 

Eq.(2.4.7),gives

$\ \ \ $

$\ \ \ \ 
\begin{array}{cc}
\begin{array}{c}
\\ 
\dfrac{1}{4}=\mathbf{cl}\left( \dfrac{1}{5},\dfrac{1}{5}+\dfrac{1}{25},%
\dfrac{1}{5}+\dfrac{1}{25}+\dfrac{1}{125},...,\right.  \\ 
\\ 
\left. \dfrac{1}{5}+\dfrac{1}{25}+\dfrac{1}{125}+...+\dfrac{1}{5^{i}}%
,...\right) -\varepsilon _{3}, \\ 
\\ 
\varepsilon _{3}\approx 0, \\ 
\\ 
\varepsilon _{3}=\dfrac{1}{4}\mathbf{cl}\left( \dfrac{1}{5^{M}},\dfrac{1}{%
5^{M+1}},...,\dfrac{1}{5^{M+i}},...\right)  \\ 
\end{array}
& \text{ }\left( 2.4.8\right) 
\end{array}%
$

\bigskip 

Finally we obtain

\bigskip

$\ \ 
\begin{array}{cc}
\begin{array}{c}
\\ 
\ \mathbf{h}-\left( 1+\dfrac{1}{2}+\dfrac{1}{4}\right) = \\ 
\mathbf{cl}\left( 1+\dfrac{1}{6},1+\dfrac{1}{6}+\dfrac{1}{7},1+\dfrac{1}{6}+%
\dfrac{1}{7}+\dfrac{1}{10},...\right) - \\ 
\\ 
-\left( \varepsilon _{1}+\varepsilon _{2}+\varepsilon _{3}\right) . \\ 
\end{array}
& \text{ \ \ }\left( 2.4.9\right)%
\end{array}%
$

\bigskip

Proceeding similarly, i.e. deleting all the all terms that remain,we get

\bigskip

$\ \ \ 
\begin{array}{cc}
\begin{array}{c}
\\ 
\ \mathbf{h}-\left[ \mathbf{\Im }_{n}\right] = \\ 
\\ 
=\mathbf{cl}\left( 1+\dfrac{1}{5},...,1+\dfrac{1}{m\left( n^{\prime }\right) 
},1+\dfrac{1}{m\left( n^{\prime }\right) }+...,...\right) - \\ 
\\ 
-\left( \#Ext-\sum_{n\in 
\mathbb{N}
}\varepsilon _{n}\right) , \\ 
\\ 
m>n^{\prime }\left( n\right) 
\end{array}
& \text{ \ }\left( 2.4.10\right) 
\end{array}%
$

\bigskip 

where $\ \ \ \ \ \ \ \ \ \ \ \ \ \ \ \ \ \ \ \ \ \ \ \ \ \ \ \ \ \ \ \ \ \ \
\ \ \ \ \ $

$\bigskip $

$\ \ \ 
\begin{array}{cc}
\begin{array}{c}
\\ 
\mathbf{\Im }_{n}=\left( 1,\dfrac{1}{2},\dfrac{1}{4},\dfrac{1}{5},...,\dfrac{%
1}{n^{\prime }\left( n\right) },...\right) \\ 
\end{array}
& \text{ \ \ \ }\left( 2.4.11\right)%
\end{array}%
$

\bigskip

whose denominators, increased by one, are all the numbers \ \ \ \ \ \ \ \ \
\ \ \ \ \ \ \ \ \ \ \ \ \ \ \ \ \ \ \ \ \ \ \ \ \ \ \ \ \ \ 

which are not powers. From Eq.(2.4.10) we obtain

\bigskip

$\ 
\begin{array}{cc}
\begin{array}{c}
\\ 
\mathbf{h}-\left[ \mathbf{\Im }_{n}\right] = \\ 
\\ 
1+\left( \#Ext-\sum_{n\in 
\mathbb{N}
}\varepsilon _{n}\right) = \\ 
=1+\epsilon \\ 
\\ 
\epsilon =\text{ }\#Ext-\sum_{n\in 
\mathbb{N}
}\varepsilon _{n}\approx 0.%
\end{array}
& \text{ \ \ }\left( 2.4.12\right)%
\end{array}%
$

\bigskip

Thus we obtain

$\ \ \ \ \ \ \ \ \ \ \ \ \ \ \ \ \ \ \ \ \ \ \ \ \ \ \ \ \ \ \ \ \ \ \ \ \ \
\ \ \ \ \ \ \ \ \ \ \ \ \ \ \ \ \ \ $

$\ \ 
\begin{array}{cc}
\begin{array}{c}
\\ 
\mathbf{h}-\left[ \mathbf{\Im }_{n}\right] =1+\epsilon , \\ 
\end{array}
& \text{ }\left( 2.4.13\right)%
\end{array}%
$

\bigskip

Substitution Eq.(2.4.3) into Eq.(2.4.13) gives

\bigskip

$\ 
\begin{array}{cc}
\begin{array}{c}
\\ 
1+\epsilon =\dfrac{1}{3}+\dfrac{1}{7}+\dfrac{1}{8}+\dfrac{1}{15}+\dfrac{1}{24%
}+\dfrac{1}{26}+... \\ 
\\ 
\epsilon \approx 0 \\ 
\end{array}
& \text{ \ }\left( 2.4.14\right)%
\end{array}%
$

\bigskip

series whose denominators, increased by one, are all the powers \ \ \ \ \ \
\ \ \ \ \ \ \ \ \ \ \ \ \ \ \ \ \ \ \ \ \ \ \ \ \ \ \ \ 

of the integers and whose sum is one.

\bigskip

\bigskip

\bigskip

\bigskip

Time passed, and more or less coherent thoughts came to him.

Well,that's it,he thought unwillingly. The road is open. He could

go down right now, but it was better, of course, to wait a while.

The meatgrinders can be tricky.\bigskip

He got up,automatically brushed off his pants, and started down into

the quarry.The sun was broiling hot, red spots floated before his eyes,

the air was quivering on the floor of the quarry, and in the shimmer it

seemed that the ball was dancing in place like a buoy on the waves.

He went past the bucket, superstitiously picking up his feet higher and

making sure not to step on the splotches. And then, sinking into the rubble,

he dragged himself across the quarry to the dancing, winking ball.

\ \ \ \ \ \ \ \ \ \ \ \ \ \ \ \ \ \ \ \ \ \ \ \ \ \ \ \ \ \ \ \ \ \ \ \ \ \
\ \ \ \ \ \ \ \ \ \ \ \ \ \ \ \ \ \ \ Arkady and Boris Strugatsky \bigskip\
\ \ \ \ \ \ \ \ \ \ \ \ \ \ \ \ \ \ \ \ \ \ \ \ \ \ \ \ \ \ \ \ \ \ \ \ \ \
\ \ \ \ \ \ \ \ \ \ \ \ \ \ \ \ \ \ \ \ \ \ \ \ \ \ \ \ \ \ \ \ \ \ \ \ \ \
\ \ \ \ 

\ \ \ \ \ \ \ \ \ \ \ \ \ \ \ \ \ \ \ \ \ \ \ \ \ \ \ \ \ \ \ \ \ \ \ \ \ \
\ \ \ \ \ \ \ \ \ \ \ \ \ \ \ \ \ \ \ \ \ \ \ \ \ \ "Roadside Picnic"\ \ \ \
\ \ \ \ \ \ \ \ \ \ \ \ \ \ \ \ \ \ \ \ \ \ \ \ \ \ \ \ \ \ \ \ \ \ \ \ \ \
\ \ \ \ \ \ \ \ \ \ \ \ \ \ \ \ \ \ \ \ \ \ \ \ \ \ \ \ \ \ \ \ \ \ \ \ \ \
\ \ 

\ \ \ \ \ \ \ \ \ \ \ 

\bigskip $\ \ \ \ \ \ \ \ \ \ \ \ \ \ \ \ \ \ \ \ \ \ \ \ $

\bigskip

\section{III.Non-archimedean analysis on the extended hyperreal line $^{\ast
}%
\mathbb{R}
_{\mathbf{d}}$ and transcendence conjectures over field $%
\mathbb{Q}
.$Proof that $\ e+\protect\pi $ and $e\cdot \protect\pi $ is irrational.\ \
\ \ \ \ \ \ \ \ }

\bigskip $\ \ \ \ \ \ \ \ $

\bigskip $\ \ \ \ \ \ \ \ \ \ \ $

\bigskip

\section{III.1.Hyperrational approximation of the irrational numbers.}

The next simple result shows that in a way the hyperrationals already

"incorporate" the real numbers (see e.g. [25] Ch.II Thm. 2).

\textbf{Theorem 3.1.1.}Let $^{\ast }%
\mathbb{Q}
_{\mathbf{fin}}$ be the ring of finite hyperrationals, and let $\Im $ be

the maximal ideal of its infinitesimals. Then $%
\mathbb{R}
$ and $^{\ast }%
\mathbb{Q}
_{\mathbf{fin}}/\Im $ are isomorphic

as ordered fields.

\textbf{Theorem 3.1.2.(Standard form of Dirichlet's Approximation }

\textbf{Theorem).}Let be $\alpha \in 
\mathbb{R}
$ positive real number and $n\in 
\mathbb{N}
$ a positive

integer. Then there is an integer $k\in 
\mathbb{N}
$ and an integer $b\in 
\mathbb{N}
$ with

$0<k<n,$ for which

\ 

\ \ \ \ \ \ \ $\ 
\begin{array}{cc}
\begin{array}{c}
\\ 
-\dfrac{1}{n}<k\cdot \alpha -b<\dfrac{1}{n}. \\ 
\end{array}
& (\mathbf{DAP})\text{\ \ }\left( 3.1.1\right) \text{\ }%
\end{array}%
$

\bigskip

\textbf{Definition 3.1.1.} A \textquotedblleft $\mathbf{D}$%
-approximation\textquotedblright\ to $\alpha $ is a rational

number $\dfrac{p}{q}\in 
\mathbb{Q}
$, whose denominator is a positive integer

$q\in 
\mathbb{N}
,$ with

\bigskip

\ \ \ \ \ \ \ $\ \ 
\begin{array}{cc}
\begin{array}{c}
\\ 
\left\vert \alpha -\dfrac{p}{q}\right\vert <\dfrac{1}{q^{2}}. \\ 
\end{array}
& \text{ }\ (\mathbf{DAP1})\text{\ \ }\left( 3.1.2\right) \text{\ }%
\end{array}%
$

\bigskip

\textbf{Theorem 3.1.3.} If $\alpha \in 
\mathbb{R}
$ is irrational it has infinitely many

$\mathbf{D}$-approximations.

\textbf{Remark} \textbf{3.1.1.}[37]. Let sequence $\dfrac{p_{n}}{q_{n}}$ be
a convergent to $\alpha $ in the

sense such that:\ 

\ \ \ \ \ \ \ \ \ \ \ \ \ \ \ \ \ \ \ \ \ \ \ \ \ \ \ \ \ \ \ \ \ \ \ $\ \ \
\ \ \ \ \ \ \ \ \ \ \ 
\begin{array}{cc}
\begin{array}{c}
\\ 
\alpha =\dfrac{p_{n}}{q_{n}}+\dfrac{\theta _{q_{n}}}{q_{n}^{2}}, \\ 
\\ 
\left( p_{n},q_{n}\right) =1,\left\vert \theta _{q_{n}}\right\vert
<1,n=0,1,2,... \\ 
\end{array}
& \text{ }(\mathbf{DAP3})\text{\ \ }\left( 3.1.3\right) \text{\ }%
\end{array}%
$

\bigskip

i.e. there is exist infinite sequence $\left( p_{n},q_{n}\right) \in 
\mathbb{Z}
\times 
\mathbb{N}
,n=0,1,2,...$

such that $q_{n+1}>q_{n}$ and

$\ \ \ \ \ \ \ \ \ \ 
\begin{array}{cc}
\begin{array}{c}
\\ 
\alpha =\dfrac{p_{n}}{q_{n}}+\dfrac{\theta _{q_{n}}}{q_{n}^{2}}, \\ 
\\ 
\left( p_{n},q_{n}\right) =1,\left\vert \theta _{q_{n}}\right\vert <1. \\ 
\end{array}
& \text{ \ \ \ \ }(\mathbf{DAP4})\text{\ \ }\left( 3.1.4\right) \text{\ }%
\end{array}%
$

Theorem 3.1.3 shows that each irrational number $\alpha $ has infinitely many

convergents of the form $\mathbf{DAP4.}$

\bigskip

\textbf{Definition 3.1.2. }(\textbf{i})\textbf{\ }Let $\alpha \in 
\mathbb{R}
$ is irrational number. A $\ast $-$\mathbf{D}$-approximation\ 

to $^{\ast }\alpha \in $ $^{\ast }%
\mathbb{R}
$ is a number $\dfrac{P}{Q}\in $ $^{\ast }%
\mathbb{Q}
,$ $P\in $ $^{\ast }%
\mathbb{Z}
_{\infty }$ whose denominator is

a positive hyperinteger $Q\in $ $^{\ast }%
\mathbb{N}
_{\infty },$with\ \ \ 

\bigskip

\ \ $\ \ \ 
\begin{array}{cc}
\begin{array}{c}
\\ 
\left\vert ^{\ast }\alpha -\dfrac{P}{Q}\right\vert <\dfrac{1}{Q^{2}}, \\ 
\\ 
\left( P,Q\right) =1. \\ 
\end{array}
& \text{ \ }(\ast \text{-}\mathbf{DAP})\text{\ \ }\left( 3.1.5\right) \text{%
\ }%
\end{array}%
$

\bigskip

\bigskip (\textbf{ii}) Let $\alpha \in 
\mathbb{R}
$ is irrational number. A \textquotedblleft $\#$-$\mathbf{D}$%
-approximation\textquotedblright\ to $\left( ^{\ast }\alpha \right) ^{\#}\in 
$ $^{\ast }%
\mathbb{R}
_{\mathbf{d}}$

is a Wattenberg hyperrational number $\dfrac{P^{\#}}{Q^{\#}}\in $ $^{\ast }%
\mathbb{Q}
_{\mathbf{d}},$ $P^{\#}\in $ $^{\ast }%
\mathbb{Z}
_{\infty ,\mathbf{d}}$ whose

denominator is a positive hyperinteger $Q^{\#}\in $ $^{\ast }%
\mathbb{N}
_{\infty ,\mathbf{d}},$with

\bigskip

\ $\ \ \ \ \ \ 
\begin{array}{cc}
\begin{array}{c}
\\ 
\left\vert \left( ^{\ast }\alpha \right) ^{\#}-\dfrac{P^{\#}}{Q^{\#}}%
\right\vert <\dfrac{1^{\#}}{Q^{\#2}}, \\ 
\\ 
\left( P^{\#},Q^{\#}\right) =1^{\#}. \\ 
\end{array}
& \text{ \ }(\#\text{-}\mathbf{DAP})\text{\ \ }\left( 3.1.6\right) \text{\ }%
\end{array}%
$

\ \ \ \ \ \ 

\textbf{Definition 3.1.3. }(\textbf{i}) Let $\alpha \in 
\mathbb{R}
$ is irrational number. A hyperrational

approximation to $\alpha $ is a $\ast $-$\mathbf{D}$-approximation\ to $%
^{\ast }\alpha \in $ $^{\ast }%
\mathbb{R}
.$

(\textbf{ii}) Let $\alpha \in 
\mathbb{R}
$ is irrational number. A Wattenberg hyperrational

approximation to $\alpha $ is a $\#$-$\mathbf{D}$-approximation%
\textquotedblright\ to $\left( ^{\ast }\alpha \right) ^{\#}\in $ $^{\ast }%
\mathbb{R}
_{\mathbf{d}}.$

\bigskip

\textbf{Theorem 3.1.3.(Nonstandard form of Dirichlet's Approximation }

\textbf{Theorem).}

(\textbf{1}) If $\alpha $ is irrational it has infinitely many $\ast $-$%
\mathbf{D}$-approximations such that

for any two $\ast $-$\mathbf{D}$-approximations $P_{1}/Q_{1}$ and $%
P_{2}/Q_{2}$ the next equality is

satisfied

\bigskip

\bigskip\ $\ \ \ \ \ \ \ \ \ 
\begin{array}{cc}
\begin{array}{c}
\\ 
\dfrac{P_{1}}{Q_{1}}\approx \dfrac{P_{2}}{Q_{2}}, \\ 
\end{array}
& \text{ \ }\left( 3.1.7\right) \text{\ \ }%
\end{array}%
$

i.e.

$\ \ \ \ \ \ \ \ \ 
\begin{array}{cc}
\begin{array}{c}
\\ 
\left( ^{\ast }\alpha \right) ^{\#}=\left( \dfrac{P_{1}}{Q_{1}}\right) ^{\#}%
\text{ }\left( \func{mod}\varepsilon _{\mathbf{d}}\right) . \\ 
\end{array}
& \text{ \ \ }\left( 3.1.8\right) \text{\ \ \ }%
\end{array}%
$

\bigskip

(\textbf{2}) If $\alpha \in 
\mathbb{R}
$ is irrational then $^{\ast }\alpha \in $ $^{\ast }%
\mathbb{R}
$ has representation

\bigskip

$\ \ \ \ \ \ 
\begin{array}{cc}
\begin{array}{c}
\\ 
^{\ast }\alpha =\dfrac{P}{Q}+\dfrac{\theta _{Q}}{Q^{2}}, \\ 
\\ 
P\in \text{ }^{\ast }%
\mathbb{Z}
_{\infty },Q\in \text{ }^{\ast }%
\mathbb{N}
_{\infty },\left( P,Q\right) =1, \\ 
\\ 
\left( P,Q\right) =1,\left\vert \theta _{Q}\right\vert <1. \\ 
\end{array}
& \text{ \ }(\ast \text{-}\mathbf{DAP1})\text{\ \ }\left( 3.1.9\right) \text{%
\ }%
\end{array}%
$

\bigskip

\bigskip (\textbf{3}) If $\alpha \in 
\mathbb{R}
$ is irrational then $\left( ^{\ast }\alpha \right) ^{\#}\in $ $^{\ast }%
\mathbb{R}
_{\mathbf{d}}$ has representation

\bigskip

$\ \ \ \ \ \ \ 
\begin{array}{cc}
\begin{array}{c}
\\ 
\left( ^{\ast }\alpha \right) ^{\#}=\dfrac{P^{\#}}{Q^{\#}}+\dfrac{\left(
\theta _{Q}\right) ^{\#}}{Q^{\#2}}, \\ 
\\ 
P^{\#}\in \text{ }^{\ast }%
\mathbb{Z}
_{\infty ,\mathbf{d}},Q^{\#}\in \text{ }^{\ast }%
\mathbb{N}
_{\infty ,\mathbf{d}},\left( P^{\#},Q^{\#}\right) =1, \\ 
\\ 
\left( P^{\#},Q^{\#}\right) =1,\left\vert \left( \theta _{Q}\right)
^{\#}\right\vert <1^{\#}. \\ 
\end{array}
& \text{ \ \ }(\#\text{-}\mathbf{DAP1})\text{\ \ }\left( 3.1.10\right) \text{%
\ }%
\end{array}%
$

\bigskip

\bigskip

\textbf{Definition} \textbf{3.1.4.} A real number $\alpha \in 
\mathbb{R}
$ is a Liouville number if for

every positive integer $m\in 
\mathbb{N}
$, there is exist infinite sequence

$\left( p_{n},q_{n}\right) \in 
\mathbb{Z}
\times 
\mathbb{N}
,n=0,1,2,...$ such that

\bigskip

\ $\ \ \ \ \ \ 
\begin{array}{cc}
\begin{array}{c}
\\ 
0<\left\vert \alpha -\dfrac{p_{n}}{q_{n}}\right\vert <\dfrac{1}{q^{m}}. \\ 
\end{array}
& \text{ \ \ }\left( 3.1.11\right) \text{\ \ }%
\end{array}%
$

\bigskip

\textbf{Remark} \textbf{3.1.2.}This is well known that all Liouville numbers
are

transcendental. From the inequality (3.1.11) one obtain directly:

\textbf{Theorem 3.1.4. }(\textbf{i}) Any Liouville number $\alpha _{l}\in 
\mathbb{R}
$ for every positive

hyperinteger $N\in $ $^{\ast }%
\mathbb{N}
_{\infty }$ has a hyperrational approximation such that\bigskip

$\bigskip $

$\ \ \ \ \ \ \ 
\begin{array}{cc}
\begin{array}{c}
\\ 
0<\left\vert ^{\ast }\alpha _{l}-\dfrac{P}{Q}\right\vert <\dfrac{1}{Q^{N}},
\\ 
\\ 
P\in \text{ }^{\ast }%
\mathbb{Z}
_{\infty },N,Q\in \text{ }^{\ast }%
\mathbb{N}
_{\infty },\left( P,Q\right) =1. \\ 
\end{array}
& \text{ }\left( 3.1.12\right) \text{\ \ }%
\end{array}%
$

\bigskip

(\textbf{ii}) Any Liouville number $\alpha _{l}\in 
\mathbb{R}
$\ for every positive hyperinteger

$N\in $ $^{\ast }%
\mathbb{N}
_{\infty }$ has a Wattenberg hyperrational approximation such that

\bigskip

\bigskip\ \ \ \ \ \ \ \ $\ \ \ 
\begin{array}{cc}
\begin{array}{c}
\\ 
0<\left\vert \left( ^{\ast }\alpha _{l}\right) ^{\#}-\dfrac{P^{\#}}{Q^{\#}}%
\right\vert <\dfrac{1^{\#}}{\left( Q^{N}\right) ^{\#}}, \\ 
\\ 
P^{\#}\in \text{ }^{\ast }%
\mathbb{Z}
_{\infty ,\mathbf{d}},Q^{\#}\in \text{ }^{\ast }%
\mathbb{N}
_{\infty ,\mathbf{d}}, \\ 
\\ 
N\in \text{ }^{\ast }%
\mathbb{N}
_{\infty },\left( P^{\#},Q^{\#}\right) =1^{\#}. \\ 
\end{array}
& \text{ }\left( 3.1.13\right) \text{\ \ }%
\end{array}%
$

\bigskip

\textbf{Theorem 3.1.5. }Every Liouville number $\alpha _{l}\in 
\mathbb{R}
$ are $\#$-transcendental

over field $%
\mathbb{Q}
,$ i.e., there is no real $%
\mathbb{Q}
$-analytic function

$g_{%
\mathbb{Q}
}\left( x\right) =\dsum\limits_{n=0}^{\infty }a_{n}x^{n}<\infty ,0\leq
\left\vert x\right\vert \leq r\leq e$ with rational coefficients

$a_{0},a_{1},...,a_{n},...\in 
\mathbb{Q}
$ such that $g_{%
\mathbb{Q}
}\left( \alpha _{l}\right) .$

\bigskip

\ \ \ \ \ 

\bigskip $\ \ \ \ \ \ \ \ \ \ \ \ $

\bigskip

\section{III.2.Proof that $e$ is \#-transcendental.}

\bigskip

\textbf{Definition 3.2.1. }Let $g\left( x\right) :%
\mathbb{R}
\rightarrow 
\mathbb{R}
$\textbf{\ }be any real analytic function

\bigskip

$\ \ \ 
\begin{array}{cc}
\begin{array}{c}
\\ 
g_{%
\mathbb{Q}
}\left( x\right) =\dsum\limits_{n=0}^{\infty }a_{n}x^{n},\left\vert
x\right\vert <r, \\ 
\\ 
\forall n\left[ a_{n}\in 
\mathbb{Q}
\right] \\ 
\end{array}
& \text{\ }\left( 3.2.1\right)%
\end{array}%
$

\bigskip

defined on an open \bigskip interval $I\subset 
\mathbb{R}
$ such that $0\in I.$ \ \ \ \ \ \ \ \ \ \ \ \ \ \ \ \ \ \ \ \ \ \ \ \ \ \ \
\ \ \ \ \ \ \ \ \ \ \ \ \ \ \ \ \ \ \ \ \ \ \ \ \ \ \ \ \ \ \ \ \ \ \ \ \ \
\ \ \ \ \ \ \ \ \ \ \ \ \ \ \ \bigskip

We call this function given by Eq.(3.2.1) $%
\mathbb{Q}
$-analytic function and denote $g_{%
\mathbb{Q}
}\left( x\right) .$

\bigskip

\textbf{Definition 3.2.2. }Arbitrary transcendental number $z\in 
\mathbb{R}
$ is called $\ \ \ \ \ $\ \ 

$\#$-transcendental number over field $%
\mathbb{Q}
$, if no exist $%
\mathbb{Q}
$-analytic function $\ \ \ \ \ \ \ \ \ \ \ \ \ \ \ \ \ $

$g_{%
\mathbb{Q}
}\left( x\right) $ such that $g_{%
\mathbb{Q}
}\left( z\right) =0,$i.e. for every $%
\mathbb{Q}
$-analytic function $g_{%
\mathbb{Q}
}\left( x\right) $ the

inequality $g_{%
\mathbb{Q}
}\left( z\right) \neq 0$ is satisfies.

\bigskip

\textbf{Definition 3.2.3.}Arbitrary transcendental number $z$ called $w$%
-transcendental \ \ \ \ \ \ \ \ \ \ \ \ \ 

number over field $%
\mathbb{Q}
$,if $z$ is not $\#$-transcendental number over field $%
\mathbb{Q}
$,i.e.

exist $%
\mathbb{Q}
$-analytic function $g_{%
\mathbb{Q}
}\left( x\right) $ such that $g_{%
\mathbb{Q}
}\left( z\right) =0.$

\textbf{Example 3.2.1. }Number $\pi $ is transcendental but number $\pi $ is
not $\ \ \ \ \ \ \ $

$\ \#$-transcendental number over field $%
\mathbb{Q}
$ as

(\textbf{1}) function $\sin x$ is a $%
\mathbb{Q}
$-analytic and \ \ \ \ \ \ \ \ \ \ \ \ \ \ \ \ \ \ \ 

(\textbf{2})$\ \sin \left( \dfrac{\pi }{2}\right) =1,$i.e.$\ \ $

$\ \ \ \ \ \ \ \ \ \ \ \ \ \ \ \ \ \ \ \ \ \ \ \ \ \ \ \ \ \ \ \ \ \ \ \ \ \
\ \ \ \ \ \ \ \ \ \ $

$\ \ \ \ \ \ \ \ \ \ \ \ \ \ \ \ \ 
\begin{array}{cc}
\begin{array}{c}
\\ 
-1+\dfrac{\pi }{2}-\dfrac{\pi ^{3}}{2^{3}3!}+\dfrac{\pi ^{5}}{2^{5}5!}- \\ 
\\ 
-\dfrac{\pi ^{7}}{2^{7}7!}+...+\dfrac{\left( -1\right) ^{2n+1}\pi ^{2n+1}}{%
2^{2n+1}\left( 2n+1\right) !}+...=0. \\ 
\end{array}
& \text{ }\left( 3.2.2\right)%
\end{array}%
$

\bigskip

\textbf{Theorem 3.2.1.}Number $e$ is\textbf{\ }$\#$-transcendental over
field $%
\mathbb{Q}
.$

\textbf{Proof I.}To prove $e$ is $\#$-transcendental number we must show
that $e$ it is

not $w$-transcendental, i.e., there is no exist real $%
\mathbb{Q}
$-analytic function

$g_{%
\mathbb{Q}
}\left( x\right) =\dsum\limits_{n=0}^{\infty }a_{n}x^{n},0\leq \left\vert
x\right\vert \leq r\leq e$ with rational coefficients

$a_{0},a_{1},...,a_{n},...\in 
\mathbb{Q}
$ such that\ $\ $

$\ \ \ \ \ \ \ \ \ \ \ \ \ \ \ \ \ \ \ \ \ \ \ \ \ \ \ \ \ \ \ \ \ \ \ \ \ \
\ \ \ \ \ \ \ \ \ \ \ \ \ \ $

$\ \ \ \ 
\begin{array}{cc}
\begin{array}{c}
\\ 
\dsum\limits_{n=0}^{\infty }a_{n}e^{n}=0. \\ 
\end{array}
& \text{ \ \ \ \ \ \ \ \ \ }\left( 3.2.3\right)%
\end{array}%
$\textbf{\ }

\bigskip

Suppose that $e$ is $w$-transcendental, i.e., there is exist an $%
\mathbb{Q}
$-analytic \ \ \ \ \ \ \ \ \ \ \ \ \ \ \ \ \ \ \ \ \ \ \ \ \ \ \ \ \ \ \ \ \ 

function $\breve{g}_{%
\mathbb{Q}
}\left( x\right) =\dsum\limits_{n=0}^{\infty }\breve{a}_{n}x^{n},$with
rational coefficients:

\bigskip

$\ \ \ 
\begin{array}{cc}
\begin{array}{c}
\\ 
\ \ \breve{a}_{0}=\dfrac{k_{0}}{m_{0}},\breve{a}_{1}=\dfrac{k_{1}}{m_{1}}%
,...,\breve{a}_{n}=\dfrac{k_{n}}{m_{n}},...\in 
\mathbb{Q}
, \\ 
\\ 
\ \breve{a}_{0}>0 \\ 
\end{array}
& \text{ \ }\left( 3.2.4\right)%
\end{array}%
$

$\ \ \ \ \ \ \ \ \ \ \ \ \ \ \ \ \ \ \ \ \ \ \ \ \ \ \ \ \ \ \ \ \ \ \ $\ \
\ \ \ \ \ \ \ \ \ \ \ \ \ \ 

such that the next equality is satisfied:

\bigskip

$\ \ \ \ \ \ \ \ \ \ \ \ \ \ \ \ \ \ \ \ \ \ \ \ \ \ \ 
\begin{array}{cc}
\begin{array}{c}
\\ 
\dsum\limits_{n=0}^{\infty }\breve{a}_{n}e^{n}=0. \\ 
\end{array}
& \text{ \ }\left( 3.2.5\right)%
\end{array}%
$

\bigskip

Hence there is exist sequences$\ \ \left\{ n_{i}\right\} _{i=0}^{\infty }$
and $\ \left\{ n_{j}\right\} _{j=1}^{\infty }$ such that$\ $

$\bigskip $

$\ \ \ \ \ \ \ \ \ \ \ \ \ \ \ \ \ \ \ \ \ \ \ \ \ \ \ \ \ \ $

$\ \ \ \ 
\begin{array}{cc}
\begin{array}{c}
\\ 
\forall k\left[ \dsum\limits_{n=0}^{n_{i}\leq k}\breve{a}_{n}e^{n}>0\right] ,%
\underset{i\rightarrow \infty }{\lim }\dsum\limits_{n=0}^{n_{i}}\breve{a}%
_{n}e^{n}=0, \\ 
\\ 
\forall k\left[ \dsum\limits_{n=1}^{n_{j}\leq k}\breve{a}_{n}e^{n}<0\right] ,%
\underset{j\rightarrow \infty }{\lim }\dsum\limits_{n=1}^{n_{j}}\breve{a}%
_{n}e^{n}=0. \\ 
\end{array}
& \text{ \ \ \ }\left( 3.2.6\right)%
\end{array}%
$ \ \ \ \ \ \ \ \ \ \ \ \ \ 

\bigskip 

From Eqs.(3.2.6) by using definitions one obtain \bigskip the next \ \ \ \ \
\ \ \ \ \ \ \ \ \ \ \ \ \ \ \ \ \ \ \ \ \ \ 

equalities:

\ \ \ \ \ \ \ \ \ \ \ \ \ \ \ \ \ \ \ \ \ \ \ \ \ \ \ $\ \ \ \ \ \ \ \ \ \ \
\ \ \ 
\begin{array}{cc}
\begin{array}{c}
\\ 
\left( ^{\ast }\breve{a}_{0}\right) ^{\#}+\left[ \overline{\overline{\#Ext%
\text{-}\dsum\limits_{n\in 
\mathbb{N}
\backslash \left\{ 0\right\} }\left( ^{\ast }\breve{a}_{n}\right)
^{\#}\times \left( ^{\ast }e^{n}\right) ^{\#}}}\right] _{\varepsilon }= \\ 
\\ 
=\left[ \left( ^{\ast }\left( \underset{i\in 
\mathbb{N}
}{\inf }\left( \dsum\limits_{n=0}^{n_{i}}\breve{a}_{n}e^{n}\right) \right)
\right) ^{\#}+\varepsilon _{\mathbf{d}}\right] _{\varepsilon }= \\ 
\\ 
\text{ }=\left( ^{\ast }\left( \underset{i\in 
\mathbb{N}
}{\inf }\left( \dsum\limits_{n=0}^{n_{i}}\breve{a}_{n}e^{n}\right) \right)
\right) ^{\#}+\varepsilon ^{\#}\times \varepsilon _{\mathbf{d}}=\varepsilon
^{\#}\times \varepsilon _{\mathbf{d}}, \\ 
\\ 
\varepsilon \approx 0,\varepsilon \in \text{ }^{\ast }%
\mathbb{R}
\\ 
\end{array}
& \text{ \ \ }\left( 3.2.7\right)%
\end{array}%
$

$\bigskip $and by similar way

$\ \ \ 
\begin{array}{cc}
\begin{array}{c}
\\ 
\left( ^{\ast }\breve{a}_{0}\right) ^{\#}+\left[ \underline{\underline{\#Ext%
\text{-}\dsum\limits_{n\in 
\mathbb{N}
\backslash \left\{ 0\right\} }\left( ^{\ast }\breve{a}_{n}\right)
^{\#}\times \left( ^{\ast }e^{n}\right) ^{\#}}}\right] _{\varepsilon }= \\ 
\\ 
\text{ }=\text{ }\underset{j\in 
\mathbb{N}
}{\sup }\left( \dsum\limits_{n=1}^{n_{j}}\breve{a}_{n}e^{n}\right)
-\varepsilon ^{\#}\times \varepsilon _{\mathbf{d}}=-\varepsilon ^{\#}\times
\varepsilon _{\mathbf{d}}, \\ 
\\ 
\varepsilon \approx 0,\varepsilon \in \text{ }^{\ast }%
\mathbb{R}
. \\ 
\end{array}
& \text{ \ }\left( 3.2.8\right)%
\end{array}%
$

\bigskip

Let us considered hypernatural number $\Im \in $ $^{\ast }%
\mathbb{N}
_{\infty }$ defined by

countable sequence

\bigskip

\ \ $\ \ 
\begin{array}{cc}
\begin{array}{c}
\\ 
\Im =\left( m_{0},m_{0}\times m_{1},...,m_{0}\times m_{1}\times ...\times
m_{n},...\right)  \\ 
\end{array}
& \text{ \ }\left( 3.2.9\right) 
\end{array}%
$\bigskip 

By using Eq.(3.2.7) and Eq.(3.2.9) one obtain

\bigskip

$\ \ \ \ \ 
\begin{array}{cc}
\begin{array}{c}
\\ 
\Im ^{\#}\times \left( ^{\ast }\breve{a}_{0}\right) ^{\#}+\Im ^{\#}\times %
\left[ \overline{\overline{\#Ext\text{-}\dsum\limits_{n\in 
\mathbb{N}
\backslash \left\{ 0\right\} }\left( ^{\ast }\breve{a}_{n}\right)
^{\#}\times \left( ^{\ast }e^{n}\right) ^{\#}}}\right] _{\varepsilon }= \\ 
\\ 
\Im ^{\#}\times \varepsilon ^{\#}\times \varepsilon _{\mathbf{d}}, \\ 
\\ 
\Im ^{\#}\times \left( ^{\ast }\breve{a}_{0}\right) ^{\#}+ \\ 
\\ 
+\left[ \overline{\overline{\#Ext\text{-}\dsum\limits_{n\in 
\mathbb{N}
\backslash \left\{ 0\right\} }\Im ^{\#}\times \left( ^{\ast }\breve{a}%
_{n}\right) ^{\#}\times \left( ^{\ast }e^{n}\right) ^{\#}}}\left\vert \Im
^{\#}\times \left( ^{\ast }c\right) ^{\#}\right. \right] _{\varepsilon }= \\ 
\\ 
=\Im _{0}^{\#}+\left[ \overline{\overline{\#\text{-}Ext\text{-}%
\dsum\limits_{n\in 
\mathbb{N}
\backslash \left\{ 0\right\} }\Im _{n}^{\#}\times \left( ^{\ast
}e^{n}\right) ^{\#}}}\left\vert \Im ^{\#}\times \left( ^{\ast }c\right)
^{\#}\right. \right] = \\ 
\\ 
=\Im ^{\#}\times \varepsilon ^{\#}\times \varepsilon _{\mathbf{d}}, \\ 
\\ 
c\in 
\mathbb{R}
, \\ 
\\ 
\Im _{n}^{\#}\triangleq \Im ^{\#}\times \left( ^{\ast }\breve{a}_{n}\right)
^{\#},n\in 
\mathbb{N}
. \\ 
\end{array}
& \text{ }\left( 3.2.10\right)%
\end{array}%
$

$\ \ \ \ \ \ \ \ \ \ \ \ \ \ \ \ \ \ $

$\ \ \ \ \ \ \ \ \ \ \ \ \ \ \ \ \ \ \ \ \ \ \ \ \ \ \ $

\bigskip

Now\ we have to pruve that Eq.(3.2.10) leads to contradiction.

$\mathbf{Proof}$\textbf{\ }$\mathbf{I.}$

\textbf{Part I. }Let be

\bigskip

\textbf{\ }$\ \ 
\begin{array}{cc}
\begin{array}{c}
\\ 
M_{0}\left( n,p\right) =\dint\limits_{0}^{+\infty }\left[ \dfrac{x^{p-1}%
\left[ \left( x-1\right) ...\left( x-n\right) \right] ^{p}e^{-x}}{\left(
p-1\right) !}\right] dx\neq 0, \\ 
\end{array}
& \text{ }\left( 3.2.11\right)%
\end{array}%
$

\bigskip

$\ \ \ 
\begin{array}{cc}
\begin{array}{c}
\\ 
M_{k}\left( n,p\right) =e^{k}\dint\limits_{k}^{+\infty }\left[ \dfrac{x\left[
^{p-1}\left( x-1\right) ...\left( x-n\right) \right] ^{p}e^{-x}}{\left(
p-1\right) !}\right] dx, \\ 
\\ 
k=1,2,... \\ 
\end{array}
& \text{ \ }\left( 3.2.12\right)%
\end{array}%
$

\bigskip\ $\ \ \ \ \ \ \ \ \ \ \ \ $

$\ 
\begin{array}{cc}
\begin{array}{c}
\\ 
\varepsilon _{k}\left( n,p\right) =e^{k}\dint\limits_{0}^{k}\left[ \dfrac{%
x^{p-1}\left[ \left( x-1\right) ...\left( x-n\right) ^{p}\right] e^{-x}}{%
\left( p-1\right) !}\right] dx, \\ 
\\ 
k=1,2,..., \\ 
\end{array}
& \text{ }\left( 3.2.13\right)%
\end{array}%
$

\bigskip

where $p\in 
\mathbb{N}
$ this is any prime number.Using Eq.(3.2.9)-Eq.(3.2.13) \ \ \ \ \ \ \ \ \ \
\ \ \ \ \ \ \ \ \ \ \ \ \ \ \ \ \ \ \ \ by simple calculation one obtain:$\ $

$\ \ \ \ \ \ \ \ \ \ \ \ \ \ \ \ \ \ \ \ \ \ \ \ \ \ \ \ \ \ \ \ \ \ \ $

$\ \ \ \ \ \ \ 
\begin{array}{cc}
\begin{array}{c}
\\ 
M_{k}\left( n,p\right) +\varepsilon _{k}\left( n,p\right) =e^{k}M_{0}\neq 0,
\\ 
\\ 
k=1,2,... \\ 
\end{array}
& \text{ \ \ }\left( 3.2.14\right)%
\end{array}%
$

\bigskip

and consequently

\bigskip

$\ \ \ \ \ \ \ \ \ \ 
\begin{array}{cc}
\begin{array}{c}
\\ 
e^{k}=\dfrac{M_{k}\left( n,p\right) +\varepsilon _{k}\left( n,p\right) }{%
M_{0}} \\ 
\\ 
k=1,2,... \\ 
\end{array}
& \text{ \ }\left( 3.2.15\right)%
\end{array}%
$

\bigskip By using equality

$\bigskip $

$\ 
\begin{array}{cc}
\begin{array}{c}
\\ 
x^{p-1}\left[ \left( x-1\right) ...\left( x-n\right) \right] ^{p}= \\ 
\\ 
\left( -1\right) ^{n}\left( n!\right) ^{n}x^{n-1}+\dsum\limits_{\mu
=p+1}^{\left( n+1\right) \times p}c_{\mu -1}x^{\mu -1}, \\ 
\\ 
c_{\mu }\in 
\mathbb{Z}
,\mu =p,p+1,...,\left[ (n+1)\times p\right] -1, \\ 
\end{array}
& \text{ \ }\left( 3.2.16\right)%
\end{array}%
$

\bigskip

\bigskip from Eq.(3.2.11) one obtain

\bigskip

\ \ \ \ \ 

\ \ \ $\ \ 
\begin{array}{cc}
\begin{array}{c}
\\ 
M_{0}\left( n,p\right) =\left( -1\right) ^{n}\left( n!\right) ^{p}\dfrac{%
\Gamma \left( p\right) }{\left( p-1\right) !}+ \\ 
\\ 
\dsum\limits_{\mu =p+1}^{\left( n+1\right) \times p}c_{\mu -1}\dfrac{\Gamma
\left( \mu \right) }{\left( p-1\right) !}= \\ 
\\ 
=\left( -1\right) ^{n}\left( n!\right) ^{p}+c_{p}p+c_{n+1}p\left( p+1\right)
+...= \\ 
\\ 
=\left( -1\right) ^{n}\left( n!\right) ^{p}+p\times \Theta _{1},\Theta
_{1}\in 
\mathbb{Z}
, \\ 
\\ 
\Gamma \left( \mu \right) =\dint\limits_{0}^{+\infty }x^{\mu -1}e^{-x}dx. \\ 
\\ 
M_{0}\left( n,p\right) =\left( -1\right) ^{n}\left( n!\right) ^{p}+p\times
\Theta _{1},\Theta _{1}\in 
\mathbb{Z}
\\ 
\end{array}
& \text{ \ }\left( 3.2.17\right)%
\end{array}%
$

Thus

\bigskip

\bigskip\ \ $\ \ \ 
\begin{array}{cc}
\begin{array}{c}
\\ 
M_{0}\left( n,p\right) =\left( -1\right) ^{n}\left( n!\right) ^{p}+ \\ 
\\ 
p\cdot \Theta _{1}\left( n,p\right) ,\Theta _{1}\left( n,p\right) \in 
\mathbb{Z}
. \\ 
\end{array}
& \text{ \ }\left( 3.2.18\right)%
\end{array}%
$

By subsitution $x=k+u\implies dx=du$ from Eq.(3.2.13.) one obtain

\bigskip

\bigskip $\ 
\begin{array}{cc}
\begin{array}{c}
\\ 
M_{k}\left( n,p\right) = \\ 
\dint\limits_{0}^{+\infty }\left[ \dfrac{\left( u+k\right) ^{p-1}\left[
\left( u+k-1\right) \times ...\times u\times ...\times \left( u+k-n\right) %
\right] ^{p}e^{-u}}{\left( p-1\right) !}\right] du \\ 
\\ 
k=1,2,3,... \\ 
\end{array}
& \text{ }\left( 3.2.19\right)%
\end{array}%
$

\bigskip

By using equality

\bigskip

$\ \ \ \ 
\begin{array}{cc}
\begin{array}{c}
\\ 
\left( u+k\right) ^{p-1}\left[ \left( u+k-1\right) \times ...\times u\times
...\times \left( u+k-n\right) \right] ^{p}= \\ 
\\ 
=\dsum\limits_{\mu =p+1}^{\left( n+1\right) \times p}d_{\mu -1}u^{\mu -1},
\\ 
\\ 
d_{\mu }\in 
\mathbb{Z}
,\mu =p,p+1,...,\left[ (n+1)\times p\right] -1, \\ 
\end{array}
& \text{ \ }\left( 3.2.20\right)%
\end{array}%
$

\bigskip

and by subsitution Eq.(3.2.20) into Eq.(3.2.19) one obtain

\bigskip

\bigskip

$\ \ \ \ \ \ 
\begin{array}{cc}
\begin{array}{c}
\\ 
M_{k}\left( n,p\right) =\dfrac{1}{\left( p-1\right) !}\dint\limits_{0}^{+%
\infty }\dsum\limits_{\mu =p+1}^{\left( n+1\right) \times p}d_{\mu -1}u^{\mu
-1}du= \\ 
\\ 
p\cdot \Theta _{2}\left( n,p\right) , \\ 
\\ 
\Theta _{2}\left( n,p\right) \in 
\mathbb{Z}
, \\ 
\\ 
k=1,2,...\text{ }. \\ 
\end{array}
& \text{ \ }\left( 3.2.21\right)%
\end{array}%
$

\bigskip

There is exists sequences $a\left( n\right) ,n\in 
\mathbb{N}
$ and $g_{k}\left( n\right) ,k\in 
\mathbb{N}
,n\in 
\mathbb{N}
$ such that

\ \ \ \ \ \ \ \ \ \ \ \ \ \ \ \ \ \ \ \ \ \ \ \ \ \ \ \ \ \ \ \ \ \ \ \ 

$\ \ \ \ 
\begin{array}{cc}
\begin{array}{c}
\\ 
\left\vert x\left( x-1\right) ...\left( x-n\right) \right\vert <a\left(
n\right) , \\ 
\\ 
0\leq x\leq n, \\ 
\\ 
\left\vert x\left( x-1\right) ...\left( x-n\right) e^{-x+k}\right\vert
<g_{k}\left( n\right) , \\ 
\\ 
0\leq x\leq n,k=1,2,...\text{ }. \\ 
\end{array}
& \text{ \ \ }\left( 3.2.22\right)%
\end{array}%
$

\bigskip

Substitution the inequalities (3.2.22.) into Eq.(3.2.13.) gives

\bigskip

$\ 
\begin{array}{cc}
\begin{array}{c}
\\ 
\varepsilon _{k}\left( n,p\right) \leq g_{k}\left( n\right) \dfrac{\left[
a\left( n\right) \right] ^{p-1}}{\left( p-1\right) !}\dint\limits_{0}^{k}dx%
\leq \\ 
\\ 
\leq \dfrac{n\cdot g_{k}\left( n\right) \cdot \left[ a\left( n\right) \right]
^{p-1}}{\left( p-1\right) !}. \\ 
\end{array}
& \text{ }\left( 3.2.23\right)%
\end{array}%
$

\bigskip

By using transfer, from Eq.(3.2.11.) and Eq.(3.2.18.) one obtain

\bigskip

$\ \ \ \ \ \ \ \ \ 
\begin{array}{cc}
\begin{array}{c}
\\ 
^{\ast }M_{0}\left( \mathbf{n,p}\right) = \\ 
\\ 
\text{ }^{\ast }\left( \dint\limits_{0}^{+\infty }\left[ \dfrac{x^{p-1}\left[
\left( x-1\right) ...\left( x-n\right) \right] ^{p}e^{-x}}{\left( p-1\right)
!}\right] dx\right) = \\ 
\\ 
=\left( -1\right) ^{\mathbf{n}}\left( \mathbf{n}!\right) ^{\mathbf{p}}+%
\mathbf{p}\times \text{ }^{\ast }\Theta _{1}\left( \mathbf{n},\mathbf{p}%
\right) , \\ 
\\ 
^{\ast }\Theta _{1}\left( \mathbf{n},\mathbf{p}\right) \in \text{ }^{\ast }%
\mathbb{Z}
_{\infty }, \\ 
\\ 
\mathbf{n},\mathbf{p\in }^{\ast }\mathbf{%
\mathbb{N}
}_{\infty }. \\ 
\end{array}
& \text{ }\left( 3.2.24\right)%
\end{array}%
$

\bigskip

From Eq.(3.2.12.) and Eq.(3.2.21) one obtain

\bigskip

$\ \ \ \ \ \ \ 
\begin{array}{cc}
\begin{array}{c}
\\ 
M_{k}\left( n,p\right) =e^{k}\dint\limits_{k}^{+\infty }\left[ \dfrac{x\left[
^{p-1}\left( x-1\right) ...\left( x-n\right) \right] ^{p}e^{-x}}{\left(
p-1\right) !}\right] dx= \\ 
\\ 
\dint\limits_{0}^{+\infty }\left[ \dfrac{\left( u+k\right) ^{p-1}\left[
\left( u+k-1\right) \times ...\times u\times ...\times \left( u+k-n\right) %
\right] ^{p}e^{-u}}{\left( p-1\right) !}\right] du= \\ 
\\ 
=p\cdot \Theta _{2}\left( n,p\right) , \\ 
\\ 
\Theta _{2}\left( n,p\right) \in 
\mathbb{Z}
, \\ 
\\ 
k\in 
\mathbb{N}
. \\ 
\end{array}
& \text{ }\left( 3.2.25\right)%
\end{array}%
$

\bigskip

Using transfer, from Eq.(3.2.25.) one obtain $\forall k\left( k\in 
\mathbb{N}
\right) :$

\bigskip

$\ \ 
\begin{array}{cc}
\begin{array}{c}
\\ 
^{\ast }M_{k}\left( \mathbf{n,p}\right) =\left( ^{\ast }e^{k}\right) \times 
\text{ }^{\ast }\left( \dint\limits_{k}^{+\infty }\left[ \dfrac{x\left[ ^{%
\mathbf{p}-1}\left( x-1\right) ...\left( x-\mathbf{n}\right) \right] ^{%
\mathbf{p}}e^{-x}}{\left( \mathbf{p}-1\right) !}\right] dx\right) = \\ 
\\ 
=\text{ }^{\ast }\left( \dint\limits_{0}^{+\infty }\left[ \dfrac{\left(
u+k\right) ^{\mathbf{p}-1}\left[ \left( u+k-1\right) \times ...\times
u\times ...\times \left( u+k-\mathbf{n}\right) \right] ^{\mathbf{p}}e^{-u}}{%
\left( \mathbf{p}-1\right) !}\right] du\right) = \\ 
\\ 
=\mathbf{p}\times \text{ }^{\ast }\Theta _{2}\left( \mathbf{n,p}\right) , \\ 
\\ 
^{\ast }\Theta _{2}\left( \mathbf{n,p}\right) \in \text{ }^{\ast }%
\mathbb{Z}
_{\infty }, \\ 
\\ 
k=1,2,3,... \\ 
\\ 
k\in 
\mathbb{N}
, \\ 
\\ 
\mathbf{n},\mathbf{p\in }^{\ast }\mathbf{%
\mathbb{N}
}_{\infty }. \\ 
\end{array}
& \text{ }\left( 3.2.26\right)%
\end{array}%
$

\bigskip

Using transfer, from inequality (3.2.23.) one obtain $\forall k\left( k\in 
\mathbb{N}
\right) :$

\bigskip

$\ \ \ \ \ 
\begin{array}{cc}
\begin{array}{c}
\\ 
^{\ast }\varepsilon _{k}\left( \mathbf{n,p}\right) \leq \text{ }^{\ast
}g_{k}\left( \mathbf{n}\right) \times \dfrac{\left[ ^{\ast }a\left( \mathbf{n%
}\right) \right] ^{\mathbf{p}-1}}{\left( \mathbf{p}-1\right) !}\times \left[ 
\text{ }^{\ast }\left( \dint\limits_{0}^{k}dx\right) \right] \leq \\ 
\\ 
\leq \dfrac{\mathbf{n}\cdot \left[ ^{\ast }g_{k}\left( \mathbf{n}\right) %
\right] \cdot \left[ ^{\ast }a\left( \mathbf{n}\right) \right] ^{\mathbf{p}%
-1}}{\left( \mathbf{p}-1\right) !}, \\ 
\\ 
k=1,2,...,k\in 
\mathbb{N}
, \\ 
\\ 
\mathbf{n},\mathbf{p\in }^{\ast }\mathbf{%
\mathbb{N}
}_{\infty }. \\ 
\end{array}
& \text{\ }\left( 3.2.27\right)%
\end{array}%
$

\bigskip

By using transfer again, from Eq.(3.2.15.) one obtain $\forall k\left( k\in 
\mathbb{N}
\right) :$

\bigskip

$\ \ \ \ \ \ \ \ \ \ \ \ \ \ 
\begin{array}{cc}
\begin{array}{c}
\\ 
\text{ }^{\ast }\left( e^{k}\right) =\left( ^{\ast }e\right) ^{k}=\dfrac{%
^{\ast }M_{k}\left( \mathbf{n,p}\right) +\text{ }^{\ast }\varepsilon
_{k}\left( \mathbf{n,p}\right) }{^{\ast }M_{0}\left( \mathbf{n,p}\right) },
\\ 
\\ 
k=1,2,..., \\ 
\\ 
k\in 
\mathbb{N}
, \\ 
\\ 
\mathbf{n},\mathbf{p\in }^{\ast }\mathbf{%
\mathbb{N}
}_{\infty }. \\ 
\end{array}
& \text{ \ }\left( 3.2.28\right)%
\end{array}%
\ \ \ \ \ \ \ \ \ \ \ \ \ \ \ \ \ \ \ \ \ \ $

$\bigskip $

\textbf{(Part II) }By using Eq.(3.2.28.) one obtain

\bigskip

$\ \ \ \ \ \ \ \ 
\begin{array}{cc}
\begin{array}{c}
\\ 
\text{ }\left[ ^{\ast }\left( e^{k}\right) \right] ^{\#}=\left[ \left(
^{\ast }e\right) ^{\#}\right] ^{k}= \\ 
\\ 
\dfrac{\left[ ^{\ast }M_{k}\left( \mathbf{n,p}\right) \right] ^{\#}+\text{ }%
\left[ ^{\ast }\varepsilon _{k}\left( \mathbf{n,p}\right) \right] ^{\#}}{%
\left[ ^{\ast }M_{0}\left( \mathbf{n,p}\right) \right] ^{\#}}, \\ 
\\ 
k=1,2,..., \\ 
\\ 
k\in 
\mathbb{N}
, \\ 
\\ 
\mathbf{n},\mathbf{p\in }^{\ast }\mathbf{%
\mathbb{N}
}_{\infty }. \\ 
\end{array}
& \text{ }\left( 3.2.29\right)%
\end{array}%
$

\bigskip

By using Eq.(3.2.24.) one obtain

\bigskip

$\ 
\begin{array}{cc}
\begin{array}{c}
\\ 
\left[ ^{\ast }M_{0}\left( \mathbf{n,p}\right) \right] ^{\#}=\left[ \left(
-1\right) ^{\mathbf{n}}\right] ^{\#}\left[ \left( \mathbf{n}!\right) ^{%
\mathbf{p}}\right] ^{\#}+ \\ 
\\ 
\mathbf{p}^{\#}\times \text{ }\left[ ^{\ast }\Theta _{1}\left( \mathbf{n},%
\mathbf{p}\right) \right] ^{\#}, \\ 
\\ 
\left[ ^{\ast }\Theta _{1}\left( \mathbf{n},\mathbf{p}\right) \right]
^{\#}\in \text{ }^{\ast }%
\mathbb{Z}
_{\infty ,\mathbf{d}}, \\ 
\\ 
\mathbf{n},\mathbf{p\in }^{\ast }\mathbf{%
\mathbb{N}
}_{\infty }. \\ 
\end{array}
& \text{ \ }\left( 3.2.30\right)%
\end{array}%
$

\bigskip

By using Eq.(3.2.26.) one obtain

\bigskip

$\ \ \ 
\begin{array}{cc}
\begin{array}{c}
\\ 
\left[ ^{\ast }M_{k}\left( \mathbf{n,p}\right) \right] ^{\#}=\mathbf{p}%
^{\#}\times \left[ \text{ }^{\ast }\Theta _{2}\left( \mathbf{n,p}\right) %
\right] ^{\#}, \\ 
\\ 
\left[ ^{\ast }\Theta _{2}\left( \mathbf{n,p}\right) \right] ^{\#}\in \text{ 
}^{\ast }%
\mathbb{Z}
_{\infty ,\mathbf{d}}, \\ 
\\ 
k=1,2,...,k\in 
\mathbb{N}
, \\ 
\\ 
\mathbf{n},\mathbf{p\in }^{\ast }\mathbf{%
\mathbb{N}
}_{\infty }. \\ 
\end{array}
& \text{ \ }\left( 3.2.31\right)%
\end{array}%
\ \ \ \ \ \ \ \ \ \ \ \ \ \ \ \ \ \ \ \ \ \ \ \ \ \ \ \ \ \ \ \ \ \ \ \ \ \
\ \ \ $

\bigskip

By using inequality (3.2.27) one obtain

\bigskip

$\ 
\begin{array}{cc}
\begin{array}{c}
\\ 
\left[ ^{\ast }\varepsilon _{k}\left( \mathbf{n,p}\right) \right] ^{\#}\leq 
\dfrac{\mathbf{n}^{\#}\cdot \left[ g_{k}\left( \mathbf{n}\right) \right]
^{\#}\cdot \left[ \left[ a\left( \mathbf{n}\right) \right] ^{\mathbf{p}-1}%
\right] ^{\#}}{\left[ \left( \mathbf{p}-1\right) !\right] ^{\#}}, \\ 
\\ 
k=1,2,...,k\in 
\mathbb{N}
, \\ 
\\ 
\mathbf{n},\mathbf{p\in }^{\ast }\mathbf{%
\mathbb{N}
}_{\infty }. \\ 
\end{array}
& \text{ }\left( 3.2.32\right)%
\end{array}%
$

\bigskip

Substitution Eq.(3.2.28) into Eq.(3.2.10) gives

\bigskip

\bigskip $\ 
\begin{array}{cc}
\begin{array}{c}
\\ 
\Im _{0}^{\#}+\left[ \overline{\overline{\#Ext\text{-}\dsum\limits_{n\in 
\mathbb{N}
\backslash \left\{ 0\right\} }\Im _{n}^{\#}\times \left( ^{\ast
}e^{n}\right) ^{\#}}}\left\vert \Im ^{\#}\times \left( ^{\ast }c\right)
^{\#}\right. \right] _{\varepsilon }= \\ 
\\ 
\Im _{0}^{\#}+\left[ \overline{\overline{\#Ext\text{-}\dsum\limits_{k=1}^{%
\infty }\Im _{k}^{\#}\times \dfrac{\left[ ^{\ast }M_{k}\left( \mathbf{n,p}%
\right) \right] ^{\#}+\text{ }\left[ ^{\ast }\varepsilon _{k}\left( \mathbf{%
n,p}\right) \right] ^{\#}}{\left[ ^{\ast }M_{0}\left( \mathbf{n,p}\right) %
\right] ^{\#}}}}\left\vert \Im ^{\#}\times \left( ^{\ast }c\right)
^{\#}\right. \right] _{\varepsilon } \\ 
\\ 
=\Im ^{\#}\times \varepsilon ^{\#}\times \varepsilon _{\mathbf{d}},. \\ 
\end{array}
& \left( 3.2.33\right)%
\end{array}%
$

\bigskip

Multiplying Eq.(3.2.33) by number $\left[ ^{\ast }M_{0}\left( \mathbf{n,p}%
\right) \right] ^{\#}\in $ $^{\ast }%
\mathbb{Z}
_{\mathbf{d}}$ one obtain

\bigskip

$\ \ \ \ 
\begin{array}{cc}
\begin{array}{c}
\\ 
\Im _{0}^{\#}\times \text{ }\left[ ^{\ast }M_{0}\left( \mathbf{n,p}\right) %
\right] ^{\#}+ \\ 
\\ 
\left[ \overline{\overline{\#Ext\text{-}\dsum\limits_{k=1}^{\infty }\text{ }%
\left\{ \Im _{k}^{\#}\times \left[ ^{\ast }M_{k}\left( \mathbf{n,p}\right) %
\right] ^{\#}+\Im _{k}^{\#}\times \left[ ^{\ast }\varepsilon _{k}\left( 
\mathbf{n,p}\right) \right] ^{\#}\right\} }}\right. \\ 
\\ 
\left. \left\vert \Im ^{\#}\times \left[ ^{\ast }M_{0}\left( \mathbf{n,p}%
\right) \right] ^{\#}\times \left( ^{\ast }c\right) ^{\#}\right. \right]
_{\varepsilon }= \\ 
\\ 
=\Im ^{\#}\times \text{ }\left[ ^{\ast }M_{0}\left( \mathbf{n,p}\right) %
\right] ^{\#}\times \varepsilon ^{\#}\times \varepsilon _{\mathbf{d}}. \\ 
\end{array}
& \text{ }\left( 3.2.34\right)%
\end{array}%
$

\bigskip

By using inequality (3.2.32) for we will choose prime hyper number $\mathbf{%
p\in }^{\ast }\mathbf{%
\mathbb{N}
}_{\infty }$

given $\varepsilon $ such that:$\ $

\bigskip

$\ 
\begin{array}{cc}
\begin{array}{c}
\\ 
\left[ \overline{\overline{\#Ext\text{-}\dsum\limits_{k=1}^{\infty }\text{ }%
\Im _{k}^{\#}\times \left[ ^{\ast }\varepsilon _{k}\left( \mathbf{n,p}%
\right) \right] ^{\#}}}\left\vert \Im ^{\#}\times \left( ^{\ast }c\right)
^{\#}\right. \right] _{\varepsilon }\in \\ 
\\ 
\in \Im ^{\#}\times \text{ }\left[ ^{\ast }M_{0}\left( \mathbf{n,p}\right) %
\right] ^{\#}\times \varepsilon \times \varepsilon _{\mathbf{d}}. \\ 
\end{array}
& \text{\ }\left( 3.2.35\right)%
\end{array}%
$

\bigskip

Hence from Eq.(3.2.34) and Eq.(3.2.35) one obtain

\bigskip

$\ \ \ \ 
\begin{array}{cc}
\begin{array}{c}
\\ 
\Im _{0}^{\#}\times \text{ }\left[ ^{\ast }M_{0}\left( \mathbf{n,p}\right) %
\right] ^{\#}+ \\ 
\\ 
+\left[ \overline{\overline{\#Ext\text{-}\dsum\limits_{k=1}^{\infty }\text{ }%
\Im _{k}^{\#}\times \left[ ^{\ast }M_{k}\left( \mathbf{n,p}\right) \right]
^{\#}}}\left\vert \Im ^{\#}\times \left[ ^{\ast }M_{0}\left( \mathbf{n,p}%
\right) \right] ^{\#}\times \left( ^{\ast }c\right) ^{\#}\right. \right]
_{\varepsilon }= \\ 
\\ 
=\Im ^{\#}\times \text{ }\left[ ^{\ast }M_{0}\left( \mathbf{n,p}\right) %
\right] ^{\#}\times \varepsilon ^{\#}\times \varepsilon _{\mathbf{d}}. \\ 
\end{array}
& \text{ }\left( 3.2.36\right)%
\end{array}%
$

\bigskip

\bigskip We will choose prime hyper number $\mathbf{p\in }^{\ast }\mathbf{%
\mathbb{N}
}_{\infty }$ such that

$\ 
\begin{array}{cc}
\begin{array}{c}
\\ 
\mathbf{p}^{\#}\mathbf{>\max }\left( \Im ^{\#},\left\vert \Im
_{0}^{\#}\right\vert ,\mathbf{n}^{\#}\mathbf{.}\right) \\ 
\end{array}
& \text{ \ }\left( 3.2.37\right)%
\end{array}%
$

\bigskip

Hence by using Eq.(3.2.20) one obtain:

\bigskip

$\ \ \ 
\begin{array}{cc}
\begin{array}{c}
\\ 
\left[ ^{\ast }M_{0}\left( \mathbf{n,p}\right) \right] ^{\#}\nmid \mathbf{p}%
^{\#} \\ 
\end{array}
& \text{ \ }\left( 3.2.38\right)%
\end{array}%
$

\bigskip

and \ consequently $\left[ ^{\ast }M_{0}\left( \mathbf{n,p}\right) \right]
^{\#}\neq 0^{\#}.$And by using (3.2.20),(3.2.28) 

one obtain:

\bigskip

$\ 
\begin{array}{cc}
\begin{array}{c}
\\ 
\left[ ^{\ast }M_{0}\left( \mathbf{n,p}\right) \right] ^{\#}\times \Im
_{0}^{\#}\nmid \mathbf{p}^{\#}\mathbf{.} \\ 
\end{array}
& \text{ \ }\left( 3.2.39\right) 
\end{array}%
$

By using Eq.(3.2.22) one obtain

\bigskip

$\ \ 
\begin{array}{cc}
\begin{array}{c}
\\ 
\left[ ^{\ast }M_{k}\left( \mathbf{n,p}\right) \right] ^{\#}\mid \mathbf{p}%
^{\#}\mathbf{,} \\ 
\\ 
k=1,2,...\mathbf{.} \\ 
\end{array}
& \text{ \ \ \ \ }\left( 3.2.40\right)%
\end{array}%
$

\bigskip

\bigskip By using Eq.(3.2.36) one obtain

$\ \ \ \ \ \ \ \ \ \ \ \ \ \ \ \ \ \ \ \ \ \ \ \ \ \ $

$\ \ \ \ 
\begin{array}{cc}
\begin{array}{c}
\\ 
\Xi \left( \mathbf{n,p,}\varepsilon \right) =\mathbf{a.p.}\left\{ \left[ \Xi
\left( \mathbf{n,p}\right) \right] _{\varepsilon }\right\} = \\ 
\\ 
\mathbf{ab.p.}\left\{ \left[ \overline{\overline{\#Ext\text{-}%
\dsum\limits_{k=1}^{\infty }\text{ }\Im _{k}^{\#}\times \left[ ^{\ast
}M_{k}\left( \mathbf{n,p}\right) \right] ^{\#}}}\left\vert \Im ^{\#}\times %
\left[ ^{\ast }M_{0}\left( \mathbf{n,p}\right) \right] ^{\#}\times \left(
^{\ast }c\right) ^{\#}\right. \right] _{\varepsilon }\right\} \\ 
\\ 
=\Im ^{\#}\times \text{ }\left[ ^{\ast }M_{0}\left( \mathbf{n,p}\right) %
\right] ^{\#}\times \varepsilon ^{\#}\times \varepsilon _{\mathbf{d}}. \\ 
\\ 
\Xi \left( \mathbf{n,p}\right) =\overline{\overline{\#Ext\text{-}%
\dsum\limits_{k=1}^{\infty }\text{ }\Im _{k}^{\#}\times \left[ ^{\ast
}M_{k}\left( \mathbf{n,p}\right) \right] ^{\#}}} \\ 
\end{array}
& \left( 3.2.41\right)%
\end{array}%
$

\bigskip

It is easy to see that Wattenberg hypernatural number $\Xi \left( \mathbf{n,p%
}\right) $ has

tipe $1\mathbf{.}$Hence Wattenberg hypernatural number $\Xi \left( \mathbf{%
n,p}\right) $ has

represantation: $\ \ \ \ \ \ \ \ \ $

$\bigskip $

$\ \ 
\begin{array}{cc}
\begin{array}{c}
\\ 
\Xi \left( \mathbf{n,p}\right) =\mathbf{p}^{\#}\mathbf{\times m+}\Im
^{\#}\times \text{ }\left[ ^{\ast }M_{0}\left( \mathbf{n,p}\right) \right]
^{\#}\times \varepsilon _{\mathbf{d}}, \\ 
\\ 
\mathbf{m\in }\text{ }^{\ast }\mathbf{%
\mathbb{Z}
}_{\mathbf{d}}. \\ 
\end{array}
& \text{ }\left( 3.2.42\right)%
\end{array}%
$

\bigskip

By using (3.2.42) one obtain represantation

\bigskip

$\ 
\begin{array}{cc}
\begin{array}{c}
\\ 
\Xi \left( \mathbf{n,p,\varepsilon }\right) =\left[ \Xi \left( \mathbf{n,p}%
\right) \right] _{\varepsilon }= \\ 
\\ 
\mathbf{p}^{\#}\mathbf{\times m+}\Im ^{\#}\times \text{ }\left[ ^{\ast
}M_{0}\left( \mathbf{n,p}\right) \right] ^{\#}\times \varepsilon ^{\#}\times
\varepsilon _{\mathbf{d}}, \\ 
\\ 
\mathbf{m\in }\text{ }^{\ast }\mathbf{%
\mathbb{Z}
}_{\mathbf{d}}. \\ 
\end{array}
& \text{ \ \ \ }\left( 3.2.43\right)%
\end{array}%
$

\bigskip

Substitution Eq.(3.2.43) into Eq.(3.2.36) gives

\bigskip

$\ 
\begin{array}{cc}
\begin{array}{c}
\\ 
\left\{ \Im _{0}^{\#}\times \text{ }\left[ ^{\ast }M_{0}\left( \mathbf{n,p}%
\right) \right] ^{\#}+\mathbf{p}^{\#}\mathbf{\times m}\right\} \mathbf{+} \\ 
\\ 
+\Im ^{\#}\times \text{ }\left[ ^{\ast }M_{0}\left( \mathbf{n,p}\right) %
\right] ^{\#}\times \varepsilon ^{\#}\times \varepsilon _{\mathbf{d}}= \\ 
\\ 
=\Im ^{\#}\times \text{ }\left[ ^{\ast }M_{0}\left( \mathbf{n,p}\right) %
\right] ^{\#}\times \varepsilon ^{\#}\times \varepsilon _{\mathbf{d}}. \\ 
\end{array}
& \text{ \ }\left( 3.2.45\right)%
\end{array}%
$

\bigskip

By using Eq.(3.2.39)-Eq.(3.2.40) one obtain:

\bigskip

$%
\begin{array}{cc}
\begin{array}{c}
\\ 
\left\{ \Im _{0}^{\#}\times \left[ ^{\ast }M_{0}\left( \mathbf{n,p}\right) %
\right] ^{\#}+\mathbf{p}^{\#}\mathbf{\times m}\right\} \nmid \mathbf{p}^{\#}
\\ 
\end{array}
& \text{ \ }\left( 3.2.46\right)%
\end{array}%
$

\bigskip

and consequently $\left\{ \Im _{0}^{\#}\times \left[ ^{\ast }M_{0}\left( 
\mathbf{n,p}\right) \right] ^{\#}+\mathbf{p}^{\#}\mathbf{\times m}\right\}
\neq 0^{\#}.$But on the

other hand, for sufficiently infinite smoll $\varepsilon \in $ $^{\ast }%
\mathbb{R}
$ idempotent

$\Im ^{\#}\times $ $\left[ ^{\ast }M_{0}\left( \mathbf{n,p}\right) \right]
^{\#}\times \varepsilon ^{\#}\times \varepsilon _{\mathbf{d}}$ does not
absorbs Wattenberg

hypernatural number $\left\{ \Im _{0}^{\#}\times \left[ ^{\ast }M_{0}\left( 
\mathbf{n,p}\right) \right] ^{\#}+\mathbf{p}^{\#}\mathbf{\times m}\right\} $
and

consequently:

\bigskip 

\bigskip $\ \ 
\begin{array}{cc}
\begin{array}{c}
\\ 
\left\{ \Im _{0}^{\#}\times \text{ }\left[ ^{\ast }M_{0}\left( \mathbf{n,p}%
\right) \right] ^{\#}+\mathbf{p}^{\#}\mathbf{\times m}\right\} \mathbf{+} \\ 
\\ 
+\Im ^{\#}\times \text{ }\left[ ^{\ast }M_{0}\left( \mathbf{n,p}\right) %
\right] ^{\#}\times \varepsilon ^{\#}\times \varepsilon _{\mathbf{d}}\neq \\ 
\\ 
\neq \Im ^{\#}\times \text{ }\left[ ^{\ast }M_{0}\left( \mathbf{n,p}\right) %
\right] ^{\#}\times \varepsilon ^{\#}\times \varepsilon _{\mathbf{d}}. \\ 
\end{array}
& \text{ \ \ }\left( 3.2.47\right)%
\end{array}%
$

\bigskip

Thus for sufficiently infinite small $\varepsilon \in $ $^{\ast }%
\mathbb{R}
$ inequality (3.2.47) in a

contradiction with Eq.(3.2.45).This contradiction proves that $e$ is not

$w$-transcendental. Hence $e$ is $\#$-transcendental.

\textbf{Proof II. (Part I) }

To prove $e$ is \#-transcendental we must show it is not $w$-transcendental,

i.e., there is no real analytic function $g_{%
\mathbb{Q}
}\left( x\right) =\dsum\limits_{n=0}^{\infty }b_{n}x^{n},e\leq \left\vert
x\right\vert \leq r$ with

rational coefficients $b_{0},b_{1},...,b_{n},...\in 
\mathbb{Q}
$ such that

\bigskip

\ $\ 
\begin{array}{cc}
\begin{array}{c}
\\ 
\dsum\limits_{n=0}^{\infty }b_{n}e^{n}=0. \\ 
\\ 
b_{0}=\dfrac{k_{0}}{m_{0}},b_{n}=\dfrac{k_{n}}{m_{n}}. \\ 
\end{array}
& \text{ \ }\left( 3.2.48\right)%
\end{array}%
$

\bigskip

\textbf{1.} Assume that $b_{0},b_{1},...,b_{n},...\in $ $%
\mathbb{Q}
,b_{0}\neq 0.$

Let $f\left( x\right) :%
\mathbb{R}
\rightarrow 
\mathbb{R}
$ be a polynomial of degree $m\in 
\mathbb{N}
.$ Then (repeated) integrations by

parts gives

\bigskip

$\ 
\begin{array}{cc}
\begin{array}{c}
\\ 
\dint\limits_{0}^{k}f\left( x\right) e^{-x}dx=\left. -f\left( x\right)
e^{-x}\right\vert _{0}^{k}+\dint\limits_{0}^{k}f^{\text{ }\prime }\left(
x\right) e^{-x}dx= \\ 
\\ 
=\left. -\left( f\left( x\right) +f^{\text{ }\prime }\left( x\right) +...+f^{%
\text{ }\left( m\right) }\left( x\right) \right) e^{-x}\right\vert _{0}^{k}
\\ 
\end{array}
& \text{ }\left( 3.2.49\right)%
\end{array}%
$

\bigskip

Multiply by $b_{k}e^{k}\in $ $^{\ast }%
\mathbb{Z}
,k=0,1,2,...,n\in 
\mathbb{N}
$ and add up: Then\bigskip

\bigskip

\ $\ 
\begin{array}{cc}
\begin{array}{c}
\\ 
\dsum\limits_{k=0}^{n}b_{k}e^{k}\dint\limits_{0}^{k}f\left( x\right)
e^{-x}dx= \\ 
=\left( f\left( 0\right) +f^{\text{ }\prime }\left( 0\right) +...+f^{\text{ }%
\left( m\right) }\left( 0\right) \right) \dsum\limits_{k=0}^{n}b_{k}e^{k}-
\\ 
\\ 
-\dsum\limits_{k=0}^{n}b_{k}\left( f\left( k\right) +f^{\text{ }\prime
}\left( k\right) +...+f^{\text{ }\left( m\right) }\left( k\right) \right) .
\\ 
\end{array}
& \text{ \ }\left( 3.2.50\right)%
\end{array}%
$

\bigskip

\bigskip By using transfer from Eq.(3.2.50) one obtain

$\ 
\begin{array}{cc}
\begin{array}{c}
\\ 
\dsum\limits_{k=0}^{n}\left( ^{\ast }b_{k}\left( ^{\ast }e^{k}\right)
\right) \times \text{ }^{\ast }\left( \dint\limits_{0}^{k}f\left( x\right)
e^{-x}dx\right) - \\ 
\\ 
-\left( ^{\ast }f\left( 0\right) +\text{ }^{\ast }f^{\text{ }\prime }\left(
0\right) +...+\text{ }^{\ast }f^{\text{ }\left( \mathbf{m}\right) }\left(
0\right) \right) \dsum\limits_{k=0}^{n}\left( ^{\ast }b_{k}^{\ast
}e^{k}\right) \\ 
\\ 
=-\dsum\limits_{k=0}^{n}\left( ^{\ast }b_{k}\right) \left( ^{\ast }f\left(
k\right) +\text{ }^{\ast }f^{\text{ }\prime }\left( k\right) +...+\text{ }%
^{\ast }f^{\text{ }\left( \mathbf{m}\right) }\left( k\right) \right) , \\ 
\\ 
\mathbf{m\in }^{\ast }\mathbf{%
\mathbb{N}
}_{\infty }. \\ 
\end{array}
& \text{ \ }\left( 3.2.51\right)%
\end{array}%
$

Hence

\ $\ \ \ \ \ \ \ \ \ \ \ \ \ \ \ \ \ \ \ \ \ \ \ $

$\ \ 
\begin{array}{cc}
\begin{array}{c}
\\ 
\dsum\limits_{k=0}^{n}\left( ^{\ast }b_{k}\left( ^{\ast }e^{k}\right)
\right) =\dsum\limits_{k=0}^{n}\left( ^{\ast }b_{k}\right) \left( \dfrac{%
\Delta _{k}}{\Delta _{0}}-\dfrac{\gamma _{k}}{\Delta _{0}}\right) , \\ 
\\ 
\gamma _{k}=\left( ^{\ast }b_{k}\left( ^{\ast }e^{k}\right) \right) \times 
\text{ }^{\ast }\left( \dint\limits_{0}^{k}f\left( x\right) e^{-x}dx\right) ,
\\ 
\\ 
\Delta _{0}=\left( ^{\ast }f\left( 0\right) +\text{ }^{\ast }f^{\text{ }%
\prime }\left( 0\right) +...+\text{ }^{\ast }f^{\text{ }\left( \mathbf{m}%
\right) }\left( 0\right) \right) , \\ 
\\ 
\Delta _{k}=\left( ^{\ast }f\left( k\right) +\text{ }^{\ast }f^{\text{ }%
\prime }\left( k\right) +...+\text{ }^{\ast }f^{\text{ }\left( \mathbf{m}%
\right) }\left( k\right) \right) . \\ 
\end{array}
& \text{ }\left( 3.2.52\right) 
\end{array}%
$

From Eq.(3.2.52) one obtain

\bigskip 

$\ \ \ 
\begin{array}{cc}
\begin{array}{c}
\\ 
\dsum\limits_{k=0}^{n}\left( ^{\ast }b_{k}\right) ^{\#}\times \left( ^{\ast
}e^{k}\right) ^{\#}=\dsum\limits_{k=0}^{n}\left( ^{\ast }b_{k}\right)
^{\#}\times \left( \dfrac{\Delta _{k}^{\#}}{\Delta _{0}^{\#}}-\dfrac{\gamma
_{k}^{\#}}{\Delta _{0}^{\#}}\right) , \\ 
\\ 
\gamma _{k}^{\#}=\left( \left( ^{\ast }b_{k}\left( ^{\ast }e^{k}\right)
\right) \times \text{ }^{\ast }\left( \dint\limits_{0}^{k}f\left( x\right)
e^{-x}dx\right) \right) ^{\#}, \\ 
\\ 
\Delta _{0}^{\#}=\left( ^{\ast }f\left( 0\right) +\text{ }^{\ast }f^{\text{ }%
\prime }\left( 0\right) +...+\text{ }^{\ast }f^{\text{ }\left( \mathbf{m}%
\right) }\left( 0\right) \right) ^{\#}, \\ 
\\ 
\Delta _{k}^{\#}=\left( ^{\ast }f\left( k\right) +\text{ }^{\ast }f^{\text{ }%
\prime }\left( k\right) +...+\text{ }^{\ast }f^{\text{ }\left( \mathbf{m}%
\right) }\left( k\right) \right) ^{\#}. \\ 
\end{array}
& \text{ \ }\left( 3.2.53\right) 
\end{array}%
$

\bigskip

\textbf{2. }We will choose $f\left( x\right) $ of the form

\bigskip 

\ $\ \ 
\begin{array}{cc}
\begin{array}{c}
\\ 
f\left( x\right) =\dfrac{1}{\left( P-1\right) !}x^{P-1}\cdot \left(
x-1\right) ^{P}\cdot \left( x-2\right) ^{P}\cdot ...\cdot \left( x-n\right)
^{P} \\ 
\end{array}
& \text{ \ }\left( 3.2.54\right) 
\end{array}%
$

\bigskip

where $P\in 
\mathbb{N}
$ is a prime number. Note that for $0\leq x\leq n\in 
\mathbb{N}
$ we have$\bigskip \ \ \ \ \ $

$\bigskip $

$\ \ 
\begin{array}{cc}
\begin{array}{c}
\\ 
\left\vert f\left( x\right) \right\vert \leq \dfrac{n^{\left( n+1\right)
\cdot P}}{\left( P-1\right) !}=\dfrac{\left[ A\left( n\right) \right] ^{P}}{%
\left( P-1\right) !}, \\ 
\\ 
A\left( n\right) =n^{n+1}. \\ 
\end{array}
& \text{ \ \ }\left( 3.2.55\right)%
\end{array}%
$

\bigskip

By using transfer from Eq.(3.2.54)-Eq.(3.2.55) one obtain

\bigskip

$%
\begin{array}{cc}
\begin{array}{c}
\\ 
^{\ast }\left( \text{ }f\left( x\right) \right) =\dfrac{1}{\left( \mathbf{P}%
-1\right) !}x^{\mathbf{P}-1}\cdot \left( x-1\right) ^{\mathbf{P}}\cdot
\left( x-2\right) ^{\mathbf{P}}\cdot ...\cdot \left( x-\mathbf{n}\right) ^{%
\mathbf{P}}, \\ 
\\ 
\left\vert ^{\ast }\left( \text{ }f\left( x\right) \right) \right\vert \leq 
\dfrac{\mathbf{n}^{\left( \mathbf{n}+1\right) \cdot \mathbf{P}}}{\left( 
\mathbf{P}-1\right) !}=\dfrac{\left[ A\left( \mathbf{n}\right) \right] ^{%
\mathbf{P}}}{\left( \mathbf{P}-1\right) !},A\left( \mathbf{n}\right) =%
\mathbf{n}^{\mathbf{n}+1} \\ 
\\ 
\mathbf{P,n\in }^{\ast }\mathbf{%
\mathbb{N}
}_{\infty }. \\ 
\end{array}
& \text{ \ }\left( 3.2.56\right) 
\end{array}%
$

Hence

$\ \ 
\begin{array}{cc}
\begin{array}{c}
\\ 
\overline{\overline{\#\text{-}\dsum\limits_{k=0}^{\infty }\left\vert \left(
\left( ^{\ast }b_{k}\right) ^{\#}\times \left( ^{\ast }e^{k}\right)
^{\#}\right) \times \text{ }\left[ ^{\ast }\left(
\dint\limits_{0}^{k}f\left( x\right) e^{-x}dx\right) \right]
^{\#}\right\vert }}\leq \\ 
\\ 
\leq \left( \overline{\overline{\#\text{-}\dsum\limits_{k=0}^{\infty
}\left\vert \left( ^{\ast }b_{k}\right) ^{\#}\right\vert \times \left(
^{\ast }e^{k}\right) ^{\#}}}\right) \times \dfrac{\left( \left[ A\left( 
\mathbf{n}\right) \right] ^{\mathbf{P}}\right) ^{\#}}{\left[ \left( \mathbf{P%
}-1\right) !\right] ^{\#}}\leq \\ 
\\ 
\leq \dfrac{\Delta _{\mathbf{d}}\mathbf{\times }\left( \left[ A\left( 
\mathbf{n}\right) \right] ^{\mathbf{P}}\right) ^{\#}}{\left[ \left( \mathbf{P%
}-1\right) !\right] ^{\#}}, \\ 
\end{array}
& \text{ \ }\left( 3.2.57\right)%
\end{array}%
$

\bigskip

It is easy to see that for $\mathbf{P\in }^{\ast }\mathbf{%
\mathbb{N}
}_{\infty }.$ large enough, this is less than $\epsilon ^{\#}$

for a given $\epsilon \approx 0,\epsilon \in $ $^{\ast }%
\mathbb{R}
$

Thus for $\mathbf{P\in }^{\ast }\mathbf{%
\mathbb{N}
}_{\infty }$ large enough by using Eq.(3.2.53) one obtain

\bigskip

\ \ \ $%
\begin{array}{cc}
\begin{array}{c}
\\ 
\overline{\overline{\#\text{-}\dsum\limits_{k=0}^{\infty }\left( ^{\ast
}b_{k}\right) ^{\#}\times \left( ^{\ast }e^{k}\right) ^{\#}}}= \\ 
\\ 
\left( ^{\ast }b_{0}\right) ^{\#}+\overline{\overline{\#\text{-}%
\dsum\limits_{k=1}^{\infty }\left( ^{\ast }b_{k}\right) ^{\#}\times \left( 
\dfrac{\Delta _{k}^{\#}}{\Delta _{0}^{\#}}-\dfrac{\gamma _{k}^{\#}}{\Delta
_{0}^{\#}}\right) }}, \\ 
\\ 
\overline{\overline{\#-\dsum\limits_{k=0}^{\infty }\left\vert \gamma
_{k}^{\#}\right\vert }}<\epsilon ^{\#},\epsilon \approx 0, \\ 
\\ 
\Delta _{0}^{\#}=\left( ^{\ast }f\left( 0\right) +\text{ }^{\ast }f^{\text{ }%
\prime }\left( 0\right) +...+\text{ }^{\ast }f^{\text{ }\left( \mathbf{m}%
\right) }\left( 0\right) \right) ^{\#}, \\ 
\\ 
\Delta _{k}^{\#}=\left( ^{\ast }f\left( k\right) +\text{ }^{\ast }f^{\text{ }%
\prime }\left( k\right) +...+\text{ }^{\ast }f^{\text{ }\left( \mathbf{m}%
\right) }\left( k\right) \right) ^{\#}. \\ 
\end{array}
& \text{ }\left( 3.2.58\right)%
\end{array}%
$

\bigskip\ Thus

$\ \ \ \ \ 
\begin{array}{cc}
\begin{array}{c}
\\ 
\left[ \overline{\overline{\#\text{-}\dsum\limits_{k=0}^{\infty }\left(
^{\ast }b_{k}\right) ^{\#}\times \left( ^{\ast }e^{k}\right) ^{\#}}}\right]
_{\varepsilon }= \\ 
\\ 
\left( ^{\ast }b_{0}\right) ^{\#}+\left[ \overline{\overline{\#\text{-}%
\dsum\limits_{k=1}^{\infty }\left( ^{\ast }b_{k}\right) ^{\#}\times \left( 
\dfrac{\Delta _{k}^{\#}}{\Delta _{0}^{\#}}-\dfrac{\gamma _{k}^{\#}}{\Delta
_{0}^{\#}}\right) }}\right] _{\varepsilon }, \\ 
\\ 
\overline{\overline{\#-\dsum\limits_{k=0}^{\infty }\left\vert \gamma
_{k}^{\#}\right\vert }}<\epsilon ^{\#},\epsilon \approx 0, \\ 
\end{array}
& \text{\ }\left( 3.2.59\right)%
\end{array}%
$

\bigskip By using \textbf{Theorem 1. }and Eq.(3.2.48) one obtain

$\ \ 
\begin{array}{cc}
\begin{array}{c}
\\ 
\overline{\overline{\#\text{-}\dsum\limits_{k=0}^{\infty }\left( ^{\ast
}b_{k}\right) ^{\#}\times \left( ^{\ast }e^{k}\right) ^{\#}}}=\varepsilon _{%
\mathbf{d}}. \\ 
\end{array}
& \text{ \ }\left( 3.2.60\right)%
\end{array}%
$

\bigskip

By using \textbf{Theorem 1.3.4} and Eq.(3.2.60) one obtain

\bigskip

$\ \ \ 
\begin{array}{cc}
\begin{array}{c}
\\ 
\left[ \overline{\overline{\#\text{-}\dsum\limits_{k=0}^{\infty }\left(
^{\ast }b_{k}\right) ^{\#}\times \left( ^{\ast }e^{k}\right) ^{\#}}}\right]
_{\varepsilon }=\left[ \varepsilon _{\mathbf{d}}\right] _{\varepsilon }. \\ 
\end{array}
& \text{ \ }\left( 3.2.61\right)%
\end{array}%
$

\ \ \ \ \ \ \ \ \ \ \ \ \ \ \ \ \ \ \ \ \ \ \ \ \ \ \ \ \ \ \ \ \ \ \ \ \ \
\ \ \ \ \ \ \ \ \ \ \ \ \ \ \ \ \ \ \ \ \ \ \ \ \ \ \ \ \ \ \ 

By using Eq.(3.2.59) and Eq.(3.2.60) one obtain

\bigskip

\bigskip $\ 
\begin{array}{cc}
\begin{array}{c}
\\ 
\left[ \varepsilon _{\mathbf{d}}\right] _{\varepsilon }=\left( ^{\ast
}b_{0}\right) ^{\#}+\left[ \overline{\overline{\#\text{-}\dsum%
\limits_{k=1}^{\infty }\left( ^{\ast }b_{k}\right) ^{\#}\times \left( \dfrac{%
\Delta _{k}^{\#}}{\Delta _{0}^{\#}}-\dfrac{\gamma _{k}^{\#}}{\Delta _{0}^{\#}%
}\right) }}\right] _{\varepsilon }, \\ 
\\ 
\overline{\overline{\#-\dsum\limits_{k=0}^{\infty }\left\vert \gamma
_{k}^{\#}\right\vert }}<\epsilon ^{\#}\left( \varepsilon \right) , \\ 
\\ 
\epsilon ,\varepsilon \approx 0. \\ 
\end{array}
& \text{ }\left( 3.2.62\right)%
\end{array}%
\ \ \ \ \ \ \ \ \ $

\bigskip

Hence

\bigskip

$\ 
\begin{array}{cc}
\begin{array}{c}
\\ 
\Delta _{0}^{\#}\times \left[ \varepsilon _{\mathbf{d}}\right] _{\varepsilon
}=\Delta _{0}^{\#}\times \left( ^{\ast }b_{0}\right) ^{\#}+ \\ 
\\ 
\left[ \overline{\overline{\#\text{-}\dsum\limits_{k=1}^{\infty }\left(
^{\ast }b_{k}\right) ^{\#}\times \left( \Delta _{k}^{\#}-\gamma
_{k}^{\#}\right) }}\left\vert \Delta _{0}^{\#}\times \left( ^{\ast }c\right)
^{\#}\right. \right] _{\varepsilon }, \\ 
\\ 
c\in 
\mathbb{R}
. \\ 
\end{array}
& \text{ \ }\left( 3.2.63\right)%
\end{array}%
$

\bigskip

Multiplying Eq.(3.2.63) by number$\ \ \Im ^{\#},$ where

$\Im =\left( m_{0},m_{0}\times m_{1},...,m_{0}\times m_{1}\times ...\times
m_{n},...\right) \ \ $gives

\bigskip $\ \ 
\begin{array}{cc}
\begin{array}{c}
\\ 
\Im ^{\#}\times \Delta _{0}^{\#}\times \varepsilon ^{\#}\times \varepsilon _{%
\mathbf{d}}= \\ 
\\ 
=\Im ^{\#}\times \Delta _{0}^{\#}\times \left( ^{\ast }b_{0}\right) ^{\#}+
\\ 
\\ 
+\Im ^{\#}\times \left[ \overline{\overline{\#\text{-}\dsum\limits_{k=1}^{n}%
\left( ^{\ast }b_{k}\right) ^{\#}\times \left( \Delta _{k}^{\#}-\gamma
_{k}^{\#}\right) }}\left\vert \Delta _{0}^{\#}\times \left( ^{\ast }c\right)
^{\#}\right. \right] _{\varepsilon }, \\ 
\end{array}
& \text{ \ }\left( 3.2.64\right)%
\end{array}%
$

$\bigskip $thus

$\ 
\begin{array}{cc}
\begin{array}{c}
\\ 
\Im ^{\#}\times \Delta _{0}^{\#}\times \varepsilon ^{\#}\times \varepsilon _{%
\mathbf{d}}= \\ 
\\ 
=\Im ^{\#}\times \Delta _{0}^{\#}\times \left( ^{\ast }b_{0}\right) ^{\#}+
\\ 
\\ 
\left[ \overline{\overline{\#\text{-}\dsum\limits_{k=1}^{n}\left( \Im
_{k}^{\#}\times \Delta _{k}^{\#}-\Im ^{\#}\gamma _{k}^{\#}\right) }}%
\left\vert \Im ^{\#}\times \Delta _{0}^{\#}\times \left( ^{\ast }c\right)
^{\#}\right. \right] _{\varepsilon }= \\ 
\\ 
=\Im ^{\#}\times \Delta _{0}^{\#}\times \left( ^{\ast }b_{0}\right) ^{\#}+
\\ 
\\ 
\left[ \overline{\overline{\#\text{-}\dsum\limits_{k=1}^{n}\left( \Im
_{k}^{\#}\times \Delta _{k}^{\#}-\Im ^{\#}\times \gamma _{k}^{\#}\right) }}%
\left\vert \Im ^{\#}\times \Delta _{0}^{\#}\times \left( ^{\ast }c\right)
^{\#}\right. \right] _{\varepsilon }, \\ 
\\ 
\Im ^{\#}\times \left( \overline{\overline{\#-\dsum\limits_{k=0}^{\infty
}\left\vert \gamma _{k}^{\#}\right\vert }}\right) <\epsilon ^{\#}\left(
\varepsilon \right) <\varepsilon ^{\#},\varepsilon \approx 0, \\ 
\\ 
\Im _{0}^{\#}=\Im ^{\#}\times \left( ^{\ast }b_{0}\right) ^{\#},\Im
_{k}^{\#}=\Im ^{\#}\times \left( ^{\ast }b_{k}\right) ^{\#},k=1,2,... \\ 
\end{array}
& \text{ \ \ }\left( 3.2.65\right)%
\end{array}%
$

$\ \ \ \ \ \ \ \ \ \ \ \ \ \ \ \ \ \ \ \ \ \ \ \ \ \ \ \ \ \ \ \ \ \ \ \ \ \
\ \ \ \ \ \ \ $

\bigskip \textbf{3.} If

$\ \ 
\begin{array}{cc}
\begin{array}{c}
\\ 
h\left( x\right) =\dfrac{g\left( x\right) \left( x-a\right) ^{\mathbf{P}}}{%
\mathbf{P}!} \\ 
\\ 
\mathbf{P\in }^{\ast }\mathbf{%
\mathbb{N}
}_{\infty }, \\ 
\end{array}
& \text{ \ \ \ }\left( 3.2.66\right)%
\end{array}%
$

\bigskip

where $g\left( x\right) $ is any hyper polynomial with hyper integer
coefficients

and $a\in $ $^{\ast }%
\mathbb{Z}
$ then the derivatives $h^{\left( j\right) }\left( a\right) =0$ for $0\leq j<%
\mathbf{P}$ and in

general $h^{\left( j\right) }\left( a\right) \in $ $^{\ast }%
\mathbb{Z}
$ for all $j\geq 0.$ Since $f\left( x\right) /\mathbf{P}$ has this form with

$a\in \left\{ 1,2,...,\mathbf{n}\right\} $ it follows that $f^{\text{ }%
\left( j\right) }\left( k\right) $ is an integer and is divisible

by $\mathbf{P},$ for all $j\geq 0$ and for $k\in \left\{ 1,2,...,\mathbf{n}%
\right\} .$ Thus $\mathbf{P}$ divides all terms

on the \textbf{RHS} of Eq.(3.2.65) having $k\neq 0.$

\textbf{4.} It remains to consider the terms with $k=0.$Note that $^{\ast
}f\left( x\right) $ has the form

\bigskip

$\ \ \ \ \ \ \ \ \ \ \ \ \ \ \ \ \ \ \ \ \ \ \ \ \ \ \ \ \ \ \ \ \ \ \ \ \ $

$\ 
\begin{array}{cc}
\begin{array}{c}
\\ 
^{\ast }f\left( x\right) =\dsum\limits_{\mathbf{j}=\mathbf{P}-1}^{\mathbf{m}}%
\dfrac{c_{\mathbf{j}}x^{\mathbf{j}}}{\left( \mathbf{P}-1\right) !} \\ 
\end{array}
& \text{ \ }\left( 3.2.67\right)%
\end{array}%
$

$\ $

where $c_{\mathbf{P}-1}=\left( \pm \text{ }\mathbf{n}!\right) ^{\mathbf{P}}$
and $c_{\mathbf{j}}\in $ $^{\ast }%
\mathbb{Z}
$ for all $\mathbf{j}\in $ $^{\ast }%
\mathbb{N}
.$Then $^{\ast }f^{\text{ }\left( \mathbf{j}\right) }\left( 0\right) =0$

for $\mathbf{j}<\mathbf{P}-1,$ $^{\ast }f^{\text{ }\left( \mathbf{P}%
-1\right) }\left( 0\right) =c_{\mathbf{P}-1}$ and $^{\ast }f^{\text{ }\left( 
\mathbf{j}\right) }\left( 0\right) =c_{\mathbf{j}}\cdot \mathbf{j}!/\left( 
\mathbf{P}-1\right) !$ for\textbf{\ }$\mathbf{j}\geq \mathbf{P}$ so

$\mathbf{P}$ divides $^{\ast }f^{\text{ }\left( \mathbf{j}\right) }\left(
0\right) $ if $\mathbf{j}\neq \mathbf{P}-1.$

\textbf{5.} The only term remaining on the \textbf{RHS} of Eq.(3.2.65) is

\bigskip

$\ \ 
\begin{array}{cc}
\begin{array}{c}
\\ 
\ \ \Im ^{\#}\times \Delta _{0}^{\#}\times \left( ^{\ast }b_{0}\right)
^{\#}\times \left( ^{\ast }f^{\text{ }\left( \mathbf{P}-1\right) }\left(
0\right) \right) ^{\#} \\ 
\\ 
=\left( \left( \pm \text{ }\mathbf{n}!\right) ^{\mathbf{P}}\right) ^{\#}. \\ 
\end{array}
& \text{ }\left( 3.2.68\right)%
\end{array}%
$

\bigskip

This term is not divisible by $\mathbf{P}^{\#}\mathbf{\in }^{\ast }\mathbf{%
\mathbb{N}
}_{\infty ,\mathbf{d}}$ if $\mathbf{P\in }^{\ast }\mathbf{%
\mathbb{N}
}_{\infty }$ is prime with

$\mathbf{P}>|^{\ast }b_{0}|\times \mathbf{n}.$Thus, we may choose $\mathbf{P}
$ so that

$\Im ^{\#}\times \left( \overline{\overline{\#-\dsum\limits_{k=0}^{\infty
}\left\vert \gamma _{k}^{\#}\right\vert }}\right) <\varepsilon ^{\#}$ and so
that in the \textbf{RHS} of Eq.(3.2.65),

$\mathbf{P}$ divides every term

\bigskip

$\ 
\begin{array}{cc}
\begin{array}{c}
\\ 
\ \ \Im ^{\#}\times \Delta _{0}^{\#}\times \left( ^{\#}b_{k}\right) \times
\left( ^{\ast }f^{\text{ }\left( \mathbf{j}\right) }\left( k\right) \right)
^{\#},\mathbf{j\in }^{\ast }\mathbf{%
\mathbb{N}
} \\ 
\end{array}
& \text{ \ }\left( 3.2.69\right)%
\end{array}%
$

$\ \ \ \ \ \ \ \ \ \ \ \ \ \ \ \ \ \ \ \ \ \ \ \ \ \ \ \ \ \ \ \ \ \ \ \ \ \ 
$

except for $\ \Im ^{\#}\times \Delta _{0}^{\#}\times \left(
^{\#}b_{0}\right) \times \left( ^{\ast }f^{\text{ }\left( \mathbf{P}%
-1\right) }\left( 0\right) \right) ^{\#}.$ Therefore the \textbf{RHS }has

representation $\Gamma ^{\#}+\Im ^{\#}\times \Delta _{0}^{\#}\times
\varepsilon ^{\#}\times \varepsilon _{\mathbf{d}}$ such that\ $\Gamma \in $ $%
^{\ast }%
\mathbb{R}
,\ \Gamma \geq 1.$

Thus one obtain

$\ 
\begin{array}{cc}
\begin{array}{c}
\\ 
\Im ^{\#}\times \Delta _{0}^{\#}\times \varepsilon ^{\#}\times \varepsilon _{%
\mathbf{d}}=\Gamma ^{\#}+\Im ^{\#}\times \Delta _{0}^{\#}\times \varepsilon
^{\#}\times \varepsilon _{\mathbf{d}} \\ 
\end{array}
& \text{ }\left( 3.2.70\right)%
\end{array}%
$\ \ \ 

\ \ \ \ \ \ \ \ \ \ \ \ \ \ \ \ \ \ \ \ \ \ \ \ \ \ \ \ \ \ \ \ \ 

This is a contradiction.This contradiction proves that $e$ is not

$w$-transcendental. Hence $e$ is $\#$-transcendental.

\bigskip

$\ \ \ \ \ \ \ \ \ \ \ \ \ \ \ \ \ \ \ \ \ \ \ \ \ \ \ \ \ \ \ \ \ \ \ \ \ \
\ \ \ \ \ \ \ \ \ \ \ \ \ \ $\ \ \ \ \ \ \ \ 

\bigskip $\ \ \ \ \ \ \ \ \ \ \ \ \ \ \ \ \ \ \ \ \ $

\bigskip

\section{III.3.Nonstandard generalization of the Lindeman Theorem.}

\textbf{Theorem 3.3.1.(Nonstandard Lindeman Theorem})The number $^{\ast }e$

cannot satisfy an equation of the next form:

\bigskip

\bigskip $\ 
\begin{array}{cc}
\begin{array}{c}
\\ 
a_{1}\cdot \left( ^{\ast }e\right) ^{\alpha _{1}}+a_{2}\cdot \left( ^{\ast
}e\right) ^{\alpha _{2}}+...+a_{N}\cdot \left( ^{\ast }e\right) ^{\alpha
_{N}}\approx 0, \\ 
\end{array}
& \text{ \ }\left( 3.3.1\right)%
\end{array}%
$

\bigskip

in which at least one coefficient $a_{n},n=1,2,...,N\in $ $^{\ast }%
\mathbb{N}
_{\infty }$ is different from

zero, no two exponents $\alpha _{n},n=1,2,...,N\in $ $^{\ast }%
\mathbb{N}
_{\infty }$ are equal, and all numbers

$\alpha _{n},n=1,2,...,N\in $ $^{\ast }%
\mathbb{N}
_{\infty }$ are hyperalgebraic.

\textbf{Proposition 3.3.1.}Let $\rho _{1},\rho _{2},...,\rho _{m},$ $m\in $ $%
^{\ast }%
\mathbb{N}
_{\infty }$ be the roots of the

hyperpolynomial equation $a\cdot z^{m}+b\cdot z^{m-1}+c\cdot z^{m-2}+...=0$
with integral

coefficients $a,b,c,...\in $ $^{\ast }%
\mathbb{Z}
.$ Then any symmetric hyperpolynomial in the

quantities $a\cdot \rho _{1},a\cdot \rho _{2},...,a\cdot \rho _{m}$ with
integral coefficients, is an hyperinteger.

\textbf{Proposition 3.3.2.}Suppose given a hyperpolynomial in $m\in $ $%
^{\ast }%
\mathbb{N}
_{\infty }$ variables $\alpha _{i_{1}},$

in $n\in $ $^{\ast }%
\mathbb{N}
_{\infty }$ variables $\beta _{i_{2}},...,$and in $k\in $ $^{\ast }%
\mathbb{N}
_{\infty }$ variables $\sigma _{i_{l}},l\in $ $^{\ast }%
\mathbb{N}
_{\infty }$ which is

symmetric in \ the $\alpha $'s in the $\beta $'s,...,and in the $\sigma $'s,
and which has hyperrational

coefficients. If the $\alpha $'s are chosen to be all the roots of a
hyperpolynomial

equation with rational coefficients, and similarly for the $\beta $'s, , and
for the $\sigma $'s

then the value of the polynomial is hyperrational.

\textbf{Definition 3.3.1.}A hyperpolynomial is said to be irreducible over
the rationale

if it cannot be factored into hyperpolynomials of lower degree with
hyperrational

coefficients.

\bigskip

\textbf{Definition 3.3.2.}If $\alpha _{1}$ is a root of an irreducible
hyperpolynomial equation with

hyperrational coefficients, whose other roots are $\alpha _{2},\alpha
_{3},,...,\alpha _{n},n\in ^{\ast }%
\mathbb{N}
$ then the

hyperalgebraic numbers $\alpha _{1},\alpha _{2},\alpha _{3},,...,\alpha _{n}$
are said to be the \textit{conjugates} of $\alpha _{1}.$

\bigskip

\textbf{Proposition 3.3.3.}Any hyperpolynomial with hyperrational
coefficients can be

factored into irreducible polynomials with hyperrational coefficients.

\textbf{Proposition 3.3.4.}Over the field $^{\ast }%
\mathbb{Q}
$, an hyperalgebraic number is a root of a

unique irreducible hyperpolynomial with hyperrational coefficients and
leading

coefficient unity. Such an equation has no multiple roots.

\bigskip

\textbf{Proposition 3.3.5.}The Van der Monde determinant $\det \left\vert
\left( \rho _{k}\right) ^{i-1}\right\vert $ vanishes only

if two or more of the $\rho $'s are equal.

\bigskip \bigskip

\section{\ III.4. The numbers $e$ and $\protect\pi $\ are analytically
independent.}

\bigskip

First of all, recall that is an entire function,in $2$ variables, with

coefficients in field $%
\mathbb{Q}
,$is a function $f(z_{1},z_{2})$ which is analytic in $G\subseteqq $\ $%
\mathbb{C}
\times 
\mathbb{C}
$\ 

\ \ \ \ \ \ \ \ \ \ \ \ \ \ \ \ \ \ \ \ \ \ \ \ \ \ \ \ \ \ \ \ \ \ \ \ \ \
\ \ \ \ \ \ \ \ \ \ \ \ \ \ \ \ \ \ \ \ \ \ \ \ \ \ \ \ \ \ \ \ \ \ \ \ \ \
\ \ \ \ \ \ $\ \ $

$\ \ \ \ \ \ 
\begin{array}{cc}
\begin{array}{c}
\\ 
f(z_{1},z_{2})=\sum_{i=0}^{\infty }\sum_{j=0}^{\infty
}c_{i,j}z_{1}^{i}z_{2}^{j}, \\ 
\\ 
c_{i,j}\in \text{ }%
\mathbb{Q}
. \\ 
\end{array}
& \text{ \ \ \ \ }\left( 3.4.1\right)%
\end{array}%
$

\bigskip

\textbf{Definition 3.3.1.}Two complex numbers $\alpha \in $ $%
\mathbb{C}
$ and $\beta \in $ $%
\mathbb{C}
$ are said to be

analytically dependent if there is a nonzero entire function $f(z_{1},z_{2})$
in $2$

variables, with hyperinteger coefficients $c_{i,j}\in $ $%
\mathbb{Q}
,$ such that $f(\alpha ,\beta )=0.$

Otherwise, $\alpha \in $ $%
\mathbb{C}
$ and $\beta \in $ $%
\mathbb{C}
$ are said to be analytically independent.

\ \ \ \ \ \ \ \ \ \ \ \ \ \ \ \ \ \ \ \ 

\bigskip \bigskip

\bigskip

\section{Appendix A.\ Hyper algebraic numbers.}

\section{\ \ \ \ \ \ \ \ \ \ \ \ \ \ \ \ \ \ \ \ \ \ \ \ \ \ \ \ \ \ \ \ \ \
\ \ \ \ \ \ \ \ \ \ \ \ \ \ \ \ \ \ \ \ \ \ \ \ \ \ \ \ \ \ \ \ \ \ \ \ \ \
\ \ \ \ \ \ \ \ \ \ \ \ \ \ \ \ \ \ \ 1.Definitions of symmetric polynomials
and symmetric functionals.}

\bigskip Consider a monic polynomial $P\left( z\right) $ in $z\in 
\mathbb{C}
$ of degree $n\in 
\mathbb{N}
$

$\ 
\begin{array}{cc}
\begin{array}{c}
\\ 
P\left( z\right) =1+a_{1}z+a_{2}z^{2}+...+a_{n-1}z^{n-1}+a_{n}z^{n} \\ 
\end{array}
& \text{ \ \ }\left( A.1.1\right)%
\end{array}%
$

\bigskip

There exist $n$ roots $z_{1},\ldots ,z_{n}$ of $P$ and that one is expressed

by the relation

\bigskip

$%
\begin{array}{cc}
\begin{array}{c}
\\ 
P\left( z\right) =1+a_{1}z+a_{2}z^{2}+...+a_{n-1}z^{n-1}+a_{n}z^{n}= \\ 
\\ 
=\left( 1-\dfrac{z}{z_{1}}\right) \left( 1-\dfrac{z}{z_{2}}\right) \cdot
\cdot \cdot \left( 1-\dfrac{z}{z_{n}}\right) = \\ 
\\ 
z_{1}^{-1}\cdot z_{2}^{-1}\cdot \cdot \cdot z_{n}^{-1}\left( z_{1}-z\right)
\cdot \left( z_{2}-z\right) \cdot \cdot \cdot \left( z_{n}-z\right) = \\ 
\\ 
=\left( -1\right) ^{n}z_{1}^{-1}\cdot z_{2}^{-1}\cdot \cdot \cdot
z_{n}^{-1}\left( z-z_{1}\right) \cdot \left( z-z_{2}\right) \cdot \cdot
\cdot \left( z-z_{n}\right)  \\ 
\\ 
\hat{P}\left( z\right) =P\left( z\right) /\left( -1\right)
^{n}z_{1}^{-1}\cdot z_{2}^{-1}\cdot \cdot \cdot z_{n}^{-1}= \\ 
\\ 
=\hat{a}_{n}z^{n}+\hat{a}_{n-1}z^{n-1}+...+\hat{a}_{2}z^{2}+\hat{a}_{1}z+%
\hat{a}_{0} \\ 
\end{array}
& \text{ \ }\left( A.1.2\right) 
\end{array}%
$

\bigskip Thus by comparison of the coefficients one finds\bigskip 

\bigskip\ $\ \ \ 
\begin{array}{cc}
\begin{array}{c}
\\ 
a_{1}=-\sum_{1\leq i\leq n}\dfrac{1}{z_{i}}, \\ 
\\ 
a_{2}=\sum_{1\leq i<j\leq n}\dfrac{1}{z_{i}z_{j}}, \\ 
\\ 
\cdot \cdot \cdot \cdot \cdot  \\ 
\\ 
a_{m}= \\ 
\\ 
\cdot \cdot \cdot \cdot \cdot  \\ 
\\ 
a_{n-1}=\left( -1\right) ^{n-1}\sum_{1\leq i_{1}<i_{2}<...<i_{n-1}\leq
n}\prod_{i\neq j}\dfrac{1}{z_{j}} \\ 
\\ 
a_{n}=\left( -1\right) ^{n}\prod_{1\leq i\leq n}\dfrac{1}{z_{i}}. \\ 
\end{array}
& \text{ \ }\left( A.1.3\right) 
\end{array}%
$

\textbf{Definition. }Let us defined $n$ polynomials expressed by the
relations

$\bigskip $

$\ $

\bigskip\ $%
\begin{array}{cc}
\begin{array}{c}
\\ 
e_{1}(z_{1},...,z_{n})=\sum_{1\leq i\leq n}\dfrac{1}{z_{i}}=-a_{1}, \\ 
\\ 
e_{2}(z_{1},...,z_{n})=\sum_{1\leq i<j\leq n}\dfrac{1}{z_{i}z_{j}}=a_{2}, \\ 
\\ 
\cdot \cdot \cdot \cdot \cdot \\ 
\\ 
e_{m}(z_{1},...,z_{n})= \\ 
\\ 
\cdot \cdot \cdot \cdot \cdot \\ 
\\ 
e_{n-1}(z_{1},...,z_{n})=a_{n-1}= \\ 
\\ 
e_{n}(z_{1},...,z_{n})=a_{n}=\left( -1\right) ^{n}\prod_{1\leq i\leq n}%
\dfrac{1}{z_{i}}. \\ 
\end{array}
& \text{ \ \ \ }\left( A.1.4\right)%
\end{array}%
$

The polynomial $e_{m}(z_{1},...,z_{n})$ is called the $m$-th symmetric
polynomial.

It has the following property:$\ \ \ \ \ \ \ \ \ \ \ \ \ \ \ \ \ \ \ \ \ \ \
\ \ \ \ \ \ \ \ \ \ \ \ \ \ \ \ \ \ \ \ $

$\bigskip $

$\ \ \ \ \ \ \ \ \ \ \ \ \ \ \ \ \ \ \ \ \ \ \ \ \ \ \ \ \ \ \ \ \ \ \ \ \ \
\ \ \ \ \ \ \ \ \ \ \ \ $

\bigskip \bigskip

\section{\ 2.Nonstandard polynomials.}

\bigskip

The set of natural,integer,rational, real, complex or any algebraic numbers
\ \ \ \ \ \ \ \ \ \ \ \ \ \ \ \ \ \ \ \ \ 

is denoted by $%
\mathbb{N}
,%
\mathbb{Z}
,%
\mathbb{Q}
,%
\mathbb{R}
,%
\mathbb{C}
,\Bbbk $ respectively, and their nonstandard extensions

$^{\ast }%
\mathbb{N}
,^{\ast }%
\mathbb{Z}
,^{\ast }%
\mathbb{Q}
,^{\ast }%
\mathbb{R}
,^{\ast }%
\mathbb{C}
,^{\ast }\Bbbk .$

\textbf{Definition A.2.1.}(\textbf{Nonstandard polynomials}) Nonstandard
polynomial of

hyper degree $\mathbf{d\in }^{\ast }%
\mathbb{N}
_{\infty }$ in $x$ with coefficients in nonstandard field $^{\ast }\Bbbk $
is an

expression defined by internal hyper finite sum of the form

\bigskip\ \ \ \ \ \ \ \ \ \ \ \ \ \ \ \ \ \ \ \ \ \ \ \ \ \ \ \ \ \ \ \ \ \
\ \ \ \ \ \ \ \ \ \ \ \ \ \ \ $\ 
\begin{array}{cc}
\begin{array}{c}
\\ 
f\left( x\right) =\sum_{j=0}^{\mathbf{d}}a_{j}x^{j}=a_{0}+a_{1}x+...+a_{%
\mathbf{d-}1}x^{\mathbf{d-}1} \\ 
\\ 
+a_{\mathbf{d}}x^{\mathbf{d}}\in \text{ }^{\ast }\Bbbk \left[ x\right] , \\ 
\\ 
\text{ \ }\forall j\left[ a_{j}\in \text{ }^{\ast }\Bbbk \right] . \\ 
\end{array}
& \text{ \ }\left( A.2.1\right) 
\end{array}%
$

\bigskip

\textbf{Definition A.2.2. }(\textbf{Algebraic hyper integers}) If $\alpha
\in $ $^{\ast }%
\mathbb{C}
$ is a root of a

monic nonstandard polynomial of hyper degree $\mathbf{d\in }^{\ast }%
\mathbb{N}
_{\infty },$namely a root of

a polynomial of the form

\ $\ 
\begin{array}{cc}
\begin{array}{c}
\\ 
f\left( x\right) =\sum_{j=0}^{\mathbf{d}}a_{j}x^{j}=a_{0}+a_{1}x+...+a_{%
\mathbf{d-}1}x^{\mathbf{d-}1} \\ 
\\ 
+a_{\mathbf{d}}x^{\mathbf{d}}\in \text{ }^{\ast }%
\mathbb{Z}
\left[ x\right] , \\ 
\\ 
\text{\ }\forall j\left[ a_{j}\in \text{ }^{\ast }%
\mathbb{Z}
\right]  \\ 
\end{array}
& \text{ \ \ \ \ \ }\left( A.2.2\right) 
\end{array}%
$\ \ \ \ \ 

\bigskip

and $\alpha $ is not the root of such a polynomial of hyper degree less then 
$\mathbf{d,}$

then $\alpha $ is colled an \textit{algebraic hyper integer} of hyper degree 
$\mathbf{d\in }^{\ast }%
\mathbb{N}
_{\infty }.$

\textbf{Definition A.2.3. }(\textbf{Nonstandard algebraic numbers}) An
nonstandard

algebraic number $\alpha $ of hyper degree $\mathbf{d\in }^{\ast }%
\mathbb{N}
_{\infty }$ is a root of a monic

nonstandard polynomial of hyper degree $\mathbf{d\in }^{\ast }%
\mathbb{N}
_{\infty },$ and not be the root of

an the nonstandard polynomial of hyper degree $\mathbf{d}_{1}\mathbf{\in }%
^{\ast }%
\mathbb{N}
_{\infty }$ less then $\mathbf{d.}$

\bigskip

\textbf{Remark.} We have to mach examles standard real numbers that are not

\textit{standard algebraic numbers}, such as $\ln 2$ and $\pi .$These are
examples of

\textit{standard transcendental} \textit{numbers}, which are not standard
algebraic

numbers

We will establish that every hyper finite extension of $^{\ast }%
\mathbb{Q}
$

\textbf{Definition A.2.4.} (\textbf{Simple} \textbf{hyper finite extentions
and nonstandard }

\textbf{polynomials}) If $\alpha \in $ $^{\ast }E$ an hyper finite extention
field of a given

nonstandard field $^{\ast }F,$then $\alpha $ is colled \textbf{hyper
algebraic }over $^{\ast }F$ if $f\left( \alpha \right) =0$

for some nonzero $f\left( x\right) \in F\left[ x\right] .$If $\alpha $

\bigskip

\bigskip

\bigskip

\section{References}

\bigskip

[1] \ \ \ Goldblatt,R.,Lectures on the Hyperreals. Springer-Verlag, New
York, \ \ \ \ \ \ \ \ \ \ \ \ \ \ \ \ \ \ \ \ \ \ 

\ \ \ \ \ \ \ \ \ NY, 1998

[2] \ \ \ Henle,J. and Kleinberg, E., Infinitesimal Calculus. Dover

\ \ \ \ \ \ \ \ Publications,\ Mineola,NY, 2003.

[3] \ \ \ Loeb, P. and Wolff, M.,Nonstandard Aalysis for the Working \ \ \ \
\ \ \ \ \ \ \ \ \ \ \ \ \ \ \ \ \ \ \ 

\ \ \ \ \ \ \ \ \ Mathematician.Kluwer Academic Publishers, Dordrecht, \ \ \
\ \ \ \ \ \ \ \ \ \ \ \ \ \ \ \ \ \ \ \ \ \ \ \ \ \ \ \ \ \ \ \ \ \ \ \ 

\ \ \ \ \ \ \ \ \ The Netherlands, 2000.

[4] \ \ \ Roberts, A. M.,Nonstandard Analysis. Dover Publications, \ \ \ \ \
\ \ \ \ \ \ \ \ \ \ \ \ \ \ \ \ \ \ \ \ \ \ \ \ \ \ \ \ 

\ \ \ \ \ \ \ \ Mineola, NY, 2003.

[5] \ \ \ Robinson, A.,Non-Standard Analysis (Rev. Ed.). Princeton

\ \ \ \ \ \ \ \ University Press, Princeton, NJ,

[6] \ \ \ Euler L.,Variae observationes circa series infinitas. 1737. 29p.

\ \ \ \ \ \ \ \ \ http://math.dartmouth.edu/\symbol{126}%
euler/docs/originals/E072.pdf

[7] \ \ \ Viader P.,Bibiloni L., Jaume P. On a series of Goldbach and
Euler.\ \ 

\ \ \ \ \ \ \ \ \ \
http://papers.ssrn.com/sol3/papers.cfm?abstract\_id=848586

[8] \ \ \ Edward C., (Author) How Euler Did it. Hardcover - Jul 3,2007.

[9] \ \ \ \ Metsnkyl T.,Catalan's conjecture: another old Diophantine

\ \ \ \ \ \ \ \ \ problem solved, Bull. Amer. Math. Soc. (N.S.) 41 (2004),

\ \ \ \ \ \ \ \ \ \ \ no. 1, 43--57.

[10] \ \ Goldblatt,R.,Topoi 2-ed.,NH,1984.

[11] \ \ Nelson E.\ Internal set theory: a new approach to nonstandard \ \ \
\ \ \ \ \ \ \ \ \ \ \ \ \ \ \ \ \ \ \ \ \ 

\ \ \ \ \ \ \ \ \ analysis. Bull. Amer. Math. Soc, 1977.

[12] \ \ Robinson A., Zakon E., Applications of Model Theory to Algebra, \ \
\ \ \ \ \ \ \ \ \ \ \ \ \ \ \ \ \ 

\ \ \ \ \ \ \ \ \ Analysis, and Probability.C.I.T.,Holt, Rinehart and
Winston, \ \ \ \ \ \ \ \ \ \ \ \ \ \ \ \ \ \ \ \ \ \ \ \ \ \ \ \ \ \ 

\ \ \ \ \ \ \ \ \ p.109 -122,1967.

[13] \ \ Hrb\u{c}ek K. Axiomatic foundations for nonstandard analysis, \ \ \
\ \ \ \ \ \ \ \ \ \ \ \ \ \ \ \ \ \ \ \ \ 

\ \ \ \ \ \ \ \ \ Fundamenta Mathematicae, vol. 98 (1978),pp.1-19.

[14] \ \ Jin R.,The sumset phenomenon.Proceedings of the AMS Volume \ \ \ \
\ \ \ \ \ \ \ \ \ \ \ \ \ \ \ \ \ \ 

\ \ \ \ \ \ \ \ \ 130,Number 3, Pages 855-861.

[15] \ \ Foukzon J., Nonstandard analysis and structural theorems of a

\ \ \ \ \ \ \ \ \ \ general nonlocallycompact Hausdorf Abelian group.

\ \ \ \ \ \ \ \ \ International Workshop on Topological Groups.Pamplona,

\ \ \ \ \ \ \ \ \ August 31st - September 2nd.

[16] \ Foukzon J., A definition of topological invariants for wild knots and
\ \ \ \ \ \ \ \ \ \ \ \ \ \ \ \ \ \ \ \ \ \ \ \ \ \ \ \ \ \ \ \ 

\ \ \ \ \ \ \ \ \ links by using non standard internal S. Albeverio
integral. \ \ \ \ \ \ \ \ \ \ \ \ \ \ \ \ \ \ \ \ \ \ \ \ \ \ \ \ \ \ \ \ \
\ \ 

\ \ \ \ \ \ \ \ \ The 22-nd Annual Geometric Topology Workshop.Colorado June
\ \ \ \ \ \ \ \ \ \ \ \ \ \ \ \ \ \ \ \ 

\ \ \ \ \ \ \ \ \ \ 9th-11th, 2005.

[17] \ \ Foukzon J., Generalized Pontryagian's duality theorem. \ \ \ \ \ \
\ \ \ \ \ \ \ \ \ \ \ \ \ \ \ \ \ \ \ \ \ \ \ \ \ \ \ \ \ \ \ \ \ \ 

\ \ \ \ \ \ \ \ \ \ 2006 International Conference on Topology and its

\ \ \ \ \ \ \ \ \ Applications,June 23-26, 2006, Aegion,Greece.

\ \ \ \ \ \ \ \ \ \ http://www.math.upatras.gr/\symbol{126}aegion/book.pdf

[18] \ \ Albeverio S., Fenstad J.E., Hoegh-Krohn R., Lindstrem T.

\ \ \ \ \ \ \ \ \ Nonstandard methods in stochastic analysis and mathematical

\ \ \ \ \ \ \ \ \ physics. Academic Press, Inc.1986, 590p.

[19] \ \ Kusraev A. G., Kutateladze S. S. Nonstandard Methods of Analysis.

\ \ \ \ \ \ \ \ \ \ Novosibirsk: Nauka, 1990; Dordrecht: Kluwer,1995.

[20] \ \ Laczkovich M. Conjecture and proof. 2001.

[21]\ \ \ Patterson E. M.The Jacobson radical of a pseudo-ring. \ \ \ \ \ \
\ \ \ \ \ \ \ \ \ \ \ \ \ \ \ \ \ \ \ \ \ \ \ \ \ \ \ \ \ \ \ \ \ \ 

\ \ \ \ \ \ \ \ \ Math. Zeitschr. 89, 348--364 (1965).\ \ 

[22] \ Nesterenko,Y.V., Philippon.Introduction to Algebraic Independence \ \
\ \ \ \ \ \ \ \ \ 

\ \ \ \ \ \ \ \ \ Theory. Series: Lecture Notes in Mathematics,Vol.1752
Patrice

\ \ \ \ \ \ \ \ \ (Eds.) 2001, XIII, 256 pp.,Softcover ISBN: 3-540-41496-7\
\ \ \ \ \ \ \ \ 

[23] \ Gonshor, H., Remarks on the Dedekind completion of a nonstandard

\ \ \ \ \ \ \ \ \ model of the reals.Pacific J. Math. Volume 118, Number1
(1985),

\ \ \ \ \ \ \ \ \ 117-132.

[24] \ Wattenberg, F. $\left[ 0,\infty \right] $-valued, translation
invariant measures on $%
\mathbb{N}
$ and

\ \ \ \ \ \ \ \ \ the Dedekind completion of $^{\ast }%
\mathbb{R}
.$Pacific J. Math. Volume 90,

\ \ \ \ \ \ \ \ \ Number 1 (1980), 223-247.

[25] \ Davis M. Applied Nonstandard Analysis.Wiley,New York,London,

\ \ \ \ \ \ \ \ \ Sydney, Toronto, 1977, xii + 181 pp.,

[26] \ Foukzon J.2006 Spring Central Sectional Meeting Notre Dame,IN,

\ \ \ \ \ \ \ \ April 8-9,2006 Meeting \#1016 The solution of one very old
problem in

\ \ \ \ \ \ \ \ \ transcendental numbers theory. Preliminary report. \ \ \ \
\ \ \ \ \ 

[27] \ Waldschmidt M., Algebraic values of analytic functions.Journal of

\ \ \ \ \ \ \ \ \ Computational and Applied Mathematics 160 (2003)
323--333.\ \ \ \ \ \ 

[28] \ 

[29]

[30] \ Laubenheimer P.,Schick T.,Stuhler U. Completions of countable

\ \ \ \ \ \ \ \ non-standard models of $%
\mathbb{Q}
.$ http://arxiv.org/abs/math/0604466v3

[31] \ Robinson A.,Nonstandard Arithmetic.Bull.Amer.Math.Soc.Volume 73,

\ \ \ \ \ \ \ \ Number 6 (1967),818-843.

[32] \ Harold G. Dales, W. Hugh Woodin Super-real fields: totally ordered

\ \ \ \ \ \ \ \ \ fields with additional structure.

[33] \ Arkeryd Leif O. Cutland, N.J. Henson C. \ Nonstandard analysis: theory

\ \ \ \ \ \ \ \ and applications.Ward (Eds.) 1997, 384 p., Hardcover

\ \ \ \ \ \ \ \ ISBN: 978-0-7923-4586-2

[34] Gonshor H., The ring of finite elements in a non-standard model of the

\ \ \ \ \ \ \ reals, J.London Math. Soc, (2) 3 (1971), 493-500.

[35] Cartier P., Functions polylogarithmes, nombres polyz\"{e}ta et groupes

\ \ \ \ \ \ \ prounipotents, S\'{e}m. Bourbaki, 53$^{\acute{e}me}$ ann\'{e}%
e, 2000--2001, no 884,

\ \ \ \ \ \ \ Mars 2001, 36 pp.

[36] \ Ap\'{e}ry R., Irrationalit\'{e} de $\zeta (2)$ et $\zeta (3),$ Ast%
\'{e}risque 61 (1979) 11--13.

[37]\ $\ $Balog,A. Perelli,A.Diophantine approximation by square-free
numbers.

\ \ \ \ \ \ \ Annali della Scuola Normale Superiore di Pisa - Classe di
Scienze,

\ \ \ \ \ \ \ S\'{e}r.4, 11 no. 3 (1984), p. 353-359. $\ \ \ \ \ \ \ \ \ \ \
\ \ \ \ \ \ \ \ \ \ \ \ \ $

\end{document}